\documentclass{smfbook}
\usepackage{sga2-smf}
\usepackage[ansinew]{inputenc}

\setboolean{original}{false}

\begin{document}
\frontmatter

\title[Cohomologie locale des faisceaux coh\'erents (SGA~2)] {S\'eminaire de G\'eom\'etrie Alg\'ebrique du~Bois~Marie \\
1962 \\
Cohomologie locale des~faisceaux coh\'erents et th\'eor\`emes de Lefschetz locaux et globaux \\(SGA~2)}

\alttitle{Cohomologie locale des faisceaux coh\'erents (SGA~2)}

\author[A.~Grothendieck]{Alexander Grothendieck\\
\smaller\smaller (r\'edig\'e par un groupe d'auditeurs)\\
\larger Augment\'e d'un expos\'e de Mme Mich\`ele Raynaud}

\begin{abstract}
Ce volume est une \'edition recompos\'ee et annot\'ee du livre \og Cohomologie locale des faisceaux coh\'erents et th\'eor\`emes de Lefschetz locaux et globaux (SGA~2)\fg, Advanced Studies in Pure Mathematics~2, North-Holland Publishing Company - Amsterdam, 1968, par A.~Grothendieck \textit{et al.}

Dans cet ouvrage, on donne des conditions n\'ecessaires et suffisantes de finitude des faisceaux de cohomologie locale d'un faisceau coh\'erent. Ces r\'esultats conduisent \`a des th\'eor\`emes d'alg\'ebrisation qui permettent en particulier d'obtenir, \`a l'aide de th\'eor\`emes de puret\'e \'egalement d\'emontr\'es dans le texte, des th\'eor\`emes de type Lefschetz pour le groupe fondamental ou de Picard.
\end{abstract}

\begin{altabstract}
This volume is a new updated edition of the book ``Cohomologie locale des faisceaux coh\'erents et th\'eor\`emes de Lefschetz locaux et globaux (SGA~2)'', Advanced Studies in Pure Mathematics~2, North-Holland Publishing Company - Amsterdam, 1968, by A.~Grothendieck \textit{et al.}

In this monograph are given necessary and sufficient conditions for the finiteness of the local cohomology sheaves of coherent sheaves. These results provide algebraization theorems leading in particular, with the help of purity results also proved in the text, to Lefschetz's theorem for both the fundamental group and the Picard group.
\end{altabstract}


\altkeywords{Algebraization, Lefschetz's theorem, local duality, local cohomology, depth, homotopical depth, fundamental group, Picard group}


\maketitle

\cleardoublepage
\renewcommand{\baselinestretch}{1.1}\normalfont
\makeschapterhead{Pr\'eface}
\thispagestyle{empty}
Le pr\'esent texte est une \'edition recompos\'ee \sisi{\ignorespaces}{et annot\'ee} du livre \og Cohomologie locale des faisceaux coh\'erents et th\'eor\`emes de Lefschetz locaux et globaux (\SGA 2)\fg, Advanced Studies in Pure Mathematics~2, North-Holland Publishing Company - Amsterdam, 1968, par A\ptbl Grothendieck \textit{et al.} C'est le deuxi\`eme volet du projet entam\'e par B\ptbl Edixhoven qui a r\'e\'edit\'e \SGA 1. Cette version reproduit le texte original \sisi{aussi fid\`element que possible.}{avec d'une part quelques modifications de forme corrigeant des erreurs typographiques et d'autre part des commentaires en bas de page, dits \og notes de l'\'editeur\fg (N.D.E.), pr\'ecisant l'\'etat actuel des connaissances sur les questions soulev\'ees dans la version originale. On a en outre donn\'e \c ca et l\`a quelques d\'etails suppl\'ementaires sur les preuves. Pour \'eviter des risques de confusion, les notes originales conservent leur syst\`eme de num\'erotation par des \'etoiles alors que les nouvelles sont num\'erot\'ees avec des nombres entiers. Les num\'eros de page de la version originale sont indiqu\'es dans la marge.

\medskip

Je remercie les math\'ematiciens qui ont assur\'e l'essentiel de la saisie initiale en \LaTeX 2e, \`a savoir L\ptbl Bayle, N\ptbl Borne, O\ptbl Brinon, J\ptbl Buresi, M\ptbl Chardin, F\ptbl Ducrot, P\ptbl Graftiaux, F\ptbl Han, P\ptbl Karwasz, L\ptbl Koelblen, D\ptbl Madore, S\ptbl Morel, D\ptbl Naie, B\ptbl Osserman, J\ptbl Riou et V\ptbl S\'echerre ainsi que C\ptbl Sabbah pour avoir mis ce texte au format de la SMF. Je remercie \'egalement J.-B\ptbl Bost, P\ptbl Colmez, O\ptbl Gabber, W\ptbl Fulton, S\ptbl Kleiman, F\ptbl Orgogozo, M\ptbl Raynaud et J.-P\ptbl Serre pour leurs commentaires et conseils.}

\bigskip
\begin{flushright}
L'\'editeur, Yves Laszlo.
\end{flushright}

\sisi{\vspace*{1cm}}{\newpage}

The present text is a new \sisi{\ignorespaces}{updated} edition of the book ``Cohomologie locale des faisceaux coh\'erents et th\'eor\`emes de Lefschetz locaux et globaux (\SGA 2)'', Advanced Studies in Pure Mathematics~2, North-Holland Publishing Company - Amsterdam, 1968, by A\ptbl Grothendieck \textit{et al.} It is the second part of the SGA project initiated by B\ptbl Edixhoven who made a new edition of \SGA 1. This version is meant to reproduce the original text \sisi{as close as possible.}{with some modifications including minor typographical corrections and footnotes from the editor (N.D.E.) explaining the current status of questions raised in the first edition. One has also given more details about some proofs. To avoid possible confusion, the original footnotes are numbered using stars whereas the new ones are numbered using integers. The page numbers of the original version are written in the margin of the text.

\medskip

Let me thank the mathematicians who have done most of the initial typeseting in \LaTeX 2e, namely L\ptbl Bayle, N\ptbl Borne, O\ptbl Brinon, J\ptbl Buresi, M\ptbl Chardin, F\ptbl Ducrot, P\ptbl Graftiaux, F\ptbl Han, P\ptbl Karwasz, L\ptbl Koelblen, D\ptbl Madore, S\ptbl Morel, D\ptbl Naie, B\ptbl Osserman, J\ptbl Riou and V\ptbl S\'echerre and also C\ptbl Sabbah for having adapted this text to the SMF layout. Let me thank also J.-B\ptbl Bost, P\ptbl Colmez, O\ptbl Gabber, W\ptbl Fulton, S\ptbl Kleiman, F\ptbl Orgogozo, M\ptbl Raynaud and J.-P\ptbl Serre for their comments and advice.}

\bigskip
\begin{flushright}
The editor, Yves Laszlo.
\end{flushright}

{\def\ndemark{}\def\sfootnotemark{}\let\\ \relax \chapterspace{2}\tableofcontents}

\mainmatter
\chapter*{Introduction} \label{I.0}

Nous\pageoriginale pr\'esentons ici, sous forme r\'evis\'ee et compl\'et\'ee, une r\'e\'edition par photo-offset du deuxi\`eme \sisi{S\'EMINAIRE DE G\'{E}OM\'{E}TRIE ALG\'EBRIQUE de l'INSTITUT des HAUTES \'ETUDES SCIENTIFIQUES}{S\'eminaire de G\'eom\'etrie Alg\'ebrique de l'Institut des Hautes \'Etudes Scientifiques} tenu en 1962 (mim\'eographi\'e).

Le lecteur se reportera \`a l'Introduction au premier de ces S\'eminaires (cit\'e \SGA 1 par la suite) pour les buts que poursuivent ces s\'eminaires, et leurs relations avec les \sisi{\'EL\'EMENTS DE G\'{E}OM\'{E}TRIE ALG\'EBRIQUE}{\'El\'ements de G\'eom\'etrie Alg\'ebrique}.

Le texte des expos\'es~\Ref{I} \`a~\Ref{XI} a \'et\'e r\'edig\'e \`a mesure, d'apr\`es mes expos\'es oraux et notes manuscrites, par un groupe d'auditeurs, comprenant I\ptbl \sisi{GIROGUITTI}{Giorgiutti}, J\ptbl \sisi{GIRAUD}{Giraud}, Mlle M\ptbl \sisi{JAFFE}{Jaffe} (devenue Mme M\ptbl \sisi{HAKIM}{Hakim}) et A\ptbl \sisi{LAUDAL}{Laudal}. Ces notes \`a l'origine \'etaient consid\'er\'ees comme devant \^etre provisoires et \`a circulation tr\`es limit\'ee, en attendant leur absorption par les EGA (absorption devenue maintenant pour le moins probl\'ematique, tout comme pour les autres parties des SGA). Comme il \'etait dit dans l'avertissement \`a l'\'edition primitive, ce caract\`ere \og confidentiel\fg des notes devait excuser certaines \og faiblesses de style\fg, sans doute plus manifestes dans le pr\'esent S\'eminaire \SGA 2 que dans les autres. J'ai essay\'e dans la mesure du possible d'y obvier dans la pr\'esente r\'e\'edition, par une r\'evision relativement serr\'ee du texte initial. J'ai notamment harmonis\'e les syst\`emes de num\'erotation des \'enonc\'es, employ\'es dans les diff\'erents expos\'es, en introduisant partout le m\^eme syst\`eme d\'ecimal, d\'ej\`a utilis\'e dans la plupart des expos\'es primitifs du pr\'esent \SGA 2, ainsi que dans toutes les autres parties des SGA. Cela m'a amen\'e en particulier \`a revoir enti\`erement la num\'erotation\nde{on a conserv\'e au maximum la num\'erotation originale, rajoutant l'adverbe bis quand des doublons ambigus se produisaient \c ca et l\`a.} des \'enonc\'es des\pageoriginale expos\'es~\Ref{III} \`a~\Ref{VIII}, (et par cons\'equent, des r\'ef\'erences auxdits expos\'es)\sfootnote{Il va sans dire que toutes les r\'ef\'erences \`a \SGA 2 qui figurent dans les parties des SGA publi\'ees dans la \sisi{SERIES IN PURE MATHEMATICS}{Series in Pure Mathematics} se rapporteront au pr\'esent volume, et non \`a l'\'edition primitive de \SGA 2 !}. J'ai essay\'e \'egalement d'extirper du texte primitif les principales erreurs de dactylographie ou de syntaxe (qui \'etaient nombreuses et g\^enantes). De plus, Mme M\ptbl\sisi{HAKIM}{Hakim} a bien voulu se charger de r\'e\'ecrire l'expos\'e~\Ref{IV} dans un style moins t\'el\'egraphique que l'expos\'e initial. Comme dans les autres r\'e\'editions des SGA, j'ai \'egalement ajout\'e un certain nombre de notes de bas de page, soit pour donner des r\'ef\'erences suppl\'ementaires, soit pour signaler l'\'etat d'une question pour laquelle des progr\`es ont \'et\'e faits depuis la r\'edaction du texte primitif. Enfin, ce S\'eminaire a \'et\'e augment\'e d'un nouvel expos\'e, \sisi{a}{\`a} savoir l'expos\'e~\Ref{XIV}, r\'edig\'e par Mme\sisi{.}{} \sisi{MICH\`ELE RAYNAUD}{Mich\`ele Raynaud} en 1967, qui reprend et compl\`ete des suggestions contenues dans les \og Commentaires \`a l'Expos\'e~\Ref{XIII}\fg (\Ref{XIII}~\Ref{XIII.6}) (r\'edig\'es en Mars~1963). Cet expos\'e reprend les th\'eor\`emes du type Lefschetz du point de vue de la cohomologie \'etale, en utilisant les r\'esultats sur la cohomologie \'etale expos\'es dans \SGA 4 et \SGA 5 (\`a para\^itre dans cette m\^eme collection \sisi{SERIES IN PURE MATHEMATICS}{Series in Pure Mathematics})\nde{en fait, ces s\'eminaires sont publi\'es chez Springer (num\'eros 269, 279, 305 et 589), mais, h\'elas, sont \'epuis\'es.}; il est donc \`a ce titre de nature moins \og \'el\'ementaire\fg que les autres expos\'es du pr\'esent volume, qui n'utilisent gu\`ere plus que la substance des chapitres~I \`a~III des EGA. Voici une esquisse du contenu du pr\'esent volume. L'expos\'e~\Ref{I} contient le sorite de la \og cohomologie \`a support dans $Y$\fg $H^*_Y(X, F)$, o\`u $Y$ est un ferm\'e d'un espace~$X$, cohomologie qui peut s'interpr\'eter comme une cohomologie de~$X$ \emph{modulo} l'ouvert $X-Y$, et qui est l'aboutissement d'une fort utile \og \emph{suite spectrale de passage du local au global}\fg \Ref{I}~\Ref{I.2.6}, faisant intervenir des \emph{faisceaux} de cohomologie \og \`a support dans~$Y$\fg $\sheaf{H}^*_Y(F)$\nde{on a conserv\'e les notations soulign\'ees pour les versions faisceautis\'ees de foncteurs, l'analogue calligraphique de $\Gamma$ n'\'etant pas clair.}. Ce formalisme peut dans de nombreuses questions\pageoriginale jouer un r\^{o}le de \og localisation\fg analogue \`a celui jou\'e par la consid\'eration de voisinages \og tubulaires\fg de $Y$ en g\'eom\'etrie diff\'erentielle. L'expos\'e~\Ref{II} \'etudie les notions pr\'ec\'edentes dans le cas des faisceaux quasi-coh\'erents sur les pr\'esch\'emas, l'expos\'e~\Ref{III} donne leur relation avec la notion classique de \emph{profondeur} (\Ref{III}~\Ref{III.3.3}).

Les expos\'es~\Ref{IV} et~\Ref{V} donnent des notions de \emph{dualit\'e locale}, qu'on peut comparer au th\'eor\`eme de dualit\'e projective de Serre (\Ref{XII}~\Ref{XII.1.1}); signalons que ces deux types de th\'eor\`emes de dualit\'e sont g\'en\'eralis\'es de fa\c con substantielle dans le s\'eminaire de \sisi{HARTSHORNE}{Hartshorne} (cit\'e en note de bas de page \`a la fin de \Exp \Ref{IV})\refstepcounter{toto}\nde{\label{noteconrad}le livre de Hartshorne comporte des erreurs de signes et surtout ne prouve pas vraiment la compatibilit\'e de la trace au changement de base. Conrad a compl\`etement repris ce travail, en prouvant cette compatibilit\'e, cruciale et hautement non triviale (\textit{Grothendieck duality and base change}, Lect. Notes in Math., vol.~1750, Springer-Verlag, Berlin, 2000). H\'elas, des erreurs subsistent (\cf deux pr\'epublications (Conrad~B., \og{Clarifications and corrections to ``\textit{Grothendieck duality and base change}''}\fg\ et \og{An addendum to Chapter 5 of ``\textit{Grothendieck duality and base change}''{\fg})}). Pour un aspect plus concret, avec une attention toute particuli\`ere \`a la notion de r\'esidu, voir les travaux de Lipman, en particulier (Lipman J., \textit{Dualizing sheaves, differentials and residues on algebraic varieties}, Ast\'erisque, vol.~117, Soci\'et\'e math\'ematique de France, 1984). Une preuve cat\'egorique du th\'eor\`eme de dualit\'e, bas\'ee sur le th\'eor\`eme de repr\'esentabilit\'e de Brown, a \'et\'e obtenue par Neeman (Neeman~A., {\og The Grothendieck duality theorem via Bousfield's techniques and Brown representability\fg}, \emph{J.~Amer. Math. Soc.} \textbf{9} (1996), \numero 1, p\ptbl 205--236).}

Les expos\'es~\Ref{VI} et~\Ref{VII} donnent des notions techniques faciles, utilis\'ees dans l'expos\'e~\Ref{VIII} pour prouver le \emph{th\'eor\`eme de finitude} (\Ref{VIII}~\Ref{VIII.2.3}), donnant des conditions n\'ecessaires et suffisantes, pour un faisceau coh\'erent $F$ sur un sch\'ema noeth\'erien $X$, pour que les faisceaux de cohomologie locale $\sheaf{H}^i_Y(F)$ soient coh\'erents pour $i\leq n$ (ou ce qui revient au m\^eme, pour que les faisceaux $\R^if_*(F|X-Y)$ soient coh\'erents pour $i\leq n-1$, o\`u $f\colon X-Y\to X$ est l'inclusion). Ce th\'eor\`eme est un des r\'esultats techniques centraux du S\'eminaire, et nous montrons dans l'expos\'e~\Ref{IX} comment un th\'eor\`eme de cette nature peut \^etre utilis\'e, pour \'etablir un \og th\'eor\`eme de comparaison\fg et un \og th\'eor\`eme d'existence\fg en g\'eom\'etrie formelle, en calquant et g\'en\'eralisant l'utilisation faite dans (\EGA III~\oldS\S 4 et~5) du th\'eor\`eme de finitude pour un morphisme propre.

On applique ces derniers r\'esultats dans~\Ref{X} et~\Ref{XI}, consacr\'es respectivement \`a des th\'eor\`emes du type Lefschetz pour le groupe fondamental, et pour\pageoriginale le groupe de \hbox{Picard}. Ces th\'eor\`emes consistent \`a comparer sous certaines conditions les invariants ($\pi_1$ ou $\Pic$) attach\'es respectivement \`a un sch\'ema $X$ et \`a un sous-sch\'ema $Y$ (jouant le r\^{o}le d'une section hyperplane), et \`a donner notamment des conditions o\`u ils sont isomorphes. Grosso modo, les hypoth\`eses faites servent \`a passer de $Y$ au compl\'et\'e formel de~$X$ le long de~$Y$, et \`a pouvoir appliquer ensuite les r\'esultats de~\Ref{IX} pour passer de l\`a \`a un \emph{voisinage} ouvert $U$ de $Y$ dans $X$. Pour pouvoir passer de $U$ \`a~$X$, il faut disposer encore de renseignements (type \og puret\'e\fg ou \og parafactorialit\'e{\fg}) pour les anneaux locaux de $X$ en les points de $Z=X-U$, (qui est un ensemble fini discret dans les cas envisag\'es). Ceci explique l'interaction dans les d\'emonstrations des expos\'es~\Ref{X}, \Ref{XI}, \Ref{XII} entre les r\'esultats locaux et globaux, notamment dans certaines r\'ecurrences. Les r\'esultats principaux obtenus dans~\Ref{X} et~\Ref{XI} sont les th\'eor\`emes \emph{de nature locale} \Ref{X}~\Ref{X.3.4} (\emph{th\'eor\`eme de puret\'e}) et \Ref{XI}~\Ref{XI.3.14} (\emph{th\'eor\`eme de parafactorialit\'e}). On notera que ces th\'eor\`emes sont d\'emontr\'es par des techniques cohomologiques, de nature essentiellement globale. Dans~\Ref{XII} on obtient en utilisant les r\'esultats locaux pr\'ec\'edents, les variantes globales de ces r\'esultats pour des sch\'emas projectifs sur un corps, ou plus g\'en\'eralement sur un sch\'ema de base plus ou moins quelconque; parmi les \'enonc\'es typiques, signalons \Ref{XII}~\Ref{XII.3.5} et \Ref{XII}~\Ref{XII.3.7}.

Dans~\Ref{XIII}, nous passons en revue quelques uns des nombreux probl\`emes et conjectures sugg\'er\'es par les r\'esultats et m\'ethodes du S\'eminaire. Les plus int\'eressants peut-\^etre concernent les th\'eor\`emes du type Lefschetz cohomologiques et homotopiques pour les espaces analytiques complexes, \cf \Ref{XIII} pages 26 et suivantes\nde{essentiellement toutes les conjectures \'enonc\'ees en~\Ref{XIII} et \Ref{XIV} sont maintenant prouv\'ees; voir les notes de bas de page de ces sections pour des r\'ef\'erences et commentaires.}. Dans le contexte de la cohomologie \'etale des sch\'emas\pageoriginale, les conjectures correspondantes sont prouv\'ees dans~\Ref{XIV} par une technique de dualit\'e qui devrait s'appliquer \'egalement dans le cas analytique complexe (\cf commentaires \Ref{XIII} p\ptbl25 et \Ref{XIV}~\Ref{XIV.6.4}). Mais les \'enonc\'es homotopiques correspondants dans le cas des espaces analytiques (et plus particuli\`erement les \'enonc\'es faisant intervenir le groupe fondamental) semblent exiger des techniques enti\`erement nouvelles (\cf \Ref{XIV}~\Ref{XIV.6.4}).

Je suis heureux de remercier tous ceux qui, \`a des titres divers, ont aid\'e \`a la parution du pr\'esent volume, dont les collaborateurs d\'ej\`a cit\'es dans cette Introduction. Plus particuli\`erement, je tiens \`a remercier Mlle \sisi{CHARDON}{Chardon} pour la bonne gr\^ace avec laquelle elle s'est acquitt\'ee de la t\^ache ingrate que constitue la pr\'eparation mat\'erielle du manuscrit d\'efinitif pour la photo-offset.

\bigskip\hfill Bures-sur-Yvette, Avril~1968
\par\smallskip \hfill A\ptbl Grothendieck.

\chapter[Les invariants cohomologiques globaux et locaux]
{Les invariants cohomologiques globaux et locaux relatifs \`a un sous-espace ferm\'e} \label{I}

\section{Les foncteurs $\Gamma_Z$, $\sheaf{\Gamma}_Z$} \label{I.1}
Soient\pageoriginale $X$ un espace topologique, $\ccat_X$ la cat\'egorie des faisceaux ab\'eliens sur $X$. Soient $\Phi$ une famille de supports au sens de Cartan on d\'efinit le foncteur $\Gamma_\Phi$ sur $\ccat_X$ par:
\begin{equation} \label{eq:I.1} { \Gamma_\Phi(F) = \text{ sous-groupe de } \Gamma(F) \text{ form\'e des sections } f \text{ telles que support } f \in \Phi}.
\end{equation}
Si $Z$ est une partie ferm\'ee de $X$, nous d\'esignons par abus de langage par $\sheaf{\Gamma}_Z$ le foncteur $\Gamma_\Phi$, o\`u $\Phi$ est l'ensemble des parties ferm\'ees de $X$ contenues dans $Z$. Donc on a:
\begin{equation} \label{eq:I.2}
{\Gamma_Z(F) = \text{ sous-groupe de\sisi{s}{} } \Gamma(F) \text{ form\'e des sections } f \text{ telles que support } f \subset Z}.
\end{equation}

Nous voulons g\'en\'eraliser cette d\'efinition au cas o\`u $Z$ est une partie \textit{localement ferm\'ee} de $X$, donc ferm\'ee dans une partie ouverte convenable $V$ de $X$. On posera dans ce cas:
\begin{equation} \label{eq:I.3} {\Gamma_Z(F)=\Gamma_Z(F\rest{V})}.
\end{equation}
Il faut v\'erifier que $\Gamma_Z(F)$ \og ne d\'epend pas\fg de l'ouvert choisi. Il suffit de montrer que si $V'$, $V\supset V' \supset Z$ est un ouvert, alors l'application $\rho_{V'}^V: F(V) \to F(V')$ applique $\Gamma_Z(F\rest{V})$ isomorphiquement sur $\Gamma_Z(F\rest{V'})$. Or
\begin{equation} \label{eq:I.4} { \Gamma_Z(F\rest{V}) = \ker \rho_{V-Z}^V}
\end{equation}
donc si $f\in \Gamma_Z(F\rest{V})$ et si $\rho_{V'}^V(f)=\rho_{V-Z}^V(f)=0$ alors $f=0$ puisque $(V', V-Z)$ est un recouvrement de $V$. De m\^eme, si $f'\in \Gamma_Z(F\rest{V'})$, alors $f'\in F(V')$ et\pageoriginale $0\in F(V-Z)$ d\'efinissent un $f\in F(V)$ tel que $\rho_{V'}^V(f)=f', f\in \Gamma_Z(F\rest{V})$, donc $\rho_{V'}^V$ induit un isomorphisme $\Gamma_Z(F\rest{V}) \to \Gamma_Z(F\rest{V'})$.

Notons que tout ouvert $W$ de $Z$ est induit par un ouvert $U$ de $X$ dans lequel $W$ est ferm\'e. Il en r\'esulte que $W \sisi{\rightsquigarrow}{\mto} \Gamma_W(F)$ d\'efinit un pr\'efaisceau sur $Z$, et on v\'erifie que c'est un faisceau que l'on notera $i^!(F)$, o\`u $i: Z \to X$ est l'immersion canonique. On trouve:
\begin{equation} \label{eq:I.5}
\Gamma_Z(F)=\Gamma(i^!(F)).
\end{equation}
Le faisceau $i^!(F)$ est un sous-faisceau de $i^*(F)$; en effet l'homomorphisme canonique:
$$
\Gamma(F\rest{U})=\Gamma(U, F) \to \Gamma(U\cap Z, i^*(F))$$
est injectif sur $\Gamma_{U\cap Z}(F\rest{U}) \subset \Gamma(F\rest{U})$. En r\'esumant, on~a le r\'esultat suivant:

\begin{proposition} \label{I.1.1} Il existe un sous-faisceau unique $i^!(F)$ de $i^*(F)$ tel que pour tout ouvert $U$ de $X$ tel que $U\cap Z$ soit ferm\'e dans $U$,
$$
\Gamma(F\rest{U})=\Gamma(U, F) \to \Gamma(U\cap Z, i^*(F))$$
induise un isomorphisme $\Gamma_{U \cap Z}(F\rest{U}) \to \Gamma(U\cap Z, i^!(F))$.
\end{proposition}

Notons que si $Z$ est un ouvert on aura simplement
\begin{equation} \label{eq:I.6} {i^!(F)=i^*(F)=F\rest{Z}, \ \Gamma_Z(F)=\Gamma(Z, F)}.
\end{equation}

Supposons \`a nouveau $Z$ quelconque. Alors pour un ouvert $U$ de $X$ variable, on voit que
$$
U \sisi{\rightsquigarrow}{\mto} \Gamma_{U\cap
Z}(F\rest{U})=\Gamma(U\cap Z, i^!(F))$$ est un faisceau sur $X$,
que nous noterons $\sheaf{\Gamma}_Z(F)$; de fa\c con pr\'ecise,
d'apr\`es la formule pr\'ec\'edente (exprimant que $i^!$ commute \`a la
restriction aux ouverts) on~a un isomorphisme
\begin{equation} \label{eq:I.7} {\sheaf{\Gamma}_Z(F)= i_*(i^!(F))}
\end{equation}
par d\'efinition\pageoriginale, on a, pour tout ouvert $U$ de $X$,
\begin{equation} \label{eq:I.8} {\Gamma(U, \sheaf{\Gamma}_Z(F)) = \Gamma_{U\cap Z}(F\rest{U})}.
\end{equation}
Notons ici une diff\'erence caract\'eristique entre le cas o\`u $Z$ est un ferm\'e, et celui o\`u $Z$ est un ouvert. Dans le premier cas, la formule (\Ref{eq:I.8}) nous montre que $\sheaf{\Gamma}_Z(F)$ peut \^etre consid\'er\'e comme un sous-faisceau de $F$, et on~a donc une \textit{immersion canonique}
\begin{equation*} \label{eq:I.8'} \tag{$8'$} {\sheaf{\Gamma}_Z(F) \hto F}.
\end{equation*}

Dans le cas o\`u $Z$ est ouvert, au contraire, on voit sur
(\Ref{eq:I.6}) que le deuxi\`eme membre de (\Ref{eq:I.8}) est
$\Gamma(U\cap Z, F)$, donc re\c coit $\Gamma(U, F)$, donc on~a un
\textit{homomorphisme canonique} en sens inverse du pr\'ec\'edent:
\begin{equation*} \label{eq:I.8''} \sisi{\tag{8}}{\tag{$8"$}} {F\to \sheaf{\Gamma}_Z(F)},\nde{la r\'ef\'erence originale \'etait (8).}
\end{equation*}
qui n'est autre d'ailleurs que l'homomorphisme canonique
$$F\to i_*i^*(F),
$$
compte tenu de l'isomorphisme
\begin{equation*} \label{eq:I.6bis} \tag{6 bis} {\sheaf{\Gamma}_Z(F)\simeq i_*i^*(F)}
\end{equation*}
d\'eduit de (\Ref{eq:I.6}) et (\Ref{eq:I.7}).

Bien entendu, pour $F$ variable, $\Gamma_Z(F)$, $\sheaf{\Gamma}_Z(F)$, $i^!(F)$ peuvent \^etre consid\'er\'es comme des foncteurs en $F$, \`a valeur respectivement dans la cat\'egorie des groupes ab\'eliens, des faisceaux ab\'eliens sur $X$, des faisceaux ab\'eliens sur $Z$. Il est parfois commode d'interpr\'eter le foncteur
$$i^!: \ccat_X \to \ccat_Z$$
comme le foncteur adjoint d'un foncteur bien connu
$$i_!: \ccat_Z \to \ccat_X$$
d\'efini par la proposition suivante:

\begin{proposition} \label{I.1.2}
Soit\pageoriginale $G$ un faisceau ab\'elien sur $Z$. Alors il existe un sous-faisceau unique de $i_*(G)$, soit $i_!(G)$ tel que pour tout ouvert $U$ de $X$ l'isomorphisme (identique)
$$
\Gamma(U\cap Z, G)= \Gamma(U, i_*(G))$$
d\'efinisse un isomorphisme
$$
\Gamma_{\Phi_{U\cap Z, U}}(U\cap Z, G)=\Gamma(U, i_! (G)),
$$
o\`u $\Phi_{U\cap Z, U}$ d\'esigne l'ensemble des parties de $U\cap Z$ qui sont ferm\'ees dans $U$.
\end{proposition}

La v\'erification se r\'eduit \`a noter que le premier membre est un faisceau pour $U$ variable, \ie que la propri\'et\'e pour une section de $i_*(G)$ sur $U$, consid\'er\'ee comme une section de $G$ sur $U\cap Z$, d'\^etre \`a support ferm\'e \textit{dans $U$} est de \textit{nature locale sur $U$}. Le faisceau $i_!(G)$ qu'on vient de d\'efinir est connu aussi sous le nom de: \textit{faisceau d\'eduit de~$G$ en prolongeant par $0$} en dehors de $Z$, \cf \sisi{Godement}{\cite{Godement}}. En particulier, si $Z$ est ferm\'e, on a
\begin{equation} \label{eq:I.9} {i_!(G)=i_*(G)};
\end{equation}
mais dans le cas g\'en\'eral, l'injection canonique $i_!(G) \to i_*(G)$ n'est pas un isomorphisme, comme il est bien connu d\'ej\`a pour $Z$ ouvert. \'Evidemment, $i_!(G)$ d\'epend fonctoriellement de $G$ (et c'est m\^eme un foncteur exact en G). Ceci dit, on a:

\begin{proposition} \label{I.1.3} Il existe un isomorphisme de bifoncteurs en $G, F$ ($G$ faisceau ab\'elien sur $Z$, $F$ faisceau ab\'elien sur $X$):
\begin{equation} \label{eq:I.10} {\Hom(i_!(G), F) = \Hom (G, i^!(F))}.
\end{equation}
\end{proposition}

Pour d\'efinir un tel isomorphisme, il revient au m\^eme de d\'efinir des homomorphismes fonctoriels
$$i_!i^!(F) \to F, \ G \to i^!i_!(G),
$$
satisfaisant aux conditions de compatibilit\'e bien connues (\cf par exemple l'expos\'e de Shih au s\'eminaire Cartan sur les op\'erations cohomologiques).

Se\pageoriginaled rappelant que $i_!$ est exact, donc transforme monomorphismes en monomorphismes, on en conclut:

\begin{corollaire} \label{I.1.4}
Si $F$ est injectif, $i^!(F)$ est injectif, donc $\sheaf{\Gamma}_Z(F) \sisi{\to}{=} i_*i^!(F)$ est \'egalement injectif.
\end{corollaire}

Rempla\c cant $X$ par un ouvert variable $U$ de $X$, on conclut
aussi de \Ref{I.1.3}

\begin{corollaire} \label{I.1.5} On a un isomorphisme fonctoriel en $F, G$:
\begin{equation} \label{eq:I.11}
\SheafHom{(i_!(G), F)} =i_*(\SheafHom{(G, i^!(F))}.
\end{equation}
\end{corollaire}

Prenant pour $G$ le faisceau constant sur $Z$ d\'efini par $\ZZ$, soit $\ZZ_Z$, \Ref{I.1.3} et \Ref{I.1.5} se sp\'ecialisent en

\begin{corollaire} \label{I.1.6} On a des isomorphismes fonctoriels en $F$:
\begin{equation} \label{eq:I.12}
\begin{array}[c]{l} {\Gamma_Z(F)=\Hom(\ZZ_{Z, X}, F)}, \\
{\sheaf{\Gamma}_Z(F)=\SheafHom(\ZZ_{Z, X}, F)},
\end{array}
\end{equation}
o\`u $\ZZ_{Z, X}= i_!(\ZZ_Z)$ est le faisceau ab\'elien sur $X$ d\'eduit du faisceau constant sur $Z$ d\'efini par $\ZZ$, en prolongeant par $0$ en dehors de $Z$.
\end{corollaire}

\begin{remarque} \label{I.1.7}
Supposons que $X$ soit un espace annel\'e, et munissons $Z$ du faisceau d'anneaux $\Oo_Z=i^{-1}(\OX)$, enfin d\'esignons par $\ccat_X$ et $\ccat_Z$ la cat\'egorie des Modules sur $X$ \resp $Z$. Alors les consid\'erations pr\'ec\'edentes s'\'etendent mot \`a mot, en prenant pour~$F$ un Module sur $X$ et pour $G$ un modules sur $Z$, et en interpr\'etant en cons\'equence les \'enonc\'es \Ref{I.1.3} \`a \Ref{I.1.6}.
\end{remarque}

Pour finir ces g\'en\'eralit\'es, examinons ce qui se passe quand on change la partie localement ferm\'ee $Z$. Soit $Z'\subset Z$ une autre partie localement ferm\'ee, et soient
$$j: Z' \to Z, i': Z' \to X, i'=ij$$
les inclusions\pageoriginale canoniques. Alors on~a des isomorphismes fonctoriels:
\begin{equation} \label{eq:I.13} {(ij)^!=j^!i^!, \quad (ij)_!=i_!j_!}.
\end{equation}
Le premier isomorphisme (\Ref{eq:I.13}) d\'efinit un isomorphisme fonctoriel
\begin{equation} \label{eq:I.14} {\Gamma_{Z'}(F)=\Gamma(Z', (ij)^!(F))\simeq \Gamma(Z', j^!(i^!(F)))=\Gamma_{Z'}(i^!(F))}.
\end{equation}

Supposons maintenant que $Z'$ soit \textit{ferm\'e} dans $Z$, et
$$Z''=Z-Z'$$
son compl\'ementaire dans $Z$, qui est ouvert dans $Z$ donc localement ferm\'e dans $X$. L'inclusion canonique (\Ref{eq:I.8'}) appliqu\'ee \`a $i^!(F)$ sur $Z$ muni de $Z'$ nous d\'efinit, gr\^ace \`a (\Ref{eq:I.14}), un homomorphisme canonique injectif fonctoriel
\begin{equation} \label{eq:I.15} {\Gamma_{Z'}(F) \to \Gamma_Z(F)}.
\end{equation}
Si on remplace dans (\Ref{eq:I.14}) $Z'$ par $Z''$ et utilise (\Ref{eq:I.8''}), on trouve un homomorphisme canonique fonctoriel:
\begin{equation*} \label{eq:I.15'} \tag{$15'$} {\Gamma_Z(F) \to \Gamma_{Z''}(F)}.
\end{equation*}

\begin{proposition} \label{I.1.8} Sous les conditions pr\'ec\'edentes, la suite d'homomorphismes fonctoriels:
\begin{equation} \label{eq:I.16} {0 \to \Gamma_{Z'}(F) \to \Gamma_Z(F) \to \Gamma_{Z''}(F)}
\end{equation}
est exacte. Si $F$ est flasque, la suite reste exacte en mettant un z\'ero \`a droite.
\end{proposition}

\begin{proof}
Rempla\c cant $X$ par un ouvert $V$ dans lequel $Z$ soit ferm\'e, on
est ramen\'e au cas o\`u $Z$ est ferm\'e, donc $Z'$ ferm\'e. Alors $Z''$
est ferm\'e dans l'ouvert $X-Z'$, et on~a une inclusion canonique
$$
\Gamma_{Z''}(F) \to \Gamma(X-Z', F),
$$
et l'exactitude de (\Ref{eq:I.16}) signifie simplement que les sections de $F$ \`a support dans $Z'$ sont celles dont la restriction \`a $X-Z'$ est nulle.

Lorsque\pageoriginale $F$ est flasque, tout \'el\'ement de $\Gamma_{Z''}(F)$, consid\'er\'e comme section de $F$ sur $X-Z'$, peut se prolonger en une section de $F$ sur $X$, et cette derni\`ere aura \'evidemment son support dans $Z$, ce qui prouve qu'alors le dernier homomorphisme dans (\Ref{eq:I.16}) est surjectif.
\skipqed
\end{proof}

\begin{corollaire} \label{I.1.9} On a une suite exacte fonctorielle
\begin{equation*} \label{eq:I.16bis} \tag{16 bis} {0 \to \sheaf{\Gamma}_{Z'}(F) \to \sheaf{\Gamma}_Z(F) \to \sheaf{\Gamma}_{Z''}(F)},
\end{equation*}
et si $F$ est flasque, cette suite reste exacte en mettant un $0$ \`a droite.
\end{corollaire}

On peut interpr\'eter (\Ref{I.1.8}) en termes de r\'esultats sur les
foncteurs $\Hom$ et $\SheafHom$ via \Ref{I.1.6}, de la fa\c con
suivante. Notons d'abord que si $G$ est un faisceau ab\'elien sur
$Z$, induisant les faisceaux $j^*(G)$ et $k^*(G)$ sur $Z'$ \resp
$Z''$ (o\`u $j: Z' \to Z$ et $k: Z'' \to Z$ sont les injections
canoniques), on~a une suite exacte canonique de faisceaux sur $X$:
\begin{equation} \label{eq:I.17} {0 \to k^*(G)_X \to G_X \to j^*(G)_X \to 0}
\end{equation}
o\`u pour simplifier les notations, l'indice $X$ d\'esigne le faisceau sur $X$ obtenu en prolongeant par $0$ dans le compl\'ementaire de l'espace de d\'efinition du faisceau envisag\'e. La suite exacte (\Ref{eq:I.17}) g\'en\'eralise une suite exacte bien connue lorsque $Z=X$ (\cf \sisi{Godement}{\cite{Godement}}), et s'en d\'eduit d'ailleurs en \'ecrivant la suite exacte en question sur~$Z$, et en appliquant le foncteur $i_!$. Prenant $G=\ZZ_{Z}$, on conclut en particulier:

\begin{proposition} \label{I.1.10}
Sous les conditions pr\'ec\'edentes, on~a une suite exacte de faisceaux ab\'eliens sur $X$:
\begin{equation} \label{eq:I.18} {0 \to \ZZ_{Z{''}, X} \to \ZZ_{Z, X} \to \ZZ_{Z', X} \to 0}.
\end{equation}
Ceci pos\'e, les deux suites exactes \Ref{I.1.8} et \Ref{I.1.9} ne sont autres que les suites exactes d\'eduites de (\Ref{eq:I.18}) par application du foncteur $\Hom(-, F)$ \resp $\SheafHom(-, F)$.
\end{proposition}

Cela redonne une d\'emonstration \'evidente du fait que les suites (\Ref{eq:I.16}) et (\Ref{eq:I.16bis}) restent exactes en mettant un z\'ero \`a droite, pourvu que $F$ soit \textit{injectif}.

\section{Les foncteurs $\H_Z^*(X, F)$ et $\SheafH_Z^*(F)$} \label{I.2}
\pageoriginale
\sisi{\ignorespaces}{\setcounter{equation}{18}}

\begin{definition} \label{I.2.1}
On d\'enote par $\H_Z^*(X, F)$ et ${\SheafH}_Z^*(F)$ les foncteurs d\'eriv\'es en $F$ des foncteurs $\Gamma_Z(F)$ \resp $\sheaf{\Gamma}_Z(F)$.
\end{definition}

Ce sont des foncteurs cohomologiques, \`a valeurs dans la cat\'egorie des groupes ab\'eliens \resp dans la cat\'egorie des faisceaux ab\'eliens sur $X$. Lorsque $Z$ est ferm\'e, $H_Z^*(X, F)$ n'est autre par d\'efinition que $H^*_{\Phi}(X, F)$ o\`u $\Phi$ d\'esigne la famille des parties ferm\'ees de $X$ contenues dans $Z$. Lorsque $Z$ est ouvert, on va voir que $H_Z^*(X, F)$ n'est autre que $H^*(Z, F)=H^*(Z, F\rest{Z})$, gr\^ace \`a la proposition suivante.

\begin{proposition}[Th\'eor\`eme d'excision] \label{I.2.2}
Soit $V$ une partie ouverte de $X$ contenant~$Z$. Alors on~a un isomorphisme de foncteurs cohomologiques en $F$:
\begin{equation} \label{eq:I.19}
\H^*_Z(X, F) \to \H^*_Z(V, F\rest{V}).
\end{equation}
\end{proposition}

En effet, on~a un isomorphisme fonctoriel $\Gamma_Z^X \simeq \Gamma_Z^V j^!$, o\`u $j: V \to X$ est l'inclusion et o\`u $j^!$ est donc le foncteur restriction (\cf (\Ref{eq:I.14})). Ce dernier est exact, et transforme injectifs en injectifs par \Ref{I.1.4}, d'o\`u aussit\^{o}t l'isomorphisme (\Ref{eq:I.19}).

Lorsque $Z$ est ouvert, on peut prendre $V=Z$ et on trouve:

\begin{corollaire} \label{I.2.3bof}
Supposons $Z$ ouvert, alors on~a un isomorphisme de foncteurs cohomologiques:
\begin{equation} \label{eq:I.20}
\H_Z^*(X, F) = \H^*(Z, F).
\end{equation}
\end{corollaire}

On conclut des isomorphismes \Ref{I.1.6} et des d\'efinitions (\cf \sisi{Tohoku}{\cite{Tohoku}}):

\label{I.2.3bis}
\begin{enonce*}{Proposition 2.3\sisi{}{~bis\ndemark}}\sisi{}{\ndetext{la proposition porte le num\'ero 2.3 dans l'\'edition originale.}}
On a\pageoriginale des isomorphismes de foncteurs cohomologiques:
\begin{align} \label{eq:I.21}
\H_Z^*(X, F)&\simeq \Ext^*(X; \ZZ_{Z, X}, F),\\
\label{eq:I.21bis}\tag{21 bis}
\SheafH_Z^*(F)&\simeq \SheafExt^*(\ZZ_{Z, X}, F).
\end{align}
\end{enonce*}

On peut donc appliquer les r\'esultats de \sisi{Tohoku}{\cite{Tohoku}} sur les $\Ext$ de Modules. Signalons d'abord l'interpr\'etation suivante des faisceaux $\SheafH_Z^*(F)$ en terme de groupes globaux $\H_Z^*(X, F)$:

\begin{corollaire} \label{I.2.4}
$\SheafH_Z^*(F)$ est canoniquement isomorphe au faisceau associ\'e au pr\'efaisceau
$$
U \sisi{\rightsquigarrow}{\mto} \H^*_{Z\cap U}(U, F\rest{U}).
$$
\end{corollaire}

En particulier, en utilisant \sisi{\ignorespaces}{le} corollaire \Ref{I.2.3bof}, on trouve:

\begin{corollaire} \label{I.2.5}
Supposons $Z$ ouvert, alors on~a un isomorphisme de foncteurs cohomologiques:
\begin{equation} \label{eq:I.22} {\SheafH_Z^*(F)=\R^*i_*i^*(F)}
\end{equation}
(o\`u $i: Z \to X$ est l'inclusion).
\end{corollaire}

La suite spectrale des $\Ext$ donne l'importante suite spectrale:

\begin{theoreme} \label{I.2.6} On a une suite spectrale fonctorielle en $F$, aboutissant \`a $\H_Z^*(X, F)$ et de terme initial
\begin{equation} \label{eq:I.23} { \E_2^{p, q}(F)=\H^p(X, \SheafH_Z^q(F))}.
\end{equation}
\end{theoreme}

\sisi{}{\setcounter{equation}{22}}

\begin{remarques} \label{I.2.7} Il r\'esulte aussit\^{o}t de \Ref{I.2.4} que les faisceaux $\SheafH_Z^q(F)$ sont nuls dans $X-\bar{Z}$, et \'egalement nuls dans l'int\'erieur de $Z$ pour $q\neq 0$ (donc pour un tel $q$, $\SheafH_Z^q(F)$ est m\^eme port\'e par la fronti\`ere \sisi{$\overset{\circ}{Z}$}{} de $Z$).

Par suite, le deuxi\`eme membre de (\Ref{eq:I.23}) peut s'interpr\'eter comme un groupe de cohomologie sur $\bar{Z}$. Nous utiliserons \Ref{I.2.6} dans le cas o\`u $Z$ est ferm\'e\pageoriginale dans $X$, et o\`u le deuxi\`eme membre de (\Ref{eq:I.23}) peut s'interpr\'eter comme un groupe de cohomologie calcul\'e sur $Z$:
\begin{equation} \label{eq:I.23bis}\sisi{\tag{\Ref{eq:I.23}}}{\tag{23 bis}} {\E_2^{p, q}(F) = \H^p(Z, \SheafH_Z^q(F))}.\nde{l'\'equation \'etait num\'erot\'ee (23) dans l'\'edition orginale.}
\end{equation}
\end{remarques}
\sisi{\ignorespaces}{\setcounter{equation}{23}}

Notons aussi que lorsque $Z$ est ouvert, la suite spectrale \Ref{I.2.6} n'est autre que la suite spectrale de Leray pour l'application continue $i: Z \to X$, compte tenu de l'interpr\'etation \Ref{I.2.5} dans le calcul du terme initial de la suite spectrale de Leray.

\bigskip Reprenons la suite exacte \sisi{(\Ref{I.1.10})}{(\Ref{eq:I.18})\nde{la r\'ef\'erence originale \'etait~(\Ref{I.1.10}).}}, elle donne naissance \`a une suite exacte des $\Ext$ (\cf \sisi{Tohoku}{\cite{Tohoku}}):

\begin{theoreme} \label{I.2.8}
Soient $Z$ une partie localement ferm\'ee de $X$, $Z'$ une partie ferm\'ee de $Z$ et $Z''=Z-Z'$. Alors on~a une suite exacte fonctorielle en $F$:
\begin{equation} \label{eq:I.24}
\begin{gathered}
0 \to \H_{Z'}^0(X, F) \to \H_Z^0(X, F) \to \H_{Z''}^0(X, F)\lto{\partial} \H^1_{Z'}(X, F) \to \H^1_Z(X, F) \dots \\
\phantom{0}\dots \H_{Z'}^i(X, F) \to \H_Z^i(X, F) \to \H_{Z''}^i(X, F)\lto{\partial} \H^{i+1}_{Z'}(X, F)\dots.
\end{gathered}
\end{equation}
\end{theoreme}

Rappelons comment on peut obtenir cette suite exacte. Soit $C(F)$ une r\'esolution injective de $F$, alors la suite exacte \sisi{(\Ref{I.1.10})}{(\Ref{eq:I.18})\nde{voir note pr\'ec\'edente.}} donne naissance \`a la suite exacte
\begin{equation} \label{eq:I.25} {0 \to \Gamma_{Z'}(C(F)) \to \Gamma(C(F)) \to \Gamma_{Z''}(C(F)) \to 0},
\end{equation}
(qui n'est autre que celle d\'efinie dans \Ref{I.1.8}). On en conclut une suite exacte de cohomologie, qui n'est autre que (\Ref{eq:I.24}).

Le cas le plus important pour nous est celui o\`u $Z$ est ferm\'e (et
on peut d'ailleurs toujours s'y ramener en rempla\c cant $X$ par
un ouvert $V$ dans lequel $Z$ est ferm\'e). Alors $Z'$ est ferm\'e,
$Z''$ est ferm\'e dans l'ouvert $X-Z'$, et on peut \'ecrire
\begin{equation} \label{eq:I.26} {\H^i_{Z''}(X, F)=\H^i_{Z''}(X-Z', F\rest{X-Z'})},
\end{equation}
ce qui nous permet d'\'ecrire la suite exacte (\Ref{eq:I.24}) en termes de cohomologies \`a support\pageoriginale dans un ferm\'e donn\'e. Le cas le plus fr\'equent est celui o\`u $Z=X$. Posant alors pour simplifier $Z'=A$, on trouve:

\begin{corollaire} \label{I.2.9}
Soit $A$ une partie ferm\'ee de $X$. Alors on~a une suite exacte fonctorielle en $F$:
\begin{equation} \label{eq:I.27}
\begin{gathered}
0 \to \H^0_A(X, F) \to \H^0(X, F) \to \H^0(X-A, F) \lto{\partial} \H^1_A(X, F) \dots\\
\phantom{0}\dots \H^i_A(X, F) \to \H^i(X, F) \to \H^i(X-A, F) \lto{\partial} \H^{i+1}_A(X, F)\dots.
\end{gathered}
\end{equation}
\end{corollaire}

Cette suite exacte montre que le groupe de cohomologie $\H^i_A(X,
F)$ joue le r\^{o}le d'un groupe de cohomologie relative de $X\mod
X-A$, \`a coefficients dans $F$. C'est \`a ce titre qu'elle
s'introduisit de fa\c con naturelle dans les applications. En \og
faisceautisant\fg (\Ref{eq:I.24}) et (\Ref{eq:I.27}), ou en
proc\'edant directement, on trouve\sisi{}{ en tenant compte de ce
que le faisceau associ\'e \`a $U\mto H^i(U, F)$ est nul si $i>0$}:

\begin{corollaire} \label{I.2.10}
Sous les conditions de \Ref{I.2.8}, on~a une suite exacte fonctorielle en~$F$:
\begin{equation*} \label{eq:I.24bis} \tag{24 bis}
\dots \SheafH^i_{Z'}(F) \to \SheafH^i_{Z}(F) \to \SheafH_{Z''}^i(F) \lto{\partial} \SheafH_{Z'}^{i+1}(F) \dots
\end{equation*}
\end{corollaire}

\begin{corollaire} \label{I.2.11}
Soit $A$ une partie ferm\'ee de $X$, alors on~a une suite exacte fonctorielle en $F$:
\begin{equation} \label{eq:I.28}
0 \to \SheafH_A^0(F) \to F \to f_*(F\rest{X-A}) \lto{\partial} \SheafH_A^1 \to 0,
\end{equation}
et des isomorphismes canoniques, pour $i\geq 2$:
\begin{equation} \label{eq:I.29}
\SheafH_A^i(F)= \SheafH_{X-A}^{i-1}(F)=\R^{i-1}f_*(F\rest{X-A}),
\end{equation}
o\`u $f: (X-A) \to X$ est l'inclusion.
\end{corollaire}

Cela d\'efinit donc $\SheafH^0_A(F)$ et $\SheafH_A^1(F)$ respectivement comme $\ker$ et $\coker$ de l'homomorphisme canonique
$$F \to f_*f^*(F)= f_*(F\rest{X-A}),
$$
et les $\SheafH_A^i(F)$ ($i\geq 2$) en terme des foncteurs d\'eriv\'es de $f_*$.

\begin{corollaire} \label{I.2.12}
Soit $F$\pageoriginale un faisceau ab\'elien sur $X$. Si $F$ est flasque, alors pour toute partie localement ferm\'ee $Z$ de $X$ et tout entier $i\neq 0$, on~a $\H^i_Z(X, F)=0$, $\SheafH^i_Z(F)=0$. Inversement, si pour toute partie ferm\'ee $Z$ de $X$ on~a $\H_Z^1(X, F)=0$, alors $F$ est flasque.
\end{corollaire}

Supposons que $F$ est flasque, alors $F$ induit un faisceau flasque sur tout ouvert, donc pour prouver $\H^i_Z(X, F)=0$ pour $i>0$, on peut supposer $Z$ ferm\'e, et alors l'assertion r\'esulte de la suite exacte \sisi{\Ref{I.2.9}}{(\Ref{eq:I.27})\nde{la r\'ef\'erence originale \'etait~(\Ref{I.2.9}).}}. On en conclut pour tout $Z$ localement ferm\'e, en \og faisceautisant\fg \ie appliquant \Ref{I.2.4}, que $\SheafH^i_Z(F)=0$ pour $i>0$. Inversement, supposons $\H^1_Z(X, F)=0$ pour tout ferm\'e $Z$, alors la suite exacte \sisi{\Ref{I.2.9}.}{(\Ref{eq:I.27})\nde{la r\'ef\'erence originale \'etait~(\Ref{I.2.9}).}} prouve que pour tout tel $Z$, $\H^0(X, F) \to \H^0(X-Z, F)$ est surjectif, ce qui signifie que $F$ est flasque.

Combinant \Ref{I.2.6} et \Ref{I.2.8}, on va en d\'eduire:

\begin{proposition} \label{I.2.13}
Soient $F$ un faisceau ab\'elien sur $X$, $Z$ une partie ferm\'ee de $X$, $U=X-Z$, $N$ un entier. Les conditions suivantes sont \'equivalentes:
\begin{enumeratei}
\item
$\SheafH^i_Z(F)=0$ pour $i\leq N$.
\item
Pour tout ouvert $V$ de $X$, consid\'erant l'homomorphisme canonique
$$
\H^i(V, F) \to \H^i(V\cap U, F),$$
cet homomorphisme est:
\begin{enumerate}
\item[\textup{a)}] bijectif pour $i<N$,
\item[\textup{b)}] injectif pour $i=N$.
\end{enumerate}
\end{enumeratei}
\noindent
(Lorsque $N>0$, on peut dans \textup{(ii)} se borner \`a exiger \textup{a)}).
\end{proposition}

Pour prouver (i) \ALORS (ii), on est ramen\'e, gr\^ace \`a la nature locale des $\SheafH_Z^i(F)$, \`a prouver~le

\begin{corollaire} \label{I.2.14}
Si la condition \Ref{I.2.13} (i) est v\'erifi\'ee, alors
$$
\H^i(X, F) \to \H^i(U, F)$$
est bijectif pour $i<N$, injectif pour $i=N$.
\end{corollaire}

En effet, en vertu de la suite exacte (\Ref{eq:I.27}), cela signifie aussi $\H^i_Z(X, F)=0$ pour $i\leq N$, et cette\pageoriginale relation est une cons\'equence imm\'ediate de la suite spectrale \sisi{\Ref{I.2.6}}{(\Ref{eq:I.23bis})\nde{voir note pr\'ec\'edente.}}.

R\'eciproquement, l'hypoth\`ese \Ref{I.2.13} (ii) signifie que pour tout ouvert $V$ de $X$, on a
$$
\H^i_{Z\cap V}(V, F\rest{V}) =0 \text{ pour } i\leq N,
$$
ce qui implique \Ref{I.2.13} (i) gr\^ace \`a \Ref{I.2.4}. \sisi{Si d'ailleurs $N>0$, on peut dire aussi que (ii) a) implique (en passant aux faisceaux associ\'es) que $F \to f_*(F\rest{U})$ est un isomorphisme, et que $\SheafH_U^i(F)=0 $ pour $1\leq i<N$.}{Si d'ailleurs $N>0$, l'hypoth\`ese b) est superflue. En effet, si $N=1$, l'hypoth\`ese a) et (\Ref{eq:I.28}) assurent la nullit\'e de $\SheafH^i_Z(F)=0$ pour $i\leq N$. Si $N>1$, l'hypoth\`ese a) pour $i=N-1>0$ et (\Ref{eq:I.29}) assurent la nullit\'e de $\SheafH^i_Z(F)$ pour $i\leq N$.}

Compte tenu du \Ref{I.2.11} cela prouve encore \Ref{I.2.13} (i)...

\begin{remarquestar}
Soit $Y \to X$ une immersion ferm\'ee, et supposons que localement elle est de la forme $\{0\}\times Y \subset R^n \times Y$. Supposons que $F$ est un faisceau localement constant sur $x$, alors on trouve
\begin{equation} \label{eq:I.30}
\SheafH_Y^i (F)\simeq
\begin{cases}
0 &\text{si } i \neq n\\
F\otimes \sheaf{T}_{Y, X}& \text{si } i=n, \text{ o\`u } \sheaf{T}_{Y, X} \simeq \SheafH_Y^n(\ZZ_X)
\end{cases}
\end{equation}
est un faisceau extension \`a $X$ d'un faisceau sur $Y$ localement isomorphe \`a $\ZZ_Y$, appel\'e \og faisceau d'orientation normale de $Y$ dans $X$\fg.
\end{remarquestar}

Utilisant la suite spectrale \sisi{\Ref{I.2.6}}{(\Ref{eq:I.23bis})\nde{la r\'ef\'erence originale \'etait~\Ref{I.2.6}.}}, on trouve dans ce cas:
\begin{equation} \label{eq:I.31}
\H^i_Y(X, F)\simeq \H^{i-n}(Y, F\otimes \sheaf{T}_{Y, X}),
\end{equation}
et on retrouve l'\og homomorphisme de Gysin\refstepcounter{toto}\label{hgysin}\fg:
\begin{equation} \label{eq:I.32} {\H^j(Y, F\otimes \sheaf{T}_{Y, X}) \to \H^{j+n}(X, F)}.
\end{equation}

\sisi{%
\section*{Bibliographie}
\noindent
R\ptbl Godement, Th\'eorie des faisceaux, Ac. Scient. et Ind.~1252 (cit\'e: Godement).
\par\smallskip\noindent
A\ptbl Grothendieck, sur quelques points d'alg\`ebre homologique, Tohoku Math. Journal, Vol.~9 p\ptbl 119--221, \S~1957, (cit\'e: Tohoku).
}{%

}

\chapter[Faisceaux quasi-coh\'erents sur les pr\'esch\'emas]{Application aux faisceaux quasi-coh\'erents sur les pr\'esch\'emas} \label{II}

\numberwithin{equation}{section}

\begin{supproposition} \label{II.1}
Soient\pageoriginale $X$ un pr\'esch\'ema, $Z$ une partie localement ferm\'ee de la forme $Z=U-V$, o\`u $U$ et $V$ sont deux parties ouvertes de $X$ telles que $V\subset U$ et que les immersions canoniques $U\to X$, $V\to X$ soient quasi-compactes. Alors pour tout Module quasi-coh\'erent $F$ sur $X$, les faisceaux ${\SheafH}_Z^i(F)$ sont quasi-coh\'erents.
\end{supproposition}

D'apr\`es (\sisi{\Ref{I}, }{\Ref{I}~}\Ref{eq:I.24}), il existe une suite exacte de cohomologie relative
$$
\to {\SheafH}_U^{i-1}(F) \to {\SheafH}_V^{i}(F) \to {\SheafH}_Z^{i}(F) \to {\SheafH}_U^{i}(F) \to {\SheafH}_V^{i+1}(F) \to.
$$

D'apr\`es (\EGA III~1.4.17), pour que les ${\SheafH}_Z^i(F)$ soient quasi-coh\'erents il suffit donc que les ${\SheafH}_U^i(F)$ et les ${\SheafH}_V^i(F)$ le soient. On peut donc supposer $Z$ ouvert et l'immersion canonique $j:Z\to X$ quasi-compacte.

Puisque $Z$ est ouvert on~a (\sisi{\Ref{I}, }{\Ref{I}~}\Ref{eq:I.22}) un isomorphisme canonique:
$$
{\SheafH}_Z^{i}(F) \simeq \R^ij_{*}(F_{|Z})
$$
mais $j$ est s\'epar\'e (\EGA I~5.5.1) et quasi-compact donc (\EGA III~1.4.10) les $\R^ij_{*}(F|Z)={\SheafH}_Z^{i}(F)$ sont quasi-coh\'erents, ce qui ach\`eve la d\'emonstration.

\begin{supcorollaire} \label{II.2}
Soit $Z$ une partie ferm\'ee de $X$ telle que l'immersion canonique \hbox{$X-Z\to X$} soit quasi-compacte, alors les Modules ${\SheafH}_Z^i(F)$ sont quasi-coh\'erents.
\end{supcorollaire}

\begin{supcorollaire} \label{II.3}
Si $X$ est localement noeth\'erien, alors pour toute partie localement ferm\'ee $Z$ de $X$, et tout Module $F$ quasi-coh\'erent sur $X$, les ${\SheafH}_Z^i(F)$ sont quasi-coh\'erents.
\end{supcorollaire}

R\'esulte imm\'ediatement du corollaire \Ref{II.2} et de (\EGA I~6.6.4).

\begin{supcorollaire} \label{II.4}
Supposons\pageoriginaled que $X$ soit le spectre d'un anneau $A$ et soient $U$ un ouvert quasi-compact de $X$, $Y=X-U$, $F$ un Module quasi-coh\'erent sur $X$, il existe un isomorphisme de foncteurs cohomologiques en $F$:
\begin{equation} \label{eq:II.4.1} {\SheafH}_Y^i(F) = \widetilde{(\H_Y^i(X, F))}.
\end{equation}

On a en outre une suite exacte fonctorielle en $F$:
\begin{equation} \label{eq:II.4.2}
0 \to \H_Y^0(X, F) \to \H^0(X, F) \to \H^0(U, F) \to \H_Y^1(X, F) \to 0
\end{equation}
et des isomorphismes fonctoriels en $F$:
\begin{equation} \label{eq:II.4.3}
\H_Y^i(X, F) \simeq \H^{i-1}(U, F), \quad i\geq 2.
\end{equation}
\end{supcorollaire}

D'apr\`es \sisi{(\Ref{II.1})}{le corollaire~\Ref{II.2}}, les ${\SheafH}_Y^i(F)$ sont quasi-coh\'erents, puisque $X$ est affine on~a donc $\H^p(X, {\SheafH}_Y^i(F))=0$ si $p>0$. La suite spectrale \sisi{(\Ref{I}, \Ref{eq:I.23})}{(\Ref{I}~\Ref{eq:I.23})} d\'eg\'en\`ere, donc
$$
\H_Y^i(X, F) = \Gamma({\SheafH}_Y^i(F)).
$$
L'\'egalit\'e (\Ref{eq:II.4.1}) r\'esulte alors de (\EGA I~1.1.3.7), (\Ref{eq:II.4.2}) et (\Ref{eq:II.4.3}) de la suite exacte de cohomologie (\sisi{I(\Ref{eq:I.27})}{\Ref{I}~\Ref{eq:I.27}}) et de ce que $\H^i(X, F)=0$ si $i>0$, puisque $X$ est affine.

Avec les hypoth\`eses \sisi{de \eqref{II.4}}{de \Ref{II.4}}\nde{dans un souci de coh\'erence et de clart\'e, on n'a num\'erot\'e entre parenth\`eses que les \'equations.}, $U$ est r\'eunion finie d'ouverts affines $X_f$, on peut donc trouver un id\'eal $I$ engendr\'e par un nombre fini d'\'el\'ements $f_{\alpha}$ et d\'efinissant~$Y$, soit $\bbf=(f_{\alpha})$. Avec les notations de (\EGA III~1)\nde{rappelons que $H^{\boule}(\bbf, M)$ est la cohomologie de Koszul $H^{\boule}(\Hom(K_{\boule}(\bbf), M))$ de $\bbf$ (\EGA III~1.1.2) \`a valeurs dans $M$ et que $H^{\boule}((\bbf), M)$ est la limite (\loccit, 1.1.6.5) $\varinjlim_n H^{\boule}(\bbf^n, M)$, les morphisme de transition \'etant induits par les morphisme naturels $K_{\boule}(\bbf^{n+1})\to K_{\boule}(\bbf^n)$ (\loccit, 1.1.6).} on a:

\begin{supproposition} \label{II.5}
Supposons que $X$ soit le spectre d'un anneau $A$, soient $\bbf=(f_{\alpha})$ une famille finie d'\'el\'ements de $A$, $Y$ la partie ferm\'ee de $X$ qu'ils d\'efinissent, $M$ un $A$-module, $F$ le faisceau associ\'e \`a $M$. On a alors des isomorphismes de $\partial$-foncteurs en $M$:
\begin{equation} \label{eq:II.5.1} \H^i((\bbf), M) \simeq \H_Y^i(X, F).
\end{equation}
\end{supproposition}

\noindent
(Nous noterons aussi $\H_J^i(M)=\H_Y^i(X, F)$, si $Y$ est la partie ferm\'ee de $X=\Spec A$ d\'efinie par un id\'eal $J$ de $A$).

Pour $i=0$\pageoriginale et $i=1$, on utilise les suites exactes (\Ref{eq:II.4.2}) et (\EGA III~1.4.3.2); si $i\geq2$, on utilise (\Ref{eq:II.4.3}) et (\EGA III~1.4.3.\sisi{2}{1}). Cela nous donne des isomorphismes fonctoriels en $M$. On v\'erifie qu'\`a un signe pr\`es ne d\'ependant que de $i$, ils sont compatibles avec l'op\'erateur bord, d'o\`u l'existence de l'isomorphisme de $\partial$-foncteurs (\Ref{eq:II.5.1}).

Soient maintenant $X$ un pr\'esch\'ema, $Y$ une partie ferm\'ee de $X$ et $f:Y\to X$ l'inclusion, $I$ un id\'eal quasi-coh\'erent d\'efinissant $Y$ dans $X$. Soit $F$ un faisceau sur $X$.

On a vu qu'il existe des isomorphismes de $\partial$-foncteurs en $F$
\begin{align}
\label{eq:II.5.*} \tag{$*$} \Ext_{\OX}^i(X;f_{*}f^{-1}(\OX), F) &\to \H_Y^i(X, F) \\
\label{eq:II.5.**} \tag{$**$} \SheafExt_{\OX}^i(f_{*}f^{-1}(\OX), F) &\to \SheafH_Y^i(F).
\end{align}

Soient $n, m$ des entiers tels que $m\geq n\geq 0$, on d\'esigne par $i_{n, m}$ l'application canonique: $\Oo_{Y_m}=\OX/I^{m+1} \to\OX/I^{n+1}=\Oo_{Y_n}$, et par $j_n$ l'application: $f_{*}f^{-1}(\OX)\to\Oo_{Y_n}$. Les $(\Oo_{Y_n}, i_{n, m})$ forment un syst\`eme projectif et les $j_n$ sont compatibles avec les $i_{n, m}$.

En appliquant le foncteur $\Ext_{\OX}^i(X;\cdot, F)$, on en d\'eduit un morphisme
$$
\varphi': \varinjlim_{n} \Ext_{\OX}^i(X;\Oo_{Y_n}, F) \to \Ext_{\OX}^i(X;f_{*}f^{-1}(\OX), F);
$$
on montre facilement que c'est un morphisme de foncteurs cohomologiques en $F$. Le morphisme
$$
\varphi: \varinjlim_{n} \Ext_{\OX}^i(X;\Oo_{Y_n}, F) \to \H_Y^i(X, F),
$$
compos\'e de $\varphi'$ et de (\Ref{eq:II.5.*}), est donc lui aussi un morphisme de foncteurs cohomologiques en $F$.

\enlargethispage{\baselineskip}%
On d\'efinit\pageoriginale de m\^eme
$$
\underline{\varphi}: \varinjlim_{n} \SheafExt_{\OX}^i(\Oo_{Y_n}, F) \to {\SheafH}_Y^i(F).
$$

On a en vue le th\'eor\`eme suivant:

\begin{suptheoreme} \label{II.6} \textup{a)}
Soient $X$ un pr\'esch\'ema localement noeth\'erien, $Y$ une partie ferm\'ee de $X$ d\'efinie par un Id\'eal coh\'erent $I$, $F$ un Module quasi-coh\'erent. Alors $\boldsymbol{\varphi}$ est un isomorphisme.

\textup{b)} Si $X$ est noeth\'erien $\varphi$ est un isomorphisme.
\end{suptheoreme}

Le \sisi{T}{t}h\'eor\`eme \sisi{(\Ref{II.6})}{\Ref{II.6}} r\'esultera de \sisi{(\Ref{II.6}).a}{\Ref{II.6}.a)} et du

\begin{suplemme} \label{II.7} Si l'espace topologique sous-jacent \`a $X$ est noeth\'erien et si $\boldsymbol{\varphi}$ est un isomorphisme, il en est de m\^eme de $\varphi$.
\end{suplemme}

On va d'abord prouver \sisi{(\Ref{II.7})}{ le lemme \Ref{II.7}}. On sait qu'il existe une suite spectrale
\begin{equation} \label{eq:II.7.1}
\H^p(X, {\SheafH}_Y^q(F)) \To \H_Y^{\sisi{m}{\ast}}(X, F).
\end{equation}
On a d'autre part un syst\`eme inductif de suites spectrales
\begin{equation*} \label{eq:II.7.2n} \tag{7.2$_n$}\stepcounter{equation}
\H^p(X, {\SheafExt}_{\OX}^q(\Oo_{Y_n}, F)) \To \Ext_{\OX}^{\sisi{m}{\ast}}(X;\Oo_{Y_n}, F).
\end{equation*}
Il r\'esulte de la d\'efinition de $\varphi$ et de $\underline{\varphi}$ que ces morphismes sont associ\'es \`a un homomorphisme $\Phi$ de suites spectrales de la limite inductive de (\Ref{eq:II.7.2n}) dans (\Ref{eq:II.7.1}). Si l'espace sous-jacent \`a $X$ est noeth\'erien, par (God. 4.12.1)\sfootnote{\Cf premi\`ere r\'ef\'erence bibliographique, \`a la fin de \Exp \Ref{I}.}
\begin{equation*}
\varinjlim_{n} \H^p(X, {\SheafExt}_{\OX}^q(\Oo_{Y_n}, F))\isomto \H^p(X, \varinjlim_{n} {\SheafExt}_{\OX}^q(\Oo_{Y_n}, F)),
\end{equation*}
alors $\Phi_2$ peut s'\'ecrire comme un morphisme:
\begin{equation*} \H^p(X, \varinjlim_{n} {\SheafExt}_{\OX}^q(\Oo_{Y_n}, F)) \to \H^p(X, {\SheafH}_Y^q(F))
\end{equation*}
qui n'est autre que celui qu'on d\'eduit de $\underline{\varphi}$.

Si $\underline{\varphi}$\pageoriginale est un isomorphisme il en est donc de m\^eme de $\Phi_2$, donc de $\varphi$ d'apr\`es ($\EGA 0_{\textup{III}}$ 11.1.5), \sisi{(\Ref{II.7})}{le lemme~\Ref{II.7}} est donc d\'emontr\'e.

On va maintenant prouver \sisi{(\Ref{II.6}.a)}{\Ref{II.6}.a)}, c'est une question locale sur $X$. D'apr\`es \sisi{\ignorespaces}{le} corollaire~\Ref{II.4} et (\EGA I~1.3.9 et 1.3.12) on peut supposer que $X$ est le spectre d'un anneau $A$. Il suffit donc de d\'emontrer que sous les hypoth\`eses du th\'eor\`eme~\Ref{II.6}.a), l'homomorphisme canonique:
\begin{equation} \label{eq:II.7.3} \varinjlim_{n} \Ext_A^i(A/I^n, M) \to \H_Y^i(X, M)
\end{equation}
est un isomorphisme.

Soient $f_{\alpha}$ un nombre fini d'\'el\'ements de $A$ engendrant $I$, $\bbf=(f_{\alpha})$; alors la suite des id\'eaux $(\bbf^n)$ est d\'ecroissante et cofinale \`a la suite des $I^n$, de sorte que (\Ref{eq:II.7.3}) est \'equivalent \`a un morphisme de $\partial$-foncteurs en $M$:
\begin{equation} \label{eq:II.7.4} \varinjlim_{n} \Ext_A^i(A/(\bbf^n), M) \to \H_Y^i(X, M).
\end{equation}

On a d'autre part des isomorphismes canoniques:
\begin{equation} \label{eq:II.7.5} \varinjlim_{n} \Hom_A(A/(\bbf^n), M) \simeq \varinjlim_{n} (m\in M\ |\ (\bbf^n)m=0) \simeq \H^0((\bbf), M).
\end{equation}
Comme $\varinjlim_{n} \Ext_A^i(A/(\bbf^n), M)$ est un $\partial$-foncteur universel en $M$, il existe un seul morphisme de $\partial$-foncteurs en $M$:
\begin{equation} \label{eq:II.7.6} \varinjlim_{n} \Ext_A^i(A/(\bbf^n), M) \to \H^i((\bbf), M),
\end{equation}
qui co\"incide en degr\'e z\'ero avec (\Ref{eq:II.7.5}).

Comme le compos\'e de (\Ref{eq:II.7.3}) et de (\Ref{eq:II.5.1}) est un morphisme de $\partial$-foncteurs en $M$ qui co\"incide avec (\Ref{eq:II.7.6}) en degr\'e $0$, il co\"incide avec (\Ref{eq:II.7.6}) en tout degr\'e. Le th\'eor\`eme \sisi{(\Ref{II.6}.a)}{\Ref{II.6}.a)} est donc une cons\'equence imm\'ediate du

\begin{suplemme} \label{II.8}
Soient $A$ un anneau noeth\'erien, $I$ un id\'eal engendr\'e par un syst\`eme fini $\bbf=(f_{\alpha})$ d'\'el\'ements, $M$ un $A$-module. Alors les homomorphismes (\Ref{eq:II.7.6}) sont des isomorphismes.
\end{suplemme}

\begin{suplemme} \label{II.9}
Soient $A$\pageoriginale un anneau, $\bbf=(f_{\alpha})$ un syst\`eme fini d'\'el\'ements de $A$, $I$ l'id\'eal engendr\'e par $\bbf$, $i$ un entier $>0$. Les conditions suivantes sont \'equivalentes:
\begin{enumeratea}
\item
L'homomorphisme (\Ref{eq:II.7.6}) est un isomorphisme pour tout $M$.
\item
$\H^i((\bbf), M)=0$ pour $M$ injectif.
\item
Le syst\`eme projectif $(\H_i(\bbf^n, A))=\H_{i, n}$ est essentiellement nul, c'est-\`a-dire: pour tout $n$, il existe $n'>n$ tel que $\H_{i, n'}\to\H_{i, n}$ soit nul.
\end{enumeratea}
\end{suplemme}

a) entra\^\i ne b) trivialement.

b) entra\^\i ne a), en effet b) entra\^\i ne que $M\mto\H^i((\bbf), M)$ est un foncteur cohomologique universel, (\Ref{eq:II.7.6}) est un alors un morphisme de foncteurs cohomologiques universels. C'est un isomorphisme en degr\'e z\'ero, donc en tout degr\'e\sisi{}{.}

c) entra\^\i ne b), en effet si $M$ est injectif, on~a pour tout $n$
\begin{equation*} \H^i(\bbf^n, M)=\Hom(\H_i(\bbf^n, A), M) = \Hom(\H_{i, n}, M),
\end{equation*}
c) entra\^\i ne donc que pour tout $i$ le syst\`eme inductif $(\H^i(\bbf^n, M))_{n\in\ZZ}$ est essentiellement nul, d'o\`u b).

b) entra\^\i ne c). Soit en effet $n>0$, et $j$ un monomorphisme de $\H_{i, n}$ dans un module injectif $M$. Soit $n'\geq n$ et soit $j_{n'}\in\Hom(\H_{i, n'}, M)$ le compos\'e de $j$ et de l'homomorphisme de transition $t_{n', n}:\H_{i, n'}\to\H_{i, n}$. Les $j_{n'}$ d\'efinissent un \'el\'ement de $\H^i((\bbf), M)$ qui est nul par hypoth\`ese. Il existe donc $n_0$ tel que $j_{n'}=0$ si $n'>n_0$. Mais puisque $j$ est un monomorphisme, $j_{n'}=0$ entra\^\i ne $\sisi{t_{n'}}{t_{n', n}}=0$, d'o\`u la proposition.

\begin{supcorollaire} \label{II.10} Supposons que l'espace sous-jacent \`a $X=\Spec(A)$ soit noeth\'erien. Pour que les conditions pr\'ec\'edentes soient v\'erifi\'ees pour toute famille finie d'\'el\'ements de $A$ et tout $i>0$ (ou encore: pour $i=1$), il faut et il suffit que pour tout $A$-module injectif $M$, le faisceau $F$ associ\'e \`a $M$ soit flasque.
\end{supcorollaire}

C'est n\'ecessaire: soient en effet $\bbf=(f_{\alpha})$ un syst\`eme fini d'\'el\'ements de $A$, $Y$ le ferm\'e d\'efini par $\bbf$ et $U=X-Y$, on~a alors la suite exacte
\begin{equation*} \H^0(X, F) \to \H^0(U, F) \to \H^1((\bbf), M) \to 0,
\end{equation*}
et gr\^ace \`a \Ref{II.9}.b, $\H^0(X, F)\to\H^0(U, F)$ est surjectif.

C'est suffisant\pageoriginale en vertu de (\Ref{eq:II.5.1}) et de ce que pour toute partie ferm\'ee $Y$ de $X$ et tout faisceau flasque $F$ sur $X$, $\H_Y^i(X, F)=0$ pour $i>0$.

\begin{suplemme} \label{II.11} Sous les hypoth\`eses du lemme~\Ref{II.9}, pour tout $A$-module noeth\'erien $N$ et pour tout $i>0$, le syst\`eme projectif $(\H_{i, n}(N))_{n\in\ZZ}$, o\`u $\H_{i, n}(N)=\H_i(\bbf^n, N)$, est essentiellement nul.
\end{suplemme}

Preuve par r\'ecurrence sur \sisi{\ignorespaces}{le} nombre $m$ d'\'el\'ements de $\bbf$.

Si $m=1$, $\bbf$ est r\'eduit \`a un seul \'el\'ement, soit $f$, $\H_{i, n}(N)$ est nul si $i>1$ et $\H_{1, n}(N)$ est canoniquement isomorphe \`a l'annulateur $N(n)$ de $f^n$ dans $N$, l'homomorphisme de transition $N(n')\to N(n)$, $n'\geq n$, \'etant la multiplication par $f^{n'-n}$. Les $N(n)$ forment une suite croissante de sous-modules de $N$, et puisque $N$ est noeth\'erien il existe $n_0$ tel que $N(n)=N(n_0)$ si $n\geq n_0$. Donc tous les $N(n)$ sont annul\'es par $f^{n_0}$ et les homomorphismes de transition $N(n')\to N(n)$ sont tous nuls si $n'\geq n+n_0$. Le lemme est donc prouv\'e pour $m=1$.

On suppose maintenant que $m>1$ et que le lemme est prouv\'e pour les entiers $m'<m$; soit alors $\bbg=(f_1, \ldots, f_{m-1})$ et $\bbh=f_m$.

Pour tout $n>0$, on~a (\EGA III~1.1.4.1) une suite exacte:
\begin{equation*}
0 \to \H_0(\bbh^n, \H_i(\bbg^n, N)) \to \H_i(\bbf^n, N) \to \H_1(\bbh^n, \H_{i-1}(\bbg^n, N)) \to 0,
\end{equation*}
et pour $n$ variable un syst\`eme projectif de suites exactes. Il r\'esulte des hypoth\`eses de r\'ecurrence que pour $i>0$ les $\H_i(\bbg^n, N)$ forment un syst\`eme projectif essentiellement nul, donc aussi les $\H_0(\bbh^n, \H_i(\bbg^n, N))$ qu'on identifie \`a des quotients de $\H_i(\bbg^n, N)$. Pour les termes de droite on va factoriser les morphismes de transition de $n'$ \`a $n$ par:
\begin{equation*}
\H_1(\bbh^{n'}, \H_{i-1}(\bbg^{n'}, N)) \to \H_1(\bbh^{n'}, \H_{i-1}(\bbg^{n}, N)) \to \H_1(\bbh^{n}, \H_{i-1}(\bbg^{n}, N)).
\end{equation*}
Puisque $\H_{i-1}(\bbg^{n}, N)$ est un module noeth\'erien il r\'esulte du cas $m=1$ qu'il existe, pour $n$ donn\'e, $n'>n$ tel que la seconde fl\`eche soit nulle. On voit\pageoriginale donc que dans ce syst\`eme projectif de suites exactes, les syst\`emes projectifs extr\^emes sont essentiellement nuls, il en est donc de m\^eme du syst\`eme projectif m\'edian.

On a donc prouv\'e le lemme \sisi{(\Ref{II.11})}{\Ref{II.11}} donc le lemme \sisi{(\Ref{II.8})}{\Ref{II.8}} et partant le th\'eor\`eme \sisi{(\Ref{II.6})}{\Ref{II.6}}.

\begin{remarquestar}
On peut aussi obtenir le th\'eor\`eme \sisi{(\Ref{II.6})}{\Ref{II.6}} en d\'emontrant la condition du \sisi{Corollaire (\Ref{II.10})}{corollaire \Ref{II.10}} \`a l'aide des th\'eor\`emes de structure des modules injectifs sur un anneau noeth\'erien (Matlis, Gabriel).
\end{remarquestar}

\chapterspace{-1}
\chapter{Invariants cohomologiques et profondeur} \label{III}

\section{Rappels} \label{III.1}

Nous\pageoriginale \'enoncerons quelques d\'efinitions et r\'esultats que le lecteur trouvera par exemple dans le chapitre I du cours profess\'e par \sisi{J.P.~SERRE}{J.-P\ptbl Serre} au Coll\`ege de France en 1957-58.\sisi{\sfootnote{\label{noteserre}\Cf Alg\`ebre locale et multiplicit\'e. Lecture Notes in Mathematics \numero 11, Springer.}}{\nde{\label{noteserre} la r\'e\'edition du texte de Serre (Serre~J.-P., \emph{Alg\`ebre locale. {M}ultiplicit\'es}, {Cours au Coll\`ege de France, 1957--1958, r\'edig\'e par Pierre Gabriel, seconde \'edition, Lect. Notes in Math.}, vol.~11, Springer-Verlag, 1965) ne contient plus les preuves de ces \'enonc\'es. On pourra se reporter \`a (Bourbaki~N., \textit{Alg\`ebre commutative}, Masson), comme le sugg\`ere Serre lui-m\^eme.}}

\begin{definition} \label{def:III.1.1}
Soit $A$ un anneau, (commutatif \`a \'el\'ement unit\'e comme dans tout ce qui suivra) et soit $M$ un $A$-module (unitaire comme dans tout ce qui suivra), on appelle:
\begin{itemize}
\item
annulateur de $M$ et on note $\Ann M$ l'ensemble des $a\in A$, tels que pour tout $m\in M$ on ait $am=0$.

\item
support de $M$ et on note $\Supp M$ l'ensemble des id\'eaux premiers $\pp $ de $A$ tels que le localis\'e $M_\pp$ soit non nul.

\item
\og assassin de $M$\fg ou \og ensemble des id\'eaux premiers associ\'es \`a $M$\fg et on note $\Ass M$ l'ensemble des id\'eaux premiers $\pp $ de $A$ tels qu'il existe un \'el\'ement non nul de~$M$ dont l'annulateur soit $\pp $.
\end{itemize}

Si $\ideal a$ est un id\'eal de $A$, nous noterons $\rr(\ideal a)$ la racine de $\ideal a$ dans $A$, \ie l'ensemble des \'el\'ements de $A$ dont un puissance est dans $\ideal a$.
\end{definition}

Les r\'esultats suivants sont valables si l'on suppose que \emph{$A$ est noeth\'erien et $M$ de type fini}.

\setcounter{subsection}{0}

\skpt

\begin{proposition} \label{III.1.1}
\begin{enumeratei}
\item
$\Ass M$ est un ensemble fini.

\item
Pour qu'un \'el\'ement de $A$ annule un \'el\'ement non nul de $M$, il faut et il suffit\pageoriginale qu'il appartienne \`a l'un des id\'eaux associ\'es \`a $M$.

\item
La racine de l'annulateur de $M$, $\rr(\Ann M)$, est l'intersection des id\'eaux associ\'es \`a $M$ qui sont minimaux (pour la relation d'inclusion dans $\Ass M$).
\end{enumeratei}
\end{proposition}

\begin{proposition} \label{III.1.2}
Soit $\pp $ un id\'eal premier de $A$, les assertions suivantes sont \'equivalentes:
\begin{enumeratei}
\item
$\pp \in \Supp M$.

\item
Il existe $\qq \in \Ass M$ tel que $\qq \subset\pp $.

\item
$\pp \supset \Ann M$.

\item[\textup{(iii bis)}]
$\pp \supset \rr(\Ann M)$.
\end{enumeratei}
\end{proposition}

\begin{proposition} \label{III.1.3}
Soit $N$ un $A$-module de type fini, on~a la formule:
\[
\Ass \Hom_A(N, M) = \Supp N \cap \Ass M.
\]
\end{proposition}

\section{Profondeur} \label{III.2}

Dans tout ce paragraphe, $A$ d\'esigne un anneau commutatif, $I$ un id\'eal de $A$, $M$~et~$N$ deux $A$-modules. On notera $X$ le spectre premier de $A$ (on ne se servira pas de son faisceau structural dans ce paragraphe) et $Y$ la vari\'et\'e de $I$, $Y = \Supp(A/I) = \{\pp \in\nobreak X, \ \pp \supset I\}$.

\begin{lemme} \label{III.2.1}
Supposons que $A$ soit noeth\'erien et que les modules $M$ et $N$ soient de type fini. Supposons de plus que $\Supp N=Y$. Alors les assertions suivantes sont \'equivalentes:
\begin{enumeratei}
\item
$\Hom_A(N, M)=0$.
\item
$\Supp N \cap \Ass M = \emptyset$.
\item
L'id\'eal $I$ n'est pas diviseur de $0$ dans $M$, ce qui signifie que pour tout $m\in M$, $Im=0$ entra\^ine $m=0$.
\item
Il existe dans $I$ un \'el\'ement $M$-r\'egulier. (Un \'el\'ement $a$ de $A$ est dit $M$-r\'egulier si l'homoth\'etie de rapport $a$ dans $M$ est injective.)

\item
Pour tout $\pp \in Y$, l'id\'eal maximal $\sisi{\mm_\pp}{\pp A_\pp}$ de l'anneau local $A_\pp$ n'est pas associ\'e \`a $M_\pp$. En formule: $\pp A_\pp \notin \Ass M_\pp$.
\end{enumeratei}
\end{lemme}

\skpt
\begin{proof}
(i) \SSI (ii) car $\Ass \Hom_A(N, M)=\emptyset$\pageoriginale est \'equivalent \`a (ii) d'apr\`es la proposition \Ref{III.1.3} et \`a (i) par une cons\'equence facile de la proposition \Ref{III.1.2}.

(iii) \ALORS (ii) par l'absurde: \og il existe $\pp \in \Supp N \cap \Ass M$\fg entra\^ine que $\pp \supset I$ et qu'il existe $m\in M$ dont l'annulateur est $\pp $, donc $Im \subset \pp m=0$, ce qui contredit (iii).

(iv) \ALORS (iii) trivialement.

(ii) \SSI (iv) car $\Supp N = Y$, donc (ii) signifie que $I$ n'est contenu dans aucun id\'eal associ\'e \`a $M$ ou encore, (car les id\'eaux associ\'es \`a $M$ sont premiers et en nombre fini), que $I$ n'est pas contenu dans la r\'eunion des id\'eaux associ\'es \`a $M$. Or, par la \sisi{Prop.}{proposition} \Ref{III.1.1}~(ii), cet ensemble est l'ensemble des \'el\'ements de $A$ qui ne sont pas $M$-r\'eguliers.

(i) \ALORS (v); en effet, si $\Hom_A(N, M)=0$ et si $\pp \in Y$, on en d\'eduit, en vertu de la formule
\[
\bigl(\Hom_A(N, M) \bigr)_\pp = \Hom_{A_\pp} \bigl(N_\pp, M_\pp \bigr),
\]
que $\Hom_{A_\pp} \bigl(N_\pp, M_\pp \bigr)=0$, donc, gr\^ace \`a la proposition \Ref{III.1.3},
\[
\Supp N_\pp \cap \Ass M_\pp = \emptyset,
\]
or $\pp A_\pp \in \Supp N_\pp, $ donc $ \pp A_\pp \notin \Ass M_\pp$.

(v) \ALORS (i); en effet, si $\pp \in \Ass M$, il existe $m\in M$ dont l'annulateur est $\pp$, donc l'image canonique de $m$ dans $M_\pp$ est non nulle, donc son annulateur est un id\'eal qui contient $\pp$, donc $\pp A_\pp$, donc lui est \'egal. L'id\'eal $\pp A_\pp$ est donc associ\'e \`a $M_\pp$, donc $\pp \notin Y$ d'apr\`es (v), d'o\`u (i).
\end{proof}

Nous allons travailler sur ces conditions en rempla\c cant le
foncteur $\Hom$ par ses d\'eriv\'es.

\begin{theoreme} \label{III.2.2}
Soit\pageoriginaled $A$ un anneau commutatif, $I$ un id\'eal de $A$, $M$ un $A$-module. Soit $n$ un entier.

\begin{enumeratea}
\item
S'il existe une suite $f_1, \ldots, f_{n+1}$, d'\'el\'ements de $I$ qui forme une suite $M$-r\'eguli\`ere (\ie si $f_1$ est $M$-r\'egulier et si $f_{i+1}$ est r\'egulier dans $M/(f_1, \ldots, f_i)M$ pour $i\leq n$), pour tout $A$-module $N$ annul\'e par une puissance de $I$, on a:
\[
\Ext_A^i(N, M)=0 \text{ pour } i\leq n.
\]

\item
Si de plus $A$ est noeth\'erien, si $M$ est de type fini, et s'il existe un $A$-module $N$ de type fini tel que $\Supp N = \variete(I)$ et tel que $\Ext_A^i(N, M)=0$ pour $i\leq n$, alors il existe une suite $f_1, \ldots, f_{n+1}$, d'\'el\'ements de $I$ qui est $M$-r\'eguli\`ere.
\end{enumeratea}
\end{theoreme}

D\'emontrons d'abord a), par r\'ecurrence. Si $n<0$ l'\'enonc\'e est vide.

Si $n\geq 0$, supposons que a) soit d\'emontr\'e pour $n'<n$; par hypoth\`ese il existe $f_1\in I$ qui est $M$-r\'egulier. D\'esignons par $f_1^i$ la multiplication par $f_1$ dans $\Ext_A^i(N, M)$ et par $f_1^M$ la multiplication par $f_1$ dans $M$. La suite
\begin{equation} \label{eq:III.2.1}
0 \to M \lto{f_1^M} M \to M/f_1M \to 0
\end{equation}
est exacte, donc aussi la suite:
\[
\Ext_A^{i-1} (N, M/f_1M) \lto{\delta}\Ext_A^i (N, M) \lto{f_1^i} \Ext_A^i (N, M).
\]
Par hypoth\`ese $I^rN=0$, donc $f_1^0$ et nilpotent; $\Ext^i$ est un foncteur universel, il en est de m\^eme de $f_1^i$ pour tout $i$. Par ailleurs, il existe une suite r\'eguli\`ere dans $M/f_1 M$ qui a $n$ \'el\'ements, donc, par hypoth\`ese de r\'ecurrence,
\[
Ext_A^{i-1}(N, M)=0 \text{ si } i\leq n-1.
\]
On en d\'eduit que si $i\leq n$, $f_1^i$ est \`a la fois nilpotent et injectif donc $\Ext_A^i(N, M)=0$.

D\'emontrons b)\pageoriginale, \'egalement par r\'ecurrence. Si $n<0$, l'\'enonc\'e est vide.

Si $n=0$, b) r\'esulte de l'assertion (i) \ALORS (iv) du lemme \Ref{III.2.1}.

Si $n>0$, d'apr\`es b) pour $n=0$, il existe un \'el\'ement $f_1\in I$ qui est $M$-r\'egulier; de la suite exacte (\Ref{eq:III.2.1}), on d\'eduit la suite exacte:
\begin{equation} \label{eq:III.2.2} \Ext_A^{i-1} (N, M) \to \Ext_A^{i-1} (N, M/f_1M) \to \Ext_A^i (N, M).
\end{equation}
On en conclut que les hypoth\`eses de b) sont v\'erifi\'ees pour le
module $M/f_1M$ et pour l'entier $n-1$. Par l'hypoth\`ese de
r\'ecurrence, il existe une suite de $n$ \'el\'ements de $I$ qui est
r\'eguli\`ere pour $M/f_1M$, ce qui entra\^ine qu'il existe une
suite de $n+1$ \'el\'ements de~$I$, commen\c cant par $f_1$, et qui
est $M$-r\'eguli\`ere.

Ce th\'eor\`eme nous invite \`a g\'en\'eraliser de la mani\`ere suivante la d\'efinition classique de la profondeur d'un module de type fini sur un anneau noeth\'erien:

\begin{definition} \label{III.2.3}
Soit $A$ un anneau commutatif \`a \'el\'ement unit\'e, soit $M$ un $A$-module, soit $I$ un id\'eal de $A$. On appelle $I$-profondeur de $M$, et on note $\prof_I M$, la borne sup\'erieure \sisi{\ignorespaces}{dans $\NN\cup\{+\infty\}$} de l'ensemble des entiers naturels $n$, qui sont tels que pour tout $A$-module de type fini $N$ annul\'e par une puissance de $I$, on ait
\[
\Ext_A^i(N, M)=0 \text{ pour tout } i< n.
\]
\end{definition}

On d\'eduit du th\'eor\`eme pr\'ec\'edent que si $n$ est la borne sup\'erieure des longueurs des suites $M$-r\'eguli\`eres d'\'el\'ements de $I$, on~a $n\leq\prof_I M$.

\noindent Plus pr\'ecis\'ement:

\begin{proposition} \label{III.2.4}
Soit $A$ un anneau commutatif, $I$ un id\'eal de $A$ et soit $M$ un $A$-module, soit enfin $n\in \NN$. Consid\'erons les assertions:
\begin{enumerate}
\item
$n\leq \prof_I M$.
\item
Pour tout $A$-module de type fini $N$ qui est annul\'e par une puissance de $I$, on~a:
\[
\Ext_A^i(N, M)=0 \text{ pour } i<n.
\]
\item
Il existe\pageoriginale un $A$-module de type fini $N$ tel que $\Supp N = \variete(I)$ et tel que $\Ext_A^i(N, M)=0$ si $i<n$.
\item
Il existe une suite $M$-r\'eguli\`ere de longueur $n$ form\'es d'\'el\'ements de $I$.
\end{enumerate}

On a les implications logiques suivantes:
\[
\begin{array}{c}\dpl \textup{(1)} \Ssi \textup{(2)} \From (4) \\
\big\Downarrow \\
\textup{(3)} \\
\end{array}
\]
De plus si $A$ est noeth\'erien et $M$ de type fini, ces conditions sont \'equivalentes.
\end{proposition}

\skpt
\begin{proof}
(1) \SSI (2) par d\'efinition et (2) \ALORS (3) en prenant $N=A/I$; De plus (4) \ALORS (2) par le th\'eor\`eme \Ref{III.2.2} a). Enfin, si $A$ est noeth\'erien et $M$ de type fini, (3) \ALORS (4) par le th\'eor\`eme \Ref{III.2.2} b).
\skipqed
\end{proof}

Nous supposons $A$ noeth\'erien et $M$ de type fini jusqu'\`a la fin de ce paragraphe.

\begin{corollaire} \label{III.2.5}
Soit $f\in I$ un \'el\'ement $M$-r\'egulier, on a:
\[
\prof_I M = \prof_I (M/fM) + 1.
\]
\end{corollaire}

En effet, si $n\leq \prof_I (M/fM)$, il existe une suite d'\'el\'ements de $I$, $f_1, \ldots, f_n$, qui est $(M/fM)$-r\'eguli\`ere; donc la suite $f, f_1, \ldots, f_n$ est $M$-r\'eguli\`ere, donc $n+1 \leq \prof_I M$, donc $\prof_I M\geq \prof_I (M/fM) + 1$. Par ailleurs, d'apr\`es la suite exacte (\Ref{eq:III.2.2}), si $i \leq \prof_I M$, on~a $\Ext_A^{i-1} (N, M/fM)=0$, donc $\prof_I M - 1 \leq \prof_I (M/fM)$.

\enlargethispage{1.2\baselineskip}%
\begin{corollaire} \label{III.2.6}
Toute suite $M$-r\'eguli\`ere finie, form\'ee d'\'el\'ements de $I$, peut \^etre prolong\'ee en une suite $M$-r\'eguli\`ere maximale, dont la longueur est n\'ecessairement \'egale \`a la $I$-profondeur de~$M$.
\end{corollaire}

\begin{remarque} \label{III.2.7}
On ne se retient qu'\`a grand peine de dire qu'un $A$-module est d'autant plus beau que sa profondeur est plus grande. Un module dont le support\pageoriginale ne rencontre pas $\variete (I)$ est des plus beaux; en effet, on peut d\'emontrer que pour que $\prof_I M$ soit fini, il est n\'ecessaire et suffisant que $\Supp M \cap \variete(I) \neq \emptyset$.
\end{remarque}

\begin{remarque} \label{III.2.8}
Si $A$ est un anneau semi-local, soit $\rr(A)$ son radical et $k=A/\rr(A)$ son anneau r\'esiduel. La notion de profondeur int\'eressante est obtenue en prenant pour $I$ le radical de~$A$. Nous conviendrons donc de noter simplement $\prof M$ la $\rr(A)$-profondeur d'un $A$-module~$M$. On retrouve dans ce cas la notion de \og codimension homologique\fg, (\cf \sisi{SERRE}{Serre}, \textit{op\ptbl cit.} note~\eqref{noteserre}, page~\pageref{noteserre}), que l'on notait $\codimh_A M$, et qui est d\'efinie comme la borne inf\'erieure des entiers $i$ tels que $\Ext_A^i(k, M)\neq 0$; en effet $\Supp k = \variete(\rr(A))$.
\end{remarque}

\begin{proposition} \label{III.2.9}
Si $A$ est noeth\'erien et $M$ de type fini, on a:
\[
\prof_I M = \inf_{\pp \in \variete (I)} \prof M_\pp.
\]
\end{proposition}

\begin{corollaire} \label{III.2.10}
Si $A$ est un anneau semi-local noeth\'erien, et si $M$ est un $A$-module de type fini, on a:
\[
\prof M = \inf_\mm \prof M_\mm,
\]
o\`u $\mm $ parcourt l'ensemble des id\'eaux maximaux de $A$.
\end{corollaire}

Le corollaire r\'esulte imm\'ediatement de le proposition \Ref{III.2.9}; en effet les id\'eaux premiers qui contiennent le radical sont les id\'eaux maximaux.

Par ailleurs, soit $f\in I$; si $f$ est $M$-r\'egulier, si $\pp\in X$ et si $\pp\supset I$, l'image $g$ de $f$ dans $A_\pp$ appartient \`a $\sisi{\mm_\pp}{\pp A_\pp}$, id\'eal maximal de $A_\pp$; de plus $g$ est $M_\pp$-r\'egulier, comme il r\'esulte de la suite exacte
\begin{equation} \label{eq:III.2.3} 0 \to M_\pp \lto{g'} M_\pp \to (M/fM)_\pp \to 0,
\end{equation}
o\`u $g'$ d\'esigne l'homoth\'etie de rapport $g$ dans $M_\pp$. Cette suite exacte donne aussi que $(M/fM)_\pp$ est isomorphe \`a $M_\pp/gM_\pp$; en appliquant le corollaire \Ref{III.2.5} \`a $M$ et \`a $M_\pp$, on en d\'eduit, par r\'ecurrence, que\pageoriginale $\prof_I M\leq \nu(M)$, o\`u l'on a pos\'e pour tout $M$:
\[
\nu(M)=\inf_{\pp\in\variete(I)} \prof M_\pp.
\]
Plus pr\'ecis\'ement, toujours par r\'ecurrence, on sait, si $f$ est $M$-r\'egulier, que $\nu(M)= \nu(M/fM)+1$; il reste donc \`a d\'emontrer que si $\nu(M)\geq 1$, il existe un \'el\'ement $M$-r\'egulier dans $I$. Or en appliquant \sisi{(\Ref{III.2.1})}{le lemme \Ref{III.2.1}} \`a $M_\pp$, $A_\pp$ et $\sisi{\mm_\pp}{\pp A_\pp}$ pour tout $\pp\in \variete(I)$, on voit que $\sisi{\mm_\pp}{\pp A_\pp} \notin \Ass M_\pp$, donc, en appliquant \sisi{(\Ref{III.2.1})}{le lemme \Ref{III.2.1}} \`a $A$, $M$ et $I$, on~a la conclusion.

\enlargethispage{\baselineskip}%
\begin{proposition} \label{III.2.11}
Soit $u:A\to B$ un homomorphisme d'anneaux noeth\'eriens. Soit~$I$ un id\'eal de $A$, $M$ un $A$-module de type fini. Posons $I_B=I\otimes_A B$ et $M_B=M\otimes_A B$. Si $B$ est $A$-plat, on~a:
\[
\prof_{I_B} M_B \geq \prof_I M;
\]
de plus si $B$ est fid\`element plat sur $A$, on~a \'egalit\'e.
\end{proposition}

En effet, soit $N=A/I$; par platitude on a: $N\otimes_A B=B/I_B$; posons $N_B=N\otimes_A B$. Toujours par platitude et hypoth\`eses noeth\'eriennes, on a:
\[
\Ext_B^i (N_B, M_B) = \Ext_A^i(N, M)\otimes_A B,
\]
donc $\Ext_A^i(N, M) = 0 $ entra\^ine $\Ext_B^i(N_B, M_B)=0$, et la r\'eciproque est vraie si $B$ est fid\`element plat sur $A$.

\section{Profondeur et propri\'et\'es topologiques} \label{III.3}

\begin{lemme} \label{III.3.1}
Soit $X$ un espace topologique, $Y$ un sous-espace ferm\'e, soit $F$ un faisceau de groupes ab\'eliens sur $X$. Posons $U=X- Y$. Si $n$ est un entier, les conditions suivantes sont \'equivalentes:
\begin{enumeratei}
\item
$\SheafH_Y^i (X, F)=0$ si $i<n$.\pageoriginale
\item
Pour tout ouvert $V$ de $X$, l'homomorphisme de groupes
\[
\H^i(V, F) \to \H^i(V\cap U, F)
\]
est bijectif si $i<n-1$ et injectif si $i=n-1$.

\item
Pour tout ouvert $V$ de $X$,
\[
\H_{Y\cap V}^i(V, F_{\mid V})=0 \text{ si } i<n.
\]
\end{enumeratei}
\end{lemme}

\skpt

\begin{proof}
(ii) \SSI (iii); en effet, soit $V$ un ouvert de $X$, posons $X'=V$, $Y'=Y\cap V$, $F'=F_{\mid V}$, $U'=X'- Y'$; $Y'$ est ferm\'e dans $X'$, on~a donc une suite exacte:
\[
\H_{Y'}^i(X', F') \to \H^i(X', F') \lto{\rho_i} \H^i(U', F') \to \H_{Y'}^{i+1}(X', F').
\]
Si les termes extr\^emes sont nuls, l'homomorphisme $\rho_i$ est bijectif, et si le terme de gauche est nul, $\rho_i$ est injectif. Donc (iii) \ALORS (ii). R\'eciproquement, si $i<n$, $\H_{Y'}^i(X', F')$ est nul car $\rho_i$ est injectif et $\rho_{i-1}$ surjectif.

(i) \ALORS (iii); en effet la suite spectrale \og de passage du local au global\fg donne:
\[
\H^p\bigl(X, \SheafH_Y^q(X, F)\bigr) \To \H_Y^*(X, F).
\]
Or, par hypoth\`ese $\SheafH_Y^q(X, F) =0$ si $q<n$, donc $\sisi{(\H_Y (X, F))^{p+q}}{\H_Y^{p+q} (X, F)} = 0$ si $p+q<n$.

(iii) \ALORS (i); en effet (iii) exprime que le pr\'efaisceau
\[
V \sisi{\rightsquigarrow}{\mto} \H_{Y\cap V}^i(V, F_{\mid V})
\]
est nul, donc aussi le faisceau associ\'e qui est $\SheafH_Y^i (X, F)$, car $Y$ est ferm\'e.
\skipqed
\end{proof}

\begin{remarque} \label{III.3.2}\pageoriginale
\sisi{Si $n\geq 2$, on peut omettre la condition que $\rho_{n-1}$ soit injectif, car le faisceau associ\'e au pr\'efaisceau
\[
V \sisi{\rightsquigarrow}{\mto} \H^{n-2}(V\cap U, F)
\]
est isomorphe \`a $\SheafH_Y^{n-1} (X, F)$, mais par l'hypoth\`ese (ii) pour $i=n-2$, ce pr\'efaisceau est isomorphe au pr\'efaisceau $V \sisi{\rightsquigarrow}{\mto} \H^i(V, F)$ dont le faisceau associ\'e est nul.}{L'\'equivalence de (i) est (ii) a \'et\'e prouv\'ee dans la proposition~\Ref{I}~\Ref{I.2.13}. Comme on~a remarqu\'e alors, si $n\geq 2$, on peut omettre la condition que $\rho_{n-1}$ soit injectif.}\nde{l'\'edition originale redonnait une preuve, pas tout \`a fait correcte.}
\end{remarque}

\begin{proposition} \label{III.3.3}
Soit $X$ un pr\'esch\'ema localement noeth\'erien, $Y$ un sous-pr\'esch\'ema ferm\'e de $X$, $F$ un $\OX$-module coh\'erent. Les conditions du lemme \Ref{III.3.1} sont \'equivalentes \`a chacune des conditions suivantes:
\begin{enumeratei}
\setcounter{enumi}{3}
\item
Pour tout $\sisi{y}{x}\in Y$, on~a $\prof F_x \geq n$;

\item
Pour tout $\OX$-module coh\'erent $G$ sur $X$, de support contenu dans $Y$ on a
\[
\SheafExt_{\OX}^i(G, F)=0 \text{ si } i<n;
\]

\item
Il existe un $\OX$-module coh\'erent $G$ dont le support est \'egal \`a $Y$ et tel que
\[
\SheafExt_{\OX}^i(G, F)=0 \text{ si } i<n.
\]
\end{enumeratei}
\end{proposition}

Si $X$ est affine, on~a fait tout ce qu'il faut \sisi{\ignorespaces}{(\cf proposition~\Ref{III.2.4})} pour d\'emontrer l'\'equivalence des trois conditions de la proposition \Ref{III.3.3}\sisi{, }{;} or elles sont \sisi{locales}{locales, mise \`a part l'implication (v)\ALORS (vi), mais on peut alors prendre $G=\Oo_Y$ et invoquer \`a nouveau proposition~\Ref{III.2.4}}\sisi{, donc il suffit}{. Il suffit donc} de prouver (i) \ALORS (vi) et (v) \ALORS (i).

Soit $J$ l'id\'eal de $Y$, c'est un faisceau coh\'erent d'id\'eaux; soit $\Oo_{\sisi{X_m}{Y_m}} = \OX/J^{m+1}$, c'est un $\OX$-module coh\'erent dont le support est \'egal \`a $Y$, et on sait \sisi{\ignorespaces}{(th\'eor\`eme \Ref{II}~\Ref{II.6}.b)} que
\[
\SheafH_Y^i(X, F) = \varinjlim_m \SheafExt_{\OX}^i(\Oo_{\sisi{X_m}{Y_m}}, F),
\]
donc (v) \ALORS (i). Par ailleurs\pageoriginale, les morphismes de transition sont des \'epimorphismes dans le syst\`eme projectif des $\Oo_{\sisi{X_m}{Y_m}}$.

Si le foncteur $\SheafExt^i$ est exact \`a gauche en sont premier argument, du moins lorsque celui-ci est dans la cat\'egorie des $\OX$-modules coh\'erents de support contenu dans $Y$, les morphismes de transition du syst\`eme inductif obtenu en appliquant $\SheafExt^i$ aux $\Oo_{\sisi{X_m}{Y_m}}$ seront injectifs, or (i) entra\^ine que la limite est nulle, donc (i) entra\^inera que les modules $\SheafExt_{\OX}^i(\Oo_{\sisi{X_m}{Y_m}}, F)$ sont nuls pour tout~$m$. Raisonnons par r\'ecurrence. L'\'enonc\'e est trivial pour $n<0$. Supposons que (i) \ALORS (vi) pour $n<q$, alors (i) \ALORS (v), donc, par la suite exacte des $\SheafExt$, $\SheafExt^q$ est exact \`a gauche en son premier argument, donc les modules $\SheafExt_{\OX}^q(\Oo_{\sisi{X_m}{Y_m}}, F)$ sont nuls pour tout $m$. Donc (i) \ALORS (vi) pour $n\leq q$.
\CQFD

\bgroup
\def\labelenumi{\theenumi)}
\begin{exemple} \label{III.3.4}
Soit $A$ un anneau local noeth\'erien, $\mm$ son id\'eal maximal, $M$ un $A$-module de type fini, soit enfin $n$ un entier. Posons $X=\Spec (A)$, $Y=\{\mm\}$, $U=X-Y$. Soit $F$ le faisceau associ\'e \`a $M$. Les conditions suivantes sont \'equivalentes:
\begin{enumerate}

\item
$\prof M \geq n$;

\item
l'homomorphisme naturel
\[
\H^i(X, F) \to \H^i(U, F)
\]
\sisi{}{est} injectif si $i=n-1$, \sisi{est}{\ignorespaces} bijectif si $i<n-1$;

\item
$\Ext_A^i(k, M) =0$ si $i<n$, o\`u $k=A/\mm$;

\item
$\H_Y^i(X, F)=0$ si $i<n$.
\end{enumerate}
\end{exemple}

\sisi{}{Tenant compte de la remarque~\Ref{III.3.2}, on obtient:}

\begin{corollaire} \label{III.3.5}
Soit $X$ un pr\'esch\'ema localement noeth\'erien, $Y$ un sous-pr\'esch\'ema ferm\'e de $X$, $F$ un $\OX$-module coh\'erent; les conditions suivantes sont \'equivalentes:
\begin{enumerate}
\item
pour tout $x\in Y$, $\prof F_x \geq 2$;
\item
pour tout ouvert\pageoriginale $V$ de $X$, l'homomorphisme naturel
\[
\H^0(V, F) \to \H^0(V\cap(X-Y), F)
\]
est bijectif.
\end{enumerate}
\end{corollaire}
\egroup

\begin{theoreme}[\sisi{HARTSHORNE}{Hartshorne}] \label{III.3.6}
Soit $X$ un pr\'esch\'ema localement noeth\'erien, $Y$ un sous-pr\'esch\'ema ferm\'e de $X$. Supposons que, pour tout $x\in Y$, $\prof \Oo_{X, x} \geq 2$; alors l'application naturelle
\[
\pi_0(X) \to \pi_0(X-Y)
\]
est bijective.
\end{theoreme}

\begin{proof}
Puisque $X$ est localement noeth\'erien, $X$ est localement connexe; il suffit donc de prouver que pour que $X$ soit connexe il est n\'ecessaire et suffisant que $X-Y$ le soit. Or, pour qu'un espace annel\'e en anneaux \vadjust{\penalty -100} locaux $(X, \OX)$ soit connexe, il est n\'ecessaire et suffisant que $\H^0(X, \OX)$ ne soit pas compos\'e direct de deux anneaux non nuls. Mais l'hypoth\`ese implique, d'apr\`es le corollaire \Ref{III.3.5} appliqu\'e \`a $F=\OX$, que l'homomorphisme
\[
\H^0(X, \OX) \to \H^0(X-Y, \OX)
\]
est un isomorphisme, d'o\`u la conclusion.
\skipqed
\end{proof}

\begin{corollaire} \label{III.3.7}
Soit $X$ un pr\'esch\'ema localement noeth\'erien. Soit $d$ un entier tel que $\dim \Oo_{X, x} \geq d$ entra\^ine $\prof \Oo_{X, x} \geq 2$. Alors, si $X$ est connexe, $X$ est connexe en codimension $d-1$, \ie si $X'$ et $X''$ sont deux composantes irr\'eductibles de $X$, il existe une suite de composantes irr\'eductibles \sisi{}{de $X$}:
\[
X'=X_0, X_1, \ldots, X_n=X''
\]
telle que pour tout $i$, \sisi{$1\leq i \leq n$, la codimension de $X_{i-1}\cap X_i$}{$0\leq i < n$, la codimension de $X_{i}\cap X_{i+1}$} dans $X$ soit inf\'erieure ou \'egale \`a $d-1$.
\end{corollaire}

Remarquons\pageoriginale d'abord que si $X$ est de \sisi{COHEN-MACAULAY}{Cohen-Macaulay}, $d=2$ jouira de la propri\'et\'e \'evoqu\'ee plus haut. \`A ce propos, rappelons que l'on d\'efinit la codimension de $Y$ dans~$X$ comme la borne inf\'erieure des dimensions des anneaux locaux dans~$X$ des points de~$Y$.

\begin{proof}
Soit $\Ff$ l'ensemble des parties ferm\'ees de $X$ dont la codimension dans $X$ est sup\'erieure ou \'egale \`a $d$. On notera que $\Ff$ est un antifiltre de parties ferm\'ees de $X$. Par ailleurs, pour qu'un ferm\'e $Y$ de $X$ soit \'el\'ement de $\Ff$, il faut et il suffit que, pour tout $y\in Y$ il existe un voisinage ouvert $V$ de $X$ tel que $Y\cap V$ soit de codimension $\geq d$ dans $V$. Enfin, si $X$ est connexe et si $Y\in \Ff$, $X-Y$ est connexe d'apr\`es le th\'eor\`eme de \sisi{HARTSHORNE}{Hartshorne}. Le corollaire r\'esulte donc du lemme suivant, qui est de nature purement topologique.

\begin{lemme} \label{III.3.8}
Soit $X$ un espace topologique connexe et localement noeth\'erien, et soit $\Ff$ un antifiltre de parties ferm\'ees de $X$. On suppose que tout ferm\'e $Y\subset X$ qui appartient localement \`a $\Ff$, (\ie pour tout point $x\in X$ il existe un voisinage ouvert $V$ de $x$ dans $X$ et un $Y'\in \Ff$ tel que $V\cap Y=V\cap Y'$), appartient \`a $\Ff$. Les conditions suivantes sont \'equivalentes:
\begin{enumeratei}

\item
pour tout $Y\in \Ff$, $X-Y$ est connexe;

\item
si $X'$ et $X''$ sont deux composantes irr\'eductibles distinctes de $X$, il existe une suite de composantes irr\'eductibles de $X$, $X_0, X_1, \ldots, X_n$, telle que $X'=X_0$, $X''=X_n$ et, pour tout $i$, $\sisi{1\leq i\leq n}{1\leq i< n}$, $\sisi{X_{i-1} \cap X_i}{X_i \cap X_{i+1}} \notin \Ff$.
\end{enumeratei}
\end{lemme}

(ii) \ALORS (i). Soit $Y\in \Ff$, il faut prouver que l'ouvert $U=X-Y$ est connexe. Or, si $U'$ et $U''$ sont deux composantes irr\'eductibles de $U$, il existe deux composantes irr\'eductibles $X'$ et $X''$ de $X$ telles que $X''\cap U=U''$ et $X'\cap U=U'$; soit $X_0, \ldots X_n$ une suite de composantes irr\'eductibles de $X$ poss\'edant la propri\'et\'e \'evoqu\'ee plus haut; si l'on pose $U_i=X_i\cap U$, $0\leq i\leq n$, les $U_i$ seront des composantes irr\'eductibles de $U$, de plus $U_{i}\cap U_{i+1}$ est non vide si \sisi{$0\leq i\leq n$}{$0\leq i< n$}, car sinon, {$X_{i}\cap X_{i+1} \subset Y$} serait \sisi{$\in$}{\'el\'ement de} $\Ff$, ce\pageoriginaled qui est contraire au choix de la suite des $X_i$. Ceci entra\^ine que $U$ est connexe.

(i) \ALORS (ii). Soit $Y=\bigcup_{X', X''} X'\cap X''$ o\`u l'on impose que $X'$ et $X''$ soient deux composantes irr\'eductibles \emph{distinctes} de $X$ telles que $X'\cap X'' \in \Ff$. La famille des $X'\cap X''$ est localement finie car $X$ est localement noeth\'erien, de plus les $X'\cap X''$ sont ferm\'es, donc $Y$ est ferm\'e. Par ailleurs, $Y$ appartient localement \`a $\Ff$, donc $Y\in \Ff$. Donc $U=X-Y$ est connexe. Soient $X'$ et $X''$ deux composantes irr\'eductibles distinctes de~$X$, soient $U'$ et $U''$ leurs traces sur $U$, qui sont non vides par construction de $Y$. Ce sont des composantes irr\'eductibles de $U$, or $U$ est connexe, donc $U$ \'etant localement noeth\'erien, il existe une suite de composantes irr\'eductibles $U_0, \ldots, U_n$ de $U$, telles que $U_0=U'$, $U_n=U''$ et $U_{i}\cap U_{i+1}\neq \emptyset$ et $\neq U_{i}, 0\leq i < n$. Soit $X_0, \ldots, X_n$ la suite de composantes irr\'eductibles de $X$ qui est telle que $X_{i}\cap U=U_{i}$; si $X_{i}\cap X_{i+1}\in \Ff$, par construction de $\Ff$, $U_{i}\cap U_{i+1}=\emptyset$ ou $U_{i}=U_{i+1}$ ce qui n'est pas possible d'apr\`es le choix des $U_{i}$.
\end{proof}

\begin{corollaire} \label{III.3.9}
Soit $A$ un anneau local noeth\'erien. On suppose que pour tout id\'eal premier $\pp$ de $A$, on a:
\[
(\dim A_\pp \geq 2) \To (\prof A_\pp \geq 2)
\]
On suppose de plus que $A$ satisfait la condition des cha\^ines\sfootnote{\Cf $\EGA 0_{\textup{IV}}$ 14.3.2}. Alors, pour tout $\pp$, id\'eal premier minimal de $A$, $\dim A/\pp =\dim A$, ou encore, toutes les composantes irr\'eductibles de $\Spec A$ ont m\^eme dimension: celle de $A$.
\end{corollaire}

Si $X'$ et $X''$ sont deux composantes irr\'eductibles de $X$, on les joint par une cha\^ine ayant les propri\'et\'es \'enum\'er\'ees dans \sisi{\ignorespaces}{le corollaire} \Ref{III.3.7}; il suffit alors de d\'emontrer que deux composantes successives ont la m\^eme dimension, ce qui r\'esulte de la deuxi\`eme hypoth\`ese.

\begin{exemple} \label{III.3.10}
Soit\pageoriginale $X$ la r\'eunion de deux sous-espaces vectoriels suppl\'ementaires de dimensions respectives $2$ et $3$ dans un espace vectoriel de dimension $5$; plus pr\'ecis\'ement, soit $X=\Spec A$, avec $A=B/{\pp \cap \qq}$, o\`u $B=k[X_1, \ldots, X_5]$, $\pp$ est l'id\'eal engendr\'e par $X_1, X_2, X_3$ et $\qq$ l'id\'eal engendr\'e par $X_4$ et $X_5$; $X$ peut \^etre disconnect\'e par le point d'intersection $x$ des deux sous-espaces vectoriel, donc la profondeur de $\Oo_{X, x}$ est \'egale \`a $1$, car elle ne peut \^etre $\geq 2$ en vertu du th\'eor\`eme \Ref{III.3.6}. Autre raison: la conclusion d'\'equidimensionnalit\'e du corollaire pr\'ec\'edent est en d\'efaut.

Plus g\'en\'eralement, prenant une r\'eunion $X$ de deux sous-espaces vectoriels de dimension $p, q\geq 2$ dans un espace vectoriel de dimension $p+q$, pour aucun plongement de $X$ dans un sch\'ema r\'egulier, $X$ n'est m\^eme ensemblistement une intersection compl\`ete \`a l'origine, car (quitte \`a le modifier sans changer l'espace topologique sous-jacent au voisinage de l'origine), $X$ serait Cohen-Macaulay donc de profondeur $\geq 2$ \`a l'origine, ce qui n'est pas le cas.
\end{exemple}

\begin{remarque} \label{III.3.11}
Soit $X$ un pr\'esch\'ema localement noeth\'erien, $Y$ un sous-pr\'esch\'ema ferm\'e de $X$, $F$ un $\OX$-module. La profondeur est une notion purement topologique qui s'exprime en termes de nullit\'e des $\SheafH_Y^i(X, F)$ pour $i<n$. On d\'esire \'egalement \'etudier ces faisceaux pour un $i$ donn\'e, ou pour $i>n$. On d\'emontre \`a ce propos le r\'esultat suivant:
\end{remarque}

\begin{lemme} \label{III.3.12}
Soit $m$ un entier, pour que $\SheafH_Y^i(X, F)=0$, $i>m$ pour tout $F$ coh\'erent, il faut et il suffit que ce soit vrai pour $F=\OX$.
\end{lemme}

Par limite inductive c'est alors vrai pour tout faisceau quasi\sisi{~}{-}coh\'erent. Par exemple, si $Y$ peut \^etre d\'ecrit localement par $m$ \'equations, ou, comme on dit, si $Y$ est localement ensemblistement une intersection compl\`ete (ce qui se produit par exemple si $X$ et $Y$ sont \emph{non singuliers}) il r\'esulte du calcul des $\SheafH_Y^i(X, F)$ par le complexe de \sisi{KOSZUL}{Koszul} que ces faisceaux sont nuls pour $i>m$. On a d'ailleurs utilis\'e ce fait implicitement dans l'exemple \Ref{III.3.10}. Cette condition cohomologique n'est cependant pas suffisante, comme le prouve l'exemple ci-apr\`es:

\begin{exemple} \label{III.3.13}
Soit\pageoriginale $X=\Spec (A)$, o\`u $A$ est un anneau local noeth\'erien, normal de dimension $2$. Soit $Y$ une courbe dans $X$. On peut d\'emontrer que le compl\'ementaire de la courbe est un ouvert affine donc\nde{Il y avait une coquille dans l'\'edition originale.} \sisi{$\SheafH_Y^i(X, \OX) \approx \SheafH_Y^i(X, \OX)$ $\H^{i-1}(X-Y, \OX)=0$}{$\SheafH_Y^i(\OX) \approx \SheafH_{X-Y}^{i-1}(\OX)=0$ pour $i>1$ car $\H^{i-1}(X-Y, \OX)=0$}. Cependant on peut construire une courbe qui n'est pas d\'ecrite par une \'equation.
\end{exemple}

Nous chercherons\sfootnote{\Cf Exp~\Ref{VIII}.} des conditions pour que les $\SheafH_Y^i(X, F)$ soient coh\'erents pour un~$i$ donn\'e, ce qui n'est pas le cas en g\'en\'eral, comme le montrent des exemples \'evidents, par exemple $\H_\mm^n(A)$ pour $A$ un anneau local noeth\'erien de dimension $n>0$; lorsque par exemple $A$ est un anneau de valuation discr\`ete de corps des fractions $K$, on trouve $\H_\mm^1(A) \simeq K/A$, qui n'est pas un module de type fini sur $A$. Pour \'eclairer la lanterne du lecteur, disons que le probl\`eme pos\'e est \'equivalent au suivant: \sisi{S}{s}oit $f:U\to X$ une immersion ouverte, soit $G$ un faisceau coh\'erent sur $U$, trouver des crit\`eres pour que les images directes sup\'erieures $\R^i f_* G$ soient des faisceaux coh\'erents sur $X$ pour un $i$ donn\'e. Ces conditions sont n\'ecessaires pour l'utilisation de la g\'eom\'etrie formelle que nous avons vu dans \sisi{Exp}{l'expos\'e} \Ref{IX} et \sisi{\ignorespaces}{les} suivants.

\chapter{Modules et foncteurs dualisants} \label{IV}

\section{G\'en\'eralit\'es sur les foncteurs de modules} \label{IV.1}
\pageoriginale\ignorespaces
\hspace*{-3mm}
\begin{tabular}{ll}
Soient &$A$ un anneau noeth\'erien commutatif,\\
&$\ccat$ la cat\'egorie des $A$-modules de type fini,\\
&$\ccat'$ la cat\'egorie des $A$-modules quelconques,\\
&$\Ab$ la cat\'egorie des groupes ab\'eliens.\refstepcounter{toto}\label{pIV.1}
\end{tabular}

\smallskip
Le but de ce paragraphe est l'\'etude de certaines propri\'et\'es des foncteurs \hbox{$T: \ccat^{\circ}\to \Ab$} (suppos\'es additifs).

Remarquons que si $M\in \Ob \ccat$, $T(M)$ peut \^etre muni de fa\c
con canonique d'une structure de $A$-module qui est la suivante:
si $f_{M}$ est l'homoth\'etie de $M$ associ\'ee \`a $f\in A$, $A$
op\`ere sur $T(M)$ par $f_{T(M)}$. En d'autres termes, $T$ se
factorise en:
$$
\xymatrix{ \ccat^{\circ}\ar^{T}[rr]\ar_{T_{\circ}}[dr]&&\Ab\\
&\ccat'\ar[ur]&}
$$
o\`u $\ccat'\to \Ab$ est le foncteur canonique.

Dans la suite, $T(M)$ sera consid\'er\'e comme muni de sa structure de $A$-module.

En composant avec l'isomorphisme \hbox{$M\isomto \Hom_{A}(A, M)$} le morphisme $\Hom_{A}(A, M)\to \Hom_{A}(T(M), T(A))$, on obtient les morphismes suivants qui se d\'eduisent l'un de l'autre de mani\`ere \'evidente:
\begin{align*}
M&\to \Hom_{A}(T(M), T(A)),\\
M\times T(M)&\to T(A),\\
T(M)&\to \Hom_{A}(M, T(A)),
\end{align*}
ce qui nous d\'efinit un morphisme $\varphi_{T}$ de foncteurs contravariants:
\begin{eqnarray*}
\varphi_{T}: T\to \Hom_{A}(M, T(A)).
\end{eqnarray*}

\begin{proposition} \label{IV.1.1}
Les\pageoriginale deux propri\'et\'es suivantes sont \'equivalentes:
\begin{enumeratei}
\item
$\varphi_{T}$ est un isomorphisme de foncteurs.

\item
$T$ est exact \`a gauche.
\end{enumeratei}
\end{proposition}

L'implication (i)\ALORS (ii) est triviale.

L'implication (ii)\ALORS (i) r\'esulte de ce que pour un morphisme $u:F\to F'$ de deux foncteurs \sisi{\ignorespaces}{additifs} exacts \`a gauche $F$ et $F'$ de $\ccat^{\circ}$ dans $\Ab$, si $u(A)$ est un isomorphisme, $u$ est un isomorphisme (on utilise le fait que $A$ est noeth\'erien et donc que tout $A$-module de type fini est de pr\'esentation finie).

\begin{remarque} \label{IV.1.2}
Ceci montre en particulier que les foncteurs $T:{\ccat'}^{\circ}\to \Ab$ qui sont repr\'esentables sont les foncteurs qui commutent aux limites projectives quelconques (sur un ensemble pr\'eordonn\'e non n\'ecessairement filtrant).
\end{remarque}

Si $\SheafHom (\ccat^{\circ}, \Ab)_{g}$ d\'esigne la sous-cat\'egorie pleine de $\SheafHom (\ccat^{\circ}, \Ab)$ dont les objets sont les foncteurs exacts \`a gauche, on~a d\'emontr\'e l'\'equivalence des cat\'egories
\[
\ccat' \sisi{\overset{\thickapprox}{\to}}{\isomto} \SheafHom (\ccat^{\circ}, \Ab)_{g}
\]
par les foncteurs quasi-inverses l'un de l'autre
\[
H \sisi{\rightsquigarrow}{\mto}\Hom_{A}(\;, H)
\]
et
\[
T(A) \mfrom T.
\]

Soient maintenant $J$ un id\'eal de $A$, $Y=V(J)\subset \Spec\, A$, et d\'esignons par $\ccat_{Y}$ la sous-cat\'egorie pleine de $\ccat$ dont les objets sont les $A$-modules de type fini $M$ tels que $\supp M\subset Y$. On a:
\[
\ccat_{Y}=\bigcup_{n}\ccat^{(n)},
\]
o\`u $\ccat^{(n)}$ est la sous-cat\'egorie pleine de $\ccat_{Y}$ des modules $M$ tels que $J^{n}M=0$.

\begin{proposition} \label{IV.1.3}
Avec\pageoriginale les m\^emes notations que pr\'ec\'edemment, soit $T:\ccat_{Y}\to \Ab$ un foncteur. \`A $H=\sisi{}{\varinjlim T(A/J^n)}$\nde{la d\'efinition de $H$ est implicite dans le texte original.} est associ\'e un morphisme naturel
\[
\varphi_{T}: T\to \Hom_{A}(\;, H),
\]
et les conditions suivantes sont \'equivalentes:
\begin{enumeratei}
\item
$\varphi_{T}$ est un isomorphisme.

\item
Le foncteur $T$ est exact \`a gauche.
\end{enumeratei}
\end{proposition}

\begin{proof}
a) D\'efinition de $\varphi_{T}$.

Soit $M\in \Ob \ccat_{Y}$. Il existe un entier $n$ tel que $J^{n}M=0$. Alors $M$ est un $A/J^{n}$-module, et si $T_{n}$ d\'esigne la restriction de $T$ \`a $\ccat^{(n)}$, on sait d\'efinir le morphisme
\begin{eqnarray*}
T_{n}\to \Hom_{A}(\;, H_{n}), \ \hbox{o\`u}\ H_{n}=T(A/J^{n})\;;
\end{eqnarray*}
\begin{equation*} {T(M)=T_{n}(M)\to \Hom_{A}(M, \varinjlim H_{n})= \Hom_{A}(M, H)}\leqno{{\text{d'o\`u}}}
\end{equation*}
\begin{equation*} \varphi_{T}: T\to \Hom_{A}(\;, H).\leqno{\text{et}}
\end{equation*}

b) \'Equivalence de (i) et (ii).

Il est clair que (i) \sisi{entraine}{entra\^ine} (ii). Supposons (ii) v\'erifi\'e et soit $M\in \Ob \ccat^{(n)}$. On a vu que $T_{n}(M)\isomto \Hom_{A}(M, H_{n})$, donc pour tout entier $n'>n$ on a
\begin{equation*} T(M)=T_{n}(M)=T_{n'}(M)=\varinjlim T_{n}(M)
\end{equation*}
\begin{equation*} T(M)=\varinjlim \Hom_{A}(M, H_{n}).\leqno\text{et}
\end{equation*}
Comme il s'agit ici de limites inductives filtrantes, on~a aussi l'isomorphisme
\begin{equation*} \varinjlim \Hom_{A}(M, H_{n})\isomto \Hom_{A}(M, \varinjlim H_{n}) = \Hom_{A}(M, H).
\end{equation*}

Si $\ccat'_{Y}$ d\'esigne la cat\'egorie des $A$-modules, de support contenu dans $Y$, mais non n\'ecessairement de type fini, on~a encore l'\'equivalence naturelle de cat\'egories: \hbox{$\ccat'_{Y}\sisi{\overset{\thickapprox}{\to}}{\isomto} \Hom (\ccat_{Y}^{\circ}, \Ab)_{g}$}.
\skipqed
\end{proof}

\subsection*{Application}
Les notations\pageoriginale \'etant les m\^emes que pr\'ec\'edemment, soit
\begin{equation*} T^{*}:\ccat_{Y}^{\circ}\to{\Ab}
\end{equation*}
un $\partial$-foncteur exact. Pour tout $i\in \ZZ$, on pose $H^{i}_{n}=T^{i}(A/J^{n})$ et $H^{i}=\varinjlim H^{i}_{n}$.

\begin{theoreme}
\label{IV.1.4} Soit $n\in \ZZ$. S'il existe $i_{0}\in \ZZ$ tel que $T^{i}=0$ pour tout $i<i_{0}$, alors les trois conditions suivantes sont \'equivalentes:
\begin{enumeratei}
\item
$T^{i}=0$ pour tout $i<n$.
\item
$H^{i}=0$ pour tout $i<n$.
\item
Il existe un module $M_{0}$ de $\ccat_{Y}$ tel que $\supp M_{0}=Y$ et $T^{i}(M_{0})=0$ pour tout $i<n$.
\end{enumeratei}
\end{theoreme}

\begin{proof}
Il est \'evident que (i) \sisi{entraine}{entra\^ine} (ii) et (iii) (on prend $M_{0}=A/J$). Montrons par r\'ecurrence sur $n$ que (ii) \sisi{entraine}{entra\^ine} (i); c'est vrai pour $n=i_{0}$, et on le suppose d\'emontr\'e jusqu'au rang $n$. Supposons alors que $H^{i}=0$ pour tout $i<n+1$; par l'hypoth\`ese de r\'ecurrence on~a $T^{i}=0$ pour $i<n$, mais $T^{n-1}=0$ \sisi{entraine}{entra\^ine} que $T^{n}$ est un foncteur exact \`a gauche et
\begin{equation*} T^{n}\isomto\Hom_{A}(\;, H^{n})=0.
\end{equation*}
Montrons alors que (iii) \sisi{entraine}{entra\^ine} (ii). C'est encore vrai pour $n=i_{0}$; supposons-le d\'emontr\'e jusqu'au rang $n$: soit $M_{0}$ un $A$-module de $\ccat_{Y}$ tel que $\supp M_{0}=Y$ et $T^{i}(M_{0})=0$ pour tout $i<n+1$; par l'hypoth\`ese de r\'ecurrence on~a alors $H^{i}=0$ pour tout $i<n$; il reste \`a montrer que $H^{n}=0$. Mais \og $H^{i}=0$ pour tout $i<n$\fg implique que $T^{n-1}=0$ et donc que $T^{n}\isomto\Hom_{A}(\;, H^{n})$. On a alors
\begin{equation*}
\Ass H^{n}=\Ass \Hom(M_{0}, H^{n})=\supp M_{0}\cap \Ass H^{n}= \Ass H^{n}
\end{equation*}
car
\begin{equation*}
\Ass H^{n}\subset \supp H^{n}\subset Y= \supp M_{0}.
\end{equation*}

D'o\`u $T^{n}(M_{0})=0\Ssi \Ass H^{n}=\emptyset \Ssi H^{n}=0$; ce qui ach\`eve la d\'emonstration.
\skipqed
\end{proof}

\section{Caract\'erisation des foncteurs exacts} \label{IV.2}

L'anneau\pageoriginale $A$ est toujours suppos\'e noeth\'erien et commutatif. Les notations sont celles de la \sisi{Proposition 2}{proposition~\Ref{IV.1.3}}:
\begin{equation*} Y=V(J), \quad T:\ccat_{Y}^{\circ}\to \Ab, \quad H=\varinjlim T(A/J^{n}),
\end{equation*}
o\`u on suppose que $T$ est un foncteur exact \`a gauche, d'o\`u:
\begin{equation*} T(M)\isomto \Hom_{A}(M, H).
\end{equation*}

\begin{proposition} \label{IV.2.1}
Les propri\'et\'es suivantes sont \'equivalentes:
\begin{enumeratei}
\item
Le foncteur $T$ est exact,

\item
$H$ est injectif dans $\ccat'$.
\end{enumeratei}
\end{proposition}

\begin{proof}
Il suffit \'evidemment de montrer que (i) implique (ii), c'est-\`a-dire de d\'emontrer que si la restriction \`a $\ccat_{Y}$ du foncteur $\Hom_{A}(\;, H)$ est un foncteur exact, $H$ est injectif. Mais comme $A$ est noeth\'erien, pour montrer que $H$ est injectif, on peut se borner \`a prouver que tout homomorphisme $f:N\to H$ de source un $A$-module $N$ \emph{de type fini}, sous-module d'un $A$-module $M$ \emph{de type fini}, se prolonge en un homomorphisme $\overline{f}:M\to H$.

La d\'efinition de $H$ et le fait que $N$ soit de type fini \sisi{entrainent}{entra\^inent} qu'il existe un entier~$n$ tel que $J^{n}.f(N)=0$. Munissons alors $M$ et $N$ de la topologie $J$-adique. La topologie $J$-adique de $N$ est \'equivalente \`a la topologie induite par la topologie $J$-adique de $M$ (th\'eor\`eme de Krull). Il existe donc $V=J^{k}\cdot M$ tel que
\begin{equation*} U=V\cap N\subset J^{n}N.
\end{equation*}
On a alors la factorisation
$$
\xymatrix{ N\ar_{f}[d]\ar[r]&N/U\ar^{u}[dl]\\
H&},
$$
$N/U$ et $M/V$ sont dans $\ccat_{Y}$. L'hypoth\`ese permet donc de prolonger $u$ en $\overline{u}$
$$
\xymatrix{ N/U\ar_{u}[d]\ar@{^{ (}->}[r]&M/V\ar^{\overline{u}}[dl]\\
H&},
$$
et\pageoriginale $M\to M/V\lto{\overline{u}} H$ donne le prolongement cherch\'e $\overline{f}$.
\skipqed
\end{proof}

\begin{corollaire} \label{IV.2.2}
Soit $K$ un $A$-module injectif, alors le sous-module $H^{0}_{J}(K)$ de $K$ form\'e des \'el\'ements annul\'es par une puissance convenable de $J$ est injectif.
\end{corollaire}

\begin{proof}
Il suffit de v\'erifier que la restriction \`a $\ccat_{Y}$ du foncteur $ \Hom_{A}(\;, H^{0}_{J}(K))$ est un foncteur exact. Or soit $M\in \Ob \ccat_{Y}$, il existe $k$ tel que $J^{k}\cdot M=0$, et l'inclusion
$$
\Hom_{A}(\;, H^{0}_{J}(K))\to \Hom_{A}(M, K)
$$
est alors un isomorphisme. D'o\`u le r\'esultat puisque $\Hom_{A}(\;, K)$ est exact.
\skipqed
\end{proof}

\section[\'Etude du cas o\`u $T$ est exact \`a gauche]
{\'Etude du cas o\`u $T$ est exact \`a gauche et $T(M)$ de type fini pour tout $M$} \label{IV.3}

Soit comme plus haut
\begin{equation*} T:\ccat_{Y}^{\circ}\to \Ab,
\end{equation*}
on suppose d\'esormais que $T$ est exact \`a gauche et qu'on a la factorisation
$$
\xymatrix@C=6mm@R=6mm{
\ccat^{\circ}_{Y}\ar[dr]\ar^{T}[rr]&&\Ab\,.\\
&\ccat_{Y}\ar[ur]&}
$$
\sisi{}{o\`u comme plus haut, $\ccat_Y\to\Ab$ est le foncteur d'oubli.} On sait donc d\'efinir $T(T(M))=T\circ T(M)$, et le morphisme canonique
$$
M\to \Hom_{A}(\Hom_{A}(M, H), H)
$$
d\'efinit un morphisme\pageoriginale
$$
M\to T\circ T(M).
$$

\begin{proposition} \label{IV.3.1}
L'anneau $A$ \'etant toujours suppos\'e noeth\'erien, si on fait l'hypoth\`ese suppl\'ementaire que $A/J$ est artinien, les conditions suivantes sont \'equivalentes:
\begin{enumeratei}
\item
$T$ est exact \`a gauche et, pour tout $M\in \Ob \ccat_{Y}$, $T(M)$ est de type fini et \hbox{$M\to T\circ T(M)$} est un isomorphisme.

\item
$T$ est exact et, pour tout corps r\'esiduel $k$ associ\'e \`a un id\'eal maximal contenant~$J$, on~a $T(k)\isomto k$.

\item
On a $T\isomto \Hom_{A}(\;, H)$ avec $H$ injectif et, pour tout $k$ comme dans \emph{(ii)}, on~a $\Hom_{A}(k, H)\isomto k$.

\item
$T$ est exact et, pour tout $M\in \Ob \ccat_{Y}$, on~a $\longueur T(M)=\longueur M$.
\end{enumeratei}
\end{proposition}

\begin{proof}
On a d\'ej\`a montr\'e l'\'equivalence de (ii) et (iii) (prop\ptbl \Ref{IV.2.1}). Montrons que (ii) \sisi{entraine}{entra\^ine} (iv): d'abord si $M\in \Ob \ccat_{Y}$, comme $M$ est un $A/J^{n}$-module avec $A/J^{n}$ artinien, $\longueur M$ est finie. Raisonnons par r\'ecurrence sur la longueur de $M$. La condition (iv) est vraie si $\longueur M=1$, parce qu'alors $M$ est un \sisi{$k$}{ corps r\'esiduel justiciable de (ii)}. Si $\longueur M>1$, il existe un sous-module $M'$ de $M$ tel que $M'\neq 0$ et que $\longueur M'<\longueur M$. Formons alors la suite exacte
$$
0\to M'\to M\to M''\to 0.
$$
Comme $T$ est exact, on~a la suite
$$
0\to T(M')\to T(M)\to T(M'')\to 0,
$$
et $\longueur T(M)=\longueur T(M')+\longueur T(M'')=\longueur M'+\longueur M''=\longueur M$.

(ii) \sisi{entraine}{entra\^ine} (i): Comme (ii) \sisi{entraine}{entra\^ine} (iv), soit $M$ un $A$-module de $\ccat_{Y}$, on~a $\longueur T(M)=\longueur M$; donc $T(M)$ est de longueur finie et par suite de type fini.

Il reste \`a montrer que $M\to T\circ T(M)$ est un isomorphisme; on raisonne encore par r\'ecurrence sur $\longueur M$. Pour $M=k$ c'est vrai. Dans le cas g\'en\'eral on \'ecrit le diagramme commutatif dont les deux lignes sont exactes\pageoriginaled:
$$
\xymatrix{ 0\ar[r]&M'\ar[r]\ar[d]&M\ar[r]\ar[d]&M''\ar[r]\ar[d]&0\\
0\ar[r]& T\circ T(M')\ar[r]& T\circ T(M)\ar[r]& T\circ T(M'')\ar[r]&0, }
$$
o\`u $M'$ est un sous-module de $M$ tel que $M'\neq 0$ et $\longueur M'<\longueur M$. Par l'hypoth\`ese de r\'ecurrence, les fl\`eches extr\^emes sont des isomorphismes, donc
$$
M\to T\circ T(M)
$$
est un isomorphisme.

(i) \sisi{entraine}{entra\^ine} (ii): soit
$$
0\to M'\to M\to M''\to 0
$$
une suite exacte de $A$-modules de $\ccat_{Y}$, et soit $Q$ le conoyau de $T(M)\to T(M')$. Appliquons $T$ \`a la suite exacte
$$
0\to T(M')\to T(M)\to T(M'')\to Q \to 0,
$$
on obtient:
$$
\xymatrix{
0\ar[r]&T(Q)\ar[r]&T\circ T(M')\ar[r]&T\circ T(M)&\\
& &M'\ar_{s}[u]\ar[r]&M\ar_{s}[u],
}
$$
donc $T(Q)=0$ et $Q\isomto T(T(Q))=0$.

D'autre part soit $k$ un corps r\'esiduel, $k=A/\m$, $J\subset \m$. Il faut montrer que $T(k)\isomto k$. Pour cela il suffit de remarquer que $T(k)$ est un $k$-espace vectoriel. On en d\'eduit:
\begin{gather*}
T(k)\simeq k\oplus V,\\
T(T(k))\simeq T(k)\oplus T(V)\simeq k\oplus V\oplus T(V)\simeq k,
\end{gather*}
d'o\`u $V=0$.

Montrons enfin que (iv) \sisi{entraine}{entra\^ine} (iii): il suffit de montrer que $T(k)\isomto k$; or $\longueur T(k)= \longueur k=1$, donc $T(k)=k'$ est un corps r\'esiduel et $\supp k'= \supp \Hom_{A}(k, H)\subset \supp k$. Donc $k\sisi{=}{\simeq}k'$.
\skipqed
\end{proof}

\begin{remarque} \label{IV.3.2}
On peut\pageoriginale montrer que la condition (iv) est \'equivalente \`a la condition

(iv)$'$ Pour tout $M\in \Ob \ccat_{Y}$, on~a $\longueur T(M)=\longueur M$.
\end{remarque}

\section{Module dualisant. Foncteur dualisant} \label{IV.4}

\begin{definition} \label{IV.4.1}
Soient $A$ un anneau local noeth\'erien et $\m$ son id\'eal maximal. On appelle foncteur dualisant pour $A$ tout foncteur
$$
T:\ccat^{\circ}_\mm\to \Ab,
$$
o\`u on note $\ccat_{\mm}$ au lieu de $\ccat_{Y}$ pour $Y=V(\mm)$, qui satisfait aux conditions \'equivalentes de la proposition \Ref{IV.3.1}. On dit qu'un $A$-module $I$ est dualisant pour $A$ si le foncteur $M\to \Hom_{A}(M, I)$ est dualisant.
\end{definition}

On peut g\'en\'eraliser la d\'efinition \Ref{IV.4.1} au cas o\`u on ne suppose plus que $A$ soit un anneau local.

\begin{definition} \label{IV.4.2}
Soient $A$ un anneau noeth\'erien et soit $\bar{\ccat}$ la sous-cat\'egorie pleine de~$\ccat$ form\'ee des $A$-modules de longueur finie; on appelle foncteur dualisant tout foncteur~$T$, $A$-lin\'eaire, de $\bar{\ccat}^{\circ}$ dans $\bar{\ccat}$, qui est exact et tel que le morphisme de foncteurs
$$
\id \to T\circ T
$$
soit un isomorphisme.
\end{definition}

On va d\'emontrer un th\'eor\`eme d'existence et aussi que le module $I$ qui repr\'esente un tel foncteur est localement artinien. On montrera aussi que, pour tout id\'eal maximal $\mm$ de $A$, la composante $\mm$-primaire du socle de $I$ est de longueur 1.

\begin{proposition} \label{IV.4.3}
Soient $A$ et $B$ deux anneaux locaux noeth\'eriens d'id\'eaux maximaux $\mm_{A}$ et $\mm_{B}$, tels que $B$ soit une $A$-alg\`ebre finie. Alors, si $I$ est un module dualisant pour~$A$, $\Hom_{A}(B, I)$ est un module dualisant pour $B$.
\end{proposition}

\begin{proof}
Soit\pageoriginale
$$
R: \ccat_{\mm_{B}}\to \ccat_{\mm_{A}}
$$
le foncteur restriction des scalaires; il est exact. Soit $T$ un foncteur dualisant pour $A$,
$$
T: \ccat_{\mm_{A}}\to \Ab;
$$
il est exact et, pour tout $M\in \Ob \ccat_{\mm_{A}}$, le morphisme naturel $M\to T\circ T(M)$ est un isomorphisme; donc $T\circ R$ est un foncteur dualisant pour $B$. Si $I$ repr\'esente $T$, d'apr\`es la formule classique $\Hom_{A}(M, I)=\Hom_{B}(M, \Hom_{A}(B, I))$, valable pour tout $B$-module $M$, on en d\'eduit que $\Hom_{A}(B, I)$ est un module dualisant pour $B$.
\skipqed
\end{proof}

\begin{corollaire} \label{IV.4.4} Soient $A$ un anneau local noeth\'erien et $\aaa$ un id\'eal de $A$; soit $B=A/\aaa$. Si $I$ est un module dualisant pour $A$, l'annulateur de $\aaa$ dans $I$ est un module dualisant pour $B$.
\end{corollaire}

\begin{lemme} \label{IV.4.5} Soient $A$ un anneau local noeth\'erien et $I$ un $A$-module localement artinien. Il existe un isomorphisme canonique
$$
I\to \hat{I}=I\otimes_{A}\hat{A}.
$$
\end{lemme}

\begin{proof}
Soit $I_{n}$ l'annulateur de $\mm^n$ dans $I$, o\`u $\mm$ est l'id\'eal maximal de $A$. Dire que $I$ est localement artinien, c'est dire que $I$ est limite inductive des $I_{n}$ et que ceux-ci sont de longueur finie. Or le produit tensoriel commute aux limites inductives, on est donc ramen\'e au cas o\`u $I$ est artinien. Dans ce cas $I$ est annul\'e par une puissance de l'id\'eal maximal, soit $\mm^{k}$; donc pour $p\geq k$, $I\isomto I\otimes_{A}A/\mm^{p}$, donc $I\isomto I\otimes_{A}\hat{A}$ car $A$ est noeth\'erien et $I$ est de type fini.

On en conclut que le foncteur restriction des scalaires de $\hat{A}$ \`a $A$ et le foncteur extension des scalaires induisent des \emph{\'equivalences} quasi-inverses l'une de l'autre de la cat\'egorie des $\hat{A}$-\emph{modules localement artiniens}\pageoriginale et de la cat\'egorie des $A$-\emph{modules localement artiniens}.
\skipqed
\end{proof}

\begin{proposition} \label{IV.4.6}
Soient $A$ un anneau local noeth\'erien, $\hat{A}$ son compl\'et\'e, $I$ un module dualisant pour $A$ (\resp pour $\hat{A}$) et $J$ le compl\'et\'e\,\nde{il faut comprendre ici le produit tensoriel $\hat{I}=I\otimes_A\hat{A}$ (\cf lemme \Ref{IV.4.5}), \`a savoir $I$ muni de sa structure de $\hat{A}$-module canonique, et non le compl\'et\'e $\m$-adique. Par exemple, si $p$ est un nombre premier et $A\!=\!\hat{A}\!=\!\ZZ_p$ est l'anneau des entiers $p$-adiques. Alors, l'enveloppe injective du corps r\'esiduel $k\!=\!\ZZ/p\ZZ$ est le $\ZZ_p$-module discret $\QQ_p/\ZZ_p$, dont le compl\'et\'e pour la topologie $p$-adique est~nul.} de $I$ (\resp le $A$-module obtenu par restriction des scalaires). Alors $J$ est un module dualisant pour $\hat{A}$ (\resp pour $A$). De plus les groupes ab\'eliens sous-jacents \`a $I$ et $J$ sont isomorphes.
\end{proposition}

\begin{proof}
On remarque simplement que l'\'equivalence entre la cat\'egorie des $A$-modules localement artiniens et la cat\'egorie des $\hat{A}$-modules localement artiniens induit un isomorphisme entre les bifoncteurs $\Hom_{A}(\;, \;)$ et $\Hom_{\hat{A}}(\;, \;)$, et que la caract\'erisation d'un foncteur ou d'un module dualisant ne fait intervenir que ces bifoncteurs.
\skipqed
\end{proof}

\begin{theoreme} \label{IV.4.7}
Soit $A$ un anneau local \sisi{neoth\'erien}{noeth\'erien}.
\begin{enumeratea}
\item
Il existe toujours un module dualisant $I$.

\item
Deux modules dualisants sont isomorphes (par un isomorphisme non canonique).

\item
Pour qu'un module $I$ soit dualisant, il faut et il suffit qu'il soit une enveloppe injective du corps r\'esiduel $k$ de $A$.
\end{enumeratea}
\end{theoreme}

\begin{remarque} \label{IV.4.8}
La proposition \Ref{IV.4.6} permet de se ramener au cas d'un anneau local noeth\'erien complet. D'apr\`es un th\'eor\`eme de structure de \sisi{COHEN}{Cohen}\nde{voir Cohen~I.S., {\og On the structure and ideal theory of complete local rings\fg}, \emph{Trans. Amer. Math. Soc.} \textbf{59} (1946), p\ptbl 54--106.}, un tel anneau est quotient d'un anneau r\'egulier. La proposition \Ref{IV.4.3} permet alors de supposer $A$ r\'egulier. Comme nous le verrons plus loin, cette remarque permet un calcul explicite du module dualisant\sfootnote{C'\'etait la m\'ethode suivie par Grothendieck (en 1957). La m\'ethode par enveloppes injectives qui va suivre est due semble-t-il \`a K\ptbl Morita, \sisi{Sc. Rep. Tokyo Kyoiku Daigaku, t\ptbl 6, 1958-59, p\ptbl 83-142}{ {\og Duality for modules and its applications to the theory of rings with minimum conditions\fg}, \emph{Sc. Rep. Tokyo Kyoiku Daigaku} \textbf{6} (1958/59), p\ptbl 83-142}. Le travail de Morita est d'ailleurs ind\'ependant de celui de Grothendieck et bien ant\'erieur au pr\'esent s\'eminaire, et ne se limite pas au cas des anneaux de base commutatifs.}; nous d\'emontrerons cependant le th\'eor\`eme \Ref{IV.4.7} par d'autres moyens.
\end{remarque}

\subsection*{Rappels}
Avant\pageoriginale de d\'emontrer le th\'eor\`eme, nous faisons quelques rappels sur la notion d'\emph{enveloppe injective}. \Cf \sisi{GABRIEL}{Gabriel}, Th\`ese, Paris 1961, \sisi{Des Cat\'egories Ab\'eliennes}{\emph{Des Cat\'egories Ab\'eliennes}}, ch\ptbl II \S 5.

Soit $\Cc$ une cat\'egorie ab\'elienne dans laquelle les $\varinjlim$ existent et sont exactes\nde{bien entendu, ce sont les petites limites inductives filtrantes qui sont suppos\'ees exactes; il faudrait aussi supposer l'existence d'un g\'en\'erateur. \Cf \sisi{Tohoku}{\cite{Tohoku}}. En ce qui concerne la cat\'egorie des modules, suffisante en ce qui nous concerne, on peut aussi se reporter au chapitre 10 de l'\textit{Alg\`ebre} de Bourbaki.} (ex. $\Cc=$cat\'egorie des modules). Tout objet $M$ se plonge dans un objet injectif et on appelle enveloppe injective de $M$ tout objet injectif, contenant $M$, minimal. On a les propri\'et\'es suivantes:
\begin{enumeratei}
\item
Tout objet $M$ a une enveloppe injective $I$.
\item
Si $I$ et $J$ sont deux enveloppes injectives de $M$, il existe entre $I$ et $J$ un isomorphisme (en g\'en\'eral non unique) qui induit l'identit\'e sur $M$.
\item
$I$ est une extension essentielle\refstepcounter{toto}\label{pIV.13} de $M$, c'est-\`a-dire que $P\subset I$ et $P\cap M=\{ 0\}$ implique ($P=\{ 0\}$). De plus si $I$ est injectif et extension essentielle de $M$, $I$ est enveloppe injective de $M$.
\end{enumeratei}

Ces r\'esultats admis, pour d\'emontrer le th\'eor\`eme \Ref{IV.4.7}, il suffit \'evidemment de prouver c).

\begin{proof}
Soit $I$ un module dualisant pour $A$. Alors $I$ est injectif et $\Hom_{A}(k, I)$ est isomorphe \`a $k$. En composant l'isomorphisme $k\simeq \Hom_{A}(k, I)$ avec l'inclusion
$$
\Hom_{A}(k, I)\hto \Hom_{A}(A, I)\simeq I,
$$
on obtient l'inclusion
$$
k\hto I.
$$
Montrons que $I$ est enveloppe injective de $k$. Soit $J$ un module injectif tel que
$$
k\subset J\subset I.
$$
Comme $J$ est injectif, il existe un sous $A$-module injectif $J'$ de $I$ tel que $I=J\oplus J'$. Montrons que $\Hom_{A}(k, J')=0$. On a
$$
\Hom_{A}(k, I)\simeq \Hom_{A}(k, J)\oplus \Hom_{A}(k, J');
$$
$\Hom_{A}(k, J)$ est\pageoriginale un sous-espace vectoriel de $\Hom_{A}(k, I)\simeq k$ non r\'eduit \`a z\'ero (puisqu'il contient l'inclusion $k\subset J$), donc $\Hom_{A}(k, J)\simeq k$ et par suite $\Hom_{A}(k, J')=0$.

En raisonnant par r\'ecurrence sur la longueur on en d\'eduit que $\Hom_{A}(M, J')=0$ pour tout $A$-module $M$ de longueur finie; \sisi{donc $J'$ est l'objet qui repr\'esente le foncteur $0:\ccat_{\mm}^{\circ}\to \Ab$}{comme $I$ est limite inductive des modules $\Hom(A/\m^n, I)$ (\cf proposition~\Ref{IV.1.3}) qui sont de longueur finie par hypoth\`ese, la projection $I\to J'$ est nulle}, et par suite $J'=0$.

R\'eciproquement, soit $I$ une enveloppe injective de $k$. Pour voir que $I$ est un module dualisant, il suffit de montrer \sisi{\ignorespaces}{d'apr\`es~\Ref{IV.2.1} et \Ref{IV.3.1} (ii)} que $V=\Hom_{A}(k, I)$ est isomorphe \`a $k$. Or on~a la double inclusion
$$
k\subset V\subset I;
$$
$V$ est un espace vectoriel sur $k$ qui se d\'ecompose en la somme directe de $k$ et d'un sous-espace vectoriel $V'$ de $I$ tel que $V'\cap k=0$. Or $I$ est une extension essentielle de~$k$, d'o\`u $V'=0$ et $V=k$.
\skipqed
\end{proof}

\begin{corollaire} \label{IV.4.9}
Soit $A$ un anneau local noeth\'erien; tout module dualisant pour $A$ est localement artinien.
\end{corollaire}

\begin{proof}
Soit $I$ un module dualisant; c'est une enveloppe injective de $k$. Utilisons les notations et le r\'esultat du corollaire 2.2. On a
$$
k\subset H^{0}_{\mm}(I)\subset I,
$$
et $H^{0}_{\mm}(I)$ est injectif. On en d\'eduit que $I=H^{0}_{\mm}(I)$ et donc que $I$ est localement artinien.\nde{comme on l'a d\'ej\`a vu, on peut aussi simplement observer que $I$ est limite inductive des modules $\Hom_A(A/\m^n, I)$.}
\skipqed
\end{proof}

\section{Cons\'equences de la th\'eorie des modules dualisants} \label{IV.5}

Le foncteur\pageoriginale
$$
T=\Hom_{A}(\;, I)\colon \ccat_{\mm}\to \ccat_{\mm}
$$
est une anti-\'equivalence. En effet $T\circ T$ est isomorphe au foncteur identique et l'argument est formel \`a partir de l\`a.

\emph{On en d\'eduit les propri\'et\'es habituelles de la notion d'orthogonalit\'e:}

Soit $M^{*}=\Hom_{A}(M, I)=T(M)$ et soit $N\subset M$ un sous-module. On appelle \emph{orthogonal de} $N$ le sous-module $N'$ de $M^{*}$ form\'e des \'el\'ements de $M^*$ nuls sur $N$. On obtient ainsi une bijection entre l'ensemble des sous-modules de $M$ et l'ensemble des sous-modules de $M^{*}$, qui renverse la notion d'ordre.

On a en particulier:
\begin{itemize}
\item
$\longueur_{M}N=\colong_{M^{*}}N'$.
\item
Aux modules monog\`enes, \ie tels que $M/\mm M$ soit de dimension 0 ou 1, correspondent dans la dualit\'e des modules dont le socle est de longueur 0 ou 1.
\item
Si $A$ est artinien, les id\'eaux de $A$ correspondent aux sous-modules de $I$.
\par\hfill etc.
\end{itemize}

\smallskip
Soient $A$ un anneau local noeth\'erien, $DA$ la cat\'egorie des $A$-modules $M$ tels que, pour tout $n\in {\bf N}$, $M^{(n)}=M/\mm^{n+1}M$ soit de longueur finie et tels que $M=\varprojlim_{n}M^{(n)}$, et soit $\hat{A}$ le compl\'et\'e de $A$. Le foncteur restriction des scalaires et le foncteur compl\'etion sont des \emph{\'equivalences} quasi-inverses entre $DA$ et $D\hat{A}$, qui commutent \`a isomorphisme pr\`es \`a la formation des groupes ab\'eliens sous-jacents aux modules consid\'er\'es. Notons $CA$ la cat\'egorie des $A$-modules localement artiniens \`a socle de dimension finie.

\begin{proposition} \label{IV.5.1}
Soit\pageoriginale $A$ un anneau local noeth\'erien et soit $I$ un module dualisant pour $A$. Les foncteurs
$$
\Hom_{A}(\;, I)\colon (CA)^{\circ}\to DA
$$
et
$$
\Hom_{\hat{A}}(\;, I)\colon DA\to (CA)^{\circ}
$$
sont des \'equivalence de cat\'egories, quasi-inverses l'une de l'autre.

De plus, si l'on transporte ces foncteurs par les \'equivalences de cat\'egories entre $DA$ et $D\hat{A}$ d'une part, et $CA$ et $C\hat{A}$ d'autre part, on trouve le foncteur $\Hom_{\hat{A}}(\;, I)$.
\end{proposition}

\begin{proof}
Soit $X\in \Ob CA$. Par d\'efinition, on a:
$$
X=\varinjlim_{k\in {\bf N}}X_{k}, \quad X_{k}=\Hom_{A}(A/\mm^{k+1}, X),
$$
donc
$$
\Hom_{A}(X, I)=\varprojlim \Hom_{A}(X_{k}, I).
$$
Donc $Y=\varprojlim X_{k}$ est un $\hat{A}$-module de type fini comme il r\'esulte de $\EGA 0_{\textup{I}}$ 7.2.9. On remarque \`a ce propos que $D\hat{A}$ est aussi la cat\'egorie des $\hat{A}$-modules de type fini ou, si l'on veut, que $DA$ est la cat\'egorie des $A$-modules complets et de type fini sur~$\hat{A}$. Soit alors $Y$ un tel module, soit $f:Y\to I$ un $\hat{A}$-homomorphisme. L'image de $f$ est un sous-module de type fini, donc est annul\'e par $\mm^{k}$ pour un certain $k$; en effet tout $x\in I$ est annul\'e par une puissance de $\mm$. Donc $f$ se factorise par $\sisi{y}{Y}/\mm^{k}Y$, d'o\`u il r\'esulte que
\begin{align*}
\Hom_{\hat{A}}(Y, I)&=\varinjlim_k\Hom_{\hat{A}}(Y^{(k)}, I)\quad\text{avec } Y^{(k)}=Y/\m^{k+1}Y\\
&=\varinjlim_k(Y^{(k)})^*
\end{align*}
appartient \`a $\Ob CA$. D'o\`u il r\'esulte imm\'ediatement que les deux foncteurs de l'\'enonc\'e sont quasi-inverses l'un de l'autre.
\skipqed
\end{proof}

Il r\'esulte des consid\'erations pr\'ec\'edentes que l'on ne change rien
aux cat\'egories ou aux foncteurs consid\'er\'es, non plus qu'aux
groupes ab\'eliens sous-jacents aux modules consid\'er\'es, en rempla\c
cant $A$ par $\hat{A}$; la proposition \Ref{IV.5.1}\pageoriginale s'\'enonce alors ainsi:

La restriction du foncteur $\Hom_{\hat{A}}(\;, I)$ \`a la cat\'egorie des $\hat{A}$-modules de type fini prend ses valeurs dans la cat\'egorie des $\hat{A}$-modules localement artiniens \`a socle de dimension finie, et admet un foncteur quasi-inverse, qui est la restriction du foncteur $\Hom_{\hat{A}}(\;, I)$. Sur l'intersection de ces deux cat\'egories, ces deux foncteurs co\"incident (\'evidemment !) et \'etablissent une auto-dualit\'e de la cat\'egorie des $\hat{A}$-modules de longueur finie.

\begin{exemple}[\sisi{MACAULAY}{Macaulay}]
\label{IV.5.2} Soit $A$ un anneau local de corps r\'esiduel $k$. Soit $k_{0}$ un sous-corps de $A$ tel que $k$ soit fini sur $k_{0}$, $[k:k_{0}]=d$. Tout $A$-module de longueur finie peut \^etre consid\'er\'e comme un $k_{0}$-espace vectoriel de dimension finie et \'egale \`a $d\cdot\longueur(M)$. Le foncteur $T$:
$$
M\to \Hom_{k_{0}}(M, k_{0})
$$
est alors exact et conserve la longueur, donc est dualisant pour $A$. Le module dualisant associ\'e est donc:
$$
A'=\varinjlim_{n}\Hom_{k_{0}}(A/\mm^{n}, k_{0}),
$$
c'est le dual topologique de $A$ muni de la topologie $\mm$-adique.
\end{exemple}

\begin{exemple} \label{IV.5.3}
Soit $A$ un anneau local noeth\'erien \emph{r\'egulier de dimension} $n$. Soit $\mm$ son id\'eal maximal, soit $k$ son corps r\'esiduel. Il existe un syst\`eme r\'egulier de param\`etres, $(x_{1}, x_{2}, \ldots, x_{n})$, qui engendre $\mm$, et qui est une \emph{suite $A$-r\'eguli\`ere}. On peut donc calculer les $\Ext^{i}_{A}(k, A)$ par le complexe de \sisi{KOSZUL}{Koszul}; on trouve:
\begin{align*}
\Ext^{i}_{A}(k, A)&=0\quad \text{si } i\neq n,\\
\Ext^{n}_{A}(k, A)&\simeq k.
\end{align*}

La profondeur\pageoriginale de $A$ \'etant $n$, pour tout $M$ annul\'e par une puissance de $\mm$, $\Ext^{i}_{A}(M, A)=0$ si $i<n$; de plus $\Ext^{i}_{A}(M, A)=0$ si $i>n$ car la dimension cohomologique globale de $A$ est \'egale \`a $n$. Donc $\Ext^{n}_{A}(\;, A)$ est exact et de plus $\Ext^{n}_{A}(k, A)\simeq k$; il en r\'esulte que:
\end{exemple}

\begin{proposition} \label{IV.5.4}
Si $A$ est un anneau local noeth\'erien r\'egulier de dimension $n$, le foncteur
$$
M\to \Ext_{A}^{n}(M, A)
$$
est dualisant. Le module dualisant associ\'e est
$$
I=\varinjlim_{r}\Ext_{A}^{n}(A/\mm^{r}, A),
$$
il est isomorphe \`a $H^{n}_{\mm}(A)$ (Expos\'e \Ref{II}, th\ptbl \Ref{II.6})\sfootnote{Soient $A$ un anneau, $J$ un id\'eal de $A$, $M$ un $A$-module, $i\in \ZZ$; on posera alors $H^{i}_{J}(M)=H^{i}_{Y}(X, F)$, o\`u $X={\Spec}(A)$, $Y=V(J)$ et $F=\tilde{M}$.}.
\end{proposition}

\begin{remarque} \label{IV.5.5}
Si $A$ v\'erifie \`a la fois les hypoth\`eses des deux exemples pr\'ec\'edents, les deux modules dualisants trouv\'es sont isomorphes. Supposons par exemple que $A$ soit r\'egulier de dimension $n$, complet et d'\'egales caract\'eristiques. Il existe alors un corps des repr\'esentants, soit $K$. Si l'on choisit un syst\`eme de param\`etres $(x_{1}, \ldots, x_{n})$ de $A$, on peut construire un isomorphisme entre $A$ et l'anneau des s\'eries formelles: $K[[T_{1}, \ldots, T_{n}]]$; d'o\`u, comme nous allons le voir, un isomorphisme \emph{explicite} entre les deux modules dualisants
$$
v\colon \sisi{H^{n}}{{H^n_\m}}(A)\to A'.
$$
On peut trouver une interpr\'etation \emph{intrins\`eque} de cet isomorphisme \`a l'aide du module $\sisi{\Omega^{n}(A/K)}{\Omega^n=\Omega^{n}(A/K)}$ des diff\'erentielles relatives compl\'et\'e de degr\'e maximum. En effet, l'on sait que $\sisi{\Omega^{n}(A/K)}{\Omega^n}$ admet une base form\'ee de l'\'el\'ement $dx_{1}\wedge dx_{2}\cdots \wedge dx_{n}$.

D'o\`u un isomorphisme
$$
u\colon \sisi{H^{n}}{{H^n_\m}}(\sisi{\Omega^{n}(A/K)}{\Omega^n})\to \sisi{H^{n}}{{H^n_\m}}(A).
$$
Un fait remarquable est alors que le compos\'e
$$
vu=w\colon \sisi{H^{n}}{{H^n_\m}}(\Omega^{n})\to A'
$$
\emph{ne\pageoriginaled d\'epend plus du choix du syst\`eme de param\`etres} et commute au changement du corps de base.

Pour construire $v$ on calcule $\sisi{H^{n}}{{H^n_\m}}(A)$ gr\^ace au complexe de \sisi{KOSZUL}{Koszul} associ\'e aux $x_{i}$, on trouve:
$$
\sisi{H^{n}}{{H^n_\m}}(A)=\varinjlim_{r}A/(x_{1}^{r}, \ldots, x_{n}^{r});
$$
o\`u les morphismes de transition sont d\'efinis comme suit: posons $I_{r}=A/(x_{1}^{r}, \ldots, x_{n}^{r})$; soit $e_{a_{1}, \ldots, a_{n}}^{r}$ l'image de $\sisi{(x_{1})^{a_{1}}(x_{2})^{a_{2}}\cdots (x_{n})^{a_{n}}}{x_{1}^{a_{1}}x_{2}^{a_{2}}\cdots x_{n}^{a_{n}}}$ dans $I_{r}$. Les $e_{a_{1}, \ldots, a_{n}}^{r}$, pour $\sisi{s}{0}\leq a_{i}<r$ forment une base de $I_{r}$.

Ceci dit, si $s$ est un entier, le morphisme de transition
$$
t_{r, r+s}\colon I_{r}\to I_{r+s}
$$
est la multiplication par $x_{1}^{s}x_{2}^{s}\cdots x_{n}^{s}$, donc:
$$
u_{r, r+s}(e_{a_{1}, \ldots, a_{n}}^{r})=e_{a_{1}+s, \ldots, a_{n}+s}^{r+s}.
$$

Notons que la donn\'ee d'un $A$-homomorphisme $w$ d'un $A$-module $M$ dans $A'$ \'equivaut \`a la donn\'ee d'une forme $K$-lin\'eaire $w':M\to K$ qui soit \emph{continue} sur les sous-modules de type fini. Dans le cas $M=\sisi{H^{n}}{{H^n_\m}}(\Omega^{n})$, la d\'efinition de $w$ \'equivaut donc \`a celle d'une forme lin\'eaire
$$
\rho \colon \sisi{H^{n}}{{H^n_\m}}(\Omega^{n})\to K,
$$
appel\'ee \emph{forme r\'esidu}\refstepcounter{toto}\label{pagecithar}\sfootnote{\label{cithar}Pour une \'etude plus d\'etaill\'ee de la notion de r\'esidu, \Cf \sisi{R.~HARTSHORNE, Lecture Notes in Math. \numero 20, (1966), Springer}{R\ptbl Hartshorne, \emph{Residues and Duality}, Lect. Notes in Math., vol.~20, Springer, 1966.}}. Pour construire $\rho$, il suffit de d\'efinir des formes $\rho_{r}:I_{r}\to K$ qui se recollent, et on prendra
$$
\rho_{r}(e_{a_{1}, \ldots, a_{n}}^{r})=
\begin{cases}
1&\text{si}\ a_{i}=r-1\ \text{pour}\ 1\leq i\leq n\\
0& \text{sinon.}
\end{cases}
$$
\end{remarque}

\chapter{Dualit\'e locale et structure des $\H^i(M)$}\label{V}

\renewcommand{\theequation}{\arabic{equation}}

\section{Complexes d'homomorphismes} \label{V.1}

\subsection{} \label{V.1.1}
Soient\pageoriginale $F^{\boule}$ et~$G^{\boule}$ deux modules gradu\'es, alors on note:
\begin{equation} \label{eq:V.1} \Hom^{\boule}(F^{\boule}, G^{\boule})
\end{equation}
le module gradu\'e des homomorphismes de modules gradu\'es de~$F^{\boule}$ dans~$G^{\boule}$. Ainsi on~a:
\begin{equation} \label{eq:V.2} \Hom^{s}(F^{\boule}, G^{\boule})=\prod_{k}\Hom(F^{k}, G^{k+s}).
\end{equation}
Soit~$F^{\boule}$ (\resp $G^{\boule}$) un complexe, et soit~$d_{1}$ (\resp $d_{2}$) sa diff\'erentielle, alors pour $h\in\Hom^{s}(F^{\boule}, G^{\boule})$ on posera\nde{la convention de signe originale \'etait diff\'erente; mais, elle n'est pas compatible avec la convention de l'expos\'e~\Ref{VIII}, qui para\^it plus raisonnable, puisque dans ce cas la cohomologie en degr\'e~$0$ est alors l'ensemble des classes d'homotopie de morphismes de $F^{\boule}$ dans $G^{\boule}$. Les calculs ont \'et\'e modifi\'es en cons\'equence dans ce qui suit.}
\begin{equation} \label{eq:V.3} d(h)=h\circ d_{1}+(-1)^{\sisi{s}{{s+1}}}d_{2}\circ h.
\end{equation}
On v\'erifie trivialement que $d\circ d=0$, donc que $\Hom^{\boule}(F^{\boule}, G^{\boule})$ muni de~$d$ est un complexe. Le groupe de cohomologie de ce complexe se note
\begin{equation} \label{eq:V.4} \uH^{\boule}(F^{\boule}, G^{\boule}).
\end{equation}
Si~$G^{\boule}$ est injectif en chaque degr\'e, alors
$$F^{\boule}\sisi{\rightsquigarrow}{\mto}\sisi{\underline\H}{\uH}^{\boule}(F^{\boule}, G^{\boule})$$
est un~$\partial$-foncteur exact. De m\^eme, pour $F^{\boule}$ quelconque,
$$G^{\boule}\sisi{\rightsquigarrow}{\mto}\sisi{\underline\H}{\uH}^{\boule}(F^{\boule}, G^{\boule})$$
est un~$\delta$-foncteur exact sur la cat\'egorie des complexes~$G^{\boule}$ injectifs en chaque degr\'e.

\begin{remarque} \label{V.1.2}
Les cycles de~$\Hom^{\boule}(F^{\boule}, G^{\boule})$ sont les homomorphismes de~$F^{\boule}$ dans~$G^{\boule}$ qui commutent ou anticommutent, suivant le degr\'e, avec les diff\'erentielles. Les bords de~$\Hom^{\boule}(F^{\boule}, G^{\boule})$ sont les homomorphismes de~$F^{\boule}$ dans~$G^{\boule}$ homotopes \`a z\'ero.
\end{remarque}

Soit $A$\pageoriginale un anneau et soient~$M$ (\resp $N$) un $A$-module, $R(M)$ (\resp $R(N)$) une r\'esolution injective de $M$ (\resp $N$), alors il existe un isomorphisme canonique\nde{on~a conserv\'e l'\'etrange num\'erotation originale.}
\begin{equation*} \label{eq:V.1.3} \tag{1.3} \uH^{s}(R(M), R(N))\simeq\Ext^{s}(M, N).
\end{equation*}
En effet, soit~$i\colon M\to R(M)$ l'augmentation canonique, et soit $h\in\Hom^{s}(R(M), R(N))$ alors on notera~$t_{s}$ l'application
$$h\sisi{\rightsquigarrow}{\mto} h_{0}\circ i$$
de~$\Hom^{s}(R(M), R(N))$ dans~$\Hom(M, R(N)^s)$. La famille $(\sisi{(-1)^{s}}{}t_{s})_{s\geq 0}$ d\'efinit un homomorphisme de complexes (ordinaires)\nde{on note encore $M$ le complexe $M[0]$ constitu\'e de $M$ plac\'e en degr\'e $0$.}
$$t\colon\Hom^{\boule}(R(M), R(N))\to\Hom\sisi{}{^\boule}(M, R(N)),
$$
\ie on~a $(dh)_{0}\circ i=\sisi{(-1)}{} d_{2}\circ h_{0}\circ i.$

On v\'erifie facilement que, passant \`a la cohomologie, $t$ donne un isomorphisme. En particulier il en r\'esulte que
$$
\uH^{\boule}(R(M), R(N))$$
ne \og d\'epend pas\fg de la r\'esolution injective~$R(M)$ (\resp $R(N)$) de~$M$ (\resp $N$) choisie.

\`A toute suite exacte de $A$-modules
\begin{equation} \label{eq:V.5} 0\to M^{\prime}\to M\to M^{\prime\prime}\to 0
\end{equation}
on fait correspondre une suite exacte de r\'esolutions injectives
\begin{equation} \label{eq:V.6} 0\to R(M^{\prime})\to R(M)\to R(M^{\prime\prime})\to 0.
\end{equation}
On v\'erifie que l'isomorphisme~(\Ref{eq:V.1.3}) commute avec les homomorphismes: \setcounter{equation}{7}
\begin{align} \label{eq:V.8}
\sisi{\H}{\uH}^{s}(R(M^{\prime}), R(N))&\to \sisi{\H}{\uH}^{s+1}(R(M^{\prime\prime}), R(N)),\\
\label{eq:V.9} \Ext^{s}(R(M^{\prime}), R(N))&\to\Ext^{s+1}(R(M^{\prime\prime}), R(N)),
\end{align}
d\'eduits de~(\Ref{eq:V.6}) et~(\Ref{eq:V.5}).

Soient $P$\pageoriginale un troisi\`eme $A$-module, $R(P)$ une r\'esolution injective de~$P$, alors la composition de morphismes gradu\'es donne un accouplement
\begin{equation} \label{eq:V.10} \Hom^{i}(R(N), R(M))\times\Hom^{j}(R(M), R(P))\to\Hom^{i+j}(R(N), R(P)),
\end{equation}
qui d\'efinit un accouplement:
\begin{equation} \label{eq:V.11} \uH^{i}(R(N), R(M))\times\uH^{j}(R(M), R(P))\to\uH^{i+j}(R(N), R(P)),
\end{equation}
donc un homomorphisme de foncteurs en~$M$:
\begin{equation*} \label{eq:V.1.4} \tag{1.4} \uH^{i}(R(N), R(M))\to\Hom(\sisi{\H}{\uH}^{j}(R(M), R(P)), \sisi{\H}{\uH}^{i+j}(R(N), R(P))).
\end{equation*}
On va voir que~\sisi{(5.4.)}{(\Ref{eq:V.1.4})} est un homomorphisme de~$\delta$-foncteurs en~$M$. Les suites exactes~(\Ref{eq:V.5}) et~(\Ref{eq:V.6}) donnent un diagramme commutatif:
$$
\xymatrix{ \Hom^{i}(R(N), R(M^{\prime}))\ar[r]\ar[d] & \Hom(\Hom^{j}(R(M^{\prime}), R(P)), \Hom^{i+j}(R(N), R(P)))\ar[d]^{\Hom(q, \id)}\\
\Hom^{i}(R(N), R(M))\ar[r]\ar[d]^{\sisi{P}{p}} & \Hom(\Hom^{j}(R(M), R(P)), \Hom^{i+j}(R(N), R(P)))\ar[d]\\
\Hom^{i}(R(N), R(M^{\prime\prime}))\ar[r] & \Hom(\Hom^{j}(R(M^{\prime\prime}), R(P)), \Hom^{i+j}(R(N), R(P))).}$$
Soit~$h\in\Hom^{i}(R(N), R(M^{\prime\prime}))$ (resp.\ $g\in\Hom^{j}(R(M^{\prime}), R(P))$) un cycle, et soit $h^{\prime}\in\Hom^{i}(R(N), R(M))$ (resp.\ $g^{\prime}\in\Hom^{j}(R(M), R(P))$) tel que~$p(h^{\prime})=h$ (resp.\ $q(g^{\prime})=g$), alors dire que~(\Ref{eq:V.1.4}) est un homomorphisme de $\delta$-foncteurs en~$M$, c'est dire que
\begin{equation} \label{eq:V.12}
g\circ dh^{\prime}-dg^{\prime}\circ h
\end{equation}
est un cobord dans~$\Hom^{\boule}(R(N), R(P))$.

Or on a:
\begin{align*}
dh^{\prime} &= h^{\prime}\circ d_{1}+(-1)^{\sisi{i}{i+1}}d_{2}\circ h^{\prime}\\
dg^{\prime} &= g^{\prime}\circ d_{2}+(-1)^{\sisi{j}{j+1}}d_{3}\circ g^{\prime}
\end{align*}
avec les notations \'evidentes. Donc~(\Ref{eq:V.12}) s'\'ecrit: \sisi{\begin{align*} g\circ h^{\prime}\circ d_{1} &+ (-1)^{i}g\circ d_{2}\circ h^{\prime}, \\
-g^{\prime}\circ d_{2}\circ h &- (-1)^{j}d_{3}\circ g^{\prime}\circ h.
\end{align*}}{\[
g\circ h^{\prime}\circ d_{1} + (-1)^{i+1}g\circ d_{2}\circ h^{\prime} -g^{\prime}\circ d_{2}\circ h - (-1)^{j+1}d_{3}\circ g^{\prime}\circ h.
\]}

D'autre part\pageoriginale, puisque $h$ et~$g$ sont des cycles on a:
\begin{align*}
g\circ d_{2} &= (-1)^{\sisi{j+1}{j}}d_{3}\circ g, \\
d_{2}\circ h &= (-1)^{\sisi{i+1}{i}}h\circ d_{1},
\end{align*}
donc, finalement,~(\Ref{eq:V.12}) s'\'ecrit:
$$d(g\circ h^{\prime}+(-1)^{\sisi{i}{i+1}}g^{\prime}\circ h),
$$
ce qui termine la d\'emonstration.

(\Ref{eq:V.1.3}) et~(\Ref{eq:V.1.4}) donnent ainsi un homomorphisme de~$\delta$-foncteurs en~$M$:
\begin{equation*} \label{eq:V.1.5} \tag{1.5} \Ext^{i}(N, M)\to\Hom(\Ext^{j}(M, P), \Ext^{i+j}(N, P)).
\end{equation*}

\section{Le th\'eor\`eme de dualit\'e locale pour un anneau local r\'egulier} \label{V.2}

\setcounter{equation}{12}

Soient~$A$ un anneau local r\'egulier de dimension~$r$, $\mm$ l'id\'eal maximal de~$A$ et~$M$ un $A$-module de type fini. \refstepcounter{toto}\label{pV2.4}On pose $\H^{i}(M)=\H^{i}_{\mm}(M)$ (donc $\H^{i}(M)=\varinjlim\Ext^{i}(A/\mm^{k}, M)$). On a vu (\Ref{IV}~\Ref{IV.5.4}) que~$I=\H^{r}(A)$ est un module dualisant pour~$A$, d\'esignons par~$D$ le foncteur dualisant associ\'e. Dans (\Ref{eq:V.1.5}) posons $N=A/\mm^{k}$, $P=A$ alors on obtient un homomorphisme de~$\delta$-foncteurs en~$M$
\begin{equation} \label{eq:V.13} \varphi_{k}\colon\Ext^{i}(A/\mm^{k}, M)\to\Hom(\Ext^{r-i}(M, A), \Ext^{r}(A/\mm^{k}, A)).
\end{equation}
Passant \`a la limite inductive suivant~$k$, on trouve un homomorphisme de~$\delta$-foncteurs
\begin{equation} \label{eq:V.14} \varphi\colon\H^{i}(M)\to D(\Ext^{r-i}(M, A)).
\end{equation}

\begin{theoreme}[\sisi{Th.\ }{Th\'eor\`eme} de dualit\'e locale] \label{V.2.1}
L'homomorphisme fonctoriel~$\varphi$ pr\'ec\'edent est un isomorphisme.
\end{theoreme}

\begin{proof}
Si $i>r$\pageoriginale, le deuxi\`eme membre de~(\Ref{eq:V.14}) est trivialement nul, et le premier membre est nul car $\H^{i}(M)=\varinjlim_{k}\Ext^{i}(A/\mm^{k}, M)$, et c'est vrai pour chaque $\Ext^{i}(A/\mm^{k}, M)$ (Th\'eor\`eme des syzygies).

Si~$i=r$, d'apr\`es ce qui pr\'ec\`ede, les deux foncteurs en~$M$, $\H^{r}(M)$ et~$D(\Hom(M, A))$ sont exacts \`a droite; puisque $A$ est noeth\'erien, et~$M$ de type fini, il suffit de v\'erifier l'isomorphisme pour~$M=A$, ce qui est imm\'ediat.

Pour montrer que~$\varphi$ est un isomorphisme fonctoriel, il suffit maintenant, en proc\'edant par r\'ecurrence descendante par rapport \`a~$i$, de remarquer que tout module de type fini admet une pr\'esentation finie, et que pour~$i<r$ les deux membres de~(\Ref{eq:V.14}) sont z\'ero si~$M$ est libre de type fini. Cela est \'evident pour le deuxi\`eme membre, et puisque $\H^{i}$ commute avec les sommes finies il suffit, quant au premier, de montrer que~$\H^{i}(A)=0$ pour~$i<r$. Or ceci r\'esulte, puisque $\sisi{\profp A}{\prof(A)}=r$, de~(\Ref{III}~\Ref{III.3.4}).
\skipqed
\end{proof}

\section{Application \`a la structure des $\H^i(M)$} \label{V.3}

\setcounter{equation}{14}

\begin{theoreme} \label{V.3.1} Soient~$A$ un anneau local noeth\'erien, $D$ un foncteur dualisant pour~$A$, $M$ un $A$-module de type fini $\neq 0$, et de dimension~$n$, alors on a:
\begin{enumeratei}
\item
$\H^{i}(M)=0$ si~$i<0$ ou si~$i>n$.
\item
$\sisi{D(\H^{i}(M))\;\widehat{\hbox{}}}{\widehat{D(\H^{i}(M))}}$ est un module de type fini sur~$\widehat{A}$, de dimension~$\leq i$.
\item
$\H^{n}(M)\neq 0$, et si~\sisi{$A=\widehat{A}$}{$A$ est complet}, $D(\H^{n}(M))$ est de dimension~$n$ et $\Ass(D(\H^{n}(M)))=\{ \pp\in\Ass(M) \mid \sisi{\dimpt}{\dim} A/\pp=n\}$.
\end{enumeratei}
\end{theoreme}

\begin{proof}
Soit~$I$ le module dualisant associ\'e \`a~$D$. On sait que~$\widehat{I}$ est un module dualisant pour~$\widehat{A}$. D'autre part, on a:
\begin{enumerate}
\item[] $\sisi{\H^{i}(M)\;\widehat{\hbox{}}}{\widehat{H^i(M)}}=\H^{i}(\widehat{M})$,
\item[] $\sisi{D(\H^{i}(M))\;\widehat{\hbox{}}}{\widehat{D(\H^{i}(M))}}=\Hom(\H^{i}(\widehat{M}), \widehat{I})$ et
\item[] $\sisi{\dimpt}{\dim}\widehat{M}=\sisi{\dimpt}{\dim} M$,
\end{enumerate}
\noindent
donc\pageoriginale on peut supposer~$A$ complet. Or, d'apr\`es un th\'eor\`eme de Cohen, tout anneau local complet est quotient d'un anneau local r\'egulier. Pour se ramener \`a ce cas, on~a besoin du lemme suivant:

\begin{lemme} \label{V.3.2}
Soient~$X$ (\resp $Y$) un espace annel\'e, $X^{\prime}$ (\resp $Y^{\prime}$) un sous-espace ferm\'e de~$X$ (\resp $Y$), et $f\colon X\to Y$ un morphisme d'espaces annel\'es tel que $f^{-1}(Y^{\prime})=X^{\prime}$. Soit~$F$ un~$\OX$-Module et d\'esignons par~$A$ (\resp $B$) l'anneau~$\Gamma(\OX)$ (\resp $\Gamma(\Oo_{Y})$) et par~$\overline{f}\colon B\to A$ l'homomorphisme d'anneaux correspondant \`a~$f$. Il existe une suite spectrale de~$B$-modules, de terme initial
\begin{equation} \label{eq:V.15} E_{2}^{p, q}=\H_{Y^{\prime}}^{p}(Y, \R^{q}f_{\ast}(F)),
\end{equation}
aboutissant au~$B$-module $\H_{X^{\prime}}^{\boule}(X, F)_{[\overline{f}]}$.
\end{lemme}

\begin{proof}
Soit~$\Oo_{Y, Y^{\prime}}$ le faisceau $\sisi{{\Oo_{Y}}{|Y^{\prime}}}{{\Oo_{Y}}_{|Y^{\prime}}}$ prolong\'e par~$0$ en dehors de~$Y^{\prime}$ (voir \Exp \Ref{I}). On a un isomorphisme de~$B$-modules:
\begin{equation} \label{eq:V.16} \Hom(\Oo_{Y, Y^{\prime}}, f_{\ast}(F))\simeq\Hom(f^{\ast}(\Oo_{Y, Y^{\prime}}), F)_{[\overline{f}]}.
\end{equation}
Or, on a:
\begin{equation} \label{eq:V.17} f^{\ast}(\Oo_{Y, Y^{\prime}})=\Oo_{X, X^{\prime}},
\end{equation}
et de plus si~$G$ est un~$\OX$-Module injectif, alors
$f_{\ast}(G)$ est un~$\Oo_{Y}$-Module injectif, du moins si $f$
est plat, cas auquel on peut se ramener ais\'ement en rempla\c cant
$\OX$ etc. par les faisceaux d'anneaux constants~$\ZZ$. Donc la
suite spectrale du foncteur compos\'e
$$F\to\Hom(\Oo_{Y, Y^{\prime}}, f_{\ast}(F)),
$$
de terme initial
$$E_{2}^{p, q}=\Ext^{p}(Y;\Oo_{Y, Y^{\prime}}, \R^{q}f_{\ast}(F)),
$$
aboutit, compte tenu de~(\Ref{eq:V.16}) et~(\Ref{eq:V.17}), \`a:
$$
\Ext^{\boule}(X;\Oo_{X, X^{\prime}}, F)_{[\overline{f}]}.
$$

Le lemme r\'esulte alors\pageoriginale de~(\sisi{\Ref{I}~\Ref{I.2.3bof}}{\Ref{I}~\Ref{I.2.3bis}~\textup{bis})}.
\end{proof}

Soit maintenant $f\colon B\to A$ un homomorphisme d'anneaux locaux surjectif. Soit
$$f\colon\Spec(A)\to\Spec(B)$$
le morphisme de sch\'emas affines correspondant. Posons~$X=\Spec(A)$ (resp.\ $X^{\prime}=\{\mm_{A}\}$), $Y=\Spec\sisi{.}{}(B)$ (\resp $Y^{\prime}=\{\mm_{B}\}$), et soient~$M$ un~$A$-module et $\widetilde{M}$ le~$\OX$-Module correspondant. Puisque $\R^{q}f_{\ast}(\widetilde{M})=0$ pour~$q>0$, la suite spectrale~(\Ref{eq:V.15}) d\'eg\'en\`ere, et on obtient d'apr\`es~(\Ref{V.3.2}) un isomorphisme de $B$\sisi{~}{-}modules:
\begin{equation} \label{eq:V.18} \H_{\{\mm_{B}\} }^{n}(Y, f_{\ast}(\widetilde{M}))\simeq\H_{\{\mm_{A}\} }^{n}(X, \widetilde{M})_{[f]},
\end{equation}
donc un isomorphisme de~$B$-modules:
\begin{equation} \label{eq:V.19} \H_{\mm_{B}}^{n}(M_{[f]})\simeq\H_{\mm_{A}}^{n}(M)_{[f]}.
\end{equation}
D'autre part si~$D_{A}$ (\resp $D_{B}$) est le foncteur dualisant pour~$A$ (\resp $B$), on a:
\begin{equation} \label{eq:V.20} D_{A}(M)_{[f]}\simeq D_{B}(M_{[f]}).
\end{equation}
Enfin, puisqu\sisi{e}{'}on a un isomorphisme d'anneaux
\begin{equation} \label{eq:V.21} B/\sisi{\Annp}{\Ann} M_{[f]}\simeq A/\sisi{\Annp}{\Ann} M,
\end{equation}
on voit que le changement d'anneaux de base envisag\'e ne change rien. Supposons donc que~$A$ est r\'egulier de dimension~$r$.

D'apr\`es (\Ref{V.2.1}) on a:
\begin{equation} \label{eq:V.22}
D(\H^{i}(M))=\Ext^{r-i}(M, A).
\end{equation}
On va\pageoriginale d\'emontrer l'\'equivalence entre les propri\'et\'es suivantes:
\begin{enumerate}
\item[(a)] $\sisi{\dimpt}{\dim} \Ext^{j}(M, A)\leq r-j$;
\item[(b)] pour tout~$\pp\in X=\Spec(A)$ tel que $\sisi{\dimpt}{\dim} A_{\pp}<j$, on a~$\Ext^{j}(M, A)_{\pp}=0$;
\item[(c)] $\codim(\sisi{\supp}{\Supp}(\Ext^{j}(M, A)), X)\geq j$.
\end{enumerate}

Pour d\'emontrer (a)\ALORS (b), soit~$\pp\in X$, $\sisi{\dimpt}{\dim} A_\pp<j$, alors $\dim A/\pp>r-j$, donc par (a) $\Ann(\Ext^{j}(M, A))\not\subset\pp$, ce qui entra\^ine $\Ext^{j}(M, A)_{\pp}=0$. Soit~$\pp\in\Supp(\Ext^{j}(M, A))$ alors $\Ext^{j}(M, A)_{\pp}\neq 0$ donc par (b) $\dim A_\pp\geq j$. Donc $\codim(\Supp(\Ext^{j}(M, A)), X)=\inf\{ \dim A_\pp\mid \pp\in\Supp(\Ext^{j}(M, A))\} \geq j$, c'est-\`a-dire (b)\ALORS (c). Enfin (c) implique~(a) trivialement.

D\'emontrons maintenant le th\'eor\`eme.
\begin{enumeratei}
\item
Soit~$x=(x_{1},\dots, x_{r})$ un syst\`eme de param\`etres pour~$A$ tel que~$x_{i}\in\sisi{\Annp}{\Ann} M$ pour~$i=1,\dots, r-n$. Soit~\sisi{$K^{\boule}((\underline{x}^{k}), M)$}{ $K^{\boule}((x^{k}), M)$} le complexe de Koszul. On voit facilement que l'application $K^{i}((x^{k}), M)\to K^{i}((x^{k^{\prime}}), M)$ pour~$k<k^{\prime}$ est z\'ero, si~$i>n$. Il en r\'esulte que $\H^{i}(M)=\varinjlim\sisi{.}{}\H^{i}((x^{k}), M)=0$ si~$i>n$. D'autre part, il est trivial que~$\H^{i}(M)=0$ si $i<0$, donc (i) est d\'emontr\'e.

\item
Puisque~$A$ est r\'egulier, $\dim A_{\pp}<j$ entra\^ine \sisi{$\dim\gl{.}A_{\pp}<j$}{que la dimension homologique globale de $A_\pp$ est strictement plus petite que $j$ et} donc $\Ext^{j}(M, A)_{\pp}=\Ext_{A_{\pp}}^{j}(M_{\pp}, A_{\pp})=\nobreak0$, donc on~a d\'emontr\'e~(b) et par suite~(a). (ii) r\'esulte alors de~(\Ref{eq:V.22}) et de~(a).

\item
Il existe un~$\pp\in\Supp(M)$ tel que $\dim A_{\pp}=r-n$ et tel que $\Supp(M_{\pp})=\{ \mm A_{\pp}\}$. Puisque $A_{\pp}$ est r\'egulier si~$A$ l'est, on trouve $\prof A_{\pp}=r-n$ donc
\begin{equation} \label{eq:V.23} \Ext_{A}^{r-n}(M, A)_{\pp}=\Ext_{A_{\pp}}^{r-n}(M_{\pp}, A_{\pp})\neq 0.
\end{equation}
Ceci implique\pageoriginale, tenant compte \sisi\`a{de}~(\Ref{eq:V.22}), que d'une part:
$$
\H^{n}(M)\neq 0,
$$
d'autre part
$$
\sisi{\dimpt}{\dim} D(\H^{n}(M))\geq n,
$$
donc d'apr\`es (ii)
$$
\sisi{\dimpt}{\dim} D(\H^{n}(M))=n.
$$

Soit maintenant $Y=\Supp(M)$. D'apr\`es~(i) on sait que $D(\H^{n}(M^{\prime}))=\Ext^{r-n}(M^{\prime}, A)$ est un foncteur en $M^{\prime}$, exact \`a gauche, sur la cat\'egorie $(\CC_{Y})^{\circ}$. Donc il existe un $A$-module~$H$ et un isomorphisme de foncteurs en~$M^{\prime}$:
$$
\Ext^{r-n}(M^{\prime}, A)=\Hom(M^{\prime}, H).
$$

Soient~$Y_{i}$, $i=1,\dots, k$ les composantes irr\'eductibles de~$Y$ de dimension maximum. On va voir que l'assertion $\Ext^{r-n}(M^{\prime}, A)\neq 0$ est \'equivalente \`a l'assertion: il existe un~$i$ tel que $\Supp M^{\prime}\supset Y_{i}$. En effet si $\Supp M^{\prime}\supset Y_{i}$ alors \sisi{$\dimpt M^{\prime}=n$}{$\dim(M^{\prime})=n$} donc $\Ext^{r-n}(M^{\prime}, A)\neq 0$.

Si~$\Supp M^{\prime}\not\supset Y_{i}$ pour tout~$i=1,\dots, k$, alors $\sisi{\dimpt}{\dim} M^{\prime}<n$
$$D(\H^{n}(M^{\prime}))=\Ext^{r-n}(M^{\prime}, A)=0.
$$

Puisque $\Ass(\Ext^{r-n}(M, A))=\Supp M\cap\Ass \sisi{H}{(H)}$ on voit que la derni\`ere assertion de~(iii) r\'esulte du lemme suivant:
\end{enumeratei}

\begin{lemme} \label{V.3.3} Soit~$X=\Spec(A)$, et soient $Y$ une partie ferm\'e de~$X$, $T\colon (\CC_{Y})^{\circ}\to\Ab$ un foncteur exact \`a gauche, et~$Y_{i}$, $i=1,\dots, k$ une famille de composantes irr\'eductibles de~$Y$ tels que l'assertion: $T(M)=0$ est \'equivalente \`a l'assertion: $\forall i\;\Supp M\not\supset Y_{i}$\sisi{,}{.} Alors $T$ est repr\'esentable par un module~$H$ tel que $\Ass\sisi{.H}{(H)}=\bigcup_{i=1}^{k}\{y_{i}\}$, o\`u~$y_{i}$ est le point g\'en\'erique de~$Y_{i}$,~$i=1,\dots, k$.
\end{lemme}

\begin{proof}
Soit $y\in Y$,\pageoriginaled on fabrique un~$A$-module~$M(y)$ tel que $\Supp(M(y))=\overline{\{ y\} }$. Supposons que $y\neq y_{i}$ pour tout~$i=1,\dots, k$, alors $Y_{i}\not\subset\Supp(M(y))$ pour tout~$i=1,\dots, k$, donc $T(M\sisi{)}{}(y))=0$. Il en r\'esulte:
$$
\Ass(T(M(y)))=\Supp(M(y))\cap\Ass \sisi{H}{(H)}=\emptyset,
$$
donc $y\not\in\Ass \sisi{H}{(H)}$. Si~$y=y_{i}$, alors $Y_{i}\subset\Supp(M(y))$, donc $T(M(y))\neq 0$, donc:
$$
\Ass(T(M(y)))=\Supp(M\sisi{}{(y)})\cap\Ass \sisi{H}{(H)}\neq\emptyset.
$$

D'apr\`es la premi\`ere partie de la d\'emonstration, ceci implique~$y\in\Ass \sisi{H}{(H)}$, d'o\`u le lemme, \hfill Q.E.D.
\skipqed
\end{proof}
\skipqed
\end{proof}

\begin{exemple} \label{V.3.4}
Soit~$A$ un anneau noeth\'erien, soient $X=\Spec(A)$, et~$Y$ une partie ferm\'ee de~$X$, tel que $X-Y$ soit affine, alors pour toute composante irr\'eductible $Y_{\alpha}$ de~$Y$ on~a $\codim\sisi{.}{}(Y_{\alpha}, X)\leq 1$.

\sisi{En effet, consid\'erons $X$ comme pr\'esch\'ema au-dessus de~$X$. Soit $y_{\alpha}\in Y_{\alpha}$ un point g\'en\'erique, et consid\'erons le morphisme $\Spec(\Oo_{X, Y_{\alpha}})\to X$. Le sch\'ema affine obtenu par extension du sch\'ema de base de~$X$ \`a~$\Spec(\Oo_{X, Y_{\alpha}})$ est canoniquement isomorphe \`a~$\Spec(\Oo_{X, Y_{\alpha}})$.

D'apr\`es~(\EGA I 3.2.7) on voit que si $y_{0}$ est l'unique point ferm\'e de~$Y_{0}=\Spec(\Oo_{X, Y_{\alpha}})$ alors $Y_{0}- y_{0}$ est affine. Par~(\EGA III 1.3.1) on trouve:
$$
\H^{i}(Y_{0}- y_{0}, \Oo_{Y_{0}})=0\quad\textrm{ si }i>0,
$$
donc par~(\Ref{I}\sisi{.}{~}\Ref{I.2.11})
$$
\H^{i}(\Oo_{X, Y_{\alpha}})=\H_{\{ y_{0}\} }^{i}(Y_{0}, \Oo_{Y_{0}})=0\quad\textrm{ si }i\geq 2.
$$

Tenant compte de~(5.7.) il vient:
$$
\sisi{\dimpt}{\dim} \Oo_{X, y_{\alpha}}\leq 1,
$$
donc $\sisi{\codimpt}{\codim} (Y_{\alpha}, X)=\underset{y\in Y_{\alpha}}{\inf}\sisi{\dimpt}{\dim} \Oo_{X, y}\leq 1$,}{En effet, consid\'erons $X$ comme pr\'esch\'ema au-dessus de~$X$. Soit $y_{\alpha}\in Y_{\alpha}$ un point g\'en\'erique, et consid\'erons le morphisme $\Spec(\Oo_{X, y_{\alpha}})\to X$. Le sch\'ema affine obtenu par extension du sch\'ema de base de~$X$ \`a~$\Spec(\Oo_{X, y_{\alpha}})$ est canoniquement isomorphe \`a~$\Spec(\Oo_{X, y_{\alpha}})$.

D'apr\`es~(\EGA I 3.2.7) on voit que si $y_{0}$ est l'unique point ferm\'e de~$Y_{0}=\Spec(\Oo_{X, y_{\alpha}})$ alors $Y_{0}- y_{0}$ est affine. Par~(\EGA III 1.3.1) on trouve:
$$
\H^{i}(Y_{0}- y_{0}, \Oo_{Y_{0}})=0\quad\textrm{ si }i>0,
$$
donc par~(\Ref{I}\sisi{.}{~}\Ref{I.2.9})
$$
\H^{i-1}(\Oo_{X, y_{\alpha}})=\H_{\{ y_{0}\} }^{i}(Y_{0}, \Oo_{Y_{0}})=0\quad\textrm{ si }i\geq 2.
$$

Tenant compte de~\Ref{V.3.1} (iii), il vient:
$$
\dim\Oo_{X, y_{\alpha}}\leq 1,
$$
donc $\sisi{\codimpt}{\codim} (Y_{\alpha}, X)=\underset{y\in Y_{\alpha}}{\inf}\sisi{\dimpt}{\dim} \Oo_{X, y}\leq 1$,}\hfill Q.E.D.
\end{exemple}

Soient $A$\pageoriginale un anneau local noeth\'erien, $\mm$ l'id\'eal maximal et~$M$ un~$A$-module de type fini. Supposons que $A$ est quotient d'un anneau local r\'egulier. Posons $X=\Spec(A)$, et pour tout $x\in X$, $\mm_{x}=\mm A_{x}$.

\begin{proposition} \label{V.3.5}
Les deux conditions suivantes sont \'equivalentes:
\begin{enumeratea}
\item
$\H^{i}(M)$ est de longueur finie,
\item
$\forall x\in X-\{ \mm\}, \, \H_{\mm_{x}}^{i-\dim\overline{\{ x\} }}(M_{x})=0$.
\end{enumeratea}
\end{proposition}

\begin{proof}
Compte tenu de~\sisi{\Ref{V.3.3})}{(\Ref{V.3.2})} nous pouvons supposer $A$ r\'egulier. D'apr\`es \sisi{la formule (\Ref{eq:V.1.3})}{(\Ref{V.2.1})} nous avons:
$$
\H^{i}(M)=D(\Ext^{r-i}(M, A)),
$$
o\`u $r=\dim\sisi{.A}{A}$. D'apr\`es (\Ref{IV}~\Ref{IV.4.7}), a) est
\'equivalent\refstepcounter{toto}\nde{\label{dualnul}en effet, il
s'agit de montrer que, $E$ \'etant est un $A$-module de type fini,
si $D(E)$ est de longueur finie alors $E$ est de longueur finie.
Soit $K$ (\resp $Q$) le noyau (\resp conoyau) du morphisme
canonique $\epsilon: E\to DD(E)$. Le compos\'e de $D(\epsilon)$ et
du morphisme canonique $\gamma:~D(E)\lto{\gamma} DDD(E)$ est
l'identit\'e de $D(E)$. Comme $D(E)$ est de longueur finie, $\gamma$
est un isomorphisme et donc $D(\epsilon)$ \'egalement. Comme $D$ est
exact, on en d\'eduit que $D(K)$ et $D(Q)$ sont nuls. Il suffit de
prouver que si $M$ est un $A$-module de dual nul, alors $M$ est
nul, car on aura alors $E=DD(E)$ de longueur finie comme $D(E)$.
En effet, soit $M_0$ un sous-module de type fini de $M$. Comme $D$
est exact, $D(M_0)$ est un quotient de $D(M)$, qui est donc nul.
Toujours par exactitude de $D$, on~a $D(M_0/\m_AM_0)=0$, et donc,
par bidualit\'e, le module de longueur finie $M_0/\m_AM_0$ est nul.
Le lemme de Nakayama assure alors la nullit\'e de $M_0$ et
finalement on obtient celle de $M$.} \`a:
\begin{equation} \label{eq:V.24} \Ext^{\sisi{n}{r}-i}(M, A)\textrm{ est de longeur finie.}
\end{equation}
Or (\Ref{eq:V.24}) est \'equivalent \`a:
\begin{equation} \label{eq:V.25} \forall x\in X-\{ \mm\}, \textrm{ on~a }\Ext^{r-i}(M, A)_{x}=0.
\end{equation}
D'autre part $A_{x}$ est r\'egulier de dimension~$r-\dim\overline{\{ x\} }$, donc d'apr\`es~(\Ref{V.2.1})
$$%\label{eq:V.26}
\H_{\mm_{x}}^{i-\dim\overline{\{ x\} }}(M_{x})\!=\!D(\Ext_{A_{x}}^{(r-\dim\overline{\{ x\} })-(i-\dim\overline{\{ x\} })}(M_{x}, A_{x}))\!=\!D(\Ext_{A_{x}}^{r-i}(M_{x}, A_{x})).\leqno(26)
$$
Puisque $M$ est de type fini on a:
$$
\Ext_{A}^{r-i}(M, A)_{x}=\Ext_{A_{x}}^{r-i}(M_{x}, A_{x})$$
d'o\`u la proposition.
\skipqed
\end{proof}

\begin{corollaire} \label{V.3.6} Pour que $\H^{i}(M)$ soit de longueur finie pour~$i\leq n$, il faut et il suffit que
$$
\sisi{\profp M_{x}}{\prof(M_x)}>n-\dim\overline{\{ x\} }$$
pour tout $x\in X-\{ \mm\}$.
\end{corollaire}

\begin{proof}
R\'esulte de~(\Ref{V.3.5}) et de~(\Ref{III}~\Ref{III.3.1}).
\skipqed
\end{proof}

\chapterspace{-4}
\chapter{Les foncteurs $\protect\Ext^{\protect\boule}_Z(X;F, G)$ et $\protect\SheafExt^{\protect\boule}_Z(F, G)$} \label{VI}

\numberwithin{equation}{subsection}

\section{G\'en\'eralit\'es} \label{VI.1}
\enlargethispage{\baselineskip}%
\subsection{} \label{VI.1.1}
Soient\pageoriginale $(X, \OX)$ un espace annel\'e et $Z$ une partie localement ferm\'ee de $X$. Soient $F$ et $G$ des $\OX$-Modules, on d\'esignera par ${\Ext_Z^i}(X;F, G)$ (\resp ${\SheafExt_Z^i}(F, G)$) le $i$-i\`eme foncteur d\'eriv\'e du foncteur $G \sisi{\rightsquigarrow}{\mto} \Gamma_Z(\SheafHom_{\OX}(F, G))$ (\resp $\SheafGamma_Z(\SheafHom_{\OX}(F, G))$).

\begin{lemme} \label{VI.1.2}
Le faisceau ${\SheafExt_Z^i}(F, G)$ est canoniquement isomorphe au faisceau associ\'e au pr\'efaisceau $U\sisi{\rightsquigarrow}{\mto} {\Ext_{Z\cap U}^i}(U;F_{|U}, G_{|U})$.
\end{lemme}

Cela r\'esulte de (\sisi{T}{\cite{Tohoku}}, 3.7.2) et de ce que $\Gamma(U;\SheafGamma_Z(\SheafHom_{\OX}(F, G)))$ est canoniquement isomorphe \`a $\Gamma_{Z\cap U}(\SheafHom_{{\OX}_{|U}}(F_{|U}, G_{|U}))$

\begin{theoreme}[Th\'eor\`eme d'excision] \label{VI.1.3}
Soit $V$ un ouvert de $X$ contenant $Z$, on~a alors un isomorphisme de foncteurs cohomologiques
\begin{equation} \label{VI.1.3.1}
\Ext_{X}^{\boule}(X;F, G) \simeq \Ext_{V}^{\boule}(V;F_{|V}, G_{|V})
\end{equation}
\end{theoreme}

En effet, si $G^\boule$ est une r\'esolution injective de $G$, alors $G^\boule_{|V}$ est une r\'esolution injective de $G_{|V}$. Le th\'eor\`eme en r\'esulte imm\'ediatement.

\subsection{} \label{VI.1.4}
Soit $\Oo_{X, Z}$ le $\OX$-Module d\'efini par les conditions suivantes (\sisi{G}{\cite{Godement}},~2.9.2): ${\Oo_{X, Z}}_{|X{-} Z}=0$ et ${\Oo_{X, Z}}_{|Z}={\OX}_{|Z}$. On a vu que pour tout $\OX$-Module $H$, il existe un isomorphisme fonctoriel: $\Gamma_Z(H)\simeq \Hom_{\OX}(\Oo_{X, Z}, H)$. On en d\'eduit donc des isomorphismes fonctoriels en $F$ et $G$: \setcounter{equation}{0}
\begin{align}
\label{VI.1.4.1}
\Gamma_Z(\SheafHom_{\OX}(F, G)) &\simeq \Hom_{\OX}(\Oo_{X, Z}, \SheafHom_{\OX}(F, G)),\\
\label{VI.1.4.2} \Gamma_Z(\SheafHom_{\OX}(F, G)) &\simeq \Hom_{\OX}(\Oo_{X, Z}\otimes_{\OX} F, G),\\
\label{VI.1.4.3} \Gamma_Z(\SheafHom_{\OX}(F, G)) &\simeq \Hom_{\OX}(F, \SheafHom_{\OX}(\Oo_{X, Z}, G))=\Hom_{\OX}(F, \SheafGamma_Z(G)).
\end{align}
Il r\'esulte\pageoriginale en particulier de~(\Ref{VI.1.4.2}) qu'il existe un isomorphisme $\partial$\sisi{~}{-}fonctoriel en $F$ et~$G$ \sisi{entre $\Ext_{Z}^i(X;F, G)$ et $\Ext_{{\OX}}^i(\Oo_{X, Z}\otimes_{\OX} F, G)$}{$$\sisi{}{\theta\colon } \Ext_{{\OX}}^i(\Oo_{X, Z}\otimes_{\OX} F, G)\isomto\Ext_{Z}^i(X;F, G).$$}

\subsection{} \label{VI.1.5}
Par d\'efinition le foncteur $G \sisi{\rightsquigarrow}{\mto} \Gamma_Z(\SheafHom_{\OX}(F, G))$ est le compos\'e du foncteur $G \sisi{\rightsquigarrow}{\mto}\SheafHom_{\OX}(F, G)$ et du foncteur $\Gamma_Z$. Comme $\Gamma_Z$ est exact \`a gauche (\Ref{I}\sisi{,}{}~\Ref{I.1.9}) et comme $\SheafHom_{\OX}(F, G)$ est flasque si $G$ est injectif, et que $\Gamma_Z$ est exact sur les flasques (\Ref{I}\sisi{,}{}~\Ref{I.2.12}), il r\'esulte de (\sisi{T}{\cite{Tohoku}}, 2.4.1) qu'il existe un foncteur spectral aboutissant \`a $\Ext_{Z}^\boule(X;F, G)$ et dont le terme initial est $H_Z^p(X, \SheafExt_{\OX}^q\sisi{{}_{(F, G))}}{(F, G))}$.

D'autre part, il r\'esulte de (\Ref{VI.1.4.3}) que $\Gamma_Z(\SheafHom_{\OX}(F, G))$ est le compos\'e de $\SheafGamma_Z$ et du foncteur $H \sisi{\rightsquigarrow}{\mto} \Hom_{\OX}(F, H)$.

Puisque le foncteur $\SheafGamma_Z$ transforme les injectifs en injectifs (\Ref{I}\sisi{,}{}~\Ref{I.1.4}), il r\'esulte de (\sisi{T}{\cite{Tohoku}},~2.4.1) qu'il existe un foncteur spectral aboutissant \`a $\Ext_{Z}^\boule(X;F, G)$ et dont le terme initial est $\Ext_{\OX}^p(X;F, \SheafH_Z^q(G))$.

Il r\'esulte enfin, de (\Ref{VI.1.4.2}) et de la suite spectrale des $\Ext$, qu'il existe un foncteur spectral aboutissant \`a $\Ext^{\boule}_Z(X;F, G)$ et dont le terme initial est $H^p(X;\SheafExt_Z^\boule(F, G))$. D'o\`u le

\setcounter{equation}{0}

\begin{theoreme} \label{VI.1.6}
Il existe trois foncteurs spectraux aboutissant \`a $\Ext_{Z}^\boule(X;F, G)$ et dont les termes initiaux sont respectivement
\begin{align}
\label{VI.1.6.1} &H_Z^p(X, \SheafExt_{\OX}^q(F, G))\\
\label{VI.1.6.2} &H^p(X, \SheafExt_Z^q(F, G))\\
\label{VI.1.6.3} &\Ext_{\OX}^p(X;F, \SheafH_Z^q(G)).
\end{align}
\end{theoreme}

\subsection{} \label{VI.1.7}
Soit maintenant\pageoriginale $Z'$ une partie ferm\'ee de $Z$ et soit $Z''=Z {-} Z'$. On a une suite exacte
\setcounter{equation}{0}
\begin{equation} \label{VI.1.7.1}
0 \to \Oo_{X, Z''} \to \Oo_{X, Z} \to \Oo_{X, Z'} \to 0
\end{equation}
qui g\'en\'eralise la suite exacte de (\sisi{G}{\cite{Godement}}, 2.9.3). Cette suite exacte splitte localement, on~a donc, pour tout $\OX$-Module $F$, une autre suite exacte:
\begin{equation} \label{VI.1.7.2} {0 \to F \otimes_{\OX} \Oo_{X, Z''} \to F \otimes_{\OX} \Oo_{X, Z} \to F \otimes_{\OX} \Oo_{X, Z'} \to 0}
\end{equation}
Soit maintenant $G$ un $\OX$-Module; si on applique le foncteur $\Hom_{\OX}(\boule, G)$ \`a la suite exacte (\Ref{VI.1.7.2}), on d\'eduit de (\Ref{VI.1.4.2}) et de la suite exacte des $\Ext$ le th\'eor\`eme suivant:

\begin{theoreme} \label{VI.1.8}
Soient $Z$ une partie localement ferm\'ee de $X$, $Z'$ une partie ferm\'ee de $Z$ et $Z''=Z{-} Z'$. On a alors une suite exacte fonctorielle en $F$ et $G$:
\begin{gather*}
0 \to \Hom_{Z'}(F, G) \to \Hom_Z(F, G) \to \Hom_{Z''}(F, G)\to \Ext_{{Z'}}^1(F, G) \to \cdots{}\\
{}\cdots \Ext_{Z}^i(F, G) \to \Ext_{{Z''}}^i(F, G) \to \Ext_{{Z'}}^{i+1}(F, G) \to \cdots{}
\end{gather*}
\end{theoreme}

\begin{corollaire} \label{VI.1.9} Soit $Y$ une partie ferm\'ee de $X$ et soit $U=X{-} Y$. On a alors une suite exacte fonctorielle en $F$ et $G$:
\begin{gather*}
0 \to \Hom_Y(F, G) \to \Hom_{\OX}(F, G) \to \Hom_{{\OX}_{|U}}(F_{|U}, G_{|U}) \to \Ext_{Y}^1(F, G) \to \cdots{}\\
{}\cdots \Ext_{{\OX}}^i(F, G) \to \Ext_{{{\OX}_{|U}}}^i(F_{|U}, G_{|U}) \to \Ext_{Y}^{i+1}(F, G) \to \cdots
\end{gather*}
\end{corollaire}

Ce corollaire est une cons\'equence imm\'ediate du th\'eor\`eme (\Ref{VI.1.3}) et du th\'eor\`eme (\Ref{VI.1.8}).

\section{Applications aux faisceaux quasi-coh\'erents sur les pr\'esch\'emas} \label{VI.2}
\pageoriginale
\begin{proposition} \label{VI.2.1}
Soit $X$ un pr\'esch\'ema localement noeth\'erien; pour toute partie localement ferm\'ee $Z$ de $X$, pour tout Module coh\'erent $F$ et tout Module quasi-coh\'erent~$G$ sur $X$, les $\SheafExt_Z^i(F, G)$ sont quasi-coh\'erents.
\end{proposition}

On montre, comme pour (\Ref{VI.1.6.3}), que les Modules $\SheafExt_Z^i(F, G)$ sont l'aboutissement d'une suite spectrale de terme initial $\SheafExt_{\OX}^p(F, \SheafH_Z^q (G))$. D'apr\`es (\Ref{II}, corol\ptbl \Ref{II.3}), les $\SheafH_Z^q(G)$ sont quasi-coh\'erents, et donc aussi les $\SheafExt_{\OX}^p (F, \SheafH_Z^q (G))$, puisque $F$ est coh\'erent. La proposition en d\'ecoule alors imm\'ediatement.

\begin{empty} \label{VI.2.2}
\end{empty}

\subsection{} Soient maintenant $Y$ un sous-pr\'esch\'ema ferm\'e de $X$ et $\SheafI$ un id\'eal de d\'efinition de $Y$. Soient $m$ et $n$ des entiers tels que $m\geq n \geq 0$; on d\'esigne par $i_{n, m}$ l'application canonique $\Oo_{Y_m} = \OX/\SheafI^{m+1} \to \OX/\SheafI^{n+1} = \Oo_{Y_n}$ et par $j_n$ l'application \hbox{$\Oo_{X, Y}\to \Oo_{Y_n}$}. Les $(\Oo_{Y_n}, i_{n, m})$ forment un syst\`eme projectif et les $j_n$ sont compatibles avec les~$i_{n, m}$.

En appliquant le foncteur $\Ext_{\OX}^i(F\otimes\boule, G)$, on en d\'eduit un morphisme
$$
\varphi'\colon \varinjlim_{n} \Ext_{\OX}^i(X;F\otimes \Oo_{Y_n}, G) \to \Ext_{\OX}^i(X;F \otimes \Oo_{X, Y}, G);
$$
c'est un morphisme de foncteurs cohomologiques en $G$. Le morphisme
$$
\varphi\colon \varinjlim_{n} \Ext_{\OX}^i(X;F\otimes \Oo_{Y_n}, G) \to \Ext_Y^i(X;F, G)
$$
compos\'e de $\varphi'$ et de $\theta$ (\cf \Ref{VI.1.4}) est donc lui aussi un foncteur de morphismes cohomologiques en $G$.

On d\'efinit de m\^eme
$$
\Sheafphi\colon \varinjlim_n \SheafExt_{\OX}^i(F \otimes \Oo_{Y_n}, G) \to \SheafExt_Y^i(X;F, G)
$$

\begin{theoreme} \label{VI.2.3}
Soient\pageoriginale $X$ un pr\'esch\'ema localement noeth\'erien, $Y$ une partie ferm\'ee de $X$ d\'efinie par un id\'eal coh\'erent $\SheafI$, $F$ un Module coh\'erent, $G$ un Module quasi-coh\'erent. Alors,
\begin{enumeratea}
\item
$\Sheafphi$ est un isomorphisme.
\item
Si $X$ est noeth\'erien, $\varphi$ est un isomorphisme.
\end{enumeratea}
\end{theoreme}

La d\'emonstration de b) \'etant presque mot \`a mot celle de (\Ref{II}\sisi{,}{}~\Ref{II.6}\sisi{,~}{}b)), gr\^ace \`a la suite spectrale~\Ref{VI.1.6.2}, nous ne la reproduirons pas.

Pour la d\'emonstration de a), on peut, d'apr\`es (\Ref{VI.2.1}), supposer $X$ affine d'anneau~$A$, $F$ (\resp $G$) d\'efini par un $A$-module $M$ (\resp $N$) et $\SheafI$ par un id\'eal $\sisi{\underline{\rm{i}}}{I}$. Il suffit de prouver que l'homomorphisme
\begin{equation} \label{VI.2.3.1} \varinjlim_{{n}} \Ext_A^i(M/\sisi{\underline{\rm{i}}}{I}^nM, N) \to \Ext_Y^i(X, F, G)
\end{equation}
d\'eduit de $\Sheafphi$ est un isomorphisme

En effet, on peut, pour $i=0$, identifier canoniquement les deux membres de (\Ref{VI.2.3.1}) au sous-module de $\Hom_A(M, N)$ d\'efini par les \'el\'ements de $\Hom_A(M, N)$ annul\'es par une puissance de $\sisi{\underline{\rm{i}}}{I}$. On voit alors que l'homomorphisme (\Ref{VI.2.3.1}) n'est autre que l'application identique.

Le foncteur $N\sisi{\rightsquigarrow}{\mto} \varinjlim_{n} \Ext_A^\boule(M/\sisi{\underline{\rm{i}}}{I}^nM, N)$ est un $\partial$-foncteur universel. On va montrer qu'il en est de m\^eme du foncteur $N \sisi{\rightsquigarrow}{\mto} \Ext_Y^\boule(M, N)$. En effet si $N$ est un module injectif, d'apr\`es (\Ref{II.9} et \Ref{II.11}), $\SheafH_Y^q(N)=0$ si $q \neq 0$; et d'apr\`es (\sisi{IV.2.3}{\Ref{IV}.\Ref{IV.2.2}}), $H_Y^{\sisi{\circ}{0}}(N)$ est injectif.

Il en r\'esulte alors que $\Ext_{\OX}^p(X;M, \SheafH_Y^q(N))=0$ si $p+q \neq 0$; on~a donc, d'apr\`es (\Ref{VI.1.6.3}), $\Ext_Y^i(M, N)=0$ pour $i\neq 0$ et $N$ injectif. Ce qui ach\`eve la d\'emonstration.

\section*{Bibliographie} M\^emes r\'ef\'erences que celles list\'ees \`a la fin de l'\Exp \Ref{I}, cit\'ees respectivement \sisi{T}{\cite{Tohoku}} et \sisi{G}{\cite{Godement}}.

\chapter[Crit\`eres de nullit\'e, conditions de coh\'erence]
{Crit\`eres de nullit\'e, conditions de coh\'erence des faisceaux $\SheafExt^{i}_{Y}(F, G)$}\label{VII}

\renewcommand{\theequation}{\thesection.\arabic{equation}}

\section{\'Etude pour $i < n$}

D\'emontrons\pageoriginale un lemme:

\begin{lemme} \label{VII.1.1}
Soient $X$ un pr\'esch\'ema localement noeth\'erien, $Y$ une partie ferm\'ee de~$X$ et $G$ un $\OX$-Module quasi-coh\'erent. Supposons que pour tout $\OX$-Module coh\'erent~$F$ de support contenu dans $Y$, on ait:
\[
\SheafExt^{n-1}(F, G) = 0.
\]
Alors pour tout $\OX$-Module coh\'erent $F$ et tout ferm\'e $Z$ de $X$ tels que $Y \supset \Supp F \cap Z$, on a
\[
\SheafExt^n_Z(F, G) \approx \SheafHom (F, \HH^n_Y(G)).
\]
\end{lemme}

On remarque d'abord que
\[
\SheafExt^i_Z(F, G) = \SheafExt^i_{Z \cap \Supp F}(F, G).
\]
(trivial, \cf Expos\'e \Ref{VI}). On fait d'abord la d\'emonstration pour $Z \!=\! X$ donc \hbox{$\Supp F \!\subset\! Y$}. Le foncteur
\[F \sisi{\rightsquigarrow}{\mto} \SheafExt^n(F, G),
\]
d\'efini sur la cat\'egorie des $\OX$-Modules coh\'erents de support contenu dans $Y$, est exact \`a gauche. En vertu de (\Ref{IV}\sisi{~\Ref{IV.1.2}}{~\Ref{IV.1.3}}), il est repr\'esentable par
\[I = \varinjlim_{k} \SheafExt^n (\OX/\mathcal{I}^{k+1}, G),
\]
o\`u $\mathcal{I}$ est l'id\'eal de d\'efinition de Y. Or, d'apr\`es (\Ref{II}~\Ref{II.6}), on sait que:
\[
\HH^n_Y(G) \approx \varinjlim_{k} \SheafExt^n (\OX/\mathcal{I}^{k+1}, G).
\]
D'o\`u\pageoriginale la conclusion si $ Z= X$. Toujours d'apr\`es (\Ref{VI}~\Ref{VI.2.3}) on sait que:
\[
\SheafExt^n_Z (F, G) \approx \varinjlim_{k} \SheafExt^n (F/\Jj^{k + 1}F, G),
\]
o\`u $\Jj$ est l'id\'eal de d\'efinition de Z. Le support de $F/\Jj^{k+1}F$ est contenu dans $Y$ si $Z \cap \Supp F \subset Y$; d'apr\`es ce que nous venons de d\'emontrer, on~a donc:
\[
\SheafExt^n_Z (F, G) \approx \varinjlim_{k} \SheafHom (F/\Jj^{k + 1}F, \HH^n_Y(G)).
\]
Il reste \`a faire voir que l'homomorphisme naturel:
\[
\varinjlim_{k} \SheafHom (F/\Jj^{k + 1}F, \HH^n_Y(G)) \to \SheafHom (F, \HH^n_Y (G)),
\]
est un isomorphisme lorsque $Z \cap \Supp F \subset Y$. Or $X$ peut \^etre recouvert par des ouverts affines noeth\'eriens; on est ainsi ramen\'e au cas o\`u $X$ est affine noeth\'erien; alors $F(X)$ est un $\OX(X)$-Module de type fini et $\Supp \HH^n_Y(G) \subset Y$. Donc tout homomorphisme \hbox{$u: F(X) \to \HH^n_Y (G) (X)$} est annul\'e par une puissance de $\Ii$ donc par une puissance de~$\Jj$, \sisi{}{\hfill}C.Q.F.D.

\begin{proposition} \label{VII.1.2}
Soient $X$ un pr\'esch\'ema localement noeth\'erien, $Y$ une partie ferm\'ee de $X$, $G$ un $\OX$-Module quasi-coh\'erent et $n$ un entier. Quelles que soient $Z$ et $S$, parties ferm\'ees de $X$ telles que $Z \cap S = Y$, les conditions suivantes sont \'equivalentes:
\begin{enumeratei}
\item
$\HH^i_Y (G) = 0$ si $i<n$;

\item
il existe un $\OX$-Module coh\'erent $F$, de support $S$, tel que:
\[
\SheafExt^i_Z (F, G) = 0 \text{ si } i<n;\]

\item
pour tout $\OX$-Module $F$ coh\'erent de support contenu dans $S$ (\ie $\Supp F \cap Z$ $=$ $\Supp F \cap Y)$, on a:
\[
\SheafExt^i_{Z}(F, G) = 0 \text{ si } i < n;
\]
\item
pour tout\pageoriginale $\OX$-Module coh\'erent $F$, on a:
\[
\SheafExt^i_Y (F, G) = 0 \text{ si } i < n.
\]
\end{enumeratei}
\noindent
De plus, si elles sont v\'erifi\'ees, alors pour tout $\OX$-Module coh\'erent F et toute partie ferm\'ee $Z'$ de $X$ telle que $Z' \cap \Supp F = Y \cap \Supp F$, on~a des isomorphismes:
\[
\SheafExt^n_{Z} (F, G) \approx \SheafExt^n_Y (F, G) \approx \SheafHom (F, \HH^n_Y (G)).
\]
\end{proposition}

\begin{proof}
Raisonnons par r\'ecurrence. La proposition est triviale pour $n < 0$. Supposons la d\'emontr\'ee pour $n < q$. Si l'une des conditions est v\'erifi\'ee pour $n = q$, et pour deux parties $Z$ et $S$ comme dit, par l'hypoth\`ese de r\'ecurrence on a, pour tout ferm\'e $Z'$ de $X$ et tout $\OX$-Module coh\'erent $F$ tels que $Z' \cap \Supp F = Y \cap \Supp F$, des isomorphismes:
\begin{equation}\label{eq:VII.1.1}\SheafExt^{q-1}_{Z'} (F, G) \approx \SheafHom (F, \HH^{q-1}_Y (G)) \approx \SheafExt^{q-1}_Y (F, G).\end{equation}

\noindent Donc:

(i) \ALORS (iv), car on prend $Z' = Y$ dans (1.1);

(iv) \ALORS (iii), car on prend $Z' = Z$ dans (1.1);

(iii) \ALORS (ii), car on prend $F = \Oo_S$;

(ii) \ALORS (i), car on prend $Z' = Z$ dans (1.1); d'o\`u $\SheafHom (F, \HH_Y^{q-1} (G)) = 0$; on remarque alors que:
\[
\Supp \HH_Y^{q-1} (G) \subset Y = Z \cap S \subset S = \Supp F,
\]
et on applique le lemme suivant:

\begin{lemme} \label{VII.1.3}
Soit $X$ un pr\'esch\'ema, soit $P$ un $\OX$-Module coh\'erent, et soit $H$ un $\OX$-Module quasi-coh\'erents tels que:
\[
\SheafHom (P, H) = 0 \text{ et } \Supp P \supset \Supp H.
\]
Alors $H = 0$.
\end{lemme}

Il suffit\pageoriginaled de d\'emontrer le lemme quand $X$ est affine, car les ouverts affines forment une base de la topologie de $X$ et les hypoth\`eses se conservent par restriction \`a un ouvert. Or, dans ce cas, on est ramen\'e \`a un probl\`eme sur les $A$-modules, o\`u $X = \Spec(A)$. On applique la formule (valable sous la seule hypoth\`ese que $M$ est de type fini):
\[
\Ass \Hom_A (P, H) = \Supp P \cap \Ass H;\]
On sait que $\Ass H \subset \Supp H \subset \Supp P$ et que $\Ass \Hom_A (P, H) = \emptyset$; donc $\Ass H = \emptyset$, donc $H = 0$.

Pour terminer la d\'emonstration de la proposition, il reste \`a remarquer que (iv) permet d'appliquer~\sisi{(\Ref{VII.1.1})}{\Ref{VII.1.1}}.
\skipqed
\end{proof}

\begin{corollaire}\label{VII.1.4}
Soit $G$ un $\OX$-Module coh\'erent et de Cohen-Macaulay, soit $n \in \ZZ$. Les conditions de \sisi{\textup{(1.2)}}{\Ref{VII.1.2}} \'equivalent \`a:
\begin{enumerate}
\item[\textup{(v)}]\hfill$\codim(Y \cap \Supp G, \Supp G) \geq n$.\hfill\null
\end{enumerate}
\end{corollaire}
\setcounter{equation}{1}
Rappelons d'abord qu'un $\OX$-module est dit de Cohen-Macaulay si, pour tout $x \in X$, la fibre $G_x$ est un $\Oo_{X, x}$-module de Cohen-Macaulay, \ie on~a pour tout $x \in S = \Supp G$:
\begin{equation}\label{eq:VII.1.2}
\prof G_x = \dim G_x = \dim \Oo_{S, x}.
\end{equation}
D'apr\`es la \sisi{(\Ref{III}~\Ref{III.3.3})}{proposition~\Ref{III}~\Ref{III.3.3}}, la condition (i) de \sisi{(\Ref{VII.1.2})}{\Ref{VII.1.2}} est \'equivalente \`a:
\begin{equation}\label{eq:VII.1.3} \prof_Y G = \inf_{x\in Y} \prof G_x \geq n,
\end{equation}
donc aussi \`a:
\[
\prof_Y G = \inf_{x \in Y \cap S} \prof G_x \geq n;
\]
car la profondeur d'un module nul est infinie.

Or, par d\'efinition:
\[
\codim (Y \cap S, S) = \inf_{x \in S \cap Y} \dim \Oo_{S, x},
\]
d'o\`u la conclusion, en appliquant la formule \eqref{eq:VII.1.2}).

Nous allons\pageoriginale maintenant d\'emontrer un r\'esultat permettant de d\'eduire les conditions de coh\'erence que nous avons en vue de certains crit\`eres de nullit\'e.

\begin{lemme}\label{VII.1.5}
Soit $X$ un pr\'esch\'ema localement noeth\'erien. Soit $T^{\ast}$ un $\partial$-foncteur contravariant exact, d\'efini sur la cat\'egorie des $\OX$-Modules coh\'erents, \`a valeurs dans la cat\'egorie des $\OX$-Modules. Soit $Y$ une partie ferm\'ee de $X$. Soit $i \in \ZZ$. Supposons que, pour tout $\OX$-Module coh\'erent de support contenu dans $Y$, $T^i F$ et $T^{i-1} F$ soient coh\'erents. Soit $F$ un $\OX$-Module coh\'erent. Pour que $T^i F$ soit coh\'erent, il faut et il suffit que $T^i F''$ le soit, o\`u l'on a pos\'e:
\[
F'' = F / \Gamma_Y (F).
\]
\end{lemme}

En effet, $F' = \Gamma_Y (F)$ est coh\'erent car $X$ est localement noeth\'erien; la suite exacte de cohomologie de $T^{\ast}$ donne alors:
\[
T^{i-1}F' \to T^i F'' \to T^i F \to T^i F'
\]
o\`u les termes extr\^emes sont coh\'erents, d'o\`u la conclusion.

\begin{lemme}\label{VII.1.6}
Si $F$ et $G$ sont coh\'erents, et si $\Supp F \subset Y$, $\SheafExt^i_Y (F, G)$ est coh\'erent.
\end{lemme}

En effet, $\SheafExt^i_Y (F, G)$ est isomorphe \`a $\SheafExt^i (F, G)$; ceci est valable sur tout espace annel\'e $X$, d'ailleurs: si $Z$ est un ferm\'e qui contient $Y \cap \Supp F$, $\SheafExt^i_Z (F, G)$ est isomorphe \`a $\SheafExt^i_Y (F, G)$ (\cf Expos\'e \Ref{VI}).

\begin{proposition}\label{VII.1.7}
\hspace*{3mm}Supposons $F$ et $G$ coh\'erents et posons $\Supp F = S$, $S' = \overline{S \cap (X-Y)}$. Supposons que, pour tout $x \in Y \cap S'$, on ait $\prof G_x \geq n$, alors $\SheafExt^i_Y (F, G)$ est coh\'erent pour $i < n$.
\end{proposition}

En effet, \sisi{(\Ref{VII.1.6})}{\Ref{VII.1.6}} permet d'appliquer \sisi{(1.5)}{\Ref{VII.1.5}} \`a $T^{\ast}(F) = \SheafExt^{\ast}_Y (F, G)$. En posant $F'' = F/ \Gamma_Y (F)$, on voit que $\Supp F'' = S'$\sisi{, donc $T^i F'' = 0$ pour $i < n$, gr\^ace \`a (1.2) et \`a (III 3.1), d'o\`u la conclusion.}{. Or, d'apr\`es~\Ref{III}~\Ref{III.3.3}, l'hypoth\`ese sur la profondeur de $G$ assure la nullit\'e de $\H_{Y\cap S'}(G)$ pour $i<n$; d'apr\`es~\Ref{VII.1.2}, on en d\'eduit la nullit\'e de $T^i F''$ pour $i < n$, d'o\`u la conclusion gr\^ace \`a~\Ref{VII.1.5}.}

\section{\'Etude pour $i>n$}\label{VII.2}
\pageoriginale
Soit $X$ un pr\'esch\'ema localement noeth\'erien \emph{r\'egulier}, c'est-\`a-dire dont tous les anneaux locaux sont r\'eguliers. Soit $Y$ une partie ferm\'ee de $X$. Soient $F$ et $G$ deux $\OX$-Modules coh\'erents. Posons $S = \Supp F$, $S' = \overline{S \cap (X - Y)}$. Posons:
\begin{align*}
m& = \sup_{x \in Y \cap S} \dim \Oo_{X, x},\\
n &= \sup_{x \in Y \cap S'} \dim \Oo_{X, x};
\end{align*}
on a $n \leq m$.

\begin{proposition}\label{VII.2.1}
Dans la situation d\'ecrite ci-dessus, on a:
\begin{enumerate}
\item
$\SheafExt^i_Y (F, G) = 0$ si $i > m$,

\item
$\SheafExt^i_Y (F, G)$ est coh\'erent si $i > n$.
\end{enumerate}
\end{proposition}

Remarquons d'abord que $\SheafExt^i_Y (F, G)$ est coh\'erent pour tout $i$ lorsque $\Supp F \subset Y$. De plus, en posant comme ci-dessus $F'' = F/\Gamma_Y(F)$, on voit que $\Supp F'' = S'$, donc (2) r\'esulte de (1) et de \sisi{(1.3)}{\Ref{VII.1.3}}.

Pour d\'emontrer (1), on remarque d'abord que
\[
\SheafExt^i_Y (F, G) \approx \varinjlim_{k} \SheafExt^i (F/\Jj^k F, G),
\]
o\`u $\Jj$ est l'id\'eal de d\'efinition de $Y$. Par ailleurs, il
r\'esulte du th\'eor\`eme 4.2.2. de ( A. Grothendieck,
\sisi{\textit{Sur quelques points d'alg\`ebre homologique}, Tohoku
Mathematical Journal 1957}{\og Sur quelques points d'alg\`ebre
homologique\fg, \textit{T\^{o}hoku Mathematical Journal}
\textbf{9}  (1957), p. 119--221.}) que les $\SheafExt$ commutent \`a
la formation des fibres, du moins lorsque $X$ est un pr\'esch\'ema
localement noeth\'erien et lorsque le premier argument est coh\'erent;
comme il en est de m\^eme des limites inductives, on trouve des
isomorphismes:
\[
(\SheafExt^i_Y (F, G))_x \approx \varinjlim_{k} \Ext_{\Oo_{X, x}} ((F/\Jj^k F)_x, G_x)
\]
pour tout $x \in X$. Comme $\Supp \Ext^i_Y (F, G) \subset S \cap Y$, il suffit, pour conclure, de remarquer\pageoriginale que $x \in Y \cap S$ entra\^ine $\dim \Oo_{X, x} \leq m$, donc:
\[
\Ext^i_{\Oo_{X, x}} ((F/\Jj^k F)_x, G_x) = 0 \text { si } i > m,
\]
car la dimension cohomologique globale d'un anneau local r\'egulier est \'egale \`a sa dimension\sfootnote{\Cf $\EGA 0_{\textup{IV}}$ 17.3.1}.

Soit $X$ un pr\'esch\'ema \sisi{loc.}{localement} noeth\'erien; pour toute partie $P$ de $X$, posons:
\[
D(P) = \left\{ \dim \Oo_{X, p} \mid p \in P \right\}.
\]

\begin{lemme}\label{VII.2.2}
Si $P$ est l'espace sous-jacent d'un sous-pr\'esch\'ema connexe de $A$, $D(P)$ est un intervalle.
\end{lemme}

En effet, soient $a$ et $B$ appartenant \`a $D(P)$, correspondant \`a des points $p$ et $q$ de $P$. Montrons qu'il existe une suite de points de $P$: $(p = p_1, \dots, p_n = q)$ telle que, pour $1 \leq i < n $, on ait $\left| \dim \Oo_{X, p_i} - \dim \Oo_{X, p_{i+1}} \right| = 1$; il en r\'esultera que D(P) contient l'intervalle $[p, q]$. Pour cela, on remarque que $p$ et $q$ peuvent \^etre joints par une suite de composantes irr\'eductibles de $P$, telle que deux composantes successives se coupent. On est ramen\'e au cas o\`u $p$ est le point g\'en\'erique d'une composante irr\'eductible $Q$ de~$P$, et o\`u $q \in Q$ et donc $q \supset p$, en tant qu'id\'eaux de $\Oo_q$, o\`u c'est trivial sur la d\'efinition de la dimension.

\begin{proposition}\label{VII.2.3} Soient $X$ un pr\'esch\'ema localement noeth\'erien r\'egulier, $Y$ une partie ferm\'ee de $X$ et $F$ un $\OX$-Module coh\'erent. Soit $P = Y \cap \Supp F \cap (X - Y)$. Soit $n \in \ZZ$, et supposons que $n \not\in D(P)$. Alors $\SheafExt^n_Y (F, \OX)$ est coh\'erent.
\end{proposition}

La conclusion est locale et les hypoth\`eses se conservent par restriction \`a un ouvert. Or $P$ est ferm\'e donc localement noeth\'erien, donc localement connexe; on peut donc supposer $X$ affine et noeth\'erien, et $P$ connexe. Posons $D(P) = [a, b[$, ce qui est licite d'apr\`es le lemme pr\'ec\'edent; si $n> b$, on conclut par \sisi{(2.1)}{\Ref{VII.2.1}}; si $n < a$, on~a \hbox{$n < \dim \Oo_{X, x} = \prof \Oo_{X, x}$} pour tout $x \in P$, et on conclut par \sisi{(1.7)}{\Ref{VII.1.7}}.

\chapterspace{-1}
\chapter{Le th\'eor\`eme de finitude} \label{VIII}

\section[Une suite spectrale de bidualit\'e]
{Une suite spectrale de bidualit\'e\protect\sfootnotemark}\label{VIII.1}

\renewcommand{\theequation}{\thesection.\arabic{equation}}

\'Enon\c cons\pageoriginale le r\'esultat\sfootnotetext{Le lecteur au courant du
langage des cat\'egories d\'eriv\'ees de Verdier reconna\^itra la suite
spectrale associ\'ee \`a un \emph{isomorphisme de bidualit\'e}. \Cf \SGA
6~I.} auquel nous d\'esirons parvenir:

\begin{proposition} \label{VIII.1.1} Soit $A$ un anneau noeth\'erien et soit $I$ un id\'eal de $A$. Posons $X=\Spec(A)$ et $Y=V(I)$. Soit $M$ un $A$-module \emph{de type fini et de dimension projective finie}. Soit ${\Fa}=\widetilde{M}$ le $\OX$-Module associ\'e \`a $M$.
\begin{enumerate}
\item[\textup{1)}] Il existe une suite spectrale:
\[
\H_Y(X, {\Fa}) \Longleftarrow \Ext^p_Y(\Ext^{-q}(M, A), A).
\]

\item[\textup{2)}] Il existe une suite spectrale:
\[
\SheafH_Y(X, {\Fa}) \Longleftarrow \SheafExt^p_Y(\SheafExt^{-q}({\Fa}, \OX), \OX).
\]
\end{enumerate}
\end{proposition}

Bien entendu, 2) se d\'eduit de 1) en remarquant que, si $M$ et $N$ sont deux $A$-modules de type fini et si l'on pose ${{\Fa}}=\widetilde{M}$ et ${{\Ga}}=\widetilde{N}$, on~a des isomorphismes:
\begin{align*}
\SheafH_Y({\Fa}) &\approx \widetilde{\H_Y(X, {\Fa})}, \\
\SheafExt_Y({\Fa}, {\Ga}) &\approx \widetilde{\Ext_Y({\Fa}, {\Ga}}), \\
\SheafExt_{\OX}({\Fa}, {\Ga}) &\approx \widetilde{\Ext_A(M, N)}.
\end{align*}

Soit $\ccat$ la cat\'egorie des $A$-modules, et $\Ab$ celle des groupes ab\'eliens. Soit ${\cF}$ le foncteur:
\begin{align*}
\cF: \ccat&\to \Ab \quad\text{d\'efini par} \\
 M &\sisi{\rightsquigarrow}{\mto} \Gamma_Y(\widetilde{M}).
\end{align*}

On sait\pageoriginale depuis l'expos\'e \Ref{II} qu'il existe un isomorphisme de $\partial$-foncteurs:
\[
\H^\ast_Y(X, \widetilde{M})\approx \R^\ast {\cF}(M).
\]
De plus, soient $\Ext^\ast_Y$ les foncteurs d\'eriv\'es droits en le deuxi\`eme argument de
\[{\cF}\circ \Hom\colon \ccat^\circ \times \ccat \to \Ab.
\]
On sait depuis l'expos\'e \Ref{VI} que l'on a un isomorphisme de $\partial$-foncteurs:
\[
\Ext^\ast_Y(M, N)\simeq \Ext^\ast_Y(\widetilde{M}, \widetilde{N}).
\]
Retenons enfin le r\'esultat suivant de l'expos\'e \Ref{VI}: si $C$ est un $A$-module injectif et si $N$ est un $A$-module de type fini, le faisceau $\SheafHom(\widetilde{N}, \widetilde{C})\approx\widetilde{\Hom(N, C)}$ est \emph{flasque}, donc $\R^1{\cF}(\Hom(N, C))=0$.

Il nous reste \`a d\'emontrer le r\'esultat suivant:

\begin{lemme} \label{VIII.1.2} Soit $A$ un anneau noeth\'erien et soit $\ccat$ la cat\'egorie des $A$-modules. Soit ${\cF}\colon \ccat\to\Ab$ un foncteur additif exact \`a gauche tel que, pour tout $A$-module $N$ \sisi{\ignorespaces}{de type fini} et tout $A$-module \emph{injectif} $C$, on ait $\R^1{\cF}(\Hom(N, C))=0$. Soit $M$ un $A$-module de type fini et de dimension projective finie. Il existe une suite spectrale:
\[
\R^\ast {\cF}(M)\Longleftarrow \Ext^p_{{\cF}}(\Ext^{-q}(M, A), A), \]
o\`u $\Ext^p_{{\cF}}$ d\'esigne le $p$-i\`eme foncteur d\'eriv\'e droit de ${\cF}\circ \Hom$.
\end{lemme}

\emph{Nous ne consid\'ererons que des complexes dont la diff\'erentielle est de degr\'e $+1$}. D'apr\`es l'hypoth\`ese faite sur $M$, il existe une r\'esolution \sisi{\emph{libre}}{\emph{projective}} de $M$ de longueur finie:
\[u\colon L^\boule \to M, \]
o\`u, de plus, les $L^p$ sont des modules de \emph{type fini} et $L^p=0$ si $p\not\in[-n, 0]$. Soit par ailleurs, $v\colon M\to I^\boule$ une r\'esolution injective de $M$. Je dis\pageoriginale que
\begin{equation} \label{eq:VIII.1.1} v\circ u\colon L^\boule\to I^\boule
\end{equation}
est une r\'esolution injective de $L^\boule$. Il convient de pr\'eciser ce que cela signifie.

\begin{definition} \label{VIII.1.3} Soit $X^\boule$ un complexe de $A$-modules; on appelle r\'esolution injective de $X^\boule$ un homomorphisme de complexes:
\[x\colon X^\boule\to CX^\boule, \]
tel que $CX^p$ soit injectif pour tout $p\in \ZZ$, et que $x$ induise un isomorphisme en homologie.
\end{definition}

\begin{proposition} \label{VIII.1.4} Tout complexe limit\'e \`a gauche, \ie tel qu'il existe $q\in\ZZ$ avec $X^p=0$ pour $p<q$, admet une r\'esolution injective. De plus, si $u\colon X^\boule\to Y^\boule$ est un homomorphisme de complexes (limit\'es \`a gauche) et si $x\colon X^\boule \to CX^\boule$ et $y\colon Y^\boule\to CY^\boule$ sont des r\'esolutions injectives de $X^\boule$ et $Y^\boule$, il existe un homomorphisme de complexes:
\[Cu\colon CX^\boule\to CY^\boule, \]
unique \`a homotopie pr\`es, tel que le diagramme:
\[
\xymatrix{ X^\boule \ar[r]^x\ar[d]_u & CX^\boule\ar[d]^{Cu} \\
Y^\boule \ar[r]^y & CY^\boule \\
}\]
soit commutatif \`a homotopie pr\`es.
\end{proposition}

La d\'emonstration est laiss\'ee au lecteur\sfootnote{\Cf aussi \sisi{H.~Cartan - S.~Eilenberg, Homological Algebra, Princeton University Press, 19}{H.~Cartan \& S.~Eilenberg, \emph{Homological Algebra}, Princeton Math. Series, vol.~19, Princeton University Press, 1956.}}.

Rappelons une notation introduite dans l'expos\'e \Ref{V}.

\subsection*{Notation}
Soient\pageoriginale $X^\boule$ et $Y^\boule$ deux complexes. On note $\Hom^\boule(X^\boule, Y^\boule)$ le \emph{complexe simple} dont la composante de degr\'e $n$ est
\[(\Hom^\boule(X^\boule, Y^\boule))^n=\prod_{-p+q=n}\Hom(X^p, Y^q)\]
not\'ee aussi $\Hom^n(X^\boule, Y^\boule)$, et dont la diff\'erentielle est donn\'ee par:
\begin{gather*}
\partial_n\colon \Hom^n(X^\boule, Y^\boule)\to \Hom^{n+1}(X^\boule, Y^\boule)\\
\partial_n=d'+(-1)^{n+1}d'',
\end{gather*}
o\`u $d'$ et $d''$ sont les diff\'erentielles (de degr\'e $+1$) induites par celles de $X^\boule$ et de $Y^\boule$.

\setcounter{equation}{1} Soit alors $A^\boule$ le complexe d\'efini par $A^p=0$ si $p\neq 0$ et $A^0=A$. Soit
\[a\colon A^\boule \to CA^\boule\]
une r\'esolution injective de $A^\boule$. Consid\'erons le double complexe:
\begin{equation} \label{eq:VIII.1.2} Q^{p, q}=\Hom(\Hom(\sisi{L^q}{L^{-q}}, A), CA^p).
\end{equation}
La premi\`ere suite spectrale du bicomplexe ${\cF}Q^{\boule\boule}$ donnera la conclusion du lemme \Ref{VIII.1.2}.

Posons
\begin{equation} \label{eq:VIII.1.3} {L'}^\boule=\Hom^\boule(L^\boule, A^\boule),
\end{equation}
et
\begin{equation} \label{eq:VIII.1.4} P^\boule=\Hom^\boule({L'}^\boule, CA^\boule).
\end{equation}

On voit facilement que $P^\boule$ est le complexe simple associ\'e \`a $Q^{\boule\boule}$. Calculons l'aboutissement de la suite spectrale \ie l'homologie de ${\cF}P^\boule$. Pour cela, utilisant le fait que $L^\boule$ est \sisi{libre}{projectif} de type fini en toute dimension, on prouve que $L^\boule$ est isomorphe \`a $\Hom^\boule({L'}^\boule, A^\boule)$. De l'homomorphisme $a\colon A^\boule\to CA^\boule$\pageoriginale, on d\'eduit un homomorphisme:
\[b\colon \Hom^\boule({L'}^\boule, A^\boule)\to \Hom^\boule({L'}^\boule, CA^\boule), \]
ou encore, un homomorphisme
\begin{equation} \label{eq:VIII.1.5} c\colon L^\boule\to P^\boule.
\end{equation}
Ceci dit, il est facile de voir, en utilisant le fait que ${L'}^\boule$ est \sisi{libre}{projectif} de type fini en toute dimension et limit\'e \`a gauche, que (\Ref{eq:VIII.1.5}) est une r\'esolution injective de $L^\boule$. Utilisant la proposition \Ref{VIII.1.4}., on en conclut que $P^\boule$ est homotopiquement \'equivalent \`a $I^\boule$, o\`u $I^\boule$ est la r\'esolution injective de $M$ introduite plus haut (\Ref{eq:VIII.1.1}). On en d\'eduit que l'\emph{aboutissement} de la premi\`ere suite spectrale du double complexe ${\cF}Q^{\boule\boule}$, qui est $\H^\ast({\cF}P^\boule)$, \emph{est isomorphe} \`a $\R^\ast {\cF}(M)$.

Le terme initial de la premi\`ere suite spectrale du bicomplexe ${\cF}Q^{\boule\boule}$ est:
\[
\E^{p, q}_2 = {'\!\H}^p({''\!\H}^q({\cF}Q^{\boule\boule})).
\]

Pour tout $p\in\ZZ$, $CA^p$ est injectif. D'apr\`es l'hypoth\`ese faite sur ${\cF}$, le foncteur\sisi{}{ (restreint \`a la cat\'egorie des modules de type fini)}:
\[
N \sisi{\rightsquigarrow}{\mto} {\cF}\Hom(N, CA^p)\]
est exact. D'o\`u l'on d\'eduit des isomorphismes:
\[
{''\!\H}^q({\cF}\Hom({L'}^\boule, CA^p))\approx {\cF}\Hom(\H^{-q}({L'}^\boule), CA^p).
\]
D'apr\`es la d\'efinition de $\Ext^\ast_{{\cF}}$ comme foncteur d\'eriv\'e de ${\cF}{\sisi{^\ast}{\circ}}\Hom$, on en d\'eduit des isomorphismes:
\[
\E^{p, q}_2 \approx \Ext^p_{{\cF}}(\H^q({L'}^\boule), A).
\]
Or ${L'}^\boule=\Hom^\boule(L^\boule, A^\boule)$, o\`u $L^\boule$ est une r\'esolution \sisi{libre}{projective} de $M$ d'o\`u des isomorphismes:
\[
\Ext^{q}(M, A)\approx \H^q({L'}^\boule),
\]
ce qui donne la conclusion.\qed

\section{Le th\'eor\`eme de finitude} \label{VIII.2}
\pageoriginale
\refstepcounter{subsection}\label{VIII.2.1}
\begin{enonce*}{Th\'eor\`eme 2.1\ndemark}
Soient $X$ un pr\'esch\'ema localement noeth\'erien, $Y$ partie ferm\'ee de $X$ et $\sisi{F}{{\Fa}}$ un $\OX$-Module coh\'erent.\ndetext{pour un \'enonc\'e analogue, mais dans une situation un peu plus g\'en\'erale, voir Mme Raynaud (Raynaud~M., {\og Th\'eor\`emes de Lefschetz en cohomologie des faisceaux coh\'erents et en cohomologie \'etale. Application au groupe fondamental\fg}, \emph{Ann. Sci. \'Ec. Norm. Sup. (4)} \textbf{7} (1974), p\ptbl 29--52, proposition II.2.1).}
Supposons que $X$ soit localement immergeable dans un pr\'esch\'ema r\'egulier\sfootnote{Cette condition peut se g\'en\'eraliser en l'hypoth\`ese d'existence localement sur $X$, d'un \og complexe dualisant\fg, au sens d\'efini dans R\ptbl Hartshorne, \sisi{Residues and duality}{\emph{Residues and duality}} (cit\'e dans \sisi{}{la note~\eqref{cithar} de l'}\Exp \Ref{IV} p\ptbl\sisi{23}{\pageref{pagecithar}}).}. Soit $i\in\ZZ$. Supposons que:
\begin{enumeratea}
\item
pour tout $x\in U=X-Y$, on a
\[
\H^{i-c(x)}({\sisi{F}{{\Fa}}}_x)=0,
\]
o\`u l'on a pos\'e\,\nde{comme dans l'expos\'e~\Ref{V}, $H^\ast({\Fa}_x)$ d\'esigne la cohomologie locale $H_{\m_x}({\Fa}_x)$.}:
\begin{equation} \label{eq:VIII.2.1} c(x)=\codim (\overline{\{x\}}\cap Y, \overline{\{x\}}).
\end{equation}
Alors:
\item
$\SheafH^i_Y(\sisi{F}{{\Fa}})$ est coh\'erent.
\end{enumeratea}
\end{enonce*}

\refstepcounter{subsection}\label{VIII.2.2}
\begin{enonce*}{Corollaire 2.2\ndemark}
Sous les hypoth\`eses du th\'eor\`eme pr\'ec\'edent, la condition \textup{a)} est \'equivalente \`a:
\ndetext{\textit{stricto sensu}, c'est un corollaire de la preuve qui suit et non de l'\'enonc\'e. L'implication c)\ALORS a) est tautologique. L'autre sens ne l'est pas, mais r\'esulte de la preuve. Pr\'ecisons. Comme plus bas, on recouvre $X$ par des ouverts immergeables dans des sch\'emas r\'eguliers, ce qui permet comme expliqu\'e plus-bas de se ramener \`a $X=\Spec(A)$ affine r\'egulier et $F=\tilde M$ o\`u $M$ est un $A$-module de dimension projective finie. Il est montr\'e dans ce cas que les conditions a) et c) \'equivalent aux conditions duales a$'$) et c$'$). On montre alors que c$'$) entra\^ine la question d) (voir \textit{infra}) qui elle-m\^eme entra\^ine a$'$). Voir les consid\'erations suivant~\Ref{VIII.2.4}.}
\begin{enumerate}
\item[\textup{c)}] pour tout $x\in U$ tel que $c(x)=1$, on~a $\H^{i-1}(\sisi{F}{{\Fa}}_x)=0$.
\end{enumerate}
\end{enonce*}

\begin{corollaire} \label{VIII.2.3} \sisi{\emph{Soient}}{Soient} $X$ un pr\'esch\'ema localement noeth\'erien et localement immergeable dans un pr\'esch\'ema r\'egulier, $Y$ une partie ferm\'ee de $X$, $\sisi{F}{{\Fa}}$ un $\OX$-Module coh\'erent, $n$ un entier. Les conditions suivantes sont \'equivalentes:
\begin{enumeratei}
\item
pour tout $x\in U$, on~a $\prof \sisi{F}{{\Fa}}_x>n-c(x)$;
\item
pour tout $x\in U$ tel que $c(x)=1$, on~a $\prof \sisi{F}{{\Fa}}_x\geq n$;
\item
pour tout $i\in \ZZ$, $\SheafH^i_Y(\sisi{F}{{\Fa}})$ est coh\'erent si $i\leq n$.
\item
$\R^ii_\ast(\sisi{F}{{\Fa}}_{|U})$ est coh\'erent pour $i<n$\ \nde{cette condition n'\'etait que dans le corps dans la preuve mais pas dans l'\'enonc\'e du corollaire; comme elle sert en~\Ref{VIII.3}, on l'a ajout\'ee.}.
\end{enumeratei}
\end{corollaire}

Supposons ces r\'esultats acquis lorsque $X$ est le spectre d'un anneau noeth\'erien r\'egulier $A$ et lorsque $\sisi{F}{{\Fa}}$ est le faisceau associ\'e \`a un $A$-module de dimension projective finie.

Remarquons\pageoriginaled d'abord que, si $(X_j)_{j\in
J}$ est un recouvrement ouvert de $X$ par des ouverts\sisi{}{
immergeables dans un sch\'ema r\'egulier}, chacune des conditions
ci-dessus est \'equivalente \`a la conjonction des conditions
analogues obtenues en rempla\c cant $X$ par~$X_j$, $Y$ par
$Y_j=Y\cap X_j$ et $\sisi{F}{{\Fa}}$ par
$\sisi{F{|X_j}}{{\Fa}_|X_j}$. En effet, seules les conditions
faisant intervenir $c(x)$ peuvent faire une difficult\'e. Soit $x\in
U$. Si $x\in X_j$, posons
\[
c_j(x)=\codim(X_j\cap \overline{\{x\}}\cap Y, X_j\cap\overline{\{x\}}),
\]
on a n\'ecessairement $c_j(x)\geq c(x)$. Soit $y\in
\overline{\{x\}}\cap Y$, qui \og donne la codimension\fg, \ie tel
que $c(x)=\dim \Oo_{\overline{\{x\}}, y}$, soit $X_j$ un ouvert du
recouvrement tel que $y\in X_j$, alors $x\in X_j$, donc
$c_j(x)=c(x)$, ce qui permet de conclure\sisi{.}{ que a) pour les
$X_j$ entra\^ine a) pour $X$. \`A ce stade, on~a seulement une
r\'eciproque partielle, \`a savoir que a) pour $X$ entra\^ine a) pour
les $X_j$ tels que $c(x)=c_j(x)$, ce qui suffit pour notre
propos.}\nde{en fait, a) pour $X$ entra\^ine a) pour tous les
$X_j$ comme affirm\'e dans le texte original, mais il faut pour cela
lire la preuve qui suit en d\'etail. Cette implication ne semble pas
formelle \`a ce stade. Notons en effet par un indice $J$ les
conjonctions d'une propri\'et\'e a), b) ou c) pour les $X_j$. Il est
d\'emontr\'e dans la preuve \textit{infra} c$_J$)\ALORS a$_J$) (c'est
la suite d'implications c$'$)\ALORS d)\ALORS a$'$)\,). Or, on~a
tautologiquement a) \ALORS c), et c)\SSI c$_J$), d'o\`u a) \ALORS
a$_J$).}

On choisit un recouvrement de $X$ par des ouverts immergeables dans un pr\'esch\'ema r\'egulier. Appliquant ce qui pr\'ec\`ede, on voit qu'on peut supposer $X$ ferm\'e dans un $X'$ r\'egulier. La r\'eduction \`a $X'$ est alors imm\'ediate.

On peut donc supposer $X$ r\'egulier, et m\^eme affine en recouvrant $X$ par des ouverts affines. Que l'on puisse supposer que $\sisi{F}{{\Fa}}=\widetilde{M}$, o\`u $M$ est de dimension projective finie r\'esultera du lemme suivant:

\begin{lemme} \label{VIII.2.4}
Soit $X$ un pr\'esch\'ema noeth\'erien r\'egulier. Soit $\sisi{F}{{\Fa}}$ un $\OX$-Module coh\'erent. La fonction qui \`a tout $x\in X$ associe la dimension projective de $\sisi{F}{{\Fa}}_x$ est born\'ee sup\'erieurement.
\end{lemme}

\setcounter{equation}{1}

Soit en effet $x\in X$ et soit $U$ un voisinage ouvert affine de $x$. Soit $L^\boule$ une r\'esolution projective du module $\sisi{F}{{\Fa}}(U)$, o\`u les $L^i$ sont de type fini. Par hypoth\`ese, l'anneau $\Oo_{X, x}$ est r\'egulier, donc la dimension projective de $\sisi{F}{{\Fa}}_x$ est finie; soit $d$ cet entier. Soit
\[
K=\ker(\sisi{L^{d-1}\to L^{d-2}}{L^{-d}\to L^{-d+1}}).
\]
Le module $K_x$ est libre, car $d$ est la dimension projective de $\sisi{F}{{\Fa}}_x$ (\sisi{[M]}{\cite{M}}, Ch\ptbl VI. Prop\ptbl 2.1). D'apr\`es ($\EGA 0_{\textup{I}}$~5.4.1 Errata), on en d\'eduit que le $\Oo_U$-Module $\widetilde{K}$ est libre sur un voisinage $U'$ de $x$, $U'\subset U$. Choisissant $f\in \OX(U)$ tel\pageoriginale que $\sisi{}{x\in}D(f)\subset U'$, on~a donc une r\'esolution projective de $M_{f'}$ ($M=\sisi{F}{{\Fa}}(U)$):
\[0\to K_f \to (L^{d-1})_f\to \cdots \to M_f\to 0, \]
ce qui prouve que la fonction \'etudi\'ee est semi-continue sup\'erieurement. Or $X$ est quasi-compact, d'o\`u la conclusion.

Nous supposons d\'esormais $X$ affine noeth\'erien r\'egulier et nous supposons que \hbox{$\sisi{F}{{\Fa}}=\widetilde{M}$}, o\`u $M$ est un $A$-module de type fini, n\'ecessairement de dimension projective finie. Nous proc\'ederons en plusieurs \'etapes. Tout d'abord, nous trouvons une condition d), \'equivalente \`a a), et prouvons qu'elle \'equivaut \'egalement \`a c). Puis, \`a l'aide de la suite spectrale du num\'ero pr\'ec\'edent, nous prouvons \text{d)}\ALORS\text{b)}. Il reste alors \`a prouver que \text{(iii)}\ALORS\text{(ii)}; en effet, \text{(i)}\SSI\text{(ii)}\ALORS\text{(iii)} r\'esulte imm\'ediatement de \text{a)}\SSI\text{c)}\ALORS\text{b)}.

Soit $x\in U$, par hypoth\`ese $\Oo_{X, x}$ est un anneau local r\'egulier; en d\'esignant par $D$ le foncteur dualisant relatif \`a l'anneau local $\Oo_{X, x}$, il r\'esulte de (\Ref{V}~\Ref{V.2.1}) que:
\[D\H^{i-c(x)}(\sisi{F}{{\Fa}}_x)\approx \Ext^{d(x)-i}_{\Oo_{X, x}}(\sisi{F}{{\Fa}}_x, \Oo_{X, x}), \]
o\`u l'on a pos\'e
\begin{equation} \label{eq:VIII.2.2} d(x)=\dim \Oo_{X, x} + c(x) = \dim \Oo_{X, x}+\codim(\overline{\{x\}}\cap Y, \overline{\{x\}}).
\end{equation}
Or $X$ est noeth\'erien et $\sisi{F}{{\Fa}}$ coh\'erent, donc:
\begin{equation} \label{eq:VIII.2.3} D\H^{i-c(x)}(\sisi{F}{{\Fa}}_x)\approx (\sisi{\Ext}{\SheafExt}^{d(x)-i}_{\OX}(\sisi{F}{{\Fa}}, \OX))_x.
\end{equation}
De plus, pour qu'un module soit nul, il faut et il suffit que son dual le soit\sisi{}{ (\cf la note de l'\'editeur~\eqref{dualnul} de la page~\pageref{dualnul})}. Pour tout $q\in\ZZ$, posons\pageoriginale:
\begin{equation} \label{eq:VIII.2.4}
\begin{cases}
S_q=\Supp \SheafExt^q_{\OX}(\sisi{F}{{\Fa}}, \OX), \\
S'_q=S_q\cap U\text{, (}U=X-Y\text{)}, \\
Z_q=\overline{S'_q}\cap Y.
\end{cases}
\end{equation}
De la formule (\Ref{eq:VIII.2.3}), il r\'esulte que a) et c) sont respectivement \'equivalentes \`a
\begin{enumerate}
\item[a$'$)]
pour tout $q\in\ZZ$ et tout $x\in S'_q$, on~a $q+i\neq d(x)$.

\item[c$'$)]
pour tout $q\in\ZZ$ et tout $x\in S'_q$, si $c(x)=1$, on~a $q+i\neq d(x)$.
\end{enumerate}
\noindent Voici la condition d) promise plus haut:
\begin{enumerate}
\item[d)]
pour tout $q\in \ZZ$ et tout $y\in Z_q$, on~a $q+i\neq \dim \Oo_{X, y}$.
\end{enumerate}
\noindent Ces conditions sont \'equivalentes:

a$'$) \ALORS c$'$) pour m\'emoire.

d) \ALORS a$'$). En effet, soit $q\in \ZZ$ et soit $x\in S'_q$; soit $y\in\overline{\{x\}}\cap Y$ qui\nde{dans tout ce qui suit, les adh\'erences de points sont munies de la structure r\'eduite.} \og donne la codimension\fg, \ie tel que:
\begin{equation} \label{eq:VIII.2.5} \dim \Oo_{\overline{\{x\}}, y}=\codim(\overline{\{x\}}\cap Y, \overline{\{x\}})=c(x).
\end{equation}
Du fait que $X$ est r\'egulier en $y$, on d\'eduit:
\begin{equation} \label{eq:VIII.2.6} \dim \Oo_{X, y}=d(x)\quad \text{(\cf (\Ref{eq:VIII.2.2})).}
\end{equation}
Mais $y\in\overline{\{x\}}$, donc $y\in Z_q$, d'o\`u la conclusion.

c$'$) \ALORS d). Soit $q\in\ZZ$ et soit $y\in Z_q$. Admettons provisoirement qu'il existe $x\in S'_q$ tel que:
\[
\sisi{}{y\in\overline{\{x\}}\textup{ et }}\dim \Oo_{\overline{\{x\}}, y}=1\, \text{;}\]
(on dit aussi que $x$ \emph{suit} $y$). Il en r\'esulte que $c(x)=1$, car $y$ \og donne la codimension de $\overline{\{x\}}\cap Y$ dans $\overline{\{x\}}$\fg, puisque $x\not\in Y$. D'apr\`es c$'$) on en tire
\[q+i\neq d(x).
\]
D'o\`u\pageoriginale la conclusion, si l'on remarque que $d(x)=\dim \Oo_{X, y}$ (\Ref{eq:VIII.2.6}). Le r\'esultat admis s'exprime dans le lemme suivant:

\begin{lemme} \label{VIII.2.5} Soit $X$ un pr\'esch\'ema localement noeth\'erien et soit $Y$ une partie ferm\'ee de $X$. Posons $U=X-Y$ et supposons que $U$ est dense dans $X$. Pour tout $y\in Y$, il existe $x\in U$ \og qui le suit\fg, \ie tel que:
\[y\in\overline{\{x\}}\text{ et } \dim\Oo_{\overline{\{x\}}, y}=1.
\]
\end{lemme}

Nous avons appliqu\'e le lemme en prenant pour $X$ le pr\'esch\'ema $\overline{S'_q}$ et pour $Y$ la partie $Y\cap \overline{S'_q}$.

\begin{proof}[D\'emonstration de~\Ref{VIII.2.5}]
Il existe $x\in U$ tel que $y\in\overline{\{x\}}$; choisissons donc un $x\in U$ tel que $y\in\overline{\{x\}}$ et tel que $\dim \Oo_{\overline{\{x\}}, y}=r$ soit minimal. Il faut prouver que $r=1$. Puisque l'on a choisi $x$ de telle sorte que tout $z\in\sisi{\overline{\{x\}}}{\Spec(\Oo_{\overline{\{x\}}, y})}$, $z\neq x$, soit dans $Y$, $\{x\}$ est ouvert dans $\Spec(\Oo_{\overline{\{x\}}, y})$. D'o\`u la conclusion.
\skipqed
\end{proof}

\emph{La deuxi\`eme \'etape consiste \`a d\'eduire \textup{b)} de \textup{d)}.}

\noindent Posons $D(Z_q)=\{\dim \Oo_{X, y} \mid y\in Z_q\}$. D'apr\`es d), on sait que, pour tout $q\in\ZZ$, on~a $q+i\not\in D(Z_q)$. On applique alors \sisi{\Ref{VIII.2.2}}{\Ref{VII}.\Ref{VII.2.3}}, et l'on voit que
\[
\SheafExt^{q+i}_Y(\SheafExt^q(\sisi{F}{{\Fa}}, \OX), \OX) \text{ est coh\'erent.}\]
Le terme initial de la suite spectrale du num\'ero pr\'ec\'edent est donn\'e par:
\[
\E^{p, q}_2=\SheafExt^p_Y(\SheafExt^{-q}(\sisi{F}{{\Fa}}, \OX), \OX).
\]
On en d\'eduit que $\E^{p, q}_2$ est coh\'erent pour tout $p\in\ZZ$ et tout $q\in\ZZ$ tels que $p+q=i$. Or il n'y a qu'un nombre fini de couples $(p, q)$ tels que $p+q=i$, et cette suite spectrale converge vers $\SheafH^\ast_Y(\sisi{F}{{\Fa}})$, d'o\`u la conclusion.

\emph{Il nous reste \`a prouver que \textup{(iii)}\ALORS\textup{(ii)}.} Notons:
\[i\colon U\to X\]
l'immersion canonique de $U$ dans $X$. Compte tenu de la suite exacte d'homologie du ferm\'e $Y$ (\Ref{I}~\Ref{I.2.11}) on voit que (iii) \emph{\'equivaut \`a}:

(iv) $\R^ii_\ast(\sisi{F}{{\Fa}}_{|U})$ est coh\'erent pour $i<n$. \refstepcounter{toto}\label{condition}

\noindent En effet\pageoriginale, on~a une suite exacte:
\[0\to \SheafH^0_Y(\sisi{F}{{\Fa}})\to \sisi{F}{{\Fa}} \to i_\ast(\sisi{F}{{\Fa}}_{|U})\to \SheafH^1_Y(\sisi{F}{{\Fa}}) \to 0.
\]
Or $\SheafH^0_Y(\sisi{F}{{\Fa}})$ est un sous-faisceau quasi-coh\'erent de $\sisi{F}{{\Fa}}$ qui est coh\'erent, donc est coh\'erent. Donc $\SheafH^1_Y(\sisi{F}{{\Fa}})$ est coh\'erent si et seulement si $i_\ast(\sisi{F}{{\Fa}}_{|U})$ l'est. Par ailleurs, si $p>0$, la suite exacte de cohomologie du ferm\'e $Y$ se r\'eduit \`a des isomorphismes:
\[
\R^pi_\ast(\sisi{F}{{\Fa}}_{|U})\overset{\approx}{\to}\SheafH^{p+1}_Y(\sisi{F}{{\Fa}}).
\]
Nous allons d\'emontrer que (iv) \ALORS (ii). Pour cela, rappelons (ii):

(ii) pour tout $x\in U$ tel que $c(x)=1$, on~a $\prof \sisi{F}{{\Fa}}_x\geq n$.

Raisonnons par r\'ecurrence sur $n$.

Si $n=0$, les deux conditions sont vides.

Si $n=1$, on suppose que $i_\ast(\sisi{F{|U}}{{\Fa}_{|U}})$ est coh\'erent. Raisonnons par l'absurde et supposons qu'il existe $x\in U$ tel que $c(x)=1$ et $\prof \sisi{F}{{\Fa}}_x=0$, \ie $x\in \Ass \sisi{F}{{\Fa}}_x$. Soit $y\in\overline{\{x\}}\cap Y$ tel que $\dim \Oo_{\overline{\{x\}}, y}=1$. Posons:
\[A=\Oo_{X, y}\text{ et } X'=\Spec(A).
\]

\setcounter{equation}{6}

Effectuons le changement de base $v\colon X'\to X$, qui est plat.
\begin{equation} \label{eq:VIII.2.7}
\begin{array}{c}
\xymatrix@C=5em{ U'=X'\times_X U\ar[r]^-{v'}\ar[d]_{i'} & U\ar[d]^i \\
X'\ar[r]^-v & \,X\,.
}
\end{array}
\end{equation}
Le morphisme $i$ est s\'epar\'e (car c'est une immersion), et de type fini (car c'est une immersion ouverte et $X$ est localement noeth\'erien), le changement de base est plat donc (\EGA III~1.4.15) on~a un isomorphisme:
\begin{equation} \label{eq:VIII.2.8}
v^\ast(i_\ast(\sisi{F|U}{{\Fa}_{|U})}) \approx i'_\ast({v'}^\ast(\sisi{F}{{\Fa}}_{|U})).
\end{equation}

Notons\pageoriginale $\underline x$ (\resp $\underline y$) l'id\'eal de $A$ correspondant \`a $x$ (\resp $y$).

Posons ${\Ga}={v'}^\ast(\sisi{F}{{\Fa}}{}_{|U})$; ${\Ga}$ est coh\'erent et $\underline{x}\in\Ass {\Ga}$, donc il existe un monomorphisme $\Oo_{\overline{\{x\}}}\to {\Ga}$, et par suite $i'_\ast(\Oo_{\overline{\{x\}}}{}_{|U'})$ est coh\'erent. D'apr\`es le choix de~$y$, $\dim A/\underline x=1$, et par suite le support de $\Oo_{\overline{\{x\}}}$ est r\'eduit \`a $\overline{\{\underline x\}}= \{\underline x\}\cup \{\underline y\}$, car $\overline{\{\underline x\}}= \Spec(A/\underline x)$ en tant que sch\'ema. Il en r\'esulte que
\[(\Oo_{\overline{\{x\}}}{}_{|U'})(U')=\Frac (A/\underline x), \]
anneau des fractions de $A/\underline x$, et
\[i'_\ast(\Oo_{\overline{\{x\}}}{}_{|U'})(X')=\Frac(A/\underline x).
\]
Mais $\Frac(A/\underline x)$ n'est pas un $A$-module de type fini car $\underline x$ est diff\'erent de l'id\'eal maximal de $A$. D'o\`u une contradiction.

Supposons $n>1$ et le r\'esultat acquis pour les $n'<n$. Par l'hypoth\`ese de r\'ecurrence, pour tout $x\in U$ tel que $c(x)=1$, on~a $x\not\in\Ass\sisi{F}{{\Fa}}_x$. Soit un tel $x$, et soit $y\in\overline{\{x\}}\cap Y$ tel que $x$ suive $y$, \ie $\dim\Oo_{\overline{\{x\}}, y}=1$. On effectue le changement de base $v\colon\Spec (\Oo_{X, y})\to X$ en conservant les notations du diagramme (\Ref{eq:VIII.2.7}). On trouve, en appliquant (\EGA III~1.4.15), des isomorphismes:
\[v^\ast(\R^pi_\ast(\sisi{F}{{\Fa}}{}_{|U})\sisi{}{)} \simeq \R^pi'_\ast({v'}^\ast(\sisi{F}{{\Fa}}{}_{|U}))\text{, }p\in\ZZ.
\]
\sisi{}{\indent}\emph{On se ram\`ene ainsi au cas o\`u $X$ est le spectre d'un anneau local $A$ dans lequel $\sisi{x}{\underline{x}}$ est un id\'eal premier de dimension $1$}, \ie $\dim A/x=1$. Posons alors $\sisi{F}{{\Fa}}'=\SheafGamma_Y(\sisi{F}{{\Fa}})$ et $\sisi{F}{{\Fa}}''=\sisi{F}{{\Fa}}/\sisi{F}{{\Fa}}'$.

\noindent On voit que $\sisi{F}{{\Fa}}_x\simeq \sisi{F}{{\Fa}}''_x$ et que $y\not\in \Ass\sisi{F}{{\Fa}}''$. Par ailleurs $\sisi{F}{{\Fa}}'{}_{|U}=0$ d'o\`u, par la suite exacte des $\R^pi_\ast$, des isomorphismes:
\[
\R^pi_\ast(\sisi{F}{{\Fa}}{}_{|U})\simeq \R^pi_\ast(\sisi{F}{{\Fa}}''{}_{|U})\text{, }p\in\ZZ.
\]
Puisque $n>1$, on en d\'eduit que ni $x$ ni $y$ n'appartiennent \`a $\Ass \sisi{F}{{\Fa}}''$. Or\sisi{ ce }{, $\underline{x}, \underline{y}$ }sont les seuls id\'eaux premiers de $A$ qui contiennent $\sisi{x}{\underline{x}}$; il en r\'esulte (\Ref{III}~\Ref{III.2.1}) qu'il existe\pageoriginale un \'el\'ement $g\in \sisi{x}{\underline{x}}$ qui est $M$-r\'egulier, o\`u l'on a pos\'e $\sisi{F}{{\Fa}}=\widetilde{M}$, $M=\sisi{F}{{\Fa}}(X)$. D'o\`u une suite exacte:
\[0\to M\lto{g'} M\to N \to 0,
\]
dans laquelle $g'$ d\'esigne l'homoth\'etie de rapport $g$ dans $M$. Par la suite exacte d'homologie, on voit que:

$\R^pi_\ast(\widetilde{N}{}_{|U})$ est coh\'erent pour $p<n-1$,

\noindent donc, par l'hypoth\`ese de r\'ecurrence, $\prof(\widetilde{N})_x\geq n-1$, donc $\prof \sisi{F}{{\Fa}}_x\geq n$, \sisi{QED.}{\qed}

\section{Applications} \label{VIII.3}

On d\'eduit de ces r\'esultats une condition de coh\'erence pour les images directes sup\'erieures d'un faisceau coh\'erent par un \emph{morphisme qui n'est pas propre}.

\begin{theoreme} \label{VIII.3.1}
Soit $f\colon X\to Y$ un morphisme de pr\'esch\'emas. Supposons que~$Y$ est \emph{localement noeth\'erien} et que $f$ est \emph{propre}. Supposons que $X$ soit \emph{localement immergeable dans un pr\'esch\'ema r\'egulier}. Soit $n\in\ZZ$. Soit $U$ un ouvert de $X$ et soit $\sisi{F}{{\Fa}}$ un $\Oo_U$-Module \emph{coh\'erent}. Supposons que, pour tout $x\in U$ tel que \hbox{$\codim(\overline{\{x\}}\cap (X-U), \overline{\{x\}})\!=\!1$}, on ait $\prof\sisi{F}{{\Fa}}_x\geq n$. \emph{Alors} les $\Oo_Y$-Modules \hbox{$\R^p(f\circ g)_\ast(F)$} sont coh\'erents pour $p<n$, o\`u $g$ est l'immersion canonique de $U$ dans~$X$.
\end{theoreme}

En effet, il existe une suite spectrale de \sisi{LERAY}{Leray} dont l'aboutissement est $\R^\ast(f\circ g)_\ast(\sisi{F}{{\Fa}})$ et dont le terme initial est donn\'e par:
\[
\E^{p, q}_2=R^pf_\ast(\R^qg_\ast(\sisi{F}{{\Fa}})).
\]
Par ailleurs, il existe un $\OX$-Module \emph{coh\'erent} ${\Ga}$ tel que ${\Ga}{}_{|U}\simeq \sisi{F}{{\Fa}}$, (\EGA I~9.4.3). Il r\'esulte alors du paragraphe pr\'ec\'edent que la condition (iv) \sisi{\ignorespaces}{de la page~\pageref{condition}} est v\'erifi\'ee\pageoriginale, \ie que $\R^qg_\ast({\Ga}{}_{|U})$ est coh\'erent pour $q<n$. On applique alors le th\'eor\`eme de finitude de \EGA III~3.2.1 \`a $f$ et aux faisceaux $\R^qg_\ast(\sisi{F}{{\Fa}})$, et on trouve que $\E^{p, q}_2$ est coh\'erent pour $q<n$, d'o\`u la conclusion.

\begin{proposition} \label{VIII.3.2} Soit $X$ un pr\'esch\'ema localement noeth\'erien et localement immergeable dans un pr\'esch\'ema r\'egulier. Soit $U$ un ouvert de $X$ et soit $i\colon U\to X$ l'immersion canonique. Soit $n\in\ZZ$. Soit enfin $\sisi{F}{{\Fa}}$ un $\Oo_U$-Module coh\'erent et de \emph{\sisi{COHEN-MACAULAY}{Cohen-Macaulay}} (sur $U$!). Les conditions suivantes sont \'equivalentes:
\begin{enumerate}
\item[(a)]
$\R^pi_\ast(\sisi{F}{{\Fa}})$ est coh\'erent pour $p<n$;

\item[(b)]
pour toute composante irr\'eductible $S'$ de l'adh\'erence $\overline S$ du support $S$ de $\sisi{F}{{\Fa}}$, on~a:
\[
\codim (S'\cap (X-U), S')>n.
\]
\end{enumerate}
\end{proposition}

Soit ${\Ga}$ un $\OX$-Module coh\'erent tel que ${\Ga}{}_{|U}\simeq \sisi{F}{{\Fa}}$ (\EGA I~9.4.3). Appliquant le corollaire \Ref{VIII.2.3} \`a ${\Ga}$, on trouve que la condition (a) \'equivaut \`a (c):
\begin{enumerate}
\item[(c)]
pour tout $x\in S$, on~a $\prof \sisi{F}{{\Fa}}_x>n-c(x)$, avec
\[c(x)=\codim(\overline{\{x\}}\cap Y, \overline{\{x\}}).
\]
\end{enumerate}

(a) \ALORS (b). En effet, soit $S'$ une composante irr\'eductible de $\overline S$ et soit $s$ son point g\'en\'erique. Puisque $\sisi{F}{{\Fa}}$ est de \sisi{COHEN-MACAULAY}{Cohen-Macaulay}, on~a $\prof \sisi{F}{{\Fa}}_s=\dim \Oo_{S, s}=0$. De plus, $\overline{\{s\}}=S'$, d'o\`u la conclusion.

(b) \ALORS (a). Soit $x\in S$ et soit $S'$ une composante irr\'eductible de $\overline{S}$ telle que $x\in S'$. Soit $s$ le point g\'en\'erique de $S'$. Puisque $\sisi{F}{{\Fa}}$ est de \sisi{COHEN-MACAULAY}{Cohen-Macaulay}, on sait que:
\[
\prof \sisi{F}{{\Fa}}_x=\dim \Oo_{\overline{\{s\}}, x}.
\]
Si\pageoriginale $c(x)=+\infty$, il n'y a rien \`a d\'emontrer. Sinon, il existe $y\in Y\cap \overline{\{x\}}$, tel que:
\[c(x)=\dim \Oo_{\overline{\{x\}}, y}.
\]
Or $\Oo_{X, y}$ est quotient d'un anneau local r\'egulier par hypoth\`ese, donc:
\[
\dim \Oo_{\overline{\{s\}}, y}=\dim \Oo_{\overline{\{s\}}, x} +\dim \Oo_{\overline{\{x\}}, y}>n.
\]
\begin{flushright}C.Q.F.D.\end{flushright}

\chapterspace{-1}
\chapter{G\'eom\'etrie alg\'ebrique et g\'eom\'etrie~formelle} \label{IX}

\newcommand{\Ib}{\Ii'}
\newcommand{\Jb}{\sisi{\mathcal J}{\mathcal{J}}}

Le but\pageoriginale de cet expos\'e est de g\'en\'eraliser au cas d'un morphisme qui n'est pas propre les th\'eor\`emes 3.4.2 et 4.1.5 de \EGA III.

\section{Le th\'eor\`eme de comparaison} \label{IX.1}

Soit $f\colon X\to X'$ un morphisme de pr\'esch\'emas, \emph{s\'epar\'e} et de \emph{type fini}. Supposons que $X'$ soit localement noeth\'erien. Soit $Y'$ une partie ferm\'ee de $X'$ et soit $Y\sisi{\subset}{=} f^{-1}(Y')$. Soient $\hat{X}$ et $\hat{X'}$ les compl\'et\'es formels de $X$ et de $X'$ le long de $Y$ et de $Y'$. Soit $\hat{f}$ le morphisme d\'eduit de $f$ par passage aux compl\'et\'es.
\begin{equation} \label{eq:IX.1.1}
\begin{array}{c}
\xymatrix{ X \ar[d]_f & \ar[l] Y \ar[d]& &\hat{X}\ar[d]_{\hat f} \ar[r]^j & X\ar[d]^f &\\
X' & \ar[l] Y' &, &\hat{X'} \ar[r]^i & X'.
}
\end{array}
\end{equation}

On d\'esigne par $j$ (\resp $i$) l'homomorphisme de $\hat{X}$ dans $X$ (\resp de $\hat{X'}$ dans $X'$). On sait que $i$ et $j$ sont \emph{plats}.

Soit $\Ii '$ un id\'eal de d\'efinition de $Y'$ et soit ${\Jb}=f^{*}(\Ii '). \OX$, c'est un id\'eal de d\'efinition de $Y$. On a donc:
\begin{equation} \label{eq:IX.1.2} {\hat{X'}=(Y', \varprojlim_{k \in \NN}\Oo_{X'}/{{\Ib}^{k+1}})}, \quad {\hat{X}=(Y, \varprojlim_{k \in \NN}\OX/{{\Jb}^{k+1}})}.
\end{equation}

Pour tout $k \in \NN$, posons:
\begin{equation} \label{eq:IX.1.3} Y'_k=(Y', \Oo_{X'}/{{\Ib}^{k+1}}), \quad {Y_k}=(Y, \OX/{{\Jb}^{k+1}}).
\end{equation}

Soit $F$ un $\OX$-Module \emph{coh\'erent}. Pour tout $k \in \NN$, on pose:
\begin{equation} \label{eq:IX.1.4} F_k=F/{\Jb}^{k+1}F, \qquad \hat{F}=j^*(F)\simeq \varprojlim F_k.
\end{equation}

Si l'on pose\pageoriginaled:
\begin{equation} \label{eq:IX.1.5} R^i f_*(F)^{\wedge}=\varprojlim_{k \in \NN} (R^i f_*(F)\sisi{^{\wedge}}{} \otimes_{{\Oo}_{X'}}{\Oo}_{Y'_k}), \quad i \in \ZZ,
\end{equation}
on a un homomorphisme naturel:
\begin{equation} \label{eq:IX.1.6} r_i: i^*(R^if_*(F)) \to R^i f_*(F)^{\wedge},
\end{equation}
qui est un \emph{isomorphisme} lorsque $ R^i f_*(F)$ est \emph{coh\'erent}.

Ainsi qu'il est expliqu\'e dans \EGA III 4.1.1, on~a un diagramme commutatif:
\begin{equation} \label{eq:IX.1.7}
\begin{array}{c}
\xymatrix{ i^*(R^i f_*(F)) \ar[rr]^-{\rho_i} \ar[d]^{r_i} & &R^i\hat{f}_*(\hat{F}) \ar[d]^{\psi_i} \\
R^i f_*(F)^{\wedge}\ar[rr]^-{\varphi_i} & &\varprojlim_{k \in \NN} R^i f_*(F_k).
}
\end{array}
\end{equation}

Dans \loccit on trouve un diagramme commutatif, car on sait que $ R^if_*(F)$ est coh\'erent, et l'on identifie la source et le but de (\Ref{eq:IX.1.6}). Dans notre cas, $ R^if_*(F)$ ne sera coh\'erent que pour certaines valeurs de $i$, pour lesquelles on \'etudiera (\Ref{eq:IX.1.7}).

Consid\'erons l'anneau gradu\'e
\begin{equation} \label{eq:IX.1.8}
\Ss = \sisi{\coprod}{\tbigoplus}_{k \in \NN} \Ii'^k,
\end{equation}
et le $\Ss$-Module gradu\'e:
\begin{equation} \label{eq:IX.1.9}
\Hh^i= \sisi{\coprod}{\tbigoplus}_{k \in \NN} R^if_*({\Jb}^k F), \quad i \in \ZZ,
\end{equation}
dont la structure de $\Ss$-Module est d\'efinie comme suit.

\sisi{Le faisceau\pageoriginale $R^if_*({\Jb}^k F)$ est associ\'e au pr\'efaisceau qui, \`a tout ouvert \emph{affine} $U$ de $X'$, associe:
\begin{equation} \label{eq:IX.1.10}
H^i(f^{-1}(U), {\Jb}^k F{|f^{-1}(U)}).
\end{equation}}
{Le faisceau\pageoriginale $R^if_*({\Jb}^k F)$ est associ\'e au pr\'efaisceau qui, \`a tout ouvert \emph{affine} $U'$ de~$X'$, associe:
\begin{equation} \label{eq:IX.1.10} H^i(f^{-1}(U'), {\Jb}^k F_{|f^{-1}(U')}).
\end{equation}
} \sisi{Soit donc $U'$ un ouvert affine de $X'$, posons
\begin{equation*} U=f^{-1}(U'),
\end{equation*}
et soit $x\in {\Ii}^m(U')$. Soit $x'$ l'image de $x$ dans ${\Jb}^{{k}} (U)$. L'homoth\'etie de rapport $x'$ dans $F_{|U'}$ applique ${\Jb}^kF{{|U'}}$ dans ${\Jb}^{k+m}F{{|U'}}$, d'o\`u, par fonctorialit\'e, un morphisme:
\begin{equation} \label{eq:IX.1.11} \mu^i_{x, k}(U):H^i(U', {\Jb}^kF{{|U'}})\to H^i(U', {\Jb}^{k+m}F{{|U'}}),
\end{equation}
d\'efini pour tout $i \in \ZZ$ et tout $k \in \NN$, qui donne, par passage au faisceau associ\'e, la structure de $\Ss$-Module gradu\'e de $\Hh^i$.} {Soit donc $U'$ un ouvert affine de $X'$, posons
\begin{equation*} U=f^{-1}(U'),
\end{equation*}
et soit $x'\in {\Ii'}^m(U')$. Soit $x$ l'image de $x'$ dans ${\Jb}^{{m}} (U)$. L'homoth\'etie de rapport $x$ dans $F_{|U}$ applique ${\Jb}^kF{_{|U}}$ dans ${\Jb}^{k+m}F{_{|U}}$, d'o\`u, par fonctorialit\'e, un morphisme:
\begin{equation} \label{eq:IX.1.11} \mu^i_{x', k}(U'):H^i(U, {\Jb}^kF{_{|U}})\to H^i(U, {\Jb}^{k+m}F{_{|U}}),
\end{equation}
d\'efini pour tout $i \in \ZZ$ et tout $k \in \NN$, qui donne, par passage au faisceau associ\'e, la structure de $\Ss$-Module gradu\'e de $\Hh^i$.}

\begin{theoreme} \label{IX.1.1} \setcounter{equation}{11}
Soit $n$ un entier. Supposons que le $\Ss$-Module gradu\'e $\Hh^i$ soit de type fini pour $i=n-1$ et $i=n$, alors:
\begin{enumerate}\setcounter{enumi}{-1}
\item
$r_n$ et $r_{n-1}$ sont des isomorphismes et $R^i\hat{f}_*(\hat{F})$ est coh\'erent pour $i=n-1$;
\item
pour $i=n-1$, $\rho_i$, $\varphi_i$ et $\psi_i$ sont des isomorphismes topologiques (en particulier la filtration d\'efinie sur $R^{n-1}f_*(F)$ par les noyaux des homomorphismes
\begin{equation} \label{eq:IX.1.12}
R^{n-1}f_*(F) \to R^{n-1}f_*(F_k)
\end{equation}
est $\Ii'$-bonne);

\item
pour $i=n$, $\rho_i$, $\varphi_i$ et $\psi_i$ sont des monomorphismes, de plus la filtration sur $R^n f_*(F)$ est $\Ii'$-bonne et $\psi_n$ est un isomorphisme;

\item
Le syst\`eme\pageoriginale projectif des $R^i f_*(F_k)$ v\'erifie, pour $i=n-2, n-1$, la condition de Mittag-Leffler uniforme, \ie il existe un entier $k\geq 0$ \emph{fixe} tel que, pour tout $p\geq 0$ et tout $p'\geq p+k$, on ait:
\begin{equation*}
\Im[R^i f_*(F_{p'})\to R^i f_*(F_p)]=\Im[R^i f_*(F_{p+k})\to R^i f_*(F_p)].
\end{equation*}
\end{enumerate}
\end{theoreme}

En proc\'edant comme dans \EGA III 4.1.8, il est facile de \emph{se ramener au cas o\`u} $X'$ \emph{est le spectre d'un anneau noeth\'erien} $A$. Dans ce cas, on sait que
\begin{equation} \label{eq:IX.1.13} R^if_*(F)=\sisi{H^i(X, F)^\thicksim}{\widetilde{H^i(X, F)}} \quad \text{(\cf \Ref{eq:IX.1.10}).}
\end{equation}

Soit $I$ l'id\'eal de $A$ tel que $\tilde{I}=\mathcal I'$ et soit
\begin{equation} \label{eq:IX.1.14}
S=\tbigoplus_{k \in \NN}I^k,
\end{equation}
\begin{equation} \label{eq:IX.1.15}
H^i=\tbigoplus_{k \in \NN}H^i(X, {\Jb}^k F), \quad i \in \ZZ,
\end{equation}
o\`u $H^i$ est muni de la structure de $S$-module gradu\'e d\'efinie par \Ref{eq:IX.1.11}, o\`u l'on a pris $U=X'$.

La d\'emonstration est calqu\'ee sur celle de \EGA III 4.1.5, donnons un r\'esum\'e.

On travaille sur $\varphi_i$ et $\psi_i$, auxquels correspondent des homomorphismes de modules:
\begin{equation} \label{eq:IX.1.16}
\begin{array}{c}
\xymatrix{ & H^i(\hat X, \hat F) \ar[d]^-{\psi_i} \\
H^i(X, F)^{\wedge}\ar[r]^-{\varphi_i}& \varprojlim_k H^i(X, F_k).
}
\end{array}
\end{equation}

\indent\textup{(a)} On suppose seulement que $H^i$ est un $S$-module gradu\'e de type fini. On en d\'eduit\pageoriginale que \emph{la filtration d\'efinie sur} $H^i(X, F)$ \emph{par les modules}:
\begin{equation} \label{eq:IX.1.17} R^i_k=\ker(H^i(X, F)\to H^i(X, F_k))
\end{equation}
\emph{est} $I$-\emph{bonne}. On utilise pour cela la suite exacte de cohomologie:
\begin{equation} \label{eq:IX.1.18} H^i(X, {\Jb}^{k+1}F)\to H^i(X, F)\to H^i(X, F_k),
\end{equation}
qui prouve que \emph{le} $S$-\emph{module gradu\'e} $\sisi{\coprod_{k \in \NN}}{\bigoplus_{k \in \NN}} R_k^i$ est \emph {quotient} du sous-$S$-module gradu\'e
\begin{equation*}
\sisi{\coprod}{\tbigoplus}_{k \in \NN}H^i(X, {\Jb}^{k+1}F)
\end{equation*}
de $H^i$, donc \emph{est de type fini}, car $S$ est noeth\'erien. D'o\`u ce premier point.

Posons:
\begin{equation} \label{eq:IX.1.19}
M^i=H^i(X, F), \quad H^i_k=H^i(X, F_k).
\end{equation}
\enlargethispage{\baselineskip}%
On a un diagramme commutatif:
\begin{equation} \label{eq:IX.1.20}
\begin{array}{c}
\xymatrix@R=6mm{ H^i(X, F)^{\wedge} \ar[r]^-{s_i}\ar[dr]_{\varphi_i} &\varprojlim_{k}(M^i/R^i_k)\ar[d]^{t_i} \\
&\varprojlim_{k}H^i_k\,,
}
\end{array}
\end{equation}
dans lequel $s_i$ est un \emph{isomorphisme}; en effet la filtration de $H^i(X, F)$ est $I$-bonne. De plus $t_i$ est un \emph{monomorphisme}; en effet le foncteur $\varprojlim$ est exact \`a gauche, et, pour tout $k\geq 0$, le morphisme naturel $M^i/R^i_k \to H^i_k$ est un monomorphisme, par d\'efinition de $R^i_k$.

Pour \'etudier la surjectivit\'e de $t_i$ on introduit:
\begin{equation} \label{eq:IX.1.21} Q_k^i = \coker (H^i(X, F) \to H^i(X, F_k)),
\end{equation}
d'o\`u\pageoriginale un syst\`eme projectif de suites exactes:
\begin{equation} \label{eq:IX.1.22} 0\to R_k^i \to M^i \to H^i_k \to Q_k^i \to 0.
\end{equation}
En utilisant la suite exacte de cohomologie:
\begin{equation} \label{eq:IX.1.23}
H^i(X, F)\to H^i(X, F_k)\to H^{i+1}(X, {\Jb}^{k+1}F),
\end{equation}
on voit que \emph{le $S$-module gradu\'e}
\begin{equation} \label{eq:IX.1.24}
Q^i=\sisi{\coprod_{k \in \NN}}{\tbigoplus_{k \in \NN}}Q^i_k
\end{equation}
\emph{est un sous-}$S$\emph{-module gradu\'e de} $H^{i+1}$. De plus, pour tout $k\geq 0$, on a:
\begin{equation} \label{eq:IX.1.25}
I^{k+1}Q^i_k=0
\end{equation}
car $Q_k^i$ est l'image de $H^i_k$.

\textup{(b)} \label{stepb} On suppose seulement que $H^{i+1}$ est de type fini et l'on s'int\'eresse \`a $t_i$ (en oubliant $s_i$). Puisque $S$ est noeth\'erien, $Q^i$ est de type fini; \sisi{en appliquant un lemme classique}{comme $I^{k+1}Q^i_k$ est nul}, on trouve qu'il existe un entier $r\geq 0$ et un entier $k_0\geq 0$ tels que
\begin{equation} \label{eq:IX.1.26}
I^rQ_k^i=0 \quad \text{pour } k\geq k_0.
\end{equation}
Il en r\'esulte que le \emph{syst\`eme projectif} $(Q_k^i)_{k\in \NN}$ \emph{est essentiellement nul}, donc que le \emph{syst\`eme projectif $(H_k^i)_{k\in \NN}$ v\'erifie la condition de Mittag-Leffler uniforme}. De la suite exacte (\Ref{eq:IX.1.22}) on d\'eduit la suite exacte
\begin{equation} \label{eq:IX.1.27}
0 \to M^i/R_k^i \to H^i_k \to Q_k^i \to 0,
\end{equation}
d'o\`u la suite exacte:
\begin{equation} \label{eq:IX.1.28}
0\to\varprojlim_{k}M^i/R_k^i \lto{t_i} \varprojlim_{k}H^i_k \to \varprojlim_{k} Q_k^i.
\end{equation}
Or\pageoriginale le syst\`eme projectif $(Q_k^i)_{k\in \NN}$ est essentiellement nul, donc $t_i$ est un \emph{isomorphisme}.

\indent\textup{(c)} Prouvons que, si $H^i$ est de type fini, $\psi_i$ est un isomorphisme. Il suffit d'appliquer $\EGA 0_{\textup{III}}$ 13.3.1 en prenant pour base d'ouverts de $X$ les ouverts affines. Ceci est licite; en effet, d'apr\`es (b), le syst\`eme projectif $(H_k^{i-1})_{k\in \NN}$ v\'erifie la condition de Mittag-Leffler.

Le th\'eor\`eme r\'esulte formellement de (a), (b) et (c). On remarquera qu'en fait la d\'emonstration n'utilise, \`a chaque pas, la finitude de $H^i$ que pour une seule valeur de~$i$.

Donnons quelques exemples o\`u l'hypoth\`ese du th\'eor\`eme \Ref{IX.1.1} est v\'erifi\'ee.

\begin{corollaire}\setcounter {equation}{28} \label{IX.1.2}
Supposons que $\Ii\sisi{}{'}$ soit engendr\'e par une section $t'$ de $\Oo_{X'}$, et notons~$t$ la section correspondante de $\OX$. Soit $F$ un $\OX$-module coh\'erent et soit $n$ un entier.

Supposons que:
\begin{enumeratei}
\item
$t$ soit $F$-r\'egulier (\ie l'homoth\'etie \sisi{$X\to tX$}{de rapport $t$} dans $F$ est un monomorphisme).
\item
$R^if_*(F)$ soit coh\'erent pour $i=n-1$ et $i=n$.
\end{enumeratei}
Alors l'hypoth\`ese du th\'eor\`eme \Ref{IX.1.1} est v\'erifi\'ee.
\end{corollaire}

En effet, on remarque que \sisi{\ignorespaces}{la multiplication par $t^k$} d\'efinit un isomorphisme \sisi{$\Jb^kF\simeq\Jb F$}{$ F\isomto\Jb^kF$} et on en d\'eduit que
\begin{equation} \label{eq:IX.1.29}
\Hh^i \simeq R^if_*(F)\otimes_{\Oo_{X'}}\Oo_{X'}[T],
\end{equation}
o\`u $T$ est une ind\'etermin\'ee. D'o\`u la conclusion.

\setcounter{equation}{29}
\refstepcounter{subsection}\label{IX.1.3}
\begin{enonce*}{Corollaire 1.3\ndemark}
Supposons\pageoriginale que $X'=\Spec(A)$, o\`u $A$ est un anneau noeth\'erien s\'epar\'e et complet pour la topologie $I$-adique. Supposons que le $S$-module $H^i$ soit de type fini pour $i=n-1$ et $i=n$. (\cf \Ref{eq:IX.1.14} et \Ref{eq:IX.1.15}). Alors les hypoth\`eses du th\'eor\`eme \Ref{IX.1.1} sont v\'erifi\'ees et on trouve un diagramme commutatif d'isomorphismes:
\begin{equation} \label{eq:IX.1.30}
\begin{array}{c}
\xymatrix{ H^i(X, F)\ar[rr]^{\rho'_i}\ar[dr]_{\varphi'_i} && H^i(\hat{X}, \hat{F}) \ar[dl]^{\psi_i}& \\
&\varprojlim_{k}H^i(X, F_k)&& \text{pour } i=n-1.}
\end{array}
\end{equation}
\ndetext{dans la m\^eme veine, voir l'article de Chow, (Chow~W.-L., {\og Formal functions on homogeneous spaces\fg}, \emph{Invent. Math.} \textbf{86} (1986), \numero 1, p\ptbl 115--130). L'auteur prouve le r\'esultat suivant. Soit $X$ une vari\'et\'e alg\'ebrique sur un corps, homog\`ene sous un groupe alg\'ebrique $G$ et soit $Z$ une sous-vari\'et\'e compl\`ete de $X$ de dimension $>0$. On suppose que $Z$ engendre $X$ au sens suivant: \'etant donn\'e $p\in Z$, soit $\Gamma_p$ l'ensemble des \'el\'ements de $G$ envoyant $p$ dans $Z$. On dit alors que $Z$ engendre si le groupe engendr\'e par la composante connexe de $1$ de $\Gamma_p$ est $G$ tout entier. Dans ce cas, toute fonction rationnelle formelle de $X$ le long de $Z$ est alg\'ebrique; \`a comparer avec les r\'esultats d'Hironaka et Matsumura cit\'es note de l'\'editeur~\eqref{noteHM}~page~\pageref{noteHM}. Dans la ligne des techniques introduites par ces auteurs, signalons le tr\`es joli r\'esultat d'alg\'ebrisation d\^u \`a Gieseker (Gieseker~D., {\og On two theorems of Griffiths about embeddings with ample normal bundle\fg}, \emph{Amer. J.~Math.} \textbf{99} (1977), \numero 6, p\ptbl 1137--1150, th\'eor\`emes 4.1 et 4.2). Soit $X$ une vari\'et\'e projective connexe de dimension $>0$ localement d'intersection compl\`ete (sur un corps alg\'ebriquement clos). On suppose qu'on a deux plongements de $X$ dans des vari\'et\'es projectives lisses $Y,W$. Alors, si les compl\'et\'es formels de $X$ dans $Y$ et $W$ sont \'equivalents, il existe un sch\'ema $U$ contenant $X$ (comme sous-sch\'ema ferm\'e) qui se plonge dans $Y$ et $W$ comme voisinage \'etale de $X$ dans $Y$ et $W$ Autrement dit, formellement \'equivalent entra\^ine \'etale-\'equivalent.
Voir aussi l'article de Faltings (Faltings~G., {\og Formale Geometrie and homogene Raüme\fg}, \emph{Invent. Math.} \textbf{64} (1981), p\ptbl 123-165).}
\end{enonce*}

On note simplement que $H^i(X, F)$ est de type fini donc isomorphe
\`a son compl\'et\'e. On obtient (\Ref{eq:IX.1.30}) en transcrivant dans
la cat\'egorie des $A$-modules le diagramme de Modules
(\Ref{eq:IX.1.7}), et en rempla\c cant la verticale de gauche par
$H^i(X, F)$.

\begin{proposition}\setcounter {equation}{30} \label{IX.1.4}
Soit $A$ un anneau noeth\'erien. Soit $t\in A$ et supposons que $A$ est s\'epar\'e et complet pour la topologie $(tA)$-adique. Posons:
\begin{equation} \label{eq:IX.1.31}
X'=\Spec (A), \quad Y'=V(t), \quad I=(tA).
\end{equation}
Soit $T$ une partie ferm\'ee de $X'$, posons
\begin{equation} \label{eq:IX.1.32}
X=X'-T, \quad Y=Y'\cap X=Y'-(Y'\cap T).
\end{equation}
Soit $F$ un $\OX$-Module coh\'erent. Soit enfin
\begin{equation} \label{eq:IX.1.33}
T'= \{ x\in X' \mid \codim(\{ \overline{x}\} \cap T, \{\overline{x}\})=1\}.
\end{equation}
Supposons que:
\begin{enumeratea}
\item
$t$ soit $F$-r\'egulier,
\item
$\prof_{T'}(F) \geq n+1$,
\item
$A$ soit quotient d'un anneau noeth\'erien r\'egulier.
\end{enumeratea}
\noindent
Alors, dans le diagramme (\Ref{eq:IX.1.30}), les morphismes $\rho'_i$, $\varphi'_i$ et $\psi_i$ sont des isomorphismes pour $i<n$ et des monomorphismes pour $i=n$. De plus $\psi_n$ est un isomorphisme.
\end{proposition}

En vertu\pageoriginale de \Ref{eq:IX.1.3} et \Ref{eq:IX.1.2}, il suffit de prouver que $R^if_*(F)$ est coh\'erent pour $i\leq n$, ce qui r\'esulte du th\'eor\`eme de finitude \Ref{VIII.2.1}.

En particulier:

\begin{exemple} \label{IX.1.5} On appliquera \Ref{eq:IX.1.4} lorsque $A$ est un anneau \emph{local} et que $t$ appartient au radical $\rr(A)$ de $A$. On prendra alors $T=\{ \rr(A)\}$. \emph{Dans ce cas}, pour $n=1$ on trouve l'\'enonc\'e suivant:

\emph{Si $A$ \sisi{\ignorespaces}{noeth\'erien} est s\'epar\'e et complet pour la topologie $t$-adique et quotient d'un anneau r\'egulier (par exemple si $A$ est complet), si de plus $t$ est $F$-r\'egulier et si $\prof F_x\geq 2$ pour tout $x\in \Spec (A)$ tel que $\dim A/x=1$, alors l'homomorphisme naturel \sisi{$$
\sisi{F}{\Gamma}(X, F)\to \Gamma(X, F)$$}{$$
\Gamma(X, F)\to\Gamma(\hat{X}, \hat{F})$$} est un isomorphisme.}

En effet, conservant les notations de \Ref{eq:IX.1.4}, on~a $T=\{ \rr(A)\}$ et la formule (\Ref{eq:IX.1.33}) dit que
$$
T'=\{x\in\Spec(A)\mid\dim A/x=1\}.
$$
\end{exemple}

\section{Th\'eor\`eme d'existence}

\'Enon\c cons d'abord \EGA III 3.4.2 sous une forme un peu plus
g\'en\'erale.

Soit $f:\formelX \to \formelX'$ un morphisme adique\sfootnote{Cette hypoth\`ese n'est pas essentielle, \cf \Ref{XII}, p\ptbl\pageref{XII.14}.} de pr\'esch\'emas formels, avec $\formelX'$ noeth\'erien. Soit $\Ii'$ un id\'eal de d\'efinition de $\formelX'$; puisque $f$ est adique, $f^*\Ii'=\Jj$ est\nde{par d\'efinition m\^eme, \cf \EGA I~10.12.1.} un id\'eal de d\'efinition de $\formelX$.

Pour tout $n\in \NN$, posons \label{IX.2}
\begin{equation} \label{eq:IX.2.1}
\formelX_n=(\formelX, \Oo_\formelX/{\Jb}^{n+1}),
\end{equation}
c'est un pr\'esch\'ema ordinaire qui a m\^eme espace topologique sous-jacent que $\formelX$.

Soit $\formelF$\pageoriginale un $\Oo_\formelX$-Module \emph{coh\'erent}. Pour tout $k\in\NN$, les $\Oo_{\formelX_k}$-Modules
\begin{equation} \label{eq:IX.2.2}
F_k=\formelF/{\Jb}^{k+1}\formelF
\end{equation}
sont \emph{coh\'erents}. Pour tout $i$, on~a un homomorphisme
\begin{equation} \label{eq:IX.2.3} \psi_i: R^if_*(\formelF) \to \varprojlim_k R^if_*(F_k),
\end{equation}
d\'eduit par fonctorialit\'e de l'homomorphisme naturel:
\begin{equation} \label{eq:IX.2.4} \formelF \to F_k.
\end{equation}

Posons
\begin{equation} \label{eq:IX.2.5}
\bS= \gr_{{\Ib}}\Oo_{\formelX'}=\sisi{\coprod_{k \in \NN}}{\tbigoplus_{k \in \NN}} {\Ib}^k/{\Ib}^{k+1},
\end{equation}
\begin{equation} \label{eq:IX.2.6} \sisi{\gr_{\Jj'}(\formelF)=\coprod_{k \in \NN} \Jj'^k\formelF/\Jj'^{k+1}\formelF}{\gr_{\Jb}(\formelF)=\tbigoplus_{k \in \NN} {\Jb}^k\formelF/{\Jb}^{k+1}\formelF},
\end{equation}
\begin{equation} \label{eq:IX.2.7} \KK^i=R^if_*(\gr_{{\Jb}}(\formelF))=\sisi{\coprod_{k \in \NN}}{\tbigoplus_{k \in \NN}}R^if_* ({\Jb}^k \formelF/{\Jb}^{k+1}\formelF\;).
\end{equation}
Il est clair que $\KK^i$ est muni d'une structure de $\bS$-Module gradu\'e.

\begin{theoreme} \label{IX.2.1}
Supposons que $\KK^i$ soit un $\bS$-Module gradu\'e \emph{de type fini} pour $i=n-1$, $i=n$, $i=n+1$, alors:
\begin{enumeratei}
\item
$R^nf_*(\formelF)$ est coh\'erent.
\item
L'homomorphisme $\psi_n$ (\Ref{eq:IX.2.3}) est un isomorphisme topologique. La filtration naturelle du deuxi\`eme membre de (\Ref{eq:IX.2.3}) est ${\sisi{\Jb}{\Ib}}$-bonne.
\item
Le syst\`eme projectif des $R^nf_*(F_k)$ v\'erifie la condition de Mittag-Leffler uniforme.
\end{enumeratei}
\end{theoreme}

\sisi{La d\'emonstration est tr\`es facile \`a partir de $\EGA 0_{\textup{III}} 13.7.7$ (\cf \EGA III $3.4.2$)\pageoriginale, \`a condition de rectifier comme suit le texte page 78\sfootnote{Comme indiqu\'e dans (EGA, $\rm{Err}_{\rm{III}}24$).}
\par
\emph{Supprimer}\begin{enumerate}
\item[ligne 6] \og et $(R^{n+1}T(A_k))_{k\in \ZZ}$\fg\
\item[ligne 13] \og et $p+q=n+1$\fg\
\item[ligne 23] \og et pour $n+1$\fg.
\end{enumerate}}{La d\'emonstration est tr\`es facile \`a partir de $\EGA 0_{\textup{III}} 13.7.7$ (\cf \EGA III $3.4.2$)\pageoriginale, \`a condition de rectifier le texte page 78 comme indiqu\'e\refstepcounter{toto}\label{pagerectif} dans (\EGA III~2, $\rm{Err}_{\rm{III}}24$).}

\setcounter{equation}{7}
\refstepcounter{subsection}\label{IX.2.2}
\begin{enonce*}{Th\'eor\`eme 2.2\ndemark}
Soit $A$ un anneau adique noeth\'erien et soit $I$ un id\'eal de d\'efinition de $A$. Soit $T$ une partie ferm\'ee de $X'=\Spec(A)$. Supposons que $I$ soit engendr\'e par un $t \in A$. Reprenons les notations \Ref{eq:IX.1.31}, \Ref{eq:IX.1.32} et \Ref{eq:IX.1.33}. Soit $\formelF$ un $\Oo_{\hat{X}}$-Module coh\'erent. Posons \ndetext{de nombreux \'enonc\'es d'alg\'ebrisation on \'et\'e obtenus depuis, sans parler de ceux cit\'es plus bas, \cf les articles de Faltings ou de Mme Raynaud cit\'es note de l'\'editeur~\eqref{HL}~p\ptbl\pageref{HL} et \eqref{MiR}~p\ptbl\pageref{MiR} respectivement. On pense notamment aux r\'esultats d'Artin (voir notamment Artin~M., {\og Algebraization of formal moduli. {I}\fg}, in \emph{Global Analysis (Papers in Honor of K\ptbl Kodaira)}, {Univ. Tokyo Press}, {Tokyo}, {1969}, p\ptbl 21-71), mais aussi aux r\'esultats r\'ecents d'alg\'ebricit\'e de feuilles de feuilletages; voir notamment Bost~J.-B., {\og Algebraic leaves of algebraic foliations over number fields\fg}, \emph{Publ. Math. Inst. Hautes \'Etudes Sci.} \textbf{93} (2001), p\ptbl 161-221 et Chambert-Loir~A. {\og Th\'eor\`emes d'alg\'ebricit\'e en g\'eom\'etrie diophantienne (d'apr\`es J.-B. Bost, Y. Andr\'e, D. \& G. Chudnovsky)\fg}, in \emph{S\'eminaire Bourbaki, Vol. 2000/2001},  Ast\'erisque~, vol.~282, Soci\'et\'e math\'ematique de France, Paris, 2002, \Exp 886, p\ptbl 175--209 et les r\'ef\'erences cit\'ees. En particulier, on trouvera dans ces deux articles des discussions sur le lien entre les questions d'alg\'ebrisation et la th\'eorie de l'approximation diophantienne.}
\begin{equation} \label{eq:IX.2.8}
F_0=\formelF / {\Jb} \formelF,
\end{equation}
o\`u ${\Jb}=t \Oo_{\hat{X}}$ est un id\'eal de d\'efinition de $\hat{X}$. \emph{Supposons que} $A$ soit quotient d'un anneau noeth\'erien r\'egulier et que:
\begin{enumerate}
\item
$t$ soit $\formelF$-r\'egulier,
\item
$\prof_{T'}\sisi{\formelF}{F}_0\geq 2$.
\end{enumerate}
\noindent
Alors il existe un $\OX$-Module coh\'erent $F$ tel que $\hat{\sisi{\formelF}{F}}\simeq \formelF$.
\end{enonce*}

Il suffit de prouver que $\hat{f}_*(\formelF)$ est un $\Oo_{\hat{X}}$-Module coh\'erent, o\`u $\hat{f}: \hat{X} \to \hat{X'}$ est le morphisme de pr\'esch\'emas formels d\'eduit de l'injection de $X$ dans $X'$ par compl\'etion par rapport \`a $t$. En effet, $A$ est s\'epar\'e et complet pour la topologie $t$-adique, il existera donc un $A$-module $F'$ dont le compl\'et\'e sera isomorphe \`a $\hat{f}_*(\formelF)$. Puisque $X$ est un ouvert de $X'$, on pourra prendre $F=\tilde{F'}_{|X}$.

Il reste \`a montrer que \Ref{IX.2.1} est applicable au morphisme de pr\'esch\'emas formels $\hat{f}$ et \`a $\formelF$. Or, d'apr\`es l'hypoth\`ese (1), pour tout $k\in \NN$ on~a un isomorphisme:
\begin{equation*} {\Jb}^k \formelF /{\Jb}^{k+1}\formelF \to \formelF / {\Jb} \formelF,
\end{equation*}
d'o\`u\pageoriginaled il r\'esulte que l'hypoth\`ese de \Ref{IX.2.1} sera v\'erifi\'ee si l'on sait que
\[
R^if_*(F_0)\text{ est coh\'erent pour }i \leq 1.
\]
Or ceci r\'esulte de (2) et du th\'eor\`eme de finitude \Ref{VIII.2.1}. D'o\`u la conclusion.

Il reste \`a sp\'ecialiser cet \'enonc\'e en supposant que $A$ est un anneau local.

\begin{corollaire}\setcounter{equation}{8} \label{IX.2.3}
Soit $A$ un anneau local noeth\'erien et soit $t \in \rr(A)$, un \'el\'ement du radical de $A$. Supposons que $A$ est s\'epar\'e et complet pour la topologie $t$-adique, et, de plus, quotient d'un anneau r\'egulier (par exemple, supposons que $A$ soit un anneau local noeth\'erien complet). Posons
\begin{equation} \label{eq:IX.2.9} X'=\Spec A, \quad T = \{ \rr(A) \},
\end{equation}
et reprenons les notations (\Ref{eq:IX.1.31}), (\Ref{eq:IX.1.32}) et (\Ref{eq:IX.1.33}). Soit $\formelF$ un $\Oo_{\hat{X}}$-Module. Supposons que:
\begin{enumerate}
\item
$t$ soit $\formelF$-r\'egulier,
\item
$\prof_{T'} \sisi{\formelF}{F}_0 \geq 2$, avec $\sisi{\formelF}{F}_0=\formelF/{\Jb}\formelF$ et ${\Jb}=t\Oo_{\hat{X}}$.

Alors il existe un $\OX$-Module coh\'erent $F$ tel que $\hat{F}\simeq \formelF$.
\end{enumerate}
\end{corollaire}

Remarquons qu'ici $T'$ est l'ensemble des id\'eaux premiers $\pp$ de $A$ tels que \hbox{$\dim A/\pp =1$}.

\chapterspace{-2}
\chapter{Application au groupe fondamental} \label{X}

\renewcommand{\theequation}{\arabic{equation}}

Dans\pageoriginale tout cet expos\'e, on d\'esignera par $X$ un pr\'esch\'ema localement noeth\'erien, par~$Y$ une partie ferm\'ee de $X$, par $U$ un voisinage ouvert variable de $Y$ dans $X$ et par~$\hat{X}$ le compl\'et\'e formel de $X$ le long de $Y$ (\EGA I 10.8). Pour tout pr\'esch\'ema $Z$, on d\'esignera par $\Et(Z)$ la cat\'egorie des rev\^etements \'etales de $Z$, et par $\LL(Z)$\refstepcounter{toto}\label{pX.1} la cat\'egorie des Modules coh\'erents localement libres sur $Z$.

\enlargethispage{3mm}%
\section{Comparaison de $\Et(\hat{X})$ et de $\Et(Y)$} \label{X.1}

Soit $\Ii$ un id\'eal de d\'efinition de $Y$ dans $X$. Posons, pour tout $n\in\NN$, $Y_n=\sisi{(Y, (\OX/\Ii^{n+1}){|Y}}{(Y, (\OX/\Ii^{n+1})_{|Y}})$. Les $Y_n$ forment un syst\`eme inductif de pr\'esch\'emas usuels, ou aussi de pr\'esch\'emas formels, en munissant les faisceaux structuraux de la topologie discr\`ete. On sait (\EGA I 10.6.2) que $\hat{X}$ est limite inductive, dans la cat\'egorie des pr\'esch\'emas formels, du syst\`eme inductif des $Y_n$. On sait aussi (\EGA I 10.13) que se donner un $\hat{X}$-pr\'esch\'ema formel de \emph{type fini} $R$, c'est se donner un syst\`eme inductif de $Y_n$-pr\'esch\'emas usuels $R_n$ de \emph{type fini}, tels que $R_n\simeq(R_{n+1})\times_{(Y_{n+1})}(Y_n)$. De plus, pour que $R$ soit un rev\^etement \'etale de $\hat{X}$, il est n\'ecessaire et suffisant que pour tout~$n$, $R_n$ soit un rev\^etement \'etale de $Y_n$. Ceci dit, il est facile de voir que les \'el\'ements nilpotents ne comptent pas pour les rev\^etements \'etales (\SGA \sisi{I}{1}~8.3), c'est-\`a-dire que le foncteur changement de base:
\[
\Et(Y_{n+1}) \to \Et(Y_n)
\]
est une \'equivalence de cat\'egories pour tout $n\in\NN$. Donc:

\begin{proposition} \label{X.1.1}
Avec les notations introduites plus haut, le foncteur naturel $\Et(\hat{X})\to\Et(Y)$ est une \'equivalence de cat\'egories (\cf \textup{\SGA \sisi{I}{1}~8.4}).
\end{proposition}

\section{Comparaison de $\Et(Y)$ et $\Et(U)$, pour $U$ variable} \label{X.2}

Nous\pageoriginale allons introduire deux conditions desquelles le th\'eor\`eme de comparaison annonc\'e r\'esultera facilement. Soit $X$ un pr\'esch\'ema localement noeth\'erien et soit $Y$ une partie ferm\'ee de $X$. On dit que le couple $(X, Y)$ v\'erifie la \emph{condition de Lefschetz}, ce qu'on \'ecrit $\Lef(X, Y)$,\refstepcounter{toto}\label{pX.2} si, pour tout ouvert $U$ de $X$ contenant $Y$ et tout faisceau coh\'erent localement libre $E$ sur $U$, l'homomorphisme naturel
\[
\Gamma(U, E) \to \Gamma(\hat{X}, \hat{E})
\]
est un isomorphisme.

On dit que le couple $(X, Y)$ v\'erifie la \emph{condition de Lefschetz effective}, ce qu'on \'ecrit $\Leff(X, Y)$, si on~a $\Lef(X, Y)$ et si de plus, pour tout faisceau coh\'erent localement libre~$\Ee$ sur $\hat{X}$, il existe un voisinage ouvert $U$ de $Y$ et un faisceau coh\'erent localement libre~$E$ sur $U$ et un isomorphisme $\hat{E}\simeq\Ee$.

Ces conditions sont v\'erifi\'ees dans deux exemples importants:

\refstepcounter{subsection}\label{X.2.1}
\begin{enonce*}[remark]{Exemple 2.1\ndemark}
\ndetext{\label{noteX.2.1}on peut l\'eg\`erement am\'eliorer (i): voir Mme Raynaud (Raynaud~M., {\og Th\'eor\`emes de Lefschetz en cohomologie des faisceaux coh\'erents et en cohomologie \'etale. Application au groupe fondamental\fg}, \emph{Ann. Sci. \'Ec. Norm. Sup. (4)} \textbf{7} (1974), p\ptbl 29--52, corollaires I.1.4 et I.5); la condition (ii) peut \^etre am\'elior\'ee pour se d\'ebarrasser des conditions de profondeur le long de $Y$ (voir le th\'eor\`eme 3.3 de \loccit pour un \'enonc\'e pr\'ecis). La preuve de ce dernier point est tr\`es technique, l'article ci-dessus ne donnant d'ailleurs que des indications de preuve, renvoyant \`a une version d\'etaill\'ee, ant\'erieure, parue au Bulletin de la Soci\'et\'e math\'ematique de France} Soit $A$ un anneau noeth\'erien et soit $t\in \rr(A)$ un \'el\'ement $A$-r\'egulier appartenant au radical $\rr(A)$ de $A$. Supposons que $A$ soit quotient d'un anneau local r\'egulier et que $A$ soit complet pour la topologie $t$-adique (par exemple $A$ complet pour la topologie $\rr(A)$-adique). Posons $X'=\Spec(A)$ et $Y'=V(t)$, par ailleurs posons $x=\rr(A)$ et $X=X'-\{x\}$, $Y=Y'-\{x\}$. Donc $X$ est ouvert dans $X'$ et $Y=X\cap Y'$. Alors:
\begin{enumeratei}
\item
Si, pour tout id\'eal premier $\pp$ de $A$ tel que $\dim A/\pp=1$ (\ie pour tout point ferm\'e de $X$) on~a $\prof A_{\pp}\geq 2$, alors on~a $\Lef(X, Y)$;

\item
si de plus, pour tout id\'eal premier $\pp$ de $A$ tel que $t\in\pp$ et $\dim A/\pp=1$ (\ie pour tout point ferm\'e de $Y$), on~a $\prof A_{\pp}\geq 3$, alors on~a $\Leff(X, Y)$.
\end{enumeratei}
\end{enonce*}

Montrons d'abord que pour tout voisinage ouvert $U$ de $Y$ dans $X$, le compl\'ementaire de $U$ dans $X$ est r\'eunion d'un nombre fini de points ferm\'es (dans $X$). Remarquons que $U$ ouvert dans $X$, donc dans $X'$, donc $Z'=X'-U$ est ferm\'e. Soit $I$\pageoriginale un id\'eal de d\'efinition de $Z'$; il suffit de prouver que $A/I$ est de dimension $1$. Or $Z'\cap Y'=\{x\}$ donc $A/(I+(t))$ est artinien, d'o\`u la conclusion gr\^ace au \og Hauptidealsatz\fg.

La \emph{premi\`ere hypoth\`ese \'equivaut \`a}: \og pour tout id\'eal premier $\pp$ de $A$, $\pp\neq \rr(A)$, on~a $\prof A_\pp\geq3-\dim A/\pp$\fg. En effet $A$ est quotient d'un anneau r\'egulier, on peut donc appliquer \Ref{VIII}\sisi{~\Ref{VIII.2.2}}{~\Ref{VIII.2.3}} au pr\'esch\'ema $X'$, \`a la partie ferm\'ee $\{x\}$ et au faisceau coh\'erent $\Oo_{X'}$\sisi{.}{, en observant que, $c(\pp)=\dim(A/\pp)$ pour $\pp\in U=X'-\{x\}$ (car $x$ est le point ferm\'e de $X'$).}

Soit $U$ un voisinage ouvert de $Y$ dans $X$ et $E$ un $\Oo_U$-module localement libre. Posons $Z=X-U$ et soit $u\colon U\to X$ l'immersion canonique. On va d'abord prouver que $u_\ast(E)$ est un $\OX$-Module \emph{coh\'erent}, ou, ce qui revient au m\^eme, que $\SheafH^i_Z(E')$ est coh\'erent pour $i=0, 1$, o\`u $E'$ est un prolongement coh\'erent de $E$ \`a $X$. Pour cela on applique le th\'eor\`eme \Ref{VIII}~\Ref{VIII.2.1} au pr\'esch\'ema $X$, \`a la partie ferm\'ee $Z$ et au faisceau coh\'erent $E'$. Il suffit de v\'erifier que pour tout point $p\in U$ tel que $c(p)=1$, on~a $\prof E'_p\geq 1$, o\`u l'on a pos\'e
\[
c(p)=\codim (\overline{\{p\}}\cap Z, \overline{\{p\}}).
\]
Or si $p\in U$ et si $c(p)=1$, notant $\pp$ l'id\'eal de $A$ correspondant \`a $p$, on voit que $\dim A/\pp=2$, car le compl\'ementaire de $U$ est r\'eunion d'un nombre fini de points ferm\'es et $A$ est quotient d'un anneau r\'egulier. De plus, $E$ est localement libre, donc, pour tout $p\in\Supp E$, on~a $\prof E_p=\prof\Oo_{U, p}$. Enfin, si $p\in U$ et si $c(p)=1$, on a
\[
\prof E'_p = \prof E_p = \prof\Oo_{U, p} = \prof A_\pp \geq 3-2=1.
\]

Il nous faut maintenant prouver que l'homomorphisme naturel
\begin{equation} \label{eq:X.1} {\Gamma(U, E) \to \Gamma(\hat{X}, \hat{E})}
\end{equation}
est un isomorphisme. Posant alors $\overline{\!E}=u_\ast(E)$, on note que $\overline{\!E}$ est coh\'erent et de profondeur $\geq 2$ en tous les points ferm\'es de $X$. Il en r\'esulte que $\R^i f_\ast(\overline{\!E})$ est coh\'erent pour $i=0, 1$, o\`u $f\colon X\to X'$ d\'esigne l'immersion canonique de $X$ dans $X'=\Spec(A)$ (\Exp \Ref{VIII}). On applique alors (\Ref{IX}~\Ref{IX.1.5}), et l'on conclut que\pageoriginale
\begin{equation} \label{eq:X.2}
\Gamma(U, \overline{\!E}) \to \Gamma(\hat{X}, \widehat{\overline{\!E}})
\end{equation}
est un isomorphisme, car $A$ est complet pour la topologie $t$-adique.

On a un diagramme commutatif:
\[
\UseTips \newdir{ >}{!/-5pt/\dir{>}} \xymatrix@=5mm{ \Gamma(X, \overline{\!E})\ar[rr]^{\simeq}\ar[dr]_{\simeq} && \Gamma(U, E)\ar[dl]\\
&\Gamma(\hat{X}, \hat{E})}
\]
d'o\`u la conclusion.

Soit maintenant $\Ee$ un faisceau coh\'erent localement libre sur $\hat{X}$. Si l'on a prouv\'e que~$\Ee$ est alg\'ebrisable, \ie est isomorphe au compl\'et\'e formel d'un $\OX$-Module coh\'erent~$E$, il est facile de voir que $E$ est localement libre au voisinage de $Y$, donc de prouver $\Leff (X, Y)$. Soit $\widehat{X'}$ le spectre formel de $A$ pour la topologie $t$-adique, qui s'identifie au compl\'et\'e formel de $X'$ le long de $Y'$. D\'esignons par $f$ l'immersion canonique de~$X$ dans $X'$, par $f'$ l'immersion canonique de $Y$ dans $Y'$, et par $\hat{f}$ le morphisme d\'eduit par passage aux compl\'et\'es. Pour que $\Ee$ soit alg\'ebrisable il suffit que $\hat{f}_\ast(\Ee)$ soit un $\Oo_{\hat{X}}$-Module coh\'erent, car $A$ est complet pour la topologie t-adique. Soit $\Ii = t\Oo_{\hat{X}}$, c'est un id\'eal de d\'efinition de $\hat{X}$.

Pour tout $n\geq 0$, posons $\E_n=\Ee/\Ii^{n+1}\sisi{}{\Ee}$. En tout point ferm\'e $y\in Y$, la profondeur de $\E_0$ est $\geq 2$; en effet $t$ est un \'el\'ement $A$-r\'egulier donc $\prof\Oo_{Y_0, y}=\prof\Oo_{X, y}-1\geq 2$. On en conclut que $\hat{f}_\ast(\Ee)$ est coh\'erent (\Ref{IX}~\Ref{IX.2.3}). \cqfd

\refstepcounter{subsection}\label{X.2.2}
\begin{enonce*}[remark]{Exemple 2.2}[{{\normalfont Permettra de comparer le groupe fondamental d'une vari\'et\'e projective et d'une section hyperplane}}]
Soit $K$ un corps et soit $X$ un $K$-pr\'esch\'ema propre. Soit $\Ll$ un $\OX$-Module inversible ample. Soit $t \in \Gamma(X, \Ll)$ un \'el\'ement $\OX$-r\'egulier\pageoriginale, ce qui signifie que, pour tout ouvert $U$ et tout isomorphisme $u\colon \Ll\sisi{~U}{_{|U}}\to\Oo_U$, \/$u(t)$ est non diviseur de z\'ero dans $\Oo_U$ (condition qui ne d\'epend pas de $u$). Soit $Y = V(t)$ le sous-sch\'ema de $X$ d'\'equation $t = 0$.\nde{la condition (ii) est superflue, voir note~page~\pageref{X.2.1}.} Alors:
\begin{enumeratei}
\item
Si, pour tout point $x$ ferm\'e dans $X$, on~a $\prof\Oo_{X, x}\geq 2$, on~a $\Lef (X, Y)$;

\item
si de plus, pour tout point ferm\'e $y\in Y$, on~a $\prof\Oo_{X, y}\geq 3$, on~a $\Leff (X, Y)$.
\end{enumeratei}
\end{enonce*}

Cet exemple sera trait\'e en d\'etail dans \Exp \Ref{XII}.

Soit S un pr\'esch\'ema, on sait (\EGA II 6.1.2) que le foncteur qui, \`a tout rev\^etement fini et plat $r\colon R\to S$, associe la $\OX$-Alg\`ebre $r_\ast(\Oo_R)$, induit une \'equivalence entre la cat\'egorie des rev\^etements finis et plats de $S$ et la cat\'egorie des $\OX$-Alg\`ebres coh\'erentes et localement libres. Soit U un voisinage ouvert de $Y$, et soit $r\colon R\to U$ un rev\^etement fini et plat de U. Soit $\hat{R}$ le rev\^etement fini et plat de $\hat{X}$ qui s'en d\'eduit par changement de base. On a $\hat{r}_\ast(\Oo_{\hat{R}})\simeq\widehat{r_\ast(\Oo_R)}$.

\emph{Supposons alors} $\Lef(X, Y)$. Ceci entra\^ine que, pour tout $U$, le foncteur image inverse:
\[
\LL(U) \to \LL(\hat{X})
\]
est \emph{pleinement fid\`ele}. En effet, soient $E$ et $F$ deux $\Oo_U$-Modules coh\'erents localement libres; $\SheafHom(E, F)$ est aussi coh\'erent et localement libre; $\SheafHom(E, F)$ est aussi coh\'erent et localement libre. Par hypoth\`ese, l'application naturelle
\[
\Gamma_U (\SheafHom(E, F)) \to \Gamma_{\hat{X}}(\widehat{\SheafHom(E, F)})
\]
est un isomorphisme, d'o\`u la conclusion, car $\SheafHom$ commute au $\, \, \hat{}\, \, $ puisque tout est localement libre. Or le $\, \, \hat{}\, \, $ commute au produit tensoriel, on en d\'eduit que le foncteur qui \`a toute $\Oo_U$-Alg\`ebre coh\'erente et localement libre $\Aa$ associe la $\Oo_{\hat{X}}$-Alg\`ebre $\hat{\Aa}$, est pleinement fid\`ele. Mieux, si $E$ est un $\Oo_U$-Module coh\'erent localement libre\pageoriginale, il y a correspondance biunivoque entre les structures de $\Oo_{\hat{X}}$-Alg\`ebre commutative sur $\hat{E}$.

\begin{proposition} \label{X.2.3}
Soit $X$ un pr\'esch\'ema localement noeth\'erien et soit $Y$ une partie ferm\'ee de $X$. Soit $\hat{X}$ le compl\'et\'e formel de $X$ le long de $Y$. Pour tout ouvert $U$ de $X$, $U\supset Y$, d\'esignons par $L_U$ (\resp $P_U$, \resp $E_U$) le foncteur qui \`a tout $\Oo_U$-Module coh\'erent localement libre (\resp \`a tout rev\^etement fini et plat de $U$, \resp \`a tout rev\^etement \'etale de $U$) associe son image inverse par $\hat{X} \to X$.
\begin{enumeratei}
\item
Si on~a $\Lef(X, Y)$, alors pour tout voisinage ouvert $U$ de $Y$, les foncteurs $L_U$, $P_U$ et $E_U$ sont pleinement fid\`eles.

\item
Si on~a $\Leff(X, Y)$, alors pour tout $\Oo_{\hat{X}}$-Module coh\'erent localement libre~$\Ee$ (\resp ...), il existe un ouvert $U$ et un $\Oo_U$-Module coh\'erent localement libre~$E$ (\resp ...), tels que $L_U(E) \simeq \Ee$ (\resp ...).
\end{enumeratei}
\end{proposition}

(i) A \'et\'e vu.

(ii) R\'esulte de (i) et de l'hypoth\`ese, du moins pour $L_U$ et $P_U$. De plus si $R$ est un rev\^etement \emph{\'etale} de $\hat{X}$, il existe un voisinage ouvert $U$ de $Y$ dans $X$ et un rev\^etement fini et plat $R'$ de $U$ tel que $\widehat{R'} \simeq R$. On en d\'eduit un rev\^etement $R''$ de $Y$ qui est \'etale d'apr\`es \Ref{X.1.1}, donc $R'$ est \'etale dans un voisinage $U'$ de $Y$. \cqfd

\begin{corollaire} \label{X.2.4}
Si on~a $\Lef (X, Y)$, pour qu'un rev\^etement fini et plat $R$ d'un voisinage ouvert $U$ de $Y$ soit connexe, il faut et il suffit que $R \times_U \hat{X}$ le soit. En particulier, pour que $Y$ soit connexe, il faut et il suffit que le voisinage ouvert $U$ de $Y$ le soit, ou encore que $X$ le soit.
\end{corollaire}

En effet, pour qu'un espace annel\'e en anneaux locaux $(X, \OX)$ soit connexe, il faut et il suffit que $\Gamma(X, \OX)$ ne soit pas compos\'e direct de deux\pageoriginale anneaux non nuls. Or on~a
\[
\Gamma(U, r_\ast(\Oo_R)) \simeq \Gamma(\hat{X}, \hat{r}_\ast(\Oo_{\hat{R}}))
\]
par $\Lef(X, Y)$.

\begin{corollaire} \label{X.2.5} Si on~a $\Lef (X, Y)$, alors pour tout $U$, le foncteur
\[
\Et(U) \to \Et(Y)
\]
est pleinement fid\`ele. Si on~a $\Leff(X, Y)$, alors pour tout rev\^etement \'etale $R$ de $Y$, il existe un voisinage ouvert $U$ de $Y$ et un rev\^etement $R'$ de $U$ tel que $R' \times_U Y \simeq R$.
\end{corollaire}

\refstepcounter{subsection}\label{X.2.6}
\begin{enonce*}{Corollaire 2.6\ndemark}
Si on~a $\Lef(X, Y)$ et si $Y$ est connexe, tout voisinage ouvert $U$ de $Y$ est connexe et l'homomorphisme naturel $\pi_1(Y)\to\pi_1(U)$ est surjectif. Si de plus on~a $\Leff(X, Y)$, l'homomorphisme naturel
\[
\pi_1(Y) \to \varprojlim_{U}\, \pi_1(U)
\]
est un isomorphisme. (N.B. On suppose choisi un \og point-base\fg dans $Y$, qu'on prend aussi pour point-base dans $X$, pour la d\'efinition des groupes fondamentaux.)
\end{enonce*}

Tout ceci r\'esulte trivialement de la prop\ptbl \Ref{X.1.1} et de la prop\ptbl \Ref{X.2.3}.
\ndetext{joint \`a~\Ref{X.3.3} et aux crit\`eres~\Ref{XII.2.4} et~\Ref{XII.3.4}, on obtient le th\'eor\`eme de Lefschetz relatif suivant. Soit $f:X\to S$ un morphisme projectif et plat de sch\'emas noeth\'eriens connexes et soit $D$ un diviseur de Cartier relatif effectif dans $X$ et relativement ample. Si, pour tout $s\in S$, la profondeur de $X_s$ en chaque point ferm\'e est $\geq 2$, alors $D$ est connexe et, pour tout ouvert $U$ de $X$ contenant $D$, la fl\`eche $i_U:\pi_1(D)\to\pi_1(U)$ est surjective. Si de plus, la profondeur de $X_s$ le long de chaque point ferm\'e de $D_s$ est $\geq 3$ et si les anneaux locaux de $X$ en ses points ferm\'es sont purs, par exemple d'intersection compl\`ete (\cf \Ref{X}~\Ref{X.3.4}), alors $i_X$ est un isomorphisme. \Cf Bost~J.-B., {\og Lefschetz theorem for Arithmetic Surfaces\fg}, \emph{Ann. Sci. \'Ec. Norm. Sup. (4)} \textbf{32} (1999), p\ptbl 241-312, th\'eor\`emes 1.1 et 2.1. Dans le cas o\`u $X$ est simplement une surface projective lisse et g\'eom\'etriquement connexe sur un corps, on~a toujours connexit\'e de $D$ et surjectivit\'e de $\pi_1(D)\to\pi_1(U)$ (o\`u $U$ ouvert contenant~$D$) pour $D$ seulement nef de carr\'e $>0$ (\cf \loccit, th\'eor\`eme 2.3 et aussi le th\'eor\`eme 2.4 pour des surfaces seulement normales et compl\`etes). Dans le cas d'une surface arithm\'etique (normale et quasi-projective) $X$ \sisi{au dessus}{au-dessus} d'un anneau d'entiers $\Oo_K$, Bost, am\'eliorant des r\'esultats de Ihara (Ihara~Y., {\og Horizontal divisors on arithmetic surfaces associated with Bely\u\i\ uniformizations\fg}, in \emph{The Grothendieck theory of dessins d'enfants (Luminy, 1993)}, London Math. Soc. Lect. Note Series, vol.~200, Cambridge Univ. Press, Cambridge, 1994, 245--254 ou \loccit, corollaire 7.2), a montr\'e que si un point $P\in X(\Oo_K)$, qui joue le r\^{o}le du diviseur $D$ dans la situation g\'eom\'etrique, v\'erifie certaines conditions de positivit\'e, alors la fl\`eche $\pi_1(X)\to\pi_1(\Spec \Oo_K)$ d\'eduite de la projection \'etait inversible d'inverse la fl\`eche $\pi_1(\Spec \Oo_K)\to\pi_1(X)$ d\'eduite de $P$ (\loccit, th\'eor\`eme 1.2).}

\section{Comparaison de $\pi_1(X)$ et de $\pi_1(U)$} \label{X.3}

\begin{definition} \label{X.3.1}
Soit $X$ un pr\'esch\'ema et $Z$ une partie ferm\'ee de $X$. Posons $U=X-Z$. On dit que le couple $(X, Z)$ est pur si, pour tout ouvert $V$ de $X$, le foncteur
\[
\begin{aligned} \Et(V) &\to \Et(V\cap U) \\
V' &\sisi{\rightsquigarrow}{\mto} V'\times_V(V\cap U)
\end{aligned}
\]
est une \'equivalence de cat\'egories\sfootnote{Pour une notion plus satisfaisante \`a certains \'egards, \cf le commentaire \Ref{XIV}~\Ref{XIV.1.6} d).}.
\end{definition}

\begin{definition} \label{X.3.2}
Soit\pageoriginale $A$ un anneau local noeth\'erien. Posons $X=\Spec A$. Soit $\rr(A)$ le radical de $A$ et soit $x=\rr(A)$ le point ferm\'e de $X$. On dit que $A$ est pur si le couple $(X, \{x\})$ l'est.
\end{definition}

Nous laissons au lecteur le soin de ne pas d\'emontrer la proposition suivante:

\begin{proposition} \label{X.3.3}
Soit $X$ un pr\'esch\'ema localement noeth\'erien et soit $Z$ une partie ferm\'e de $X$. Pour que le couple $(X, Z)$ soit pur il est n\'ecessaire et suffisant que, pour tout $z\in Z$, l'anneau $\Oo_{X, z}$ soit pur\sfootnote{Comparer le cas non commutatif de \Ref{XIV}~\Ref{XIV.1.8}, dont la d\'emonstration est essentiellement la m\^eme que celle de \Ref{X.3.3}.}.
\end{proposition}

Ceci dit le th\'eor\`eme suivant est le r\'esultat essentiel de ce num\'ero:

\refstepcounter{subsection}\label{X.3.4}
\begin{enonce*}{Th\'eor\`eme 3.4}[Th\'eor\`eme de puret\'e\ndemark]
\ndetext{pour l'historique des m\'ethodes employ\'ees, voir la lettre du 1\ier octobre 1961 de Grothendieck \`a Serre, \emph{Correspondance Grothendieck-Serre}, \'edit\'ee par Pierre Colmez et Jean-Pierre Serre, Documents Math\'ematiques, vol.~2, Soci\'et\'e Math\'ematique de France, Paris, 2001.}
\textup{(i)}\enspace
Un anneau local noeth\'erien r\'egulier de dimension $\geq 2$ est pur (th\'eor\`eme de puret\'e de Zariski-Nagata).

\textup{(ii)}\enspace
Un anneau local noeth\'erien de dimension $\geq 3$ qui est une intersection compl\`ete est pur.
\end{enonce*}

Rappelons qu'on dit qu'un anneau local est une \emph{intersection compl\`ete} s'il existe un anneau local noeth\'erien \emph{r\'egulier} $B$ et une suite $(t_1, \ldots, t_k)$ $B$-r\'eguli\`ere d'\'el\'ements du radical $\rr(B)$ de $B$ tels que
\[
A \simeq B/(t_1, \ldots, t_k).
\]

\sisi{A}À ce propos remarquons qu'il serait moins ambigu de dire que $A$ est une intersection compl\`ete \emph{absolue}, par opposition \`a la situation, que nous avons d\'ej\`a rencontr\'ee, o\`u $X$ est un pr\'esch\'ema localement noeth\'erien (qui n'a pas besoin d'\^etre r\'egulier) et o\`u $Y$ est une partie ferm\'ee de $X$, dont on dit qu'elle est \og localement ensemblistement une intersection compl\`ete \emph{dans $X$}\fg.

Prouvons d'abord quelques lemmes.

\begin{lemme} \label{X.3.5}
Soit\pageoriginale $X$ un pr\'esch\'ema localement noeth\'erien et soit $U$ une partie ouverte de $X$. Posons $Z=X-U$. Soit $i\colon U\to X$ l'immersion canonique de $U$ dans $X$. Les conditions suivantes sont \'equivalentes:
\begin{enumeratei}
\item
Pour tout ouvert $V$ de $X$, si on pose $V'=V\cap U$, le foncteur $F\sisi{\rightsquigarrow}{\mto} \sisi{F|{|V'}}{F_{|V'}}$ de la cat\'egorie des $\Oo_V$-Modules coh\'erents localement libres dans la cat\'egorie des $\Oo_{{V'}}$-Modules coh\'erents localement libres est pleinement fid\`ele;

\item
l'homomorphisme naturel $\OX \to i_\ast(\Oo_{U})$ est un isomorphisme;

\item
pour tout $z\in Z$, on~a $\prof\Oo_{X, z}\geq 2$.
\end{enumeratei}
\end{lemme}

On a d\'ej\`a vu (\Ref{III}~\Ref{III.3.3}) l'\'equivalence de (ii) et (iii). Montrons que (ii) entra\^ine (i). Soient $F$ et $G$ deux $\Oo_V$-Modules coh\'erents localement libres, $\SheafHom(F, G)$ l'est aussi, donc $\SheafHom(F, G)\to i_\ast(\SheafHom(\sisi{F|{|V'}}{F_{|V'}}, \sisi{G|{|V'}}{G_{|V'}}))$ est un isomorphisme donc $\Hom(F, G)\simeq\Hom(\sisi{F|{|V'}}{F_{|V'}}, \sisi{G|{|V'}}{G_{|V'}})$. Inversement on prend $F=G=\OX$ et on applique (i) \`a tout ouvert~$V$ de~$X$.

Voici un \og lemme de descente\fg utile:

\begin{lemme} \label{X.3.6} Soit $X$ un pr\'esch\'ema localement noeth\'erien et soit $Z$ une partie ferm\'ee de $X$. Posons $U=X-Z$. Supposons que l'homomorphisme $\OX \to i_\ast(\Oo_U)$ soit un isomorphisme. Soit $f\colon X_1\to X$ un morphisme fid\`element plat et quasi-compact. Posons $Z_1=f^{-1}(Z)$. Si le couple $(X_1, Z_1)$ est pur, il en est de m\^eme de $(X, Z)$.
\end{lemme}

Remarquons que l'hypoth\`ese $\OX\simeq i_\ast(\Oo_U)$ se conserve par extension plate de la base, car $i$ est un morphisme quasi-compact et, dans ce cas, l'image directe commute \`a l'image inverse. Or cette hypoth\`ese entra\^ine que le foncteur
\[
\Et(V) \to \Et(U\cap V)
\]
d\'efini par
\[
V'\sisi{\rightsquigarrow}{\mto} V'\times_V(V\cap U)
\]
est pleinement fid\`ele, comme le montre l'interpr\'etation d'un rev\^etement \'etale en termes\pageoriginale d'Alg\`ebre coh\'erente localement libre. Il reste \`a prouver l'effectivit\'e. On peut par exemple introduire le carr\'e $X_2$ et le cube $X_3$ de $X_1$ sur $X$ et on remarque qu'un morphisme fid\`element plat et quasi compact est un morphisme de \emph{descente effective universelle} pour la cat\'egorie fibr\'ee des rev\^etements \'etales, au-dessus de la cat\'egorie des pr\'esch\'emas. La conclusion est formelle \`a partir de l\`a\sfootnote{\Cf J.~Giraud, \sisi{M\'ethode de la descente}{\emph{M\'ethode de la descente}}, M\'emoire \numero 2 du Bulletin de la Soci\'et\'e Math\'ematique de France (1964).}.

\begin{remarque} \label{X.3.7}
On a prouv\'e chemin faisant que si $\OX\to i_\ast(\Oo_U)$ est un isomorphisme, $X$ est connexe si et seulement si $U$ l'est, et alors $\pi_1(U)\to\pi_1(X)$ est surjectif.
\end{remarque}

\begin{corollaire} \label{X.3.8}
Soit $A$ un anneau local noeth\'erien. Supposons que $\prof A\geq 2$. Alors si $\hat{A}$ est pur, $A$ est pur.
\end{corollaire}

R\'esulte du lemme \Ref{X.3.5} et du lemme \Ref{X.3.6}.

Le lemme suivant est le point essentiel de la d\'emonstration du th\'eor\`eme de puret\'e:

\begin{lemme} \label{X.3.9}
Soit $A$ un anneau local noeth\'erien et soit $t\in\rr(A)$ un \'el\'ement $A$-r\'egulier. Supposons que $A$ soit complet pour la topologie $t$-adique et de plus quotient d'un anneau local r\'egulier (par exemple $A$ complet). Posons $B=A/tA$.
\begin{enumeratei}
\item
Si pour tout id\'eal premier $\pp$ de $A$ tel que $\dim A/\pp=1$, on~a $\prof A_{\pp}\geq 2$, alors $B$ pur entra\^ine $A$ pur.

\item
Si pour tout id\'eal premier $\pp$ de $A$ tel que $\dim A/\pp=1$, on~a $\prof A_{\pp}\geq 2$, si $A_\pp$ pur lorsque $t\notin\pp$, et \sisi{\ignorespaces}{si}\nde{cette derni\`ere condition peut \^etre am\'elior\'ee, \cf la note de l'\'editeur~\eqref{noteX.2.1} de la page~\pageref{noteX.2.1}.} $\prof A_\pp\geq3$ lorsque $t\in\pp$, alors $A$ pur entra\^ine $B$ pur.
\end{enumeratei}
\end{lemme}

Soit $X'=\Spec(A)$ et soit $Y'=V(t)$, que l'on identifie au spectre de $B$. Soit $x=\rr(A)$, et posons $X=X'-\{x\}$ et $Y=Y'-\{x\}=X\cap Y'$. D\'esignons par $\sisi{X'}{\widehat{X'}}$ le spectre formel de $A$ pour la topologie $t$-adique, qui s'identifie au compl\'et\'e formel de $X'$ le long de $Y'$.

Puisque $A$ est complet pour la topologie $t$-adique, on remarque que $\Et(X')\to\Et(\widehat{X'})$\pageoriginale est une \'equivalence de cat\'egories. De m\^eme $\Et(\widehat{X'})\to\Et(Y')$ par la prop\ptbl \Ref{X.1.1}, donc $\Et(X')\to\Et(Y')$ est une \'equivalence de cat\'egories.

Montrons (i). Consid\'erons les diagrammes
\[
\UseTips \newdir{ >}{!/-5pt/\dir{>}} \xymatrix{ X' & X\ar[l] && \Et(X')\ar[r]^a \ar[d]_c & \Et(X)\ar[d]^b\\
Y'\ar[u] & Y\ar[l] \ar[u] && \Et(Y')\ar[r]^d & \Et(Y)}
\]
On vient de voir que $c$ est une \'equivalence, $d$ l'est aussi d'apr\`es l'hypoth\`ese que $B$ est pur et enfin $b$ est pleinement fid\`ele comme on l'a vu dans l'exemple \Ref{X.2.1}, \cf \Ref{X.2.3} (i).

Montrons (ii). Cette fois on suppose que $A$ est pur donc $a$ est une \'equivalence; de m\^eme $c$. Voyons que $b$ est une \'equivalence. D'apr\`es l'exemple \Ref{X.2.1} on sait que l'on~a $\Leff(X, Y)$, donc $b$ est d\'ej\`a pleinement fid\`ele, prouvons qu'il est essentiellement surjectif. On utilise \Ref{X.2.3} (ii) en notant que, si $U$ est un voisinage ouvert de $Y$ dans $X$, le compl\'ementaire de $U$ dans $X$ est une r\'eunion d'un nombre fini de points ferm\'es; le couple $(X,X-U)$ est donc pur d'apr\`es la prop\ptbl \Ref{X.3.3}, car en un tel point $\pp$, $\Oo_{X,\pp}=A_\pp$ est pur par hypoth\`ese. D'o\`u la conclusion.

\begin{proof}[D\'emonstration du th\'eor\`eme de puret\'e]
Montrons d'abord (i) par r\'ecurrence sur la dimension. Soit $A$ un anneau local noeth\'erien de \emph{dimension} $2$. Posons $X'=\Spec(A)$, $x=\rr(A)$, $X=X'-\{x\}$. On a $\prof A=2$. On peut donc appliquer le lemme \Ref{X.3.5} au couple $(X', \{x\})$ et donc $\Et(X')\to\Et(X)$ est pleinement fid\`ele. Soit maintenant $r\colon R\to X$ un rev\^etement \'etale d\'efini par une $\OX$-Alg\`ebre coh\'erente localement libre et \'etale $\Aa=r_\ast(\Oo_R)$. D\'esignons par $i\colon X\to X'$ l'immersion canonique de $X$ dans $X'$. Je dis que $i_\ast(\Aa)=\Bb$ est une $\Oo_{X'}$-Alg\`ebre \emph{coh\'erente}. En effet, il suffit d'appliquer le \og th\'eor\`eme de finitude\fg \Ref{VIII}~\Ref{VIII.2.3}. Je dis que cette alg\`ebre est de profondeur $\geq 2$ en $x$. En effet, c'est\pageoriginale l'image directe d'un $\OX$-Module, avec $X=X'-\{x\}$. Puisque $A$ est un anneau \emph{r\'egulier} de dimension $2$ on a: $\dimp \Bb + \prof\Bb=\dim A=2$, o\`u $\dimp\Bb$ d\'esigne la dimension projective de $\Bb$. Donc $\dimp\Bb=0$, donc $\Bb$ est projectif, donc libre. Il en r\'esulte que $\Bb$ d\'efinit un rev\^etement fini et plat de $X'=\Spec(A)$. L'ensemble des points de $X'$ o\`u ce rev\^etement n'est pas \'etale est une partie ferm\'ee de $X'$ dont l'\'equation est un id\'eal principal: l'id\'eal discriminant de $\Bb/A$. Or, par construction, ce ferm\'e est contenu dans $x=\rr(A)$ donc est vide car $\dim A=2$.

Soit $A$ un anneau local noeth\'erien r\'egulier, $\dim A=n\geq3$. Supposons (i) d\'emontr\'e pour les anneaux de dimension $<n$. Pour prouver que $A$ est pur, on peut supposer $A$ complet par \Ref{X.3.8}. Soit $t\in\rr(A)$ dont l'image dans $\rr(A)/\rr(A)^2$ soit non nulle. Alors $B=A/tA$ est un anneau local noeth\'erien \emph{r\'egulier} de dimension $n-1$ donc est pur, car $n-1\geq2$. On conclut par le lemme \Ref{X.3.9}~(i), qui est applicable car $A$ est complet.

\emph{Montrons} (ii). Soit $A$ un anneau local noeth\'erien de dimension $\geq 3$. Supposons qu'il existe un anneau local noeth\'erien \emph{r\'egulier} $B$ et une $B$-suite $(t_1, \ldots, t_k)$ tels que $A\simeq B/(t_1, \ldots, t_k)$. Pouvons que $A$ est pur, par r\'ecurrence sur $k$. Si $k=0$, on le sait par (i). \emph{Supposons $k\geq1$ et le r\'esultat acquis pour $k'<k$.} D'apr\`es le corollaire \Ref{X.3.9} on peut supposer que $A$ (donc aussi $B$) est complet. Posons $C=B/(t_1, \ldots, t_{k-1})$, donc $A\simeq C/t_kC$ et $t_k$ est $C$-r\'egulier. Par l'hypoth\`ese de r\'ecurrence on sait que $C$ est pur, il suffit de prouver que le lemme \Ref{X.3.9}~(ii) est applicable. Notation: $A$ et $B$ du lemme deviennent $C$ et $A$. On a $\dim C\geq4$, donc pour tout id\'eal premier $\pp$ de $C$ tel que $\dim C/\pp=1$, on~a $\prof C_\pp\geq3$. De plus, $C_\pp$ est un intersection compl\`ete avec $k'\leq k-1$, donc est pur par l'hypoth\`ese de r\'ecurrence.
\end{proof}

\begin{theoreme} \label{X.3.10}
Soit\pageoriginale $X$ un pr\'esch\'ema localement noeth\'erien et soit $Y$ une partie ferm\'ee de $X$. Supposons que l'on ait $\Leff(X, Y)$ (\cf Exemples \Ref{X.2.1} et \Ref{X.2.2}). Supposons de plus que, pour tout voisinage ouvert $U$ de $Y$ et tout $x\in X-U$, l'anneau local $\Oo_{X, x}$ soit r\'egulier de dimension $\geq2$ ou bien une intersection compl\`ete de dimension $\geq3$. Alors
\[
\pi_0(Y) \to \pi_0(X)
\]
est une bijection, et si $X$ est connexe
\[
\pi_1(Y) \to \pi_1(X)
\]
est un isomorphisme.
\end{theoreme}

Il n'y a plus rien \`a d\'emontrer. On remarque que, dans les deux exemples cit\'es \Ref{X.2.1} et \Ref{X.2.2}, le compl\'ementaire de $U$ est une r\'eunion d'un nombre fini de points ferm\'es, d'o\`u il r\'esulte que l'hypoth\`ese sur la dimension de $\Oo_{X, x}$ n'est pas canularesque.

\numberwithin{equation}{section}

\chapter{Application au groupe de Picard} \label{XI}

Cet expos\'e est calqu\'e sur le pr\'ec\'edent, mais cette fois le r\'esultat du \numero \Ref{XI.1} est moins fort.

Dans\pageoriginale tout cet expos\'e, $X$ d\'esignera un pr\'esch\'ema localement noeth\'erien, $\Ii$ un id\'eal quasi-coh\'erent de $\OX$, ($Y=V(\Ii)$ est donc une partie ferm\'ee de $X$), $U$ un voisinage ouvert variable de $Y$ dans $X$, et $\hat X$ le compl\'et\'e formel de $X$ le long de $Y$. Pour tout espace annel\'e $(Z, \Oo_Z)$, on d\'esigne par $\PP(Z)$\refstepcounter{toto}\label{pXI.1} la cat\'egorie des $\Oo_Z$-Modules inversibles, autrement dit localement libres de rang~$1$, et par $\Pic(Z)$ le groupe des classes \`a isomorphisme pr\`es de Modules inversibles sur $Z$.

\section{Comparaison de $\Pic(\hat X)$ et de $\Pic(Y)$} \label{XI.1}

Pour tout $n\in\NN$, posons $X_n = (Y, \OX/\Ii^{n+1})$ et $P_n = \Ii^{n+1}/\Ii^{n+2}$. La suite de faisceaux de groupes ab\'eliens sur $Y$
\begin{equation} \label{eq:XI.1.1} 0\to P_n\lto{u}\Oo_{X_{n+1}}^*\lto{v}\Oo_{X_n}^*\to 1
\end{equation}
est \emph{exacte}. Pr\'ecisons que la structure de groupe de $P_n$ est la structure additive, que $u(x)=1+x$ pour tout $x\in P_n$, et que $v$ est l'homomorphisme d\'eduit de l'injection $\Ii^{n+2}\to\Ii^{n+1}$. On voit que $v$ est surjectif en remarquant que, pour tout $y \in Y$, $\Oo_{X_n, y}$ est un anneau local, quotient de $\Oo_{X_{n+1}, y}$ par un id\'eal nilpotent; le reste est tout aussi trivial. On d\'eduit de (\Ref{eq:XI.1.1}) une suite exacte de cohomologie:
\begin{equation*} \label{eq:XI.1.*} \tag{$*$}
\H^1(Y, P_n)\lto{u^1}\H^1(Y, \Oo_{X_{n+1}}^*)\lto{v^1} \H^1(Y, \Oo_{X_n}^*)\lto{d}\H^2(Y, P_n){\quoi}.
\end{equation*}
Par ailleurs, pour tout $n\in\NN$, on sait identifier $\Pic(X_n)$ et $\H^1(Y, \Oo_{X_n}^*)$; de plus, si~$E$ est un $\Oo_{X_{n+1}}$-Module inversible, correspondant \`a une classe de cohomologie $c(E)$\pageoriginale, la classe de cohomologie correspondant \`a l'image inverse de $E$ sur $X_n$ est \'egale \`a $v^1(c(E))$. D'o\`u la proposition suivante:

\begin{proposition} \label{XI.1.1} Conservons les notations introduites ci-dessus. Soit $p\in\NN$. L'application $\Pic(\hat X)\to \Pic(Y_n)$:
\begin{enumeratei}
\item
est injective pour $n\geq p$, si $\H^1(Y, P_n)=0$ pour $n\geq p$;
\item
est un isomorphisme pour $n\geq p$, si $\H^i(Y, P_n)=0$ pour $n\geq p$ et $i=1, 2$.
\end{enumeratei}
\end{proposition}

Bien entendu, la suite exacte~\eqref{eq:XI.1.*} contient plus d'information que la proposition ci-dessus. Le lecteur aura remarqu\'e que l'on n'a rien dit du foncteur $\PP(\hat X)\to\PP(Y)$. \'Etant donn\'es deux $\Oo_{\hat X}$-Modules inversibles\sisi{}{ $E, F$}, $H=\SheafHom(E, F)$ est aussi inversible. Si l'on indique par un indice $n$ la r\'eduction modulo $\Ii^{n+1}$, on trouve une suite exacte:
\[
0\to H_0\otimes P_n \to \SheafHom(E_{n+1}, F_{n+1}) \to \SheafHom(E_n, F_n)\to 0{\quoi}.
\]
D'o\`u une suite exacte de cohomologie que nous n'\'ecrirons pas et dont l'interpr\'etation est \'evidente; on peut utiliser cette remarque pour \'etudier le foncteur~$\PP$.

\section{Comparaison de $\Pic(X)$ et de $\Pic(\hat X)$} \label{XI.2}

Le lecteur trouvera dans l'expos\'e~\Ref{X}, \numero \Ref{X.2}, la d\'emonstration de ce qui suit:

\begin{proposition} \label{XI.2.1} Supposons que l'on ait $\Lef(X, Y)$; alors pour tout voisinage ouvert $U$ de $Y$ dans $X$, le foncteur
\begin{equation} \label{eq:XI.2.1} \PP(U)\to\PP(\hat X)
\end{equation}
est pleinement fid\`ele, l'application
\begin{equation} \label{eq:XI.2.2} \Pic(U)\to\Pic(\hat X)
\end{equation}
est\pageoriginale donc injective. Si l'on a $\Leff(X, Y)$, l'application~(\Ref{eq:XI.2.3}) est un isomorphisme:
\begin{equation} \label{eq:XI.2.3}
\varinjlim_U\Pic(U)\to\Pic(\hat X){\quoi}.
\end{equation}
\end{proposition}

\begin{corollaire} \label{XI.2.2} Supposons que l'on ait $\Lef(X, Y)$ et que pour tout entier $n\geq p$, on ait $\H^1(Y, P_n)=0$; alors pour tout ouvert $U\supset Y$, les applications
\[
\Pic(X)\to\Pic(U)\to\Pic(Y_n)
\]
sont injectives si $n\geq p$. Si l'on a $\Leff(X, Y)$ et si de plus, pour tout entier $n\geq p$, on~a $\H^i(Y, P_n)=0$ pour $i=1$ et $i=2$, alors l'application
\[
\varinjlim_U\Pic(U)\to\Pic(Y_n)
\]
est un isomorphisme pour $n\geq p$.
\end{corollaire}

\section{Comparaison de $\PP(X)$ et de $\PP(U)$} \label{XI.3}

Une d\'efinition:
\refstepcounter{subsection}\label{XI.3.1}
\begin{enonce*}[remark]{D\'efinition 3.1\sfootnotemark}\sfootnotetext{Pour une \'etude plus d\'etaill\'ee de la notion de parafactorialit\'e, et la d\'emonstration de~\Ref{XI.3.3}, \cf \EGA IV 21.13, 21.14.}
Soit $X$ un pr\'esch\'ema et soit $Z$ une partie ferm\'ee de~$X$. Posons $U=X-Z$. On dit que $X$ est \emph{parafactoriel} aux points de $Z$ si, pour tout ouvert $V$ de~$X$, le foncteur $\PP(V)\to \PP(V\cap U)$ est une \'equivalence de cat\'egories. On dit aussi que le couple $(X, Z)$ est parafactoriel.
\end{enonce*}

Rappelons que $\PP(Z)$ d\'esigne la cat\'egorie des Modules localement libres de rang~$1$ sur~$Z$.

\begin{definition} \label{XI.3.2}
Un anneau local noeth\'erien est dit \emph{parafactoriel} si le couple $(\Spec(A), \{\rr(A)\})$ est parafactoriel.
\end{definition}

On d\'emontre la proposition suivante qui prouve que la notion est \og ponctuelle\fg:

\begin{proposition} \label{XI.3.3}
Supposons\pageoriginale $X$ localement noeth\'erien. Pour que le couple $(X, Z)$ soit parafactoriel, il faut et il suffit que, pour tout $z\in Z$, l'anneau local $\Oo_{X, z}$ le soit.
\end{proposition}

Remarquons que dans parafactoriel il y a \og pleinement fid\`ele\fg. On d\'emontre comme dans le lemme~\Ref{X.3.5} de l'expos\'e~\Ref{X} le:

\begin{lemme} \label{XI.3.4}
Si $X$ est un pr\'esch\'ema localement noeth\'erien et si $Z=X-U$ est une partie ferm\'ee de $X$, les conditions suivantes sont \'equivalentes:
\begin{enumeratei}
\item
pour tout ouvert $V$ de $X$, le foncteur $\PP(V)\to \PP(V\cap U)$ est pleinement fid\`ele;
\item
l'homomorphisme $\OX \to i_*(\Oo_U)$ est un isomorphisme;
\item
pour tout $z\in Z$, on~a $\prof(\Oo_{X, z})\geq 2$.
\end{enumeratei}
\end{lemme}

Donc \og parafactoriel\fg signifie que les conditions de~\Ref{XI.3.4} sont v\'erifi\'ees et que, pour tout ouvert $V$ de $X$, l'homomorphisme $\Pic(V)\to\Pic(V\cap U)$, est surjectif. En particulier, si $X$ est le spectre d'un anneau local noeth\'erien, on trouve:

\begin{proposition} \label{XI.3.5}
Soit $A$ un anneau local noeth\'erien; pour qu'il soit parafactoriel il est n\'ecessaire et suffisant que $\prof A\geq 2$ et $\Pic(X'-\{x\})=0$, o\`u l'on a pos\'e $X'=\Spec(A)$ et o\`u $x$ est l'unique point ferm\'e de $X'$.
\end{proposition}

On notera qu'un anneau local de dimension $\leq 1$ n'est jamais parafactoriel car sa profondeur est $\leq 1$. Donc \emph{factoriel n'entra\^ine pas parafactoriel}; c'est cependant vrai pour les anneaux locaux noeth\'eriens de dimension $\geq 2$, comme nous le verrons plus bas.

\begin{lemme} \label{XI.3.6} Soit $X$ un pr\'esch\'ema localement noeth\'erien et soit $Z$ une partie ferm\'ee de $X$. Soit $f\colon X_1\to X$ un morphisme fid\`element plat et quasi-compact. Posons $Z_1=f^{-1}(Z)$. Si $(X_1, Z_1)$ est parafactoriel, alors $(X, Z)$ l'est.
\end{lemme}

On remarque\pageoriginale d'abord que, si $i\colon(X-Z)\to X$ d\'esigne l'immersion canonique de $U=X-Z$ dans $X$, la formation de l'image directe par $i$ d'un $\Oo_U$-Module quasi-coh\'erent commute au changement de base $f$, car celui-ci est plat. Il est donc \'equivalent de supposer les conditions \'equivalentes du lemme~\Ref{XI.3.5} pour $(X, Z)$ ou pour $(X_1, Z_1)$, car $f$ est un morphisme de descente pour la cat\'egorie des faisceaux quasi-coh\'erents. Il reste \`a prouver que, pour tout ouvert $V$ de $X$, $\Pic(V)\to\Pic(V\cap U)$ est surjectif. On fait le changement de base $V\to X$, qui ne change rien (\sisi{sic}{\textit{sic}}), et on est ramen\'e \`a $V=X$. On remarque alors que, si ${\Lcal}$ est un $\Oo_U$-Module inversible et si ${\Lcal}$ admet un plongement localement libre, ce prolongement est isomorphe \`a $i_*({\Lcal})$, \`a cause de ce qui vient d'\^etre vu. Reste \`a prouver que $i_*({\Lcal})$ est inversible. Utilisant \`a nouveau le fait que l'image directe par $i$ commute au changement de base plat, et que \og localement libre de rang~$1$\fg est une propri\'et\'e qui se descend par morphisme fid\`element plat et quasi-compact, on~a gagn\'e.

\begin{corollaire} \label{XI.3.7} Soit $A$ un anneau local noeth\'erien; si $\hat A$ est parafactoriel $A$ l'est.
\end{corollaire}

Ne pas croire que, si $A$ est parafactoriel, $\hat A$ l'est\nde{pour une \'etude pr\'ecise du lien entre la factorialit\'e de $A$ et de son compl\'et\'e, voir (Heitmann~R., {\og Characterization of completions of unique factorization domains\fg}, \emph{Trans. Amer. Math. Soc.} \textbf{337} (1993), \numero 1, p\ptbl 379--387).}.

\danger

Avant d'aborder l'\'enonc\'e du th\'eor\`eme principal de ce \no, faisons le lien avec la th\'eorie des diviseurs et la notion d'anneau factoriel\sfootnote{Pour les g\'en\'eralit\'es qui suivent, \cf aussi \EGA IV~21.}.

Soit $X$ un pr\'esch\'ema \emph{noeth\'erien} et \emph{normal}. Soit $Z^1(X)$ le groupe ab\'elien libre engendr\'e par les $x\in X$ tels que $\dim\OXx=1$. L'anneau local d'un tel point est un anneau de valuation discr\`ete. On notera $v_x$ la valuation norm\'ee correspondante. Soit $K(X)$ l'anneau des fractions rationnelles de $X$ et soit
\[
p\colon K(X)^*\to Z^1(X)
\]
l'application qui \`a tout $f\in K(X)^*$ associe le cycle 1-codimensionnel:
\[
(f) = \sum_{x\in X, \;\dim\OXx=1} v_x(f)\cdot x{\quoi}.
\]
L'image\pageoriginale de $p$ est not\'ee $P(X)$ et on appelle ses \'el\'ements \emph{diviseurs principaux}\sfootnote{Conform\'ement \`a la terminologie de \EGA IV~21, nous pr\'ef\'erons maintenant r\'eserver le nom de \og diviseurs\fg aux \og diviseurs localement principaux\fg ou \og diviseurs de Cartier\fg.}. On pose
\[
\Cl(X) = Z^1(X)/P(X){\quoi}.
\]
Soit $\sisi{Z}{Z^{\prime1}(X)}$ le sous-groupe de $Z^1(X)$ dont les \'el\'ements sont les \emph{diviseurs localement principaux}. On sait que
\[
\Pic(X)\simeq Z^{\prime1}(X)/P(X){\quoi},
\]
et par suite $\Pic(X)$ s'identifie \`a un sous-groupe de $\Cl(X)$.

Remarquons que si $U$ est un ouvert dense de $X$, $K(X)\to K(U)$ est un isomorphisme, et que si $\codim(X-U, X)\geq 2$, \ie si tout $x\in X$ tel que $\dim\OXx\leq 1$ appartient \`a $U$, l'homomorphisme $Z^1(X)\to Z^1(U)$ et par suite $\Cl(X)\to\Cl(U)$ est aussi un isomorphisme. Enfin, si tout $x\in U$ est factoriel, \ie $\OXx$ l'est, alors $Z^1(U) = Z^{\prime1}(U)$ donc $\Pic(U)\simeq\Cl(U)$.

\begin{propositionblah} \label{XI.3.7.1}
Soit $X$ un pr\'esch\'ema noeth\'erien et normal. Soit $(U_i)_{i\in I}$ une famille d'ouverts de $X$ telle que:
\begin{enumeratea}
\item
les $U_i$ forment une base de filtre\nde{\ie une famille filtrante d\'ecroissante.};
\item
si on pose $Y_i=X-U_i$, on~a $\codim(Y_i, X)\geq 2$ pour tout $i$,
\item
si $x\in U_i$ pour tout $i\in I$, alors $\OXx$ est factoriel.
\end{enumeratea}
\noindent
Alors on~a un isomorphisme:
\[
\varinjlim_{i\in I}\Pic(U_i)\lto{\approx}\Cl(X){\quoi}.
\]
\end{propositionblah} Remarquons que b)~entra\^ine que tout $x\in X$ tel que $\dim\OXx\leq 1$ appartient \`a $U_i$ pour tout $i$. Donc les $U_i$ sont denses et de plus l'homomorphisme $Z^1(U_i)\to Z^1(X)$ est un isomorphisme, de m\^eme que $K(U_i)\to K(X)$. Donc $\Pic(U)\subset \Cl(U_i) \simeq \Cl(X)$. Pour prouver ce que l'on d\'esire, il suffit donc de montrer que tout $D\in Z^1(X)$ appartient \`a $Z^{\prime1}(U_i)$ pour un $i$ convenable. Il suffit de le faire pour les \og diviseurs\fg irr\'eductibles et positifs. Soit donc $x\in X$ tel\pageoriginaled que $\dim\OXx=1$. Il suffit de prouver qu'il existe un $i\in I$ tel que $\overline{\{x\}}$ soit localement principal en les points de $U_i$. Soit $\Ii$ le plus grand id\'eal de d\'efinition du ferm\'e $\overline{\{x\}}$. L'ensemble des points au voisinage desquels $\Ii$ est libre est un ouvert $U$. Or $U\supset \bigcap_{i\in I}U_i$ d'apr\`es~c). Si on pose $Y=X-U$, on~a $Y\subset \bigcup_{i\in I}Y_i$\sisi{}{ avec $Y_i=X-U_i$}, or $Y$ est ferm\'e, donc admet un nombre \emph{fini} de points g\'en\'eriques, donc est contenu dans la r\'eunion d'un nombre fini de $Y_i$, donc dans $Y_j$ pour un $j\in I$, car les $U_i$ forment une base de filtre. Donc $U\supset U_j$.\qed

\begin{corollaire} \label{XI.3.8}
Soit $X$ un pr\'esch\'ema noeth\'erien et normal et soit $Y$ une partie ferm\'ee de codimension $\geq 2$. Supposons que, pour tout $p\in X-Y$, $\Oo_{X, p}$ soit factoriel, alors
\[
\Pic(X-Y)\to\Cl(X-Y)\to\Cl(X)
\]
sont des isomorphismes.
\end{corollaire}

\begin{corollaire} \label{XI.3.9}
Soit $A$ un anneau local noeth\'erien et normal. Posons $X'=\Spec(A)$ et $x = \rr(A)$. Pour que $A$ soit \emph{factoriel} il faut et il suffit que $\Pic(X'-\{x\})=0$ et que $\pp\in X'-\{x\}$ implique que $A_\pp$ est factoriel.
\end{corollaire}
En effet, pour que $A$ soit factoriel il est n\'ecessaire et suffisant que $\Cl(X')=0$\,\nde{voir Bourbaki, \emph{Alg\`ebre commutative} VII.1.4, cor\ptbl du th\ptbl 2 et VII.3.2, th\ptbl 1}.

\begin{corollaire} \label{XI.3.10}
Soit $A$ un anneau local noeth\'erien de dimension $\geq 2$. Posons $X'=\Spec(A)$ et soit $x = \rr(A)$. Posons $X=X'-\{x\}$. Les conditions suivantes sont \'equivalentes:
\begin{enumeratei}
\item
$A$ est factoriel;
\item
\textup{a)}\enspace pour tout $y\in X$, $\Oo_{X, p}$ est factoriel, et
\par\hspace*{6.3mm}\textup{b)}\enspace $A$ est parafactoriel, \ie $\prof A\geq 2$ et $\Pic(X)=0$.
\end{enumeratei}
\end{corollaire}

Avant\pageoriginale de d\'emontrer ce corollaire, \'enon\c cons le
\refstepcounter{subsection}\label{XI.3.11}
\begin{enonce*}{Crit\`ere de normalit\'e de Serre 3.11\sfootnotemark}\sfootnotetext{\Cf \EGA IV~5.8.6.} Soit $A$ un anneau local noeth\'erien. Pour que $A$ soit normal, il et n\'ecessaire et suffisant que
\begin{enumeratei}
\item
pour tout id\'eal premier ${\pp}$ de $A$ tel que $\dim A_{{\pp}}\leq 1$, $A_{{\pp}}$ soit normal,
\item
pour tout id\'eal premier ${\pp}$ de $A$ tel que $\dim A_{{\pp}}\geq 2$, on ait $\prof A_{{\pp}}\geq 2$.
\end{enumeratei}
\end{enonce*}

Prouvons~\Ref{XI.3.10}.

(i)\ALORS (ii). Sachant qu'un localis\'e de factoriel l'est aussi, on~a (ii)~a). De plus $A$ est normal donc $\prof A\geq 2$ car $\dim A\geq 2$ (\Ref{XI.3.11}~(ii)). Enfin $A$ est parafactoriel; en effet $\Pic(X)\simeq\Cl(X')=0$ (\cf \Ref{XI.3.9}).

(ii)\ALORS (i). Prouvons d'abord que $A$ est \emph{normal} en appliquant le crit\`ere de \sisi{SERRE}{Serre}. Puisque $\dim A\geq 2$, la condition~(i) du crit\`ere est dans les hypoth\`eses. De plus, pour tout ${\pp}\in X$, $A_{{\pp}}$ est factoriel, donc normal, donc de profondeur $\geq 2$, du moins si $\dim A_{{\pp}}\geq 2$. Enfin $\prof A\geq 2$ d'apr\`es (ii)~b). Il reste \`a appliquer~\Ref{XI.3.9}.

\smallskip
R\'esumons ce qui pr\'ec\`ede:

\begin{proposition} \label{XI.3.12}
Soit $X$ un pr\'esch\'ema localement noeth\'erien et soit $\Ii$ un id\'eal quasi-coh\'erent de $X$. Soit $Y=V(\Ii)$. Soit $p\in\NN$. Supposons que:
\begin{enumerate}
\item[1)] On \sisi{a}{ait} $\Leff(X, Y)$ (\textup{Expos\'e} \Ref{X});
\item[2)] $\H^i(X, \Ii^{n+1}/\Ii^{n+2})=0$ si $i=1$ ou $2$ et si $n\geq p$;
\item[3)] pour tout voisinage ouvert $U$ de $Y$ dans $X$ et tout $x\in X-U$, l'anneau $\OXx$ soit parafactoriel.
\end{enumerate}
\noindent
Alors, pour tout $n\geq p$, et tout voisinage ouvert $U$ de $Y$, les homomorphismes
\[
\xymatrix{\Pic(X)\ar[rr]\ar[dr]&&\Pic(X_n)\\&\Pic(U)\ar[ur]&}
\]
sont des isomorphismes.
\end{proposition}

On\pageoriginale conna\^it des anneaux parafactoriels:

\skpt
\begin{theoreme} \label{XI.3.13}
\begin{enumeratei}
\item
(Auslander-Buchsbaum)\nde{\`a comparer avec le r\'esultat de puret\'e suivant, d\^u \`a Gabber. Soit $X$ le spectre d'un anneau $A$ local r\'egulier de dimension $3$, $a$ de diff\'erentielle non nulle, \ie $a\in\mm-\mm^2$ et $U$ le compl\'ementaire de $V(a)$. Alors, un fibr\'e vectoriel sur $U$ est libre (pour une d\'emonstration simple, voir (Swan~R.G., {\og A simple proof of Gabber's theorem on projective modules over a localized local ring\fg} \emph{Proc. Amer. Math. Soc.} \textbf{103} (1988), \numero 4, p\ptbl 1025-1030). Le cas du rang $1$ est un cas particulier du th\'eor\`eme~\Ref{XI.3.13}. Pour des r\'esultats de puret\'e concernant des fibr\'es vectoriels de rang arbitraire, que ce soit dans le cadre analytique ou alg\'ebrique, voir (Gabber~O., {\og On purity theorems for vector bundles\fg}, \emph{Internat. Math. Res. Notices} (2002), \numero 15, p\ptbl 783--788).} Un anneau local noeth\'erien r\'egulier est factoriel (donc parafactoriel si sa dimension est $\geq 2$).
\item
Un anneau local noeth\'erien de dimension $\geq 4$ et qui est une \emph{intersection compl\`ete} est parafactoriel.
\end{enumeratei}
\end{theoreme}

\begin{corollaire}[Conjecture de\label{XI.3.14} \sisi{SAMUEL}{Samuel}\ndemark]\ndetext{pour une preuve dans la m\^eme veine, mais plus \'el\'ementaire, voir Call~F.\ \&\ Lyubeznik~G., {\og A simple proof of Grothendieck's theorem on the parafactoriality of local rings\fg}, in \emph{Commutative algebra: syzygies, multiplicities, and birational algebra (South Hadley, MA, 1992)}, Contemp. Math., vol.~159, American Mathematical Society, Providence, RI, 1994, p\ptbl 15--18.} Un anneau local noeth\'erien $A$ qui est une intersection compl\`ete et qui est factoriel en codimension $\geq 3$ (\ie $\dim A_{{\pp}}\leq 3$ entra\^ine que $A_{{\pp}}$ est factoriel) est factoriel.
\end{corollaire}

\emph{Prouvons le corollaire}

On raisonne par r\'ecurrence sur la dimension de $A$.

Si $\dim A\leq 3$, $A$ est factoriel par hypoth\`ese.

Si $\dim A>3$, par l'hypoth\`ese de r\'ecurrence, en remarquant qu'un localis\'e d'une intersection compl\`ete l'est aussi, tous les localis\'es de $A$ autres que $A$ sont factoriels. Par le th\'eor\`eme~\Ref{XI.3.13}~(ii), $A$ est parafactoriel, donc factoriel par~\Ref{XI.3.10}.

D\'emontrons~\Ref{XI.3.13}~(i) (suivant \sisi{KAPLANSKY}{Kaplansky})\sfootnote{C'est la d\'emonstration reproduite dans \EGA IV~21.11.1.}.

Soit $A$ un anneau local noeth\'erien r\'egulier, posons $\dim A=n$. Si $n=0$ ou $1$, le r\'esultat est connu. Supposons $n\geq 2$, raisonnons par r\'ecurrence sur $n$ et supposons $n\geq 2$ et le th\'eor\`eme d\'emontr\'e pour les anneaux de dimension $<n$. Posons $X'=\Spec(A)$ et $X=X'-\{x\}$, o\`u $x=\rr(A)$. Les localis\'es de $A$ autres que $A$ sont r\'eguliers et de dimension $<n$, donc factoriels. De plus $\prof A = \dim A\geq 2$. Il suffit donc de prouver que $\Pic(X)=0$ (cor\ptbl\Ref{XI.3.10}). Soit donc ${\Lcal}$ un $\OX$-Module inversible, on sait qu'on peut le prolonger en un $\Oo_{X'}$-Module coh\'erent ${\Lcal}'$. Il existe une r\'esolution de ${\Lcal}'$ par des $\OX$-Modules libres:
\[
0\from {\Lcal}'\from {\Lcal}'_1\from\cdots\from{\Lcal}'_n\from 0{\quoi},
\]
car\pageoriginale la dimension cohomologique de $A$ est finie. Par restriction \`a $X'$ on obtient une r\'esolution libre finie. Il suffit donc de prouver le lemme suivant:

\begin{lemme} \label{XI.3.15} Soit $(X, \OX)$ un espace annel\'e et soit ${\Lcal}$ un $\OX$-Module localement libre qui admet une r\'esolution finie par des modules \emph{libres de type fini}. Alors $\det({\Lcal}) \simeq \OX$.
\end{lemme}

Rappelons que l'on d\'efinit $\det({\Lcal})$ comme la puissance ext\'erieure maximum de ${\Lcal}$. Dans le cas envisag\'e, $\det({\Lcal}) \simeq{\Lcal}$ car ${\Lcal}$ est inversible, donc le lemme permet de conclure. D\'emontrons ce lemme. Soit
\[
0\from {\Lcal}_0\from {\Lcal}_1\from {\Lcal}_2\from\cdots\from{\Lcal}_n\from 0{\quoi},
\]
la suite exacte annonc\'ee, o\`u ${\Lcal}_0={\Lcal}$. Puisque tout est localement libre, on a:
\[
\tbigotimes_{0\leq i\leq n}(\det({\Lcal}_i))^{(-1)^i} \simeq \OX
\]
or tous les ${\Lcal}_i$, $i>0$, sont libres, donc aussi leurs d\'eterminants, donc aussi celui de ${\Lcal}_0 = {\Lcal}$.\qed

\medskip
Il nous faut encore d\'emontrer le (ii) du th\'eor\`eme. Auparavant, prouvons un lemme qui permettra de proc\'eder par r\'ecurrence:

\begin{lemme} \label{XI.3.16}
Soit $A$ un anneau local noeth\'erien quotient d'un r\'egulier. Soit $t\in\rr(A)$ un \'el\'ement $A$-r\'egulier. Supposons que $A$ soit complet pour la topologie $t$-adique. Posons $X'=\Spec(A)$, $Y'=V(t)\simeq \Spec(B)$, $B=A/tA$, $X=X'-\{x\}$, $Y=Y'-\{x\}$, $x=\rr(A)$. Supposons que:
\begin{enumeratea}
\item
pour tout $y\in X$ ferm\'e dans $X$ on ait $\prof\Oo_{X, y}\geq 2$,
\item
$\prof A/tA\geq 3$,
\end{enumeratea}
\noindent
alors\pageoriginale l'application $\Pic(X)\to\Pic(Y)$ est injective. En particulier, si $B$ est parafactoriel, $A$ l'est aussi.
\end{lemme}

On sait que a)~entra\^ine $\Lef(X, Y)$ gr\^ace \`a \Ref{X}~\Ref{X.2.1}. Si on prouve que $\H^1(Y, P_n)=0$ pour tout $n\geq 0$, on saura \sisi{\ignorespaces}{gr\^ace \`a} (\Ref{XI.2.2}) que $\Pic(X)\to\Pic(Y)$ est injectif. Si, de plus, $B$ est parafactoriel, on saura que $\Pic(Y)=0$ (\Ref{XI.3.5}) donc $\Pic(X)=0$, or $\prof(A)=3+1\geq 2$ car $t$ est $A$-r\'egulier, donc $A$ sera parafactoriel par~\Ref{XI.3.5}.

Soit $\Ii=\widetilde{(tA)}$ le $\Oo_{X'}$-Module associ\'e \`a l'id\'eal $tA$. Au \numero \Ref{XI.1}, on~a pos\'e $P_n = \sisi{(\Ii^{n+1}/\Ii^{n+2}){|Y}}{(\Ii^{n+1}/\Ii^{n+2})_{|Y}}$, pour tout $n\geq 0$. Or $t$ est $A$-r\'egulier donc $P_n\simeq \Oo_Y$. \emph{Il nous reste donc \`a prouver que $\H^1(Y, \Oo_Y)=0$.} Or $Y = Y'-\{x\}$ est un ouvert de $Y'$, on~a donc une suite exacte~(\Ref{I}~(\Ref{eq:I.27})):
\[
\H^1(Y', \Oo_{Y'})\to\H^1(Y, \Oo_{Y})\to\H^2_x(Y', \Oo_{Y'}){\quoi},
\]
dont le terme de droite est nul en vertu de l'hypoth\`ese~b), et celui de gauche parce que $Y'$ est affine.\qed

\begin{lemme} \label{XI.3.17}
Conservant les hypoth\`eses de~\Ref{XI.3.16}, supposons de plus que:
\begin{itemize}
\item[c)] pour tout $y$ ferm\'e dans $Y$, on $\prof\Oo_{X, y}\geq 3$,
\item[d)] $\prof A/tA\geq 4$ (plus fort que b),
\item[e)] pour tout $y$ ferm\'e de $X$, $y\in Y$, l'anneau $\Oo_{X, y}$ est parafactoriel.
\end{itemize}
Alors l'application $\Pic(X)\to\Pic(Y)$ est un isomorphisme, en particulier pour que $A$ soit parafactoriel, il faut et il suffit que $B$ le soit.
\end{lemme}

\enlargethispage{5mm}%
On sait~(\Ref{X}~\Ref{X.2.1}) que a)~et~c) entra\^inent $\Leff(X, Y)$. De plus, par le raisonnement qui vient d'\^etre fait, d)~entra\^ine que $\H^1(Y, P_n)=0$ pour tout $n\geq 0$ et $i=1$ ou $i=2$. De plus, pour tout voisinage ouvert $U$ de $Y$ dans $X$, le compl\'ementaire de $U$ dans $X$ est form\'e d'un nombre fini de points ferm\'es. Gr\^ace \`a~e) et au th\'eor\`eme~\Ref{XI.3.12}, on en d\'eduit que $\Pic(X)\to \Pic(Y)$ est un isomorphisme. D'autre part $\prof A\geq \prof B\geq 2$; par le crit\`ere~\Ref{XI.3.5} on en d\'eduit\pageoriginale que $A$ est parafactoriel si et seulement si $B$ l'est.

D\'emontrons maintenant~\Ref{XI.3.13}~(ii). Soit $R$ un anneau local noeth\'erien r\'egulier. Soit $(t_1, \ldots, t_k)$ une $R$-suite. Posons $B=R/(t_1, \ldots, t_k)$ et supposons que $\dim B\geq 4$. Il faut prouver que $B$ est parafactoriel. Raisonnons par r\'ecurrence sur~$k$. Si $k=0$, $B$ est r\'egulier, donc factoriel par \Ref{XI.3.13}~(i), donc parafactoriel par~\Ref{XI.3.10}. Supposons $k\geq 1$ et le th\'eor\`eme d\'emontr\'e pour $k'<k$. Posons $A=R/(t_1, \ldots, t_{k-1})$ donc $B=A/t_k A$. On peut supposer $B$ complet d'apr\`es~\Ref{XI.3.7}. Par l'hypoth\`ese de r\'ecurrence, $A$ est parafactoriel. Prouvons que l'on peut appliquer le lemme~\Ref{XI.3.17}. On a suppos\'e $B$ complet donc aussi $A$, et donc $A$ est complet pour la topologie $t$-adique. Si $x\in X$, et si $x$ est ferm\'e dans $X$, $A_x$ est une intersection compl\`ete de dimension $\geq 4$, avec $k'<k$. D'apr\`es l'hypoth\`ese de r\'ecurrence $A_x$ est parafactoriel, et par ailleurs de profondeur $\geq 4$. Ce qui donne a), c)~et~e). De plus, $\dim A\geq 5$, d'o\`u~d).\qed

\begin{theoreme} \label{XI.3.18}
Soit $X$ un pr\'esch\'ema localement noeth\'erien et soit $\Ii$ un faisceau coh\'erent d'id\'eaux de $X$. Posons $Y=V(\Ii)$. Soit $n$ un entier. Supposons que:
\begin{enumeratei}
\item
on ait $\sisi{\leff}{\Leff}(X, Y)$ (\cf exemples \Ref{X}~\Ref{X.2.1}~et~\Ref{X}~\Ref{X.2.2}),
\item
pour tout $p\geq n$, on ait $\H^i(Y, \Ii^{p+1}/\Ii^{p+2})=0$ pour $i=1$ et $i=2$,
\item
pour tout ouvert $U\supset Y$ et tout $x\in X-U$, l'anneau $\OXx$ soit r\'egulier de dimension $\geq 2$ ou une intersection compl\`ete de dimension $\geq 4$.
\end{enumeratei}
\noindent
Alors, pour tout ouvert $U\supset Y$ et tout entier $p\geq n$, les homomorphismes
\[
\Pic(X)\to\Pic(U)\to\Pic(Y_p)
\]
sont des isomorphismes, o\`u $Y_p$ d\'esigne le pr\'esch\'ema $(Y, \OX/\Ii^{p+1})$.
\end{theoreme}

Il suffit de conjuguer \Ref{XI.3.12}~et~\Ref{XI.3.13}.

\chapter[Applications aux sch\'emas alg\'ebriques projectifs]
{Applications\\ aux sch\'emas alg\'ebriques projectifs} \label{XII}

\renewcommand{\theequation}{\arabic{equation}}

\section[Th\'eor\`eme de dualit\'e projective et th\'eor\`eme de finitude]
{Th\'eor\`eme de dualit\'e projective et th\'eor\`eme de finitude\protect\sfootnotemark}\label{XII.1}

Le\pageoriginale th\'eor\`eme suivant,\sfootnotetext{Le pr\'esent expos\'e, r\'edig\'e en Janvier~1963, est nettement plus d\'etaill\'e que n'\'etait l'expos\'e oral, en Juin~1962.} essentiellement contenu dans FAC\sfootnote{J.-P.~Serre, \sisi{Faisceaux alg\'ebriques Coh\'erents, Ann. of Math. t.61 (1955) 197-278}{{\og Faisceaux alg\'ebriques coh\'erents\fg}, \emph{Ann. of Math.} \textbf{61} (1955), p\ptbl197-278}.} (\`a cela pr\`es que dans le temps, Serre ne disposait pas du langage des Ext de faisceaux de modules\nde{le lecteur friand d'Histoire des Math\'ematiques consultera avec int\'er\^et la lettre que Grothendieck \'ecrivit \`a Serre le 15 d\'ecembre 1955 et la r\'eponse de celui-ci du 22 d\'ecembre de la m\^eme ann\'ee; voir \emph{Correspondance Grothendieck-Serre}, \'edit\'ee par Pierre Colmez et Jean-Pierre Serre, Documents Math\'ematiques, vol.~2, Soci\'et\'e Math\'ematique de France, Paris, 2001.}), est l'analogue global du th\'eor\`eme de dualit\'e locale (Expos\'e \Ref{IV}), qui a \'et\'e calqu\'e sur lui.

\begin{theoreme} \label{XII.1.1}
Soit $k$ un corps, $X=\PP_k^r$ \sisi{le sch\'ema projectif type}{l'espace projectif} de dimension $r$ sur $k$, $F$~un module coh\'erent variable sur $X$; alors on~a un isomorphisme de $\partial$-foncteurs en~$F$:
\steco
\begin{equation} \label{eq:XII.1} {}\H^i(X, F)' \isomfrom \Ext^{r-i}(X;F, \omx^r),
\end{equation}
o\`u on pose
\steco
\begin{equation} \label{eq:XII.2} {}\Omega_{X/k}^r= \Oo_{\PP_k^r}(-r-1).
\end{equation}
\end{theoreme}

\begin{remarquestar}
Bien entendu, ce Module est aussi le Module des diff\'erentielles relatives de degr\'e $r$ de $X$ sur $k$. Sous cette forme, le th\'eor\`eme reste vrai si $X$ est un sch\'ema propre et lisse sur $k$ (pour le cas projectif, voir A. Grothendieck, \og \sisi{Th\'eor\`emes de dualit\'e pour les faisceaux alg\'ebriques coh\'erents}{\emph{Th\'eor\`emes de dualit\'e pour les faisceaux alg\'ebriques coh\'erents}}\fg, S\'eminaire Bourbaki Mai 1957)\sfootnote{Pour un th\'eor\`eme de dualit\'e plus g\'en\'eral, \cf le S\'eminaire Hartshorne cit\'e \`a la fin de \Exp \Ref{IV}.}. Lorsque $F$ est localement libre, on retrouve le th\'eor\`eme de dualit\'e de Serre $\H^i(X, F)' \isomto \H^{r-i}(\SheafHom_{\OX}(F, \omx^r))$. Le th\'eor\`eme \Ref{XII.1.1} (qui redonne d'ailleurs le cas d'un $X$ projectif sur $k$, comme dans {\loccit}) suffira pour notre propos.
\end{remarquestar} L'homomorphisme\pageoriginale (\Ref{eq:XII.1}) se d\'eduit de l'accouplement de Yoneda
\steco
\begin{equation} \label{eq:XII.3} {}\H^i(X, F) \times \Ext^{r-i}(X;F, \omx^r) \to \H^r(X, \omx^r),
\end{equation}
et d'un isomorphisme bien connu (\cf FAC, ou \EGA III 2.1.12):
\steco
\begin{equation} \label{eq:XII.4} {}\H^r(X, \omx^r) = \H^r(\PP_k^r, \Oo_{\PP_k^r}(-r-1)) \isomto k.
\end{equation}
Pour montrer que (\Ref{eq:XII.1}) est un isomorphisme, on proc\`ede comme dans le cas du th\'eor\`eme de dualit\'e locale, en notant que $\H^r(X, F)$ comme foncteur en $F$ est exact \`a droite (puisque $\H^n(X, F)=0$ pour $n>r$), et que tout Module coh\'erent est isomorphe \`a un quotient d'une somme directe de Modules de la forme $\Oo(-m)$, avec $m$ grand. Cela nous ram\`ene par r\'ecurrence descendante sur $i$ \`a faire la v\'erification pour un faisceau de la forme $\Oo(-m)$, o\`u cela est contenu dans les calculs explicites bien connus (FAC, ou \EGA III 2.1.12). On peut d'ailleurs supposer $-m \leq -r-1$, auquel cas $\H^i(X, \Oo(-m))=0$ pour $i\neq r$.

\begin{corollaire} \label{XII.1.2} Pour $F$ coh\'erent donn\'e, et $m$ assez grand, on~a un isomorphisme canonique
\steco
\begin{equation} \label{eq:XII.5} {}\H^i(X, F(-m))'\isomto\H^0(X, \SheafExt_{\OX}^{r-i}(F, \omx^r)(m))
\end{equation}
(o\`u le $'$ d\'esigne le vectoriel dual).
\end{corollaire}
En effet, sur \sisi{un sch\'ema}{l'espace} projectif $X$, on~a pour tout couple de faisceaux coh\'erents $F, G$ et pour $n$ assez grand un isomorphisme canonique:
\steco
\begin{equation} \label{eq:XII.6} {}\Ext^n(X;F(-m), G)\simeq \Ext^n(X;F, G(m)) \isomto \H^0(X, \SheafExt_{\OX}^{n}(F, G)(m)),
\end{equation}
(l'isomorphisme des deux premiers membres \'etant trivialement vrai pour tout $m$), comme il r\'esulte de la suite spectrale des $\Ext$ globaux
\[
\H^p(X, \SheafExt^q_{\OX}(F, G(m))) \To \Ext^{\boule}(X;F, G(m)),
\]
qui\pageoriginale d\'eg\'en\`ere pour $m$ assez grand gr\^ace au fait que
\[
\SheafExt^q_{\OX}(F, G(m)) \simeq \SheafExt^q_{\OX}(F, G)(m),
\]
et que les $\SheafExt^q_{\OX}(F, G)$ sont des faisceaux coh\'erents. Donc (\Ref{eq:XII.5}) r\'esulte de (\Ref{eq:XII.6}) et (\Ref{eq:XII.1}).

\begin{corollaire} \label{XII.1.3}
\sisi{Conditions \'equivalentes pour $i, F$ donn\'es}{Pour $i, F$ donn\'es, les conditions suivantes sont \'equivalentes}:
\begin{enumeratei}
\item
$\H^i(X, F(-m))=0$ pour $m$ grand.
\item[\textup{(i bis)}] $\H^i(X, F(\cdot))= \sisi{\coprod}{\bigoplus}_{m\in\ZZ}\H^i(X, F(-m))$ est un $S$-module de type fini, o\`u $S=k[t_o, \ldots, t_r]$.
\item
$\SheafExt^{r-i}_{\OX}(F, \omx^r)=0$.
\item[\textup{(ii bis)}] $\SheafExt^{r-i}_{\OX}(F, \OX)=0$.
\item
$\H^i_x(F_x)=0$ pour tout point ferm\'e $x$ de $X$.
\item
$\H^{i+1}_x(\widetilde{F}_x)=0$ pour tout point ferm\'e $x$ du c\^{o}ne projetant \'epoint\'e $\widetilde{X}=\Spec(S)-\Spec(k)$ de $X$, o\`u $\widetilde{F}$ d\'esigne l'image inverse de $F$ par le morphisme canonique $\widetilde{X}\to X$.
\end{enumeratei}
\end{corollaire}

\begin{proof}
(i) \SSI (i bis) puisque le sous-module de $\H^i(X, F(\cdot))$ somme des composantes homog\`enes de degr\'e $\geq \nu$ est de type fini sur $S$ (en fait, pour $i\neq 0$, il est m\^eme de type fini sur $k$), (\cf FAC ou \EGA III 2.2.1 et 2.3.2).

(i) \SSI (ii) en vertu de corollaire \Ref{XII.1.2}.

(ii) \SSI (ii bis) car $\Omega_{X/k}^r$ est localement isomorphe \`a $\OX$.

(ii bis) \SSI (iii) en vertu du th\'eor\`eme de dualit\'e locale pour $\OXx$ (qui est un anneau\pageoriginale local r\'egulier de dimension $r$), d'apr\`es lequel le \og dual\fg de $\SheafExt^{r-i}_{\OX}(F, \OX)_x$ s'identifie \`a $\H^i_x(F_x)$ (V 2.1).

(ii bis) \'equivaut \`a la relation analogue
\[
\SheafExt^{r-i}_{\Oo_{\widetilde{X}}}(\widetilde{F}, \Oo_{\widetilde{X}})=0
\]
(gr\^ace au fait que $\widetilde{X}\to X$ est fid\`element plat, par suite l'image inverse de $\SheafExt^{r-i}_{\OX}(F, \OX)$ est isomorphe \`a $\SheafExt^{r-i}_{\Oo_{\widetilde{X}}}(\widetilde{F}, \Oo_{\widetilde{X}})$) et cette derni\`ere relation \'equivaut \`a (iv) par le th\'eor\`eme de dualit\'e locale pour l'anneau local $\OXx$, qui est r\'egulier de dimension $r+1$.
\skipqed
\end{proof}

En particulier, appliquant ceci \`a tous les $i\leq n$, on trouve:

\begin{corollaire} \label{XII.1.4}
Conditions \'equivalentes pour $n, F$ donn\'es:
\begin{enumeratei}
\item
$\H^i(X, F(-m))=0$ pour $i\leq n$ et $m$ grand.
\item[\textup{(i bis)}] $\H^i(X, F(\cdot))$ est un $S$-module de type fini pour $i\leq n$.
\item
$\prof(F_x)>n$ pour tout point ferm\'e $x$ de $X$.
\item
$\prof(\widetilde{F}_x)>n+1$ pour tout point ferm\'e $x$ de $\widetilde{X}$.
\end{enumeratei}
\end{corollaire}

L'int\'er\^et des corollaires \Ref{XII.1.3} et \Ref{XII.1.4} est d'exprimer une condition globale (i) ou~(i~bis) en termes de conditions locales, \sisi{\ignorespaces}{\`a} savoir l'annulation d'invariants locaux comme $\H_x^i(X, F_x)$ ou $\H_x^i(X, \widetilde{F}_x)$, ou une in\'egalit\'e sur la profondeur. Sous cette forme, ces r\'esultats restent trivialement valables pour un sch\'ema projectif quelconque $X$, et un faisceau inversible tr\`es ample $\OX(1)$ sur $X$, comme on voit en induisant ce dernier \`a l'aide d'une immersion projective $X\to\PP_k^r$ convenable. (Bien entendu, les conditions \Ref{XII.1.3}~(ii) et \Ref{XII.1.3}~(ii bis) ne sont plus \'equivalentes aux autres dans ce cas g\'en\'eral, sauf si on suppose que $X$ est r\'egulier par exemple). On peut d'ailleurs g\'en\'eraliser au cas d'un morphisme projectif $X\to S$ de la fa\c con suivante:

\begin{proposition} \label{XII.1.5}
Soient\pageoriginaled $f:X\to S$ un morphisme projectif, avec $S$ noeth\'erien, $\OX(1)$ un Module inversible sur $X$ tr\`es ample relativement \`a $S$, $F$ un Module coh\'erent sur $X$, plat par rapport \`a $S$, $s$ un \'el\'ement de $S$, $X_s$ la fibre de $X$ en $s$ (consid\'er\'e comme sch\'ema projectif sur $k(s)$), $F_s$ le faisceau induit sur $X_s$ par $F$, enfin $i$ un entier. Supposons que pour tout point ferm\'e $x$ de $X_s$, on ait $\H^i_x(F_{s, x})=0$ (par exemple $\prof(F_{s, x})>i$). Alors il existe un voisinage ouvert $U$ de $s$, tel que la m\^eme condition soit v\'erifi\'ee pour les $s'\in U$. De plus, pour un tel $U$, on a
\[
\R^if_\ast(F(-m))= 0 \; \text{pour}\; m \; \text{grand, }
\]
et si $\Ss$ est une alg\`ebre gradu\'ee quasi-coh\'erente sur $S$, engendr\'ee par $\Ss^1$, et qui d\'efinit~$X$ avec $\OX(1)$ comme $X=\Proj(\Ss)$, $\OX(1)=\Proj(\Ss(1))$, alors le $S$-module
\[
\R^if_\ast(F(\cdot)) = \sisi{\coprod}{\tbigoplus}_{m\in\ZZ}\R^if_\ast(F(m))
\]
est de type fini sur $U$.
\end{proposition}

Plongeons $X$ dans un $X'=\PP_S^r$ de fa\c con que $\OX(1)$ soit induit par $\Oo_{\sisi{P}{{X'}}}(1)$ (c'est possible \`a condition de remplacer $S$ par un voisinage affine de $s$). Posons\refstepcounter{toto}\label{P5}\sfootnote{La premi\`ere partie de \Ref{XII.1.5} peut s'obtenir aussit\^{o}t en appliquant l'\'enonc\'e purement local \EGA IV 12.3.4 aux $\underline{E}^j$ pr\'ec\'edents, ce qui court-circuite la plus grande partie de la d\'emonstration qui suit.} pour tout entier $j$, et tout $t\in S$:
\steco
\begin{equation} \label{eq:XII.7} {}\underline{E}^j(t)= \SheafExt^j_{\Oo_{X'_t}} (F_{t'}, \Oo_{X'_t}(-r-1))
\end{equation}
Donc $\underline{E}^j(t)$ est un Module coh\'erent sur $X_t$, je dis que pour $t$ variable, la famille de ces Modules est \og constructible\fg au sens suivant: il existe pour tout $t\in S$ un ouvert non vide $V$ de l'adh\'erence $\overline{t}$, qu'on munit de la structure r\'eduite induite, et un Module coh\'erent $\underline{E}^j(V)$ sur $X_V=X\times_SV$, \emph{plat} relativement \`a $V$, tel que pour tout $t'\in V$, $\underline{E}^j(t)$ soit isomorphe au Module induit par $\underline{E}^j(V)$ sur $X_t$. Pour v\'erifier cette assertion, posant $Z=\overline{t}$ muni de la structure induite, on consid\`ere les Modules coh\'erents\pageoriginale
\[
\underline{E}^j(Z)= \SheafExt^j_{\Oo_{X'_Z}} (F_Z, \Oo_{X'_Z}(-r-1))
\]
(o\`u l'indice $Z$ signifie encore qu'on s'induit \sisi{au dessus}{au-dessus} de $Z$), et on prend pour $V$ un ouvert non vide de $Z$ tel que les Modules $\underline{E}^j(Z)$ soient plats \sisi{au dessus}{au-dessus} de $V$: c'est possible, car on v\'erifie tout de suite que $\underline{E}^j(Z)=0$ pour $j$ non compris dans l'intervalle $[0, r]$, et on peut appliquer \SGA 1 IV 6.11. On prend alors pour $\underline{E}^j(V)=\underline{E}^j(Z)|_{X_V}$, et on v\'erifie facilement qu'il r\'epond \`a la question.

De la remarque pr\'ec\'edente r\'esulte qu'il existe une partition finie de $S$ form\'ee par des ensemble$V_\alpha$ de la forme $V=V(t)$ comme dessus, (r\'ecurrence noeth\'erienne), et appliquant le th\'eor\`eme de Serre \EGA III 2.2.1 aux $\underline{E}^j(V_i)$, on voit qu'il existe un entier $m_0$ tel que
\[
\R^if_{V_\alpha\ast}(\underline{E}^j(V_\alpha))=0\textup{ pour } i\neq 0, m\geq m_0, \textup{ pour }j,
\]
d'o\`u il r\'esulte, utilisant la platitude de $\underline{E}^j(V_\alpha)$ par rapport \`a $V_\alpha$ et des relations \`a la Künneth faciles (\cf \EGA III par\ptbl7), que
\[
\H^i(X_t, \underline{E}^j(t)(m))=0\textup{ pour } i\neq 0, m\geq m_0, \textup{ pour }j,
\]
pour tout $t\in V_\alpha$, donc pour tout $t$ puisque les $ V_\alpha$ recouvrent $S$. De ceci et de la suite spectrale des Ext globaux r\'esulte, gr\^ace \`a \Ref{XII.1.1}, comme dans la d\'emonstration de \Ref{XII.1.2}, un isomorphisme
\steco
\begin{equation} \label{eq:XII.8} {}\H^i(X_t, F_t(-m))' \isomfrom \H^0(X_t, \underline{E}^{r-i}(t)(m))\textup{ pour } m\geq m_0,
\end{equation}
tout entier $i$, et tout $t\in S$.

Utilisons maintenant l'hypoth\`ese faite sur $F_s$, qui s'\'ecrit
\steco
\begin{equation} \label{eq:XII.9} {}\underline{E}^{r-i}(s)=0,
\end{equation}
et gr\^ace \`a (\Ref{eq:XII.8}) \'equivaut \`a
\steco
\begin{equation} \label{eq:XII.10} {}\H^i(X_s, F_s(-m))=0 \; \text{pour}\; m\geq m_0.
\end{equation}
Comme $F$ donc $F(-m)$ est plat par rapport \`a $S$, il s'ensuit par
les relations \`a la Künneth d\'ej\`a invoqu\'ees, que (pour $m$ donn\'e) la
m\^eme relation (\Ref{eq:XII.10}) vaut en rempla\c cant\pageoriginale $s$ par un
$t$ voisin de $s$, en particulier pour toute g\'en\'erisation $t$ de
$s$. En vertu de (\Ref{eq:XII.8}), on aura donc, pour une telle
g\'en\'erisation \steco
\begin{equation} \label{eq:XII.11} {}\underline{E}^{r-i}(t)=0,
\end{equation}
or l'ensemble des $t\in s$ pour lesquels cette relation est vraie est \'evidemment un ensemble constructible (car il induit sur chaque $V_\alpha$ un ouvert); comme il contient les g\'en\'erisations de $s$, il contient un voisinage ouvert $U$ de $s$. Cela prouve la premi\`ere assertion de \Ref{XII.1.5}. De plus, pour $t\in U$, on conclut de (\Ref{eq:XII.11}) et (\Ref{eq:XII.8}) que
\steco
\begin{equation} \label{eq:XII.12} {}\H^i(X_t, F_t(-m)=0 \; \text{pour}\; m\geq m_0, \; t\in U,
\end{equation}
ce qui en vertu des relations \`a la Künneth, implique (en fait, est bien plus fort que)
\steco
\begin{equation} \label{eq:XII.13} {}\R^if_\ast(F(-m))=0 \; \text{sur}\; U, \; \text{pour} \; m\geq m_0.
\end{equation}
Cela prouve la deuxi\`eme assertion de \Ref{XII.1.5}. Enfin la derni\`ere en r\'esulte aussit\^{o}t, en proc\'edant comme au d\'ebut de la d\'emonstration de \Ref{XII.1.3}.

\refstepcounter{subsection}\label{XII.1.6}
\begin{enonce*}[remark]{Remarque 1.6\sfootnotemark} \sfootnotetext{Cette remarque se pr\'ecise par la note de bas de page \sisi{5.}{\pageref{P5}.}} La d\'emonstration se simplifie notablement (en \'eliminant toute consid\'eration de constructibilit\'e) lorsqu'on suppose d\'ej\`a que l'hypoth\`ese faite pour $s\in S$ est v\'erifi\'ee en tous les $s'\in S$. En fait, lorsqu'on fait l'hypoth\`ese que $F_s$ est de profondeur $>i$ aux points ferm\'es de $X$, on dispose d'un \'enonc\'e g\'en\'eral, \emph{de nature locale sur} $X$, qui dit que la m\^eme condition est v\'erifi\'ee pour tous les $X_t$, \`a condition de remplacer $X$ par un voisinage ouvert convenable de la fibre $X_s$ (en d'autres termes, une certaine partie de $X$, d\'efinie par des conditions sur les Modules induits par $F$ sur les fibres, est ouverte, \cf \EGA IV). Comme $f$ est ici propre, on pourra donc prendre ce voisinage de la forme $f^{-1}(U)$, ce qui redonne la premi\`ere assertion de \Ref{XII.1.5} sans p\'enible d\'evissage. Dans ce cas g\'en\'eral, on peut encore prouver par la m\'ethode de \loccit que la premi\`ere assertion de \Ref{XII.1.5} (d\'emontr\'ee ici par voie globale, en utilisant que $X$ est projectif sur $S$) r\'esulte encore d'un \'enonc\'e purement local sur $X$ (que le lecteur explicitera s'il le juge utile).
\end{enonce*}

\section[Th\'eorie de Lefschetz pour un morphisme projectif: comparaison]
{Th\'eorie de Lefschetz pour un morphisme projectif: th\'eor\`eme de comparaison de Grauert} \label{XII.2}

\setcounter{equation}{13}
C'est\pageoriginale le th\'eor\`eme suivant:

\begin{theoreme} \label{XII.2.1}
Soient $f:X\to S$ un morphisme projectif, avec $S$ noeth\'erien, $\OX(1)$ un Module inversible sur $X$, ample relativement \`a $S$, $Y$ le pr\'esch\'ema des z\'eros d'une section $t$ de $\OX(1)$, $J$ l'id\'eal d\'efinissant $Y$, $X_n$ le sous-pr\'esch\'ema de $X$ d\'efini par $J^{n+1}$, $\widehat{X}$ le compl\'et\'e formel de $X$ le long de $Y$, $\widehat{f}:\widehat{X}\to S$ le morphisme compos\'e $\widehat{X}\to X\to S$, $F$ un Module coh\'erent sur $X$, plat relativement \`a $S$. On suppose de plus que pour tout $s\in S$, le Module $F_s$ induit sur la fibre $X_s$ est de profondeur $>n$ en les points de ladite fibre, et que $t$ est $F$-r\'egulier. Sous ces conditions:
\begin{enumeratei}
\item
L'homomorphisme canonique
\[
\R^if_\ast(F) \to \R^i\widehat{f}_\ast(\widehat{F})
\]
est un isomorphisme pour $i<n$, un monomorphisme pour $i=n$.
\item
L'homomorphisme canonique
\[
\R^i\widehat{f}_\ast(\widehat{F}) \to \varprojlim_m \R^if_\ast(F_m)
\]
est un isomorphisme pour $i\leq n$.
\end{enumeratei}
\end{theoreme}

\begin{proof}
On se ram\`ene aussit\^{o}t au cas o\`u $S$ est affine, et \`a prouver alors le

\begin{corollaire} \label{XII.2.2}
Sous les conditions de \Ref{XII.2.1} supposons de plus $S$ affine. alors:
\begin{enumeratei}
\item
L'homomorphisme canonique
\[
\H^i(X, F) \to \H^i(\widehat{X}, \widehat{F})
\]
est\pageoriginale un isomorphisme pour $i<n$, un monomorphisme pour $i=n$.
\item
L'homomorphisme canonique
\[
\H^i(\widehat{X}, \widehat{F}) \to \varprojlim_m \H^i(X_m, F_m)
\]
est un isomorphisme pour $i\leq n$.
\end{enumeratei}
\end{corollaire}

Quitte \`a remplacer $\OX(1)$ par une puissance tensorielle, et $t$ par une puissance de~$t$, on peut supposer $\OX(1)$ tr\`es ample relativement \`a $S$. D'autre part $t$ donc $t^m$ \'etant $F$-r\'egulier, la multiplication par $t^m$, consid\'er\'ee comme homomorphisme de $F(-m)$ dans $F$, est injectif, donc on~a pour tout $m\geq 0$ une suite exacte:
\steco
\begin{equation} \label{eq:XII.14}
0\to F(-m) \lto{t^m} F \to F_m \to0,
\end{equation}
d'o\`u une suite exacte de cohomologie
\[
\H^i(X, F(-m)) \to \H^i(X, F) \to \H^i(X, F_m) \to \H^{i+1}(X, F(-m)).
\]
Or en vertu de \Ref{XII.1.5} on~a $\H^i(X, F(-m))=0$ pour $i\leq n$, et $m$ assez grand, ce qui prouve~le

\begin{lemme} \label{XII.2.3}
Pour $m$ grand, l'homomorphisme canonique
\[
\H^i(X, F) \to \H^i(X, F_m)
\]
est bijectif si $i<n$, injectif si $i=n$.
\end{lemme}

\enlargethispage{-2\baselineskip}%
Cela montre que pour $i<n$, le syst\`eme projectif $(\H^i(X_m, F_m))_{m\geq 0}$ est essentiellement constant, \afortiori satisfait la condition de Mittag-Leffler, donc (compte tenu \sisi{que}{de} $\widehat{F}= \varprojlim F_m$) on conclut (ii) par $\EGA 0_{\textup{III}}$ 13.3. D'autre part (i) en r\'esulte trivialement, compte tenu de \Ref{XII.2.3}.

\refstepcounter{subsection}\label{XII.2.4}
\begin{enonce*}{Corollaire 2.4\ndemark}\ndetext{voir le corollaire I.1.4 de l'article de Mme Raynaud (Raynaud~M., {\og Th\'eor\`emes de Lefschetz en cohomologie des faisceaux coh\'erents et en cohomologie \'etale. Application au groupe fondamental\fg}, \emph{Ann. Sci. \'Ec. Norm. Sup. (4)} \textbf{7} (1974), p\ptbl 29--52).}
Soient\pageoriginale $f:X\to S$ un morphisme projectif et plat, avec $S$ localement noeth\'erien, $\OX(1)$ un Module inversible sur $X$, ample relativement \`a $S$, $t$ une section de ce Module qui soit $\OX$-r\'eguli\`ere, $Y$ le sous-pr\'esch\'ema des z\'eros de $t$, $\widehat{X}$ le compl\'et\'e formel de $X$ le long de $Y$. On suppose que pour tout $s\in S$, $X_s$ est de profondeur $\geq 1$ (\resp de profondeur $\geq 2$) en ses points ferm\'es. Alors pour tout voisinage ouvert $U$ de $Y$, le foncteur
\[
F \mto \widehat{F}
\]
de la cat\'egorie des Modules coh\'erents localement libres sur $U$, dans la cat\'egorie des Modules coh\'erents localement libres sur $\widehat{X}$, est fid\`ele (\resp pleinement fid\`ele, \ie la condition de Lefschetz ($\Lef$) de \Ref{X}~\Ref{X.2} est v\'erifi\'ee).
\end{enonce*}
Introduisons pour deux Modules localement libres $F$ et $G$ sur $U$ le Module
\[
H=\SheafHom_{\Oo_U}(F, G)
\]
on est ramen\'e \`a prouver que l'homomorphisme canonique
\steco
\begin{equation} \label{eq:XII.15} {}\H^0(U, H) \to \H^0(\widehat{U}, \widehat{H})
\end{equation}
est injectif (\resp bijectif). Or les Modules $H_t$ sont de profondeur $\geq 1$ (\resp $\geq 2$) aux points ferm\'es de $X_t$, on peut donc appliquer \Ref{XII.2.1}, qui implique la conclusion \Ref{XII.2.4} dans le cas o\`u $U=X$. Dans le cas d'un $U$ quelconque, on note que la question est locale sur $S$, donc on peut supposer $S$ affine. Alors tout Module coh\'erent sur $X$ est quotient d'un Module coh\'erent localement libre, (puisque $\OX(1)$ est un Module inversible \sisi{absolument}{relativement} ample sur $X$). Comme le Module dual $H'=\SheafHom(H, \Oo_U)$ se prolonge en un Module coh\'erent sur $X$, qui est donc isomorphe \`a un conoyau d'un homomorphisme de Modules localement libres sur $X$, il s'ensuit par transposition qu'on peut trouver un homomorphisme
\[
u':{L'}^0 \to {L'}^1
\]
de\pageoriginale Modules localement libres sur $X$, induisant un homomorphisme
\[
u:L^0 \to L^1
\]
de Modules localement libres sur $U$, tel qu'on ait une suite exacte
\[
0\to H \to L^0 \lto{u} L^1.
\]
Utilisant le lemme des cinq (qui devient le lemme des trois), et l'exactitude \`a gauche du foncteur $H^0$, on est ramen\'e \`a prouver que (\Ref{eq:XII.15}) est injectif (\resp bijectif) lorsqu'on y remplace $H$ par $L^0$, $L^1$, ce qui nous ram\`ene au cas o\`u $H$ est induit par un Module localement libre $H'$ sur $X$. D'ailleurs dans le cas non resp\'e cette r\'eduction est m\^eme inutile, car le noyau de (\Ref{eq:XII.15}) est en tout cas form\'e des sections de $H$ sur~$U$ qui s'annulent dans un voisinage ouvert convenable $V$ de $Y$, or l'homomorphisme de restriction $\H^0(U, H)\to\H^0(V, H)$ est injectif, car $H$ est de profondeur $\geq 1$ en les points de toute partie ferm\'ee $Z$ de $X$ ne rencontrant pas $Y$ (\cf lemme plus bas). Dans le cas resp\'e, on est ramen\'e \`a prouver que
\[
H^0(X, H')\to\H^0(U, H')
\]
est bijectif, ce qui r\'esulte du fait que $H'$ est de profondeur $\geq 2$ en tous les points d'un ferm\'e $Z=X- U$ de $X$ ne rencontrant pas $Y$. Il faut donc simplement prouver le

\begin{lemme} \label{XII.2.5} Soit $F$ un Module coh\'erent sur $X$, plat par rapport \`a $S$, tel que pour tout $s\in S$, $F_s$ soit de profondeur $\geq n$ en tout point ferm\'e de $X_s$. Alors pour toute partie ferm\'ee $Z$ de $X$ ne rencontrant pas $Y$, $F$ est de profondeur $\geq n$ en tous les points de $Z$.
\end{lemme}

En effet, pour tout $x\in X$, posant $s=f(x)$, on~a
\steco
\begin{equation} \label{eq:XII.16} {}\prof(F_x)\geq \prof(F_{s, x}),
\end{equation}
comme on voit en relevant n'importe comment une suite $F_{s, x}$ r\'eguli\`ere d'\'el\'ements de\pageoriginale $\rr(\Oo_{X_s, x})$ maximale, ce qui donne une suite $F_x$-r\'eguli\`ere en vertu de \SGA 1 IV 5.7. Or si $x$ appartient \`a un $Z$ comme dans le lemme \Ref{XII.2.5}, alors $x$ est n\'ecessairement ferm\'e dans $X_s$, en d'autres termes, $Z$ est \emph{fini} sur $S$. En effet $Z$ (muni d'une structure induite par $X$) est projectif sur $S$ comme sous-pr\'esch\'ema ferm\'e de $X$ qui l'est, et $Z$ est affine sur $S$ comme sous-pr\'esch\'ema ferm\'e de $X- Y$, qui l'est.
\skipqed
\end{proof}

\begin{remarque} \label{XII.2.6}
Supposons que pour tout $s\in S$, la section $t_s$ de $\Oo_{X_s}(1)$ induite par $t$ soit $\Oo_{X_s}$-r\'eguli\`ere (ce qui implique par \SGA 1 IV 5.7 que $t$ est $\OX$-r\'eguli\`ere). Alors les hypoth\`eses faites sont stables par extension de la base $S'\to S$ ($S'$ localement noeth\'erien). Donc la conclusion reste valable apr\`es tout changement de base.
\end{remarque}

\section[Th\'eorie de Lefschetz pour un morphisme projectif: existence]{Th\'eorie de Lefschetz pour un morphisme projectif: th\'eor\`eme d'existence}\label{XII.3}%

\setcounter{equation}{16}%
\refstepcounter{subsection}\label{XII.3.1}%
\begin{enonce*}{Th\'eor\`eme 3.1\ndemark}
\ndetext{pour une version sans hypoth\`ese de platitude, voir (Raynaud~M., {\og Th\'eor\`emes de Lefschetz en cohomologie des faisceaux coh\'erents et en cohomologie \'etale. Application au groupe fondamental\fg}, \emph{Ann. Sci. \'Ec. Norm. Sup. (4)} \textbf{7} (1974), p\ptbl 29--52, th\'eor\`eme II.3.3).}
Soit $f:X\to S$ un morphisme projectif, avec $S$ noeth\'erien, $\OX(1)$ un Module inversible sur $X$ ample relativement \`a $S$, $X_0$ le sous-pr\'esch\'ema des z\'eros d'une section $t$ de $\OX(1)$, $\widehat{X}$ le compl\'et\'e formel de $X$ le long de $X_0$, $\formelF$ un Module coh\'erent sur $\widehat{X}$, $F_0$ le Module qu'il induit sur $X_0$. On suppose de plus:
\begin{enumeratea}
\item\label{XII.3.1a}
$\formelF$ est plat par rapport \`a $S$.
\item\label{XII.3.1b}
Pour tout $s\in S$, la section $t_s$ induite par $t$ sur $X_s$ est $\formelF_s$-r\'eguli\`ere (ce qui implique que $F_0$ est \'egalement plat par rapport \`a $S$, \cf \textup{\SGA 1 IV 5.7}).
\item\label{XII.3.1c}
Pour tout $s\in S$, $F_{0s}$ est de profondeur $\geq 2$ en les points ferm\'es de $X_{0s}$.
\end{enumeratea}
On suppose de plus que $S$ admet un faisceau inversible ample. Sous ces conditions, il existe un Module coh\'erent $F$ sur $X$, et un isomorphisme de son compl\'et\'e formel $\widehat{F}$ avec $\formelF$.
\end{enonce*}

Cet \'enonc\'e va r\'esulter du suivant:

\begin{corollaire} \label{XII.3.2}
Sous les conditions \Ref{XII.3.1a}), \Ref{XII.3.1b}), \Ref{XII.3.1c}) ci-dessus on~a ce qui suit:
\begin{enumeratei}
\item
Le\pageoriginale Module $\widehat{f}_\ast(\formelF)$ sur $S$ est coh\'erent, donc pour tout $n$, le Module $\widehat{f}_\ast(\formelF(n))$ sur~$S$ est coh\'erent.
\item
Pour $n$ grand, l'homomorphisme canonique $\widehat{f}^\ast\widehat{f}_\ast(\formelF(n))\to \formelF(n)$ est surjectif.
\end{enumeratei}
\end{corollaire}

Admettons le corollaire pour l'instant, et prouvons \Ref{XII.3.1}. Gr\^ace \`a la derni\`ere hypoth\`ese faite dans \Ref{XII.3.1}, on peut se ramener au cas o\`u $X=\PP^r_S$, quitte \`a remplacer $\OX(1)$, $t$ par une puissance convenable. Je dis qu'on peut supposer de plus que pour tout~$s$, on~a $t_s\neq 0$. Sinon, on~a en effet $\formelF_s=0$ par b), ou ce qui revient au m\^eme par Nakayama, $F_{0s}=0$ \ie $s$ n'appartient pas \`a l'image de $\Supp F_0$ par le morphisme $f_0:X_0\to S$ induit par $f$. Or cette image $S'$ est ouverte en vertu de a), b) puisque $F_0$ est plat par rapport \`a $S$, or il est \'evident qu'il suffit de prouver la conclusion de \Ref{XII.3.1} dans la situation obtenue en se restreignant au-dessus de $S'$, car le Module coh\'erent $F'$ sur $X|_{S'}$ obtenu par sera restriction d'un Module coh\'erent $F$ sur $X$, qui r\'epondra \`a la question. On peut donc supposer qu'en plus des hypoth\`eses a), b), c), les hypoth\`eses suivantes sont \'egalement v\'erifi\'ees:
\begin{enumerate}
\item[\textup{a$'$)}] $\OX$ est plat par rapport \`a $S$.
\item[\textup{b$'$)}] Pour tout $s\in S$, la section $t_s$ est $\Oo_{X_s}$-r\'eguli\`ere.
\item[\textup{c$'$)}] Pour tout $s\in S$, $\Oo_{X_{0s}}$ est de profondeur $\geq 2$ en les points ferm\'es de $X_{0s}$. (Il~suffit de choisir $X=\PP_S^r$ avec $r\geq 3$, ce qui est loisible).
\end{enumerate}
Or \Ref{XII.3.2} implique que l'on peut trouver un \'epimorphisme
\steco
\begin{equation} \label{eq:XII.17}
\widehat{L} \to \formelF \to 0,
\end{equation}
o\`u $L$ est un Module sur $X$ de la forme $f_\ast(G)(-n)$, $G$ \'etant un Module coh\'erent localement libre sur $S$: il suffit en effet; pour $n$ grand, de repr\'esenter le Module coh\'erent $\widehat{f}_\ast(\formelF)$ sur $S$ comme quotient d'un tel $G$. D'autre part, les\pageoriginale hypoth\`eses a), b), c) sur $f$, $t$ impliquent que $\widehat{L}$ satisfait aux m\^emes conditions a), b), c) que $\formelF$. On en conclut facilement qu'il en est de m\^eme du noyau de (\Ref{eq:XII.17}), auquel on peut donc appliquer le m\^eme argument, de sorte que $\formelF$ est repr\'esent\'e comme conoyau d'un homomorphisme
\steco
\begin{equation} \label{eq:XII.18}
\widehat{L'} \to \widehat{L},
\end{equation}
o\`u $L$, $L'$ sont des Modules localement libres sur $X$. Or en vertu de a$'$) et la deuxi\`eme partie de c$'$), et de \Ref{XII.2.1} ou \Ref{XII.2.4} au choix, l'homomorphisme (\Ref{eq:XII.18}) provient d'un homomorphisme $L'\to L$ de Modules sur $X$. Il suffit maintenant de prendre pour $F$ le conoyau de $L'\to L$, et on gagne.

Reste \`a prouver \Ref{XII.3.2}. Cela avait \'et\'e fait dans le
s\'eminaire par un exp\'edient un peu p\'enible, consistant \`a tout
interpr\'eter en termes de cohomologie sur le c\^{o}ne projetant
\'epoint\'e de $X$ relativement \`a $S$, pour pouvoir se ramener au
th\'eor\`eme \Ref{IX.2.1}. Une fa\c con plus directe et plus
satisfaisante (bien que substantiellement identique), me semble
maintenant la suivante. Elle consiste \`a noter que dans \Ref{IX},
\numero \Ref{IX.2}\refstepcounter{toto}\label{XII.14} (et avec les
notations de cet expos\'e) l'hypoth\`ese que le morphisme
$f:\Xx\to\Xx'$ soit adique n'intervient nulle part dans la
d\'emonstration de \Ref{IX.2.1}, via $\EGA 0_{\textup{III}}$ 13.7.7;
il suffit de supposer \`a la place que $\Xx$ est \'egalement adique,
et de choisir deux id\'eaux de d\'efinition $\Jj$ pour $\Xx'$, $\Ii$
pour $\Xx$, tels que $\Jj\Oo_\Xx\subset\Ii$, et de d\'efinir
$\Ss=\gr_\Jj(\Oo_{\Xx'})$, et consid\'erer $\gr_\Ii(\formelF)$. De
toute fa\c con, \Ref{IX.2.1} peut s'appliquer directement au
morphisme $\widehat{f}:\widehat{X}\to S$ consid\'er\'e dans le pr\'esent
num\'ero, o\`u on prend simplement $\Jj=0$. Ainsi pour v\'erifier que
$\widehat{f}_\ast(\formelF)$ est coh\'erent, il suffit en vertu de
\loccit de v\'erifier que $\R^if_{0\ast}(\gr_\Ii(\formelF))$ est
coh\'erent sur $S$ pour $i=0, 1$; pour ceci on note qu'en vertu de
a) et b), le Module envisag\'e n'est autre que
$\sisi{\coprod}{\bigoplus}_{m\geq 0}\R^if_{0\ast}(F_0(-m))$, qui
est bien coh\'erent en vertu de l'hypoth\`ese c) et de
\Ref{XII.1.5}.

Cela prouve \Ref{XII.3.2} (i). Pour \Ref{XII.3.2} (ii), nous aurons besoin du

\begin{lemme} \label{XII.3.3}
Sous\pageoriginaled les conditions de \Ref{XII.3.1a}), \Ref{XII.3.1b}), \Ref{XII.3.1c}) de \Ref{XII.3.1}, posons
\[
G_m= \widehat{f}_\ast(\formelF(\cdot)_m) = \sisi{\coprod}{\tbigoplus}_n\widehat{f}_\ast(\formelF_m(n))
\]
Alors le syst\`eme projectif $(G_m)$ satisfait la condition de Mittag-Leffler.
\end{lemme}

On peut supposer $S$ affine, d'anneau $A$. Soit alors $\Ss$ une $A$-alg\`ebre gradu\'ee de type fini \`a degr\'es positifs, et $t'\in\Ss_1$, tels que $X$ s'immerge dans $\Proj(\Ss)$, $\OX(1)$ \'etant induit par \sisi{$\SheafProj(\Ss(1))$}{$\Oo(1)$} et la section $t$ \'etant l'image de $t'$. Munissons $\Ss$ de la filtration $\Jj$-adique, o\`u $\Jj=t'\Ss$, et consid\'erons le syst\`eme projectif des $\formelF(\cdot)_m$ dans la cat\'egorie des faisceaux ab\'eliens sur $X_0$. On est encore sous les conditions pr\'eliminaires de $\EGA 0_{\textup{III}}$ 13.7.7\sfootnote{Rectifi\'e comme indiqu\'e dans \Ref{IX} p\ptbl\sisi{11}{\pageref{pagerectif}.}} et de plus $\H^i(X_0, \gr(\formelF(\cdot))$ est un Module de type fini sur $\gr_\Jj(\Ss)$ pour $i=0, 1$. En effet, comme $t'$ est $\formelF$-r\'egulier, on constate aussit\^{o}t qu'en tant que Module sur $(\Ss /t'\Ss)[T]$ (dont $\gr_\Jj(\Ss)$ est un quotient), le Module envisag\'e s'identifie \`a $\H^i(X_0, F_0(\cdot))\otimes_{\Ss /t'\Ss}(\Ss /t'\Ss)[T]$, or en vertu de \Ref{XII.1.5}, $\H^i(X_0, F_0(\cdot))$ est de type fini sur $\Ss$, donc sur $\Ss /t'\Ss$, pour $i=0, 1$, ce qui prouve notre assertion. Par suite on est sous les conditions d'application de $0_{\text{III}}$ 13.7.7 avec $n=1$, ce qui prouve \Ref{XII.3.3}.

Ce point acquis (et supposant toujours $S$ affine, ce qui est loisible pour prouver \Ref{XII.3.2} (ii)) soit $m_0$ tel que $m\geq m_0$ implique $\Im(G_m\to G_0)= \Im(G_{m_0}\to G_0)$, de sorte que les deux membres sont aussi \'egaux \`a $\Im(\varprojlim G_m\to G_0) = \Im(\widehat{f}_\ast(\formelF(\cdot))\to f_\ast\formelF_0(\cdot))$. Remarquons maintenant que pour $n$ grand, $\formelF_{m_0}(n)$ est engendr\'e par ses sections, donc $\formelF_0(n)$ est engendr\'e par des sections qui se remontent \`a $\formelF_{m_0}$, donc (gr\^ace au choix de~$m_0$) qui se remontent \`a $\formelF$. Donc les sections de $\formelF$ engendrent $\formelF_0$, donc aussi~$\formelF$ gr\^ace \`a Nakayama. Cela prouve \Ref{XII.3.2} (ii), donc \Ref{XII.3.1}.

\begin{corollaire} \label{XII.3.4}
Soient $f:X\to S$ un morphisme projectif et plat, avec $S$ localement noeth\'erien, $\OX(1)$ un Module inversible sur $X$, ample relativement \`a $S$, $t$ une section de ce Module telle que pour tout $s\in S$, la section $t_s$ induite sur\pageoriginale la fibre $X_s$ soit $\Oo_{X_s}$-r\'eguli\`ere, $X_0$ le sous-pr\'esch\'ema des z\'eros de $t$, $\widehat{X}$ le compl\'et\'e formel de $X$ le long de~$X_0$. On suppose que pour tout $s\in S$, $X_{0s}$ est de profondeur $\geq 2$ en ses points ferm\'es, (\ie $X_s$ est de profondeur $\geq 3$ aux points ferm\'es de $X_{0s}$) et $X_s$ est de profondeur $\geq 2$ en ses points ferm\'es. Sous ces conditions, le couple $(X, X_0)$ satisfait la condition de Lefschetz effective ($\Leff$) \sisi{de \Ref{X}~\Ref{X.2}}{du paragraphe~\Ref{X.2} de l'expos\'e~\emph{\Ref{X}}}, \ie:
\begin{enumeratea}
\item
Pour tout voisinage ouvert $U$ de $X_0$ dans $X$, le foncteur
\[
F \sisi{\rightsquigarrow}{\mto} \widehat{F}
\]
de la cat\'egorie des Modules coh\'erents localement libres sur $U$ dans la cat\'egorie des Modules coh\'erents localement libres sur $\widehat{X}$ est pleinement fid\`ele.
\item
Pour tout Module coh\'erent localement libre $\formelF$ sur $\widehat{X}$, il existe un voisinage ouvert~$U$ de $X_0$, et un Module coh\'erent localement libre $F$ sur $U$ tel que $\formelF$ soit isomorphe \`a~$\widehat{F}$.
\end{enumeratea}
\end{corollaire}

En effet, a) a d\'ej\`a \'et\'e not\'e dans \Ref{XII.2.4} sous des conditions plus faibles. Pour b), on applique \Ref{XII.3.1} qui donne la conclusion, du moins si $S$ est noeth\'erien et admet un Module inversible absolument ample, en particulier si $S$ est affine. En effet, si $F$ est un module coh\'erent sur $X$ tel que $\widehat{F}$ soit isomorphe \`a $\formelF$ donc localement libre, il s'ensuit que $F$ est localement libre sur un voisinage $U$ de $X_0$, et $\sisi{F|U}{F|_U}$ satisfera la condition voulue. Mais notons maintenant qu'en vertu de \Ref{XII.2.5}, pour un tel $F$, son image par l'immersion $U\to X$ est coh\'erente, et d'ailleurs ind\'ependante de la solution $(U, F)$ choisie (compte tenu \sisi{\ignorespaces}{de ce} que deux solutions co\"incident au voisinage de $X_0$, en vertu de a)). De fa\c con pr\'ecise, on peut trouver un Module coh\'erent $F$ sur $X$ et un isomorphisme $\widehat{F}\isomto \formelF$, tel que $F$ soit de profondeur $\geq 2$ en tous les points de $X$ qui sont ferm\'es dans leur fibre, et ceci d\'etermine $F$ \`a un isomorphisme unique pr\`es. Gr\^ace \`a cette propri\'et\'e d'unicit\'e, les solutions du probl\`eme qu'on trouve en s'induisant au-dessus des ouverts affines de $S$ se recollent, d'o\`u un $F$ coh\'erent sur $X$ tout entier et un isomorphisme $\widehat{F}\simeq \formelF$. Restreignant $F$ \`a l'ouvert $U$ des points en lesquels il est libre, on trouve ce qu'on a cherch\'e.

Gr\^ace\pageoriginale \`a \Ref{XII.2.4} et \Ref{XII.3.4}, on peut exploiter, dans la situation d'un sch\'ema alg\'ebrique projectif et d'une \og section hyperplane\fg dudit, les faits g\'en\'eraux \'etablis dans les expos\'es \Ref{X} et \Ref{XI} concernant les conditions (Lef) et (Leff). Ainsi:

\begin{corollaire} \label{XII.3.5} Soient $X$ un sch\'ema alg\'ebrique projectif muni d'un Module inversible ample $\OX(1)$, soit $t$ une section de ce Module qui soit $\OX$-r\'eguli\`ere, et soit $X_0$ le sous-sch\'ema des z\'eros de $t$. Supposons que $X$ soit de profondeur $\geq 2$ en ses points ferm\'es (\resp et de profondeur $\geq 3$ en les points ferm\'es de $X_0$). Alors $\pi_0(X_0)\to \pi_0(X)$ est bijectif, en particulier $X$ est connexe si et seulement si $X_0$ l'est, et choisissant un point-base g\'eom\'etrique dans $X_0$, $\pi_1(X_0)\to \pi_1(X)$ est surjectif, et plus g\'en\'eralement pour tout ouvert $U\supset X_0$, l'homomorphisme $\pi_1(X_0)\to \pi_1(U)$ est surjectif (\resp l'homomorphisme $\pi_1(X_0)\to \varprojlim_U\pi_1(U)$ est bijectif). Dans le cas resp\'e, si on suppose de plus que l'anneau local de tout point ferm\'e de $X$ non dans $X_0$ est pur (\Ref{X.3.2}) (par exemple est r\'egulier, ou seulement une intersection compl\`ete) alors $\pi_1(X_0)\to \pi_1(X)$ est un isomorphisme.
\end{corollaire}

On applique \Ref{X.2.4} et \Ref{XII.3.4}. On notera que dans le cas resp\'e, l'hypoth\`ese que $X$ soit de profondeur $\geq 3$ en les points ferm\'es de $X_0$, implique que toutes les composantes irr\'eductibles de dimension $\neq 0$ de $X$ sont de dimension $\geq 3$ (comme on voit en notant qu'une telle composante rencontre n\'ecessairement $X_0$, et en regardant en un point ferm\'e de l'intersection).

\begin{remarquestar}
Lorsque $X$ est normal, de dimension $\geq 2$ en tous ses points, il est de profondeur $\geq 2$ en ses points ferm\'es et on est sous les conditions non resp\'ees de \Ref{XII.3.5}. Dans ce cas, on~a une d\'emonstration plus \'el\'ementaire de la surjectivit\'e de $\pi_1(X_0)\to \pi_1(X)$ \`a l'aide du th\'eor\`eme de Bertini (\cf \SGA 1 X.2.10). Lorsqu'on suppose de plus $X_0$ normal, et $X$ de dimension $\geq 3$ en tous ses points, alors on est sous les conditions resp\'ees de \Ref{XII.3.5}. Dans ce cas, \Ref{XII.3.5} a \'et\'e \'etabli par Grauert (il se trouve en effet que gr\^ace \`a l'hypoth\`ese de normalit\'e, on arrive alors \`a se passer du th\'eor\`eme d'existence \Ref{XII.3.1} par des exp\'edients)\pageoriginale. C'est cette d\'emonstration de Grauert qui a \'et\'e le point de d\'epart de la \og th\'eorie de Lefschetz\fg qui fait l'objet du pr\'esent s\'eminaire.
\end{remarquestar}

\begin{corollaire} \label{XII.3.6} Soient $X$, $\OX(1)$, $t$, $X_0$ comme dans \Ref{XII.3.5}. Supposons que $X$ soit de profondeur $\geq 2$ en ses points ferm\'es, et que
\[
\H^i(X_0, \Oo_{X_0}(-n))=0
\]
pour $n>0$ et pour $i=1$ (\resp pour $i=1$ et pour $i=2$), ce qui implique en vertu de \Ref{XII.1.4} que $X_0$ est de profondeur $\geq 2$ (\resp $\geq 3$) en ses points ferm\'es, \ie que $X$ est de profondeur $\geq 3$ (\resp $\geq 4$) en les points ferm\'es de $X_0$. Sous ces conditions, pour tout voisinage ouvert $U$ de $X_0$, $\Pic(U)\to\Pic(X_0)$ est injectif, en particulier $\Pic(X)\to\Pic(X_0)$ est injectif, (\resp $\varprojlim_U \Pic(U)\to\Pic(X_0)$ est bijectif). Dans le cas resp\'e, si on suppose de plus que l'anneau local de $X$ en tout point ferm\'e non dans $X_0$ est parafactoriel (\Ref{XI.3.1}) (par exemple est r\'egulier, ou plus g\'en\'eralement une intersection compl\`ete), alors $\Pic(X)\to\Pic(X_0)$ est bijectif.
\end{corollaire}
On applique \Ref{XI}~\Ref{XI.3.12} et \Ref{XI.3.13}, en notant que l'hypoth\`ese resp\'ee implique que les composantes irr\'eductibles de dimension $\neq 0$ de $X$ sont de dimension $\geq 4$. On trouve en particulier, en appliquant ceci au cas o\`u $X$ est une intersection compl\`ete globale de dimension $\geq 4$ dans l'espace projectif:

\begin{corollaire} \label{XII.3.7} Soit $X$ un sch\'ema alg\'ebrique de dimension $\geq 3$, qui soit une intersection compl\`ete dans un sch\'ema $\PP_k^r$. Alors $\Pic(X)$ est le groupe libre engendr\'e par la classe du faisceau $\OX(1)$.
\end{corollaire}
On raisonne par r\'ecurrence sur le nombre d'hypersurfaces dont $X$ est l'intersection, en appliquant \Ref{XII.3.6} et notant que pour une intersection compl\`ete $X$ de dimension $\geq 3$, on $\H^i(X, \OX(n))=0$ pour $i=1, 2$, et tout $n$.

\begin{remarque} \label{XII.3.8}
Dans\pageoriginale le cas o\`u $X$ est une hypersurface non singuli\`ere, \Ref{XII.3.7} est d\^u \`a Andreotti. Le r\'esultat \Ref{XII.3.7} s'exprime aussi (lorsque $X$ est non singuli\`ere) en disant que l'anneau de coordonn\'ees homog\`enes de $X$ est factoriel, et sous cette forme est contenu dans \Ref{XI}~\Ref{XI.3.13} (ii). Signalons aussi que Serre avait donn\'e une d\'emonstration de \Ref{XII.3.7} dans le cas non singulier, par voie transcendante, en utilisant un argument de sp\'ecialisation pour se ramener au cas de la caract\'eristique 0, o\`u on dispose du th\'eor\`eme de Lefschetz sous sa forme classique. Bien entendu, le fait que la d\'emonstration purement alg\'ebrique donn\'ee ici permette de se d\'ebarrasser d'hypoth\`eses de non singularit\'e dans l'\'enonc\'e du th\'eor\`eme de Lefschetz, invite \`a reconsid\'erer celui-ci \'egalement dans le cas classique. \Cf l'expos\'e suivant qui propose des conjectures dans ce sens.

Dans les corollaires \Ref{XII.3.5} et \Ref{XII.3.7} nous nous sommes plac\'es sur un corps de base, alors que les th\'eor\`emes-clefs \Ref{XII.2.4} et \Ref{XII.3.4} sont valables sur une base quelconque. Pour g\'en\'eraliser \`a un $S$ g\'en\'eral ces corollaires \Ref{XII.3.5} et \Ref{XII.3.6}, nous devons donner des crit\`eres serviables pour qu'un point de $X$ (plat sur $S$) ait un anneau local \og pur\fg \resp parafactoriel. Ce sera l'objet du n\textsuperscript{o} suivant.
\end{remarque}

\section{Compl\'etion formelle et platitude normale} \label{XII.4}

\setcounter{equation}{18}

\begin{theoreme} \label{XII.4.1} Soit $X$ un pr\'esch\'ema localement noeth\'erien, localement immergeable dans un sch\'ema r\'egulier, $Y$ une partie ferm\'ee de $X$, $U=X- Y$, $X_0$ un sous-pr\'esch\'ema ferm\'e de $X$ d\'efini par un id\'eal $\Jj$, $\widehat{X}$ le compl\'et\'e formel de $X$ le long de $X_0$, $U_0$ la trace de $X_0$ sur $U$, $\widehat{U}$ le compl\'et\'e formel de $U$ le long de $U_0$, $i:U\to X$ et $\widehat{i}:\widehat{U}\to \widehat{X}$ les immersions canoniques, $n$ un entier. On suppose
\begin{enumeratea}
\item
$X$ est normalement plat le long de $X_0$ aux points de $Y\cap X_0$, \ie en ces points les \sisi{m}{M}odules $\Jj^n/\Jj^{n+1}$ sur $X_0$ sont plats \ie \sisi{\ignorespaces}{localement} libres.
\item
Pour tout $x\in Y\cap X_0$, on~a $\prof \Oo_{X_0, x}\geq n+2$.
\end{enumeratea}
Sous\pageoriginale ces conditions, on~a ce qui suit:
\begin{enumerate}
\sisi{\item[$1^\circ)$]}{\item} Soit $F$ un Module coh\'erent sur $U$, supposons que l'on ait:
\begin{enumerate}
\item[\textup{c)}] Pour tout $x\in Y- Y\cap X_0$, on~a $\prof \Oo_{X, x}\geq n+2$.

\item[\textup{d)}] $F$ est libre en les points de $U_0$, et de profondeur $\geq n+1$ en tous les points de $U$ o\`u il n'est pas libre.
\end{enumerate}
Alors le Module gradu\'e
\[
\sisi{\coprod}{\tbigoplus}_{m\geq 0}\R^pi_\ast(\Jj^mF)
\]
sur $\sisi{\coprod}{\bigoplus}_{m\geq 0}\Jj^m$ est de type fini pour $p\leq n$. \sisi{\item[$2^\circ$]}{\item} Soit $\formelF$ un Module coh\'erent sur $\widehat{U}$, alors le Module gradu\'e
\[
\sisi{\coprod}{\tbigoplus}_{m\geq 0}\R^pi_\ast(\Jj^m\formelF/\Jj^{m+1}\formelF)
\]
sur $\sisi{\coprod}{\bigoplus}_{m\geq 0}\Jj^m/\Jj^{m+1}=\gr_{\Jj}(\OX)$ est de type fini pour $p\leq n$.
\end{enumerate}
\end{theoreme}

\begin{proof}
1) Soit $X'\!=\!\Spec(\sisi{\coprod}{\bigoplus}_{m\geq 0}\Jj^m)$; le changement de base \hbox{$f:X'\!\to\! X$} d\'efinit alors $U'=X'- Y'$, $X_0'$, $U_0'=X_0'\sisi{\cap}{} U'$, et des immersions $i':U'\to X'$, \hbox{$i_0':U_0'\to X_0'$}. On a donc un carr\'e cart\'esien
\[
\xymatrix{
X'\ar[d]_f&\ar[l]_-{i'} U'\ar[d]^g \\
X&\ar[l]_-{i} U
}
\]
et on~a
\steco
\begin{equation} \label{eq:XII.19}
\sisi{\coprod}{\tbigoplus}_{m\geq 0}\R^pi_\ast(\Jj^mF) = \R^pi_\ast \Big(\sisi{\coprod}{\tbigoplus}_{m\geq 0} \Jj^mF \Big) =\R^pi_\ast (g_\ast(F')),
\end{equation}
o\`u\pageoriginale $F'=g^\ast F$, de sorte que l'on a bien un isomorphisme canonique
\[
g_\ast(F') \isomto \sisi{\coprod}{\tbigoplus}_{m\geq 0}\Jj^mF,
\]
car c'est vrai en les points de $U_0$, du fait que $F$ y est libre en vertu de d), et \'egalement en les points de $U_0$, du fait qu'on y a $\Jj^m=\Oo_U$, (de sorte que dans les deux cas, $\Jj^m\otimes_{\Oo_U}F\to \Jj^mF$ est un isomorphisme).

D'autre part, comme $f$ et par suite $g$ sont affines, on~a
\steco
\begin{equation} \label{eq:XII.20} {}\R^pi_\ast(g_\ast(F')) = \R^p(ig)_\ast(F') = \R^p(fi')_\ast(F') = f_\ast(R^pi'_\ast(F'))
\end{equation}
donc comparant (\Ref{eq:XII.19}) et (\Ref{eq:XII.20}), on voit que l'assertion 1\sisi{${}^\circ$}{}) est \'equivalente \`a la suivante: $R^pi'_\ast(F')$ est un Module de type fini \ie coh\'erent sur $X'$, pour tout $p\leq n$. Or comme $X$ est localement immergeable dans un sch\'ema r\'egulier, il en est de m\^eme de $X'$ qui est de type fini sur $X$, et on peut appliquer le crit\`ere de coh\'erence \sisi{\ignorespaces}{\Ref{VIII}}~\Ref{VIII.2.3} \`a un prolongement coh\'erent $F^{''}$ de $F'$: on veut exprimer que $\SheafH^p_Y(F^{''})$ est coh\'erent pour $p\leq n+1$, et cela \'equivaut aussi \`a dire que pour tout $x'\in U'$ tel que
\begin{equation*} \label{eq:XII.20bis} \tag{20 bis} {}\codim (\overline{x'}\cap Y', \overline{x'})=1,
\end{equation*}
on ait
\steco
\begin{equation} \label{eq:XII.21} {} \prof F'_{x'}\geq n+1.
\end{equation}
Or cette condition est v\'erifi\'ee en les points $x'$ o\`u $F'$ n'est pas libre, car pour un tel $x'$ on~a $x'\notin U_0'$ en vertu de d), donc $g$ y est un isomorphisme, et en vertu de d) encore, $F$ est de profondeur $\geq n+1$ en $g(x')$, donc $F'$ est de profondeur $\geq n+1$ en $x'$. Il suffit donc de v\'erifier la condition (\Ref{eq:XII.21}) en les $x'\in U'$ satisfaisant (\Ref{eq:XII.21}), et en lesquels $F'$ est libre. Or il suffit pour cela de prouver que l'on a\pageoriginale
\begin{equation*} \label{eq:XII.21bis} \tag{21 bis}
\prof \Oo_{X', x'}\geq n+1
\end{equation*}
\enlargethispage{\baselineskip}%
en ces points, \afortiori il suffit d'\'etablir que l'on a cette relation en tous les points~$x$ de $U'$ satisfaisant (\Ref{eq:XII.20}\sisi{}{ bis}). Or\sisi{ en vertu du crit\`ere cit\'e}{, \`a nouveau en vertu du crit\`ere~\Ref{VIII.2.3}} de l'expos\'e \Ref{VIII}, ceci \'equivaut \`a l'assertion que les Modules
\[
\SheafH^p_{Y'}(\Oo_{X'}) \hspace{5mm} \text{pour} \; p\leq n+1
\]
sont coh\'erents. En fait, nous allons prouver qu'ils sont m\^eme nuls, ou ce qui revient au m\^eme en vertu de l'expos\'e~\Ref{III}, que l'on a
\steco
\begin{equation} \label{eq:XII.22} {} \prof \Oo_{X', x'}\geq n+2 \hspace{5mm} \text{pour tout }\; x'\in Y'.
\end{equation}
Pour ceci, nous distinguons deux cas. Si $x'\notin X'_0$, alors $f$ est une immersion en $x'$, et il faut v\'erifier que $F$ est de profondeur \sisi{dans}{$\geq n+2$ en } l'image $x=f(x')$, ce qui n'est autre que la condition c). Si par contre $x'\in X'_0$, \ie $x=f(x')\in X_0$ donc $x\in Y\cap X_0$, on applique les conditions a) et b) gr\^ace au

\begin{lemme} \label{XII.4.2}
Soient $X$ un pr\'esch\'ema localement noeth\'erien, $X_0$ un sous-pr\'esch\'ema ferm\'e de $X$ d\'efini par un id\'eal $\Jj$, $X'=\Spec(\sisi{\coprod}{\bigoplus}_{m\geq 0} \Jj^m)$, $X_0'=\Spec(\sisi{\coprod}{\bigoplus}_{m\geq 0} \Jj^m/\Jj^{m+1}) = X'\times_XX_0$, $x$ un point de $X_0$ en lequel $X$ est normalement plat le long de $X_0$ \ie tel que $\gr_{\Jj}(\OX)= \sisi{\coprod}{\bigoplus}_{m\geq 0}\Jj^m/\Jj^{m+1}$ y soit plat comme Module sur $X_0$. Alors pour toute suite d'\'el\'ements $f_i$ ($1\leq i\leq m$) de $\Oo_{X, x}$ dont les images dans $\Oo_{X_0, x}$ forment une suite $\Oo_{X_0, x}$-r\'eguli\`ere, et tout $x'\in X'$ au-dessus de $x$, les images des $f_i$ dans $\Oo_{X', x'}$ (\resp dans $\Oo_{X'_0, x'}$) forment \'egalement une suite $\Oo_{X', x'}$-r\'eguli\`ere (\resp $\Oo_{X'_0, x'}$-r\'eguli\`ere), en particulier on~a
\steco
\begin{equation} \label{eq:XII.23}
\prof \Oo_{X', x'} \geq \prof \Oo_{X_0, x}, \quad \prof \Oo_{X'_0, x'} \geq \prof \Oo_{X_0, x}.
\end{equation}
\end{lemme}

Pour\pageoriginale le d\'emontrer, on peut supposer que $X$ est local de point ferm\'e $x$, donc affine d'anneau $A=\Oo_{X, x}$, $\Jj$ \'etant d\'efini par Id\'eal $J$, et il suffit de prouver que pour toute suite $f_i$ ($1\leq i\leq m$) d'\'el\'ements de $A$ dont les images dans $A/J$ forment une suite $A/J$-r\'eguli\`ere, les $f_i$ forment \'egalement une suite $(\sisi{\coprod}{\bigoplus}_{m\geq 0} J^m)$-r\'eguli\`ere et une suite $(\sisi{\coprod}{\bigoplus}_{m\geq 0} J^m/J^{m+1})$-r\'eguli\`ere, \ie pour tout $m$, elle forme une suite $J^m$-r\'eguli\`ere et $J^m/J^{m+1}$-r\'eguli\`ere. La deuxi\`eme assertion est triviale, puisque $J^m/J^{m+1}$ est un module libre sur $A/J$. La premi\`ere s'ensuit, en regardant la filtration $J$-adique de $J^m$, et notant que pour le module gradu\'e associ\'e \`a $J^m$ pour ladite filtration, la suite des~$f_i$ est r\'eguli\`ere.

Cela prouve \Ref{XII.4.2} et par suite \Ref{XII.4.1}~1\sisi{\noo}{}).

Prouvons \Ref{XII.4.1} 2\sisi{\noo}{}). Pour ceci, utilisons le carr\'e cart\'esien
\[
\xymatrix{
X'_0\ar[d]_{f_0}&\ar[l]_-{i'_0} U'_0\ar[d]^{g_0} \\
X_0&\ar[l]_-{i_0} U_0
}
\]
et proc\'edant comme dans le d\'ebut de la d\'emonstration de 1\sisi{\noo}{}), on trouve que l'on a
\steco
\begin{equation} \label{eq:XII.24}
\sisi{\coprod}{\tbigoplus}_{m\geq 0} \R^pi_\ast(\Jj^m\formelF/\Jj^{m+1}\formelF) \simeq {f_0}_\ast (\R^p{\sisi{i}{i'}_0}_\ast(\formelF_0')),
\end{equation}
o\`u $\formelF_0=\formelF/\Jj\formelF$ et o\`u $\formelF_0'=g_0^\ast(\formelF_0)$, (en utilisant le fait que $\formelF$ est localement libre). Donc la conclusion de 2\sisi{\noo}{}) \'equivaut \`a dire que pour $p\leq n$, $\R^p{i'_0}_\ast(\formelF_0')$ est un Module coh\'erent.

\enlargethispage{\baselineskip}%
Ici encore, compte tenu que $\formelF_0'$ est localement libre, le crit\`ere \Ref{VIII}~\Ref{VIII.2.3} nous permet de nous ramener \`a prouver qu'il en est ainsi si on remplace $\formelF_0'$ par $\Oo_{U'_0}$, \ie \`a prouver que les Modules\pageoriginale
\[
\SheafH^p_{Y_0'}(\Oo_{X_0'}) \hspace{5mm} \text{pour} \; p\leq n+1 \hspace{5mm} \text{(o\`u} \; Y'_0=Y'\cap X_0'=X_0'- U'_0 \text{)}
\]
sont coh\'erents. On prouve encore qu'ils sont en fait nuls, i.e; que l'on a
\steco
\begin{equation} \label{eq:XII.25} {}\prof \Oo_{X'_0, x'} \geq n+2 \hspace{5mm} \text{pour tout} \; x'\in Y'_0.
\end{equation}
Or ceci r\'esulte en effet des conditions a) et b), compte tenu de \Ref{XII.4.2}. Cela ach\`eve la d\'emonstration de \Ref{XII.4.1}.
\skipqed
\end{proof}

\begin{remarque} \label{XII.4.3}
On voit tout de suite, par descente, que l'hypoth\`ese: $X$ localement immergeable dans un sch\'ema r\'egulier, peut \^etre remplac\'ee par la suivante plus faible: il existe un morphisme $\overline{X}\to X$, fid\`element plat et quasi-compact, tel que $\overline{X}$ soit localement immergeable dans un sch\'ema r\'egulier.
\end{remarque}

Le th\'eor\`eme \Ref{XII.4.1} nous met en mesure d'appliquer les r\'esultats de l'expos\'e \Ref{IX} (th\'eor\`emes de comparaison et d'existence). Nous nous int\'eresserons notamment au

\begin{corollaire} \label{XII.4.4}
Supposons les conditions a), b), c) du th\'eor\`eme \Ref{XII.4.1} v\'erifi\'ees, avec $n=1$, et $X=\Spec(A)$, $A$ \'etant s\'epar\'e complet pour la topologie $J$-adique. Alors
\begin{enumerate}
\sisi{\item[\textup{1\noo})]}{\item} Le foncteur $F\mto \widehat{F}$ de la cat\'egorie des Modules coh\'erents localement libres sur~$U$, dans la cat\'egorie des Modules coh\'erents localement libres sur $\widehat{U}$, est pleinement fid\`ele. \sisi{\item[\textup{2\noo})]}{\item} Pour tout Module coh\'erent localement libre $\formelF$ sur $\widehat{U}$, il existe un Module coh\'erent $F$ sur $U$ et un isomorphisme $\widehat{F}\isomto \formelF$.
\end{enumerate}
En particulier, si pour tout $x\in U$ dont l'adh\'erence dans $U$ ne rencontre pas $U_0$ \ie tel que $\overline{x}\cap X_0\subset Y$, on~a $\prof \Oo_{U, x}\geq 2$, alors le couple $(U, U_0)$ satisfait la condition de Lefschetz effective ($\Leff$) de l'expos\'e \Ref{X}.
\end{corollaire}
\noindent
(Pour la derni\`ere assertion, on proc\`ede comme dans \Ref{X}~\Ref{X.2.1})

\pagebreak[2]
Cas particulier de \Ref{XII.4.4}:

\begin{corollaire} \label{XII.4.5}
Soient\pageoriginaled $A$ un anneau noeth\'erien, $J$ un id\'eal de $A$ contenu dans le radical, $A_0=A/J$. On suppose
\begin{enumeratei}
\item
$\prof A_0\geq 3$.
\item
$\gr_J(A)$ est un $A_0$-module libre.
\item
$A$ est complet pour la topologie $J$-adique.
\end{enumeratei}

Soient $X=\Spec(A)$, $X_0=\Spec(A_0)=V(J)$, $a$ le point ferm\'e de $X$, $U=X-\{a\}$, $U_0=X_0-\{a\}$, $\widehat{U}$ le compl\'et\'e formel de $U$ le long de $U_0$. Alors le foncteur $F\mto \widehat{F}$ de la cat\'egorie des Modules coh\'erents localement libres sur $U$ dans la cat\'egorie des Modules coh\'erents localement libres sur $\widehat{U}$, est pleinement fid\`ele. De plus, pour tout Module coh\'erent localement libre $\formelF$ sur $\widehat{U}$, il existe un Module coh\'erent (pas n\'ecessairement localement libre !) $F$ sur $U$, et un isomorphisme $\widehat{F}\simeq\formelF$.
\end{corollaire}

On notera que gr\^ace \`a \Ref{XII.4.3}, nous n'avons pas eu \`a supposer que $A$ est quotient d'un anneau r\'egulier, car le compl\'et\'e de $A$ pour la topologie $\rr(A)$-adique satisfait en tout cas \`a cette condition.

Proc\'edant comme dans les expos\'es \Ref{X} et \Ref{XI}, on conclut de \Ref{XII.4.5}

\begin{corollaire} \label{XII.4.6}
Sous les conditions de \Ref{XII.4.5}, on~a ceci:
\begin{enumeratea}
\item
$U$ et $U_0$ sont connexes (\Ref{III}~\Ref{III.3.1}).

Choisissant un point-base g\'eom\'etrique dans $U_0$, l'homomorphisme
\[
\pi_1(U_0)\to \pi_1(U)
\]
est surjectif.
\item
L'homomorphisme
\[
\Pic(U)\to \Pic(U_0)
\]
est injectif.
\end{enumeratea}
\end{corollaire}

Pour\pageoriginale prouver b), compte tenu de \Ref{XII.4.5}, cela revient \`a v\'erifier que tout isomorphisme $L_0'\isomto L_0$ se remonte en un isomorphisme $\widehat{L'}\isomto \widehat{L}$. Or pour ceci on remonte de proche en proche en des isomorphismes $L_n'\isomto L_n$, les obstructions se trouvent dans $\H^1(U_0, \Jj^n/\Jj^{n+1})$, or ces modules sont nuls du fait que $J^n/J^{n+1}$ est libre et $\prof A_0\geq 3$.

Nous sommes maintenant en mesure de prouver le
\enlargethispage{-1\baselineskip}%
\begin{theoreme} \label{XII.4.7}
Soient $A$ un anneau local noeth\'erien, $J$ un id\'eal de $A$ contenu dans son radical, $A_0=A/J$. On suppose
\begin{enumeratei}
\item
$\prof A_0\geq 3$.
\item
$\gr_J(A)$ est un module libre sur $A_0$.
\end{enumeratei}
\noindent
Alors, si $A_0$ est \og pur\fg (\Ref{X}~\Ref{X.3.1}) (\resp parafactoriel (\Ref{XI}~\Ref{XI.3.1})), il en est de m\^eme de~$A$.
\end{theoreme}

\begin{proof}
Par descente, on peut supposer qu'on a aussi
\begin{enumeratei}\setcounter{enumi}{2}
\item
$A$ est complet pour la topologie $J$-adique.
\end{enumeratei}

En effet, en vertu de (i) et (ii), on~a $\prof(A)\geq 3$ donc $\prof(\widehat{A})\geq 3$, o\`u $\widehat{A}$ est le compl\'et\'e de $A$ pour la topologie $J$-adique, et on applique \Ref{X}~\Ref{X.3.6} et \Ref{XI}~\Ref{XI.3.6}. On est donc sous les conditions de \Ref{XII.4.5}. Comme $\prof(A)\geq 3\geq 2$, dire que $A$ est parafactoriel signifie simplement que $\Pic(U)=0$, et en vertu de \Ref{XII.4.6} b) il suffit pour ceci que $\Pic(U_0)=0$, \ie que $A_0$ soit parafactoriel. Pour prouver que $A$ est \og pur\fg si $A_0$ l'est, il faut prouver que si $V$ est un rev\^etement \'etale de $U$, d\'efini par une alg\`ebre $\Bb$ sur~$U$, alors $\H^0(U, \Bb)$ est une alg\`ebre finie \'etale sur $A$. Or $A_0$ \'etant pur, il en est de m\^eme des $A_n$ (qui n'en diff\`erent que par des \'el\'ements nilpotents), donc pour tout $n$, $\sisi{}{B_n=}\H^0(U, \Bb_n)$ est une alg\`ebre \'etale sur $A/J^{n+1}$, et bien entendu ces alg\`ebres se recollent, de sorte que $\varprojlim B_n$ est une alg\`ebre \'etale sur $A$. Or en vertu de \Ref{XII.4.5}, cette alg\`ebre n'est autre que $\H^0(U, \Bb)$, ce qui \'etablit notre assertion.
\skipqed
\end{proof}

\begin{corollaire} \label{XII.4.8}
Soient\pageoriginale $f:X\to Y$ un morphisme plat de pr\'esch\'emas localement noeth\'eriens, $x\in X$, $y=f(x)$, on suppose que $\Oo_{X_y, x}$ est un anneau local \og pur\fg (\resp parafactoriel) de profondeur $\geq 3$, alors il en est de m\^eme de $\Oo_{X, x}$.
\end{corollaire}

C'est le r\'esultat du type promis \`a la fin du \sisi{N\textsuperscript{o}}{\numero} pr\'ec\'edent, pour g\'en\'eraliser les corollaires \Ref{XII.3.5} et suivants. On trouve donc, en utilisant \Ref{XII.3.4}, le

\begin{corollaire} \label{XII.4.9}
Soient $f:X\to S$ un morphisme projectif et plat, avec $S$ localement noeth\'erien, $\OX(1)$ un Module inversible sur $X$ ample par rapport \`a $S$, $t$ une section de $\OX(1)$ telle que pour tout $s\in S$, la section $t_s$ induite sur $X_s$ soit $\Oo_{X_s}$-r\'eguli\`ere, $X_0$ le sous-sch\'ema des z\'eros de $t$, $X_m$ le sous-\sisi{pr\'e}{}sch\'ema des z\'eros de $t^{m+1}$. On suppose que pour tout $s\in S$, $X_s$ est de profondeur $\geq 3$ en tous ses points ferm\'es. Alors:
\begin{enumeratea}
\item
Si les anneaux locaux des points ferm\'es de $X_s- X_{0, s}$ ($s\in S$) sont \og purs\fg\, par exemple sont des intersections compl\`etes, alors le foncteur $X'\mto X_0'=X'\times_X X_0$ de la cat\'egorie des rev\^etements \'etales de $X$ dans la cat\'egorie des rev\^etements \'etales de $X_0$ est une \'equivalence de cat\'egories; en particulier, choisissant un point base g\'eom\'etrique dans $X_0$, l'homomorphisme
\[
\pi_1(X_0)\to \pi_1(X)
\]
est un isomorphisme\refstepcounter{toto}\nde{\label{FH}signalons
le spectaculaire r\'esultat de connexit\'e obtenu depuis par Fulton et
Hansen, dans le cas o\`u $S=\Spec(k)$ ($k$ corps alg\'ebriquement
clos). Soit $g: X\to {\PP^m_k}\times{\PP^m_k}$ tel que $\dim
g(X)>m$; alors, l'image inverse de la diagonale est connexe. Ceci
permet entre autres choses de g\'en\'eraliser le
corollaire~\Ref{XII.4.9} lorsque $f$ est le morphisme structural
du projectif $\PP^m_k$ sur $S=\Spec(k)$: pr\'ecis\'ement, une
sous-vari\'et\'e irr\'eductible $X$ de ${\PP^m_k}$ de dimension $>m/2$ a
un groupe fondamental trivial ! (\cf Fulton~W. \& Hansen J., {\og
A connectedness theorem for projective varieties, with
applications to intersections and singularities of mappings\fg},
\emph{Ann. of Math. (2)} \textbf{110} (1979), \numero 1, p\ptbl
159-166). Pour des g\'en\'eralisations au cas des grassmanniennes ou
des vari\'et\'es ab\'eliennes, voir {Debarre~O.}, {\og Th\'eor\`emes de
connexit\'e pour les produits d'espaces projectifs et les
grassmanniennes\fg}, \emph{Amer. J.~Math.} \textbf{118} (1996),
\numero 6, p\ptbl 1347-1367 et {\og Th\'eor\`emes de connexit\'e et
vari\'et\'es ab\'eliennes\fg}, \emph{Amer. J.~Math.} \textbf{117}
(1995), \numero 3, p\ptbl 787--805. Le r\'esultat de trivialit\'e du
groupe fondamental de $X$ comme plus haut a \'et\'e obtenu
ind\'ependamment par Faltings, qui prouve en outre que
$\text{Pic}(X)$ n'a pas de torsion premi\`ere \`a la caract\'eristique
de $k$, ce par des m\'ethodes d'alg\'ebrisation de fibr\'es formels,
plus dans la lign\'ee des techniques de Grothendieck, \cf
(Faltings~G., \og{Algebraization of some formal vector bundles\fg}
\emph{Ann. of Math. (2)} \textbf{110} (1979), \numero 3, p\ptbl
501--514).}

\item
Si les anneaux locaux des points ferm\'es des $X_s- X_{0, s}$ ($s\in S$) sont \og parafactoriels\fg, par exemple r\'eguliers, ou des intersections compl\`etes de dimension $\geq 4$, alors pour tout entier $m$ tel que $\R^if_{0\ast}(\Oo_{X_0}(-n))=0$ pour $n>m$ et $i=1, 2$, l'application $\Pic(X)\to \Pic(X_m)$ est bijective.

D'ailleurs si $S$ est noeth\'erien, et les $ X_{0, s}$ sont de profondeur $\geq 3$ en leurs points ferm\'es, il existe de tels $m$ (\cf \Ref{XII.1.5}).
\end{enumeratea}
\end{corollaire}

\begin{remarque} \label{XII.4.10}
Sous\pageoriginale les conditions de la derni\`ere assertion de \Ref{XII.4.9} b), on~a vu dans \Ref{XII.1.5} qu'il existe un $m$ tel que $n>m$ implique m\^eme $\H^i(X_0, \Oo_{X_0}(-n))=0$ pour $i=1, 2$, et m\^eme pour $i\leq 2$). Cette condition est plus forte que $\R^if_{\ast}(\Oo_{X_0}(-n))=0$ pour $i=1, 2$, et elle a de plus l'avantage d'\^etre stable par changement de base. Il en est de m\^eme des hypoth\`eses de profondeur qu'on a faites dans \Ref{XII.4.9}, et \'egalement d'une hypoth\`eses du type \og les $X_s$ sont localement des intersections compl\`etes\fg. Il s'ensuit alors, sous ces conditions que \Ref{XII.4.9} b) implique \'egalement que le morphisme de foncteurs
\[
\SheafPic_{X/S}\to \SheafPic_{X_n/S}
\]
dans $\catSch_{/S}$ est un isomorphisme, donc aussi le morphisme pour les sch\'emas de Picard relatifs, lorsque ceux-ci existent:
\[
\Pic_{X/S}\to \Pic_{X_n/S}.
\]
M\^eme dans le cas o\`u $S$ est le spectre d'un corps alg\'ebriquement clos, cet \'enonc\'e est nettement plus pr\'ecis que l'\'enonc\'e consistant \`a dire seulement que $\Pic(X)\to \Pic(X_n)$ est bijectif.

On peut se demander si on peut toujours prendre $n=0$ dans les conclusions pr\'ec\'edentes (supposant donc les $X_{0, s}$ de profondeur $\geq 3$ en leurs points ferm\'es). Lorsque $X_0$ est lisse sur $S$ et les caract\'eristiques r\'esiduelles de $S$ sont nulles, il en bien ainsi, en vertu du \og vanishing theorem\fg de Kodaira, (d\'emontr\'e par voie transcendante, utilisant une m\'etrique k\"ahl\'erienne) qui implique que pour tout sch\'ema projectif lisse connexe de dimension $n$ sur un corps $k$ de caract\'eristique nulle, et tout Module inversible ample $L$ sur $X$, on~a $\H^i(X, L^{-1})=0$ pour $i\neq n$. Il n'est pas connu\nde{comme l'a remarqu\'e Raynaud, le r\'esultat de d\'ecomposition du complexe de de Rham de Deligne et Illusie entra\^ine facilement la nullit\'e du groupe $H^i(X, L^{-1})$ (avec $L$ ample sur $X$ projective lisse sur $k$ de caract\'eristique $p>0$) pour $i<\inf(p, \dim(X))$ d\`es lors qu'on suppose que $X$ est relevable en un sch\'ema plat sur $W_2(k)$ (\cf Deligne~P. \& Illusie L., {\og Rel\`evements modulo~$p^2$ et d\'ecomposition du complexe de de Rham\fg}, \emph{Invent. Math.} \textbf{89} (1987), \numero 2, p\ptbl 247--270); ceci donne un preuve purement alg\'ebrique du r\'esultat de Kodaira pour les vari\'et\'es projectives en caract\'eristique nulle. Si $X$ n'est pas relevable, il est bien connu que le \og vanishing theorem\fg est faux; \cf l'exemple dans (Raynaud~M., {\og Contre-exemple au "vanishing theorem" en caract\'eristique $p>0$\fg}, in \emph{C.P.~Ramanujam---a tribute}, Tata Inst. Fund. Res. Studies in Math., vol.~8, Springer, Berlin-New York, 1978, p\ptbl 273--278); voir aussi les tr\`es jolis exemples dans (Haboush~W. \& Lauritzen~N., {\og Varieties of unseparated flags\fg}, in \emph{Linear algebraic groups and their representations (Los Angeles, CA, 1992)}, Contemp. Math., vol.~153, American Mathematical Society, Providence, RI, 1993, p\ptbl 35--57), simplifi\'es dans (Lauritzen~N. \& Rao~A.P., {\og Elementary counterexamples to Kodaira vanishing in prime characteristic\fg}, \emph{Proc. Indian Acad. Sci. Math. Sci.} \textbf{107} (1997), \numero 1, p\ptbl 21--25). En revanche, je ne connais pas d'exemple o\`u la fl\`eche $\Pic(X_{n+1})\to\Pic(X_n)$ n'est pas surjective pour $n>1$ en caract\'eristique positive, o\`u $X_n$ d\'esigne une section hyperplane \'epaissie de $X$ projective lisse comme plus haut.} \`a l'heure actuelle si ce th\'eor\`eme peut \^etre remplac\'e par un g\'en\'eralisation en caract\'eristique $p>0$, et si l'hypoth\`ese de lissit\'e peut \^etre remplac\'ee par une hypoth\`ese de nature plus g\'en\'erale (portant sur la profondeur, ou du type \og intersection compl\`ete\fg...).
\end{remarque}

\section[Conditions de finitude universelles]
{Conditions de finitude universelles pour un morphisme non propre} \label{XII.5} Rappelons\pageoriginale pour m\'emoire la \setcounter{equation}{25}

\begin{proposition} \label{XII.5.1}
Soit $f:X\to S$ un morphisme propre de pr\'esch\'emas, avec $S$ localement noeth\'erien, $U$ une partie ouverte de $X$, $g:U\to X$ l'immersion canonique, $h=fg:U\to S$, $F$ un Module sur $U$. Supposons que les Modules $\R^ig_\ast(F)$ soient coh\'erents pour $i\leq n$ (hypoth\`ese de nature locale sur $X$, qui se v\'erifie pratiquement \`a l'aide du crit\`ere \Ref{VIII}~\Ref{VII.2.3}). Alors $\R^ih_\ast(F)$ est coh\'erent pour $i\leq n$.
\end{proposition}

Cela r\'esulte aussit\^{o}t de la suite spectrale de Leray
\steco
\begin{equation} \label{eq:XII.26} {} \E_2^{p, q}=\R^pf_\ast(\R^qg_\ast(F)) \To \R^\ast h_\ast(F),
\end{equation}
et du fait que les images directes sup\'erieures par $f$ d'un Module coh\'erent sur $X$ sont coh\'erentes (\EGA III 3.2.1).

\begin{proposition} \label{XII.5.2}
Soit $S$ un pr\'esch\'ema localement noeth\'erien, $\Ss$ une Alg\`ebre gradu\'ee quasi-coh\'erente, de type fini sur $S$, engendr\'ee par $\Ss_1$, $X$ un sous-pr\'esch\'ema de $\Proj(\Ss)$, $\OX(1)$ le Module inversible sur $X$ tr\`es ample relativement \`a $S$ induit par $\Proj(\Ss(1))$, $U$ une partie ouverte de $X$, $g:U\to X$ l'immersion canonique, $h=fg:U\to S$, $F$ un Module quasi-coh\'erent sur $U$, d'o\`u des Modules tordus $F(m)=F\otimes \OX(m)$) ($m\in\ZZ$), $n$ un entier, $m_0$ un entier. Les conditions suivantes sont \'equivalentes:
\begin{enumeratei}
\item
$\R^ig_\ast(F)$ est coh\'erent pour $i\leq n$.
\item
$\sisi{\coprod}{\bigoplus}_{m>m_0}\R^ih_\ast(F(m))$ est un $\Ss$-Module de type fini pour $i\leq n$.
\end{enumeratei}
\end{proposition}

\begin{proof}
Rempla\c cant dans la suite spectrale ci-dessus $F$ par $F(m)$ on
trouve une suite spectrale de $\Ss$-Modules gradu\'es\pageoriginale
\[
\E_2^{p, q}= \sisi{\coprod}{\tbigoplus}_{m\geq m_0} \R^pf_\ast(\R^qg_\ast(F(m)) \To \sisi{\coprod}{\tbigoplus}_{m\geq m_0} \R^\ast h_\ast(F(m)).
\]
Comme on a
\[
\R^qg_\ast(F(m)) \simeq \R^qg_\ast(F) (m),
\]
on voit que si les $\R^ig_\ast(F)$ sont coh\'erents, $\E_2^{p, q}$ est de type fini sur $\Ss$ pour $q\leq n$, gr\^ace \`a la partie a) du lemme \Ref{XII.5.3} ci-dessous, ce qui implique que l'aboutissement est de type fini sur $\Ss$ en degr\'e $i\leq n$. Cela prouve (i) \ALORS (ii). De plus, raisonnant dans la cat\'egorie ab\'elienne des $\Ss$-Modules gradu\'es modulo la sous-cat\'egorie \'epaisse $\mathcal{C}$ de ceux qui sont quasi-coh\'erents de type fini, on trouve par la suite spectrale pr\'ec\'edente
\[
\sisi{\coprod}{\tbigoplus}_{m\geq m_0} \R^{n+1} h_\ast(F(m)) \simeq \sisi{\coprod}{\tbigoplus}_{m\geq m_0} f_\ast (\R^{n+1} g_\ast(F)(m)) \mod \mathcal{C},
\]
ce qui prouve que si le membre de gauche est un $\Ss$-Module de type fini, alors $\R^{n+1} g_\ast(F)$ est coh\'erent, en vertu de la partie b) du lemme \Ref{XII.5.3}. Ceci prouve l'implication (ii) \ALORS (i) par r\'ecurrence sur $n$. Reste \`a prouver:

\begin{lemme} \label{XII.5.3} Soient $S$, $\Ss$, $X$, $f$ comme dans \Ref{XII.5.2}, et $G$ un Module quasi-coh\'erent sur $X$, $m_0$ un entier. Alors
\begin{enumeratea}
\item
Si $G$ est coh\'erent, alors pour tout entier $i$, le Module gradu\'e
\[
\sisi{\coprod}{\tbigoplus}_{m\geq m_0} \R^i f_\ast(G(m))
\]
sur $\Ss$ est de type fini.
\item
Inversement, supposons que le Module $\sisi{\coprod}{\bigoplus}_{m\geq m_0} \R^i f_\ast(G(m))$ sur $\Ss$ soit de type fini, alors $G$ est coh\'erent.
\end{enumeratea}
\end{lemme}

\begin{proof}[D\'emonstration de \Ref{XII.5.3}]
Pour a), le cas $i=0$ est donn\'e dans \EGA III 2.3.2, le cas $i>0$ dans \EGA III 2.2.1 (i)(ii) qui dit que les $ \R^i f_\ast G(m)$ sont coh\'erents\pageoriginale, et nuls pour $m$ grand (si on suppose $S$ noeth\'erien, ce qui est loisible). Pour b), on note que $G$ est isomorphe \`a $\SheafProj (\sisi{\coprod}{\bigoplus}_{m\geq m_0}f_\ast(G(m))$ (\EGA II 3.4.4 et 3.4.2), ce qui prouve que $G$ est coh\'erent si $\sisi{\coprod}{\bigoplus}_{m\geq m_0}f_\ast(G(m))$ est de type fini sur $\Ss$, en vertu de \loccit 3.4.4.
\skipqed
\end{proof}
\skipqed
\end{proof}

\begin{corollaire}[{{\normalfont$S$ noeth\'erien}}] \label{XII.5.4}
Supposons que $\R^ig_\ast(F)$ soit coh\'erent pour $i\leq n$, alors pour $i\leq n+1$, et $m$ grand, on~a un isomorphisme canonique:
\[
\R^ih_\ast(F(m)) \simeq f_\ast(\R^ig_\ast(G)(m)).
\]
\end{corollaire}

En effet la suite spectrale (\Ref{eq:XII.26}) pour $F(m)$ d\'eg\'en\`ere alors en degr\'e $\leq n$, par \EGA III 2.2.1~(ii), d'o\`u aussit\^{o}t le r\'esultat (qui redonne d'ailleurs l'implication (ii) \ALORS (i) de~5.2).

\enlargethispage{1.2\baselineskip}%
\begin{corollaire} \label{XII.5.5}
Sous les conditions pr\'eliminaires de \Ref{XII.5.2}, $S$ noeth\'erien, les conditions suivantes sont \'equivalentes:
\begin{enumeratei}
\item
$\sisi{\coprod}{\bigoplus}_{m\geq m_0}h_\ast(F(m)$ est de type fini sur $S$, et $\R^ih_\ast(F(m))=0$ pour $0<i\leq n$ et $m$ grand.
\item
$g_\ast(F)$ est coh\'erent, et $\R^ig_\ast(F)=0$ pour $0<i\leq n$.
\item[\textup{(ii bis)}] $g_\ast(F)$ est coh\'erent, et $\prof_Y g_\ast(F)>n+1$.
\end{enumeratei}
\end{corollaire}

L'\'equivalence de (ii) et (ii bis) est contenue dans \Ref{III}~\Ref{III.3.3}. D'ailleurs, en vertu de \Ref{XII.5.2} les conditions (i) et (ii) impliquent toutes deux que les $\R^ig_\ast(F)$ ($i\leq n$) sont coh\'erents. L'\'equivalence de (i) et (ii) r\'esulte alors de \Ref{XII.5.4}, compte tenu du fait que tout un Module coh\'erent $G$ sur $X$, on~a $G=0$ si et seulement si $f_\ast(G(m))=0$ pour $m$ grand, par exemple en vertu de \EGA III 2.2.1 (iii).

\begin{remarque} \label{XII.5.6}
On\pageoriginale peut interpr\'eter les crit\`eres \Ref{XII.5.2} et \Ref{XII.5.5} en disant que la \og condition de finitude simultan\'ee\fg \Ref{XII.5.2} (ii) s'exprime par des propri\'et\'es de r\'egularit\'e locale (en termes de profondeur gr\^ace \`a \Ref{VIII}~\Ref{VIII.2.1}) de $F$ en les points de $U$ voisins de $Y=X- U$, alors que la \og condition de nullit\'e asymptotique\fg \Ref{XII.5.5} (i) est de nature nettement plus forte, et s'exprime par des conditions de r\'egularit\'e locale de $g_\ast(F)$ en les points de~$Y$ lui-m\^eme. Il serait int\'eressant, pour g\'en\'eraliser les th\'eor\`emes \`a la Lefschetz pour les morphismes projectifs aux morphismes quasi-projectifs, de trouver des crit\`eres locaux sur $X$ n\'ecessaires et suffisants pour que les $\Ss$-Modules $\sisi{\coprod}{\bigoplus}_{m\geq 0} \R^ih_\ast(F(m))$ pour $i\leq n$ soient de type fini. Lorsque $S$ est le spectre d'un corps (et sans doute plus g\'en\'eralement, si c'est le spectre d'un anneau artinien) et $Y=X- U$ est fini, on peut montrer qu'il est n\'ecessaire et suffisant que les conditions suivantes soient v\'erifi\'ees:
\begin{enumerate}
\item[\textup{1\noo)}] $\prof F_x >n$ pour tout point ferm\'e $x$ de $U$ (comparer \Ref{XII.1.4}).
\item[\textup{2\noo)}] $\R^ig_\ast(F)$ est coh\'erent pour $i\leq n$, ou ce qui revient au m\^eme, il existe un voisinage ouvert $V$ de $Y$ tel que pour tout point ferm\'e $x$ de $U\cap V$, on ait $\prof F_x >n+1$.
\end{enumerate}
\end{remarque}

\begin{proposition} \label{XII.5.7}
Soient $S$ un pr\'esch\'ema localement noeth\'erien, $g:U\to X$ un morphisme de pr\'esch\'emas de type fini sur $S$\sfootnote{Il suffit en fait que $g$ soit quasi-compact et quasi-s\'epar\'e (\textup{\EGA IV1.2.1}), sans condition sur $U, X$.}, de morphismes structuraux $h$ et $f$, $F$ un Module quasi-coh\'erent sur $U$, $n$ un entier. Les conditions suivantes sont \'equivalentes:
\begin{enumeratei}
\item
Pour tout changement de base $S'\to S$, avec $S'$ noeth\'erien, le module $\R^ng'_\ast(F')$ sur $X'$ est coh\'erent.
\item
Pour tout changement de base comme ci-dessus, et tout Id\'eal coh\'erent $\Jj$ sur~$S'$, d\'esignant par $\Ii$ l'Id\'eal $\Jj\Oo_{X'}$ sur $X'$, le Module gradu\'e
\[
\sisi{\coprod}{\tbigoplus}_{m\geq 0} \R^ng'_\ast(\Ii^mF')
\]
sur $\sisi{\coprod}{\bigoplus}_{m\geq 0} \Ii^m$ est de type fini.
\item
Pour\pageoriginale tout changement de base $S'\to S$, et $\Jj$ comme ci-dessus, le Module gradu\'e
\[
\sisi{\coprod}{\tbigoplus}_{m\geq 0} \R^ng'_\ast(\Ii^mF'/\Ii^{m+1}F')
\]
sur $\gr_I(\Oo_{X'})=\sisi{\coprod}{\bigoplus}_{m\geq 0} \Ii^m/\Ii^{m+1}$ est de type fini.
\end{enumeratei}
\end{proposition}

\'Evidemment (ii) \ALORS (i) et (iii) \ALORS (i), comme on voit en faisant $\Ii=0$ dans les conditions (ii) et (iii). Les implications inverses s'obtiennent en appliquant (i) au changement de base compos\'e $S^{''}\to S'\to S$, o\`u $S^{''}$ est \'egal \`a $\Spec \sisi{\coprod}{\bigoplus}_{m\geq 0} \Jj^m$ \resp $\Spec \sisi{\coprod}{\bigoplus}_{m\geq 0} \Jj^m/\Jj^{m+1}$.

L'int\'er\^et de cette proposition est que les conditions de la forme (ii) sont celles qui interviennent dans les \og th\'eor\`emes de comparaison alg\'ebrico-formels\fg, alors que les conditions de la forme (iii) interviennent dans les \og th\'eor\`emes d'existence\fg qui les compl\`etent, \cf expos\'e \Ref{IX}. Un premier cas int\'eressant est celui o\`u $f:X\to S$ est l'identit\'e, et o\`u il s'agit donc de conditions sur un morphisme $h:U\to S$ localement de type fini et un Module $F$ quasi-coh\'erent sur $U$ plat par rapport \`a $S$. Pour obtenir des conditions suffisantes, nous allons supposer que $U$ se plonge par $g:U\to X$, comme sous-pr\'esch\'ema ouvert d'un $X$ \emph{propre} sur $S$. Appliquant \Ref{XII.5.1}, on voit donc:

\begin{corollaire} \label{XII.5.8}
Soient $f:X\to S$ un morphisme propre, avec $S$ localement noeth\'erien, $U$ un ouvert de $X$, $g:U\to X$ l'immersion canonique, $h=fg:U\to S$, $F$ un Module quasi-coh\'erent sur $X$, plat par rapport \`a $S$. Supposons que pour tout changement de base $S'\to S$, avec $S'$ localement noeth\'erien, on ait $\R^ig'_\ast(F')$ coh\'erent sur $X'$ pour $i\leq n$. Alors on~a ce qui suit:
\begin{enumeratei}
\item
Pour tout changement de base $S'\to S$, avec $S'$ localement noeth\'erien, $\R^ih'_\ast(F')$ est coh\'erent sur $S'$ pour $i\leq n$.
\item
Pour tout $S'\to S$ comme ci-dessus, et tout Id\'eal coh\'erent $\Jj$ sur $S'$, les Modules gradu\'es\pageoriginale
\[
\sisi{\coprod}{\tbigoplus}_{m\geq 0} \R^ih'_\ast(\Jj^mF')
\]
sur $\sisi{\coprod}{\bigoplus}_{m\geq 0}\Jj^m$ sont de type fini pour $i\leq n$.
\item
Pour tout $S'\to S$ et $\Jj$ comme ci-dessus, les Modules gradu\'es
\[
\sisi{\coprod}{\tbigoplus}_{m\geq 0} \R^ih'_\ast(\Jj^mF'/\Jj^{m+1}F')
\]
sur $\gr_\Jj(\Oo_{S'})=\sisi{\coprod}{\bigoplus}_{m\geq 0}\Jj^m/\Jj^{m+1}$ sont de type fini pour $i\leq n$.

De plus, sous les conditions de (ii), et en vertu du th\'eor\`eme de comparaison \Ref{IX.1.1}, d\'esignant par $\widehat{S'}$ le compl\'et\'e formel de $S'$ par rapport \`a $\Jj$, et par $\widehat{U'}$ celui de $U'$ par rapport \`a $\Jj\Oo_{U'}$, les homomorphismes canoniques
\[
\widehat{\R^ih'_\ast(F')} \to \R^i\widehat{h'}_\ast(\widehat{F'}) \to \varprojlim_k\R^ih'_\ast(F'_k)
\]
sont des isomorphismes pour $i\leq n-1$.
\end{enumeratei}
\end{corollaire}

\begin{remarque} \label{XII.5.9}
Supposons qu'on soit de plus sous les conditions de \Ref{XII.5.8} avec $F$ \emph{coh\'erent}, et consid\'erons un changement de base $S'\to S$ comme dans \Ref{XII.5.9} (i). Supposons de plus que $S'$ soit localement immergeable dans un sch\'ema r\'egulier, ou plus g\'en\'eralement, qu'il existe un morphisme $S^{''}\to S$ fid\`element plat et quasi-compact, tel que $S^{''}$ soit localement immergeable dans un sch\'ema r\'egulier; cette condition est v\'erifi\'ee en particulier si $S'$ est local. Alors la conclusion de \Ref{XII.5.8} (i) et (ii) reste valable lorsqu'on y remplace $F'$ par un Module $G'$ sur $U'$, tel que tout point de $U'$ ait un voisinage ouvert sur lequel $G'$ soit isomorphe \`a un Module de la forme ${F'}^n$. En effet, on est ramen\'e au cas $S'$ lui-m\^eme est localement immergeable dans un sch\'ema r\'egulier, de sorte qu'il en est de m\^eme de $S^{''}=\Spec\big(\sisi{\coprod}{\bigoplus}_{m\geq 0}\Jj^m\big)$ et de $X^{''}=X'\times_{S'} S^{''}= X\times_{S} S^{''}$, qui sont de type fini dessus. On applique alors le crit\`ere de finitude \Ref{VIII}~\Ref{VIII.2.3} pour les images directes pour $i\leq n$ de $G^{''}$ sous l'immersion $U^{''}\to X^{''}$, en notant qu'elles sont satisfaites par hypoth\`eses pour $F^{''}$\pageoriginaled, donc aussi pour $G^{''}$ puisqu'elles s'expriment en termes de profondeur et que $G^{''}$ est localement isomorphe \`a un ${F^{''}}^n$. Le m\^eme argument montre que si $\Gg'$ est un Module coh\'erent sur $\widehat{U'}$ (compl\'et\'e de $U'$ par rapport \`a l'Id\'eal $\Jj\Oo_{U'}$), tel que $\Gg'_0=\Gg'/\Jj\Gg'$ soit localement de la forme $F_0^{'n}$, alors la conclusion de (iii) reste valable en y rempla\c cant $F'$ par $\Gg'$. On obtient ainsi le r\'esultat suivant, en utilisant les r\'esultats de l'expos\'e \Ref{IX}:
\end{remarque}

\begin{corollaire} \label{XII.5.10}
Soient $f:X\to S$ un morphisme propre, avec $S$ localement noeth\'erien, $U$ une partie ouverte de $X$; on suppose $U$ plat par rapport \`a $S$, et que pour tout changement de base $S'\to S$, avec $S'$ localement noeth\'erien, on ait $R^ig'_\ast(\Oo_{U'})$ coh\'erent sur $X'$ pour $i=0, 1$. Supposons alors que $S'$ soit de la forme $\Spec(A)$, o\`u $A$ est anneau noeth\'erien muni d'un id\'eal $J$ tel que $A$ soit s\'epar\'e et complet pour la topologie $J$-adique. Sous ces conditions:
\begin{enumeratei}
\item
Le foncteur $F\mto\widehat{F}$ de la cat\'egorie des Modules localement libres sur $U'$ dans la cat\'egorie des Modules localement libres sur $\widehat{U'}$ est pleinement fid\`ele.
\item
Pour tout Module localement libre $\formelF$ sur $\widehat{U'}$, il existe un Module coh\'erent $F$ sur~$U'$ (pas n\'ecessairement localement libre, h\'elas), et un isomorphisme $\widehat{F}\simeq \formelF$.
\end{enumeratei}
\end{corollaire}

Il reste seulement \`a prouver (ii), gr\^ace \`a \Ref{XII.5.9}. Or d'apr\`es cette remarque et \Ref{IX.2.1} il s'ensuit que $\formelF$ est induit par un Module coh\'erent $\formelG$ sur $\widehat{X'}$. D'apr\`es le th\'eor\`eme d'existence \EGA III 5.1.4, $\formelG$ est de la forme $\widehat{F}$, o\`u $F$ est coh\'erent sur $X$, d'o\`u la conclusion.

\begin{remarques} \label{XII.5.11}
1\textsuperscript{o}) Utilisant \Ref{XII.5.10}, \Ref{XII.4.7} et une hypoth\`ese convenable, disant que certains anneaux locaux des fibres g\'eom\'etriques de $X'\to S'$ sont \og purs\fg \resp parafactoriels, on doit pouvoir obtenir des \'enonc\'es disant que le foncteur $Z'\mto \widehat{Z'}$ de la cat\'egorie des rev\^etements \'etales de $X'$ dans la cat\'egorie des rev\^etements\pageoriginale \'etales de $\widehat{X'}$ (ou ce qui revient au m\^eme, de $X_0'$) est une \'equivalence de cat\'egories, \resp que le foncteur $L\mto \widehat{L}$ de la cat\'egorie des Modules inversibles sur $X'$ dans la cat\'egorie des Modules inversibles sur $\widehat{X'}$ est une \'equivalence. Utilisant des r\'esultats r\'ecents de Murre, il est probable que l'on doit pouvoir d\'eduire des th\'eor\`emes d'existence des sch\'emas de Picard pour certains sch\'emas alg\'ebriques \emph{non propres}\nde{bien entendu, on se reportera dans le cas projectif aux th\'eor\`emes d'existence de Grothendieck de FGA; \cf Grothendieck~A., {\og Technique de descente et th\'eor\`emes d'existence en g\'eom\'etrie alg\'ebrique. VI. Les sch\'emas de Picard: propri\'et\'es g\'en\'erales\fg}, in \emph{S\'eminaire Bourbaki}, vol.~7, Soci\'et\'e math\'ematique de France, Paris, 1995, \Exp 236, p\ptbl 221--243 et {\og Technique de descente et th\'eor\`emes d'existence en g\'eom\'etrie alg\'ebrique. V. Les sch\'emas de Picard: th\'eor\`emes d'existence\fg}, in \emph{S\'eminaire Bourbaki}, vol.~7, Soci\'et\'e math\'ematique de France, Paris, 1995, \Exp 232, 143--161. Les neuf conjectures de finitude qui s'y trouvent sont d\'emontr\'ees dans les expos\'es \Ref{XII} et \Ref{XIII} de Mme Raynaud et de Kleiman de \SGA 6. Pour un excellent texte \'el\'ementaire sur le sujet, voir l'article d'exposition de Kleiman (Kleiman~S., {\og The Picard scheme\fg}, \`a para\^itre dans Contemp. math.). Pour une application de ces techniques aux jacobiennes g\'en\'eralis\'ees globales d'une courbe lisse relative, voir (Contou-Carr\`ere~C., {\og La jacobienne g\'en\'eralis\'ee d'une courbe relative; construction et propri\'et\'e universelle de factorisation\fg}, \textit{C.~R. Acad. Sci. Paris S\'er.~A-B} \textbf{289}, (1979), \numero3, A203--A206 et {\og Jacobiennes g\'en\'eralis\'ees globales relatives\fg}, in \emph{The Grothendieck Festschrift, Vol. II}, Progr. Math., vol.~87, Birkh\"auser, Boston, 1990, p\ptbl 69--109). Voir aussi du même auteur, dans le cadre purement local, la construction et l'\'etude du foncteur \og jacobienne g\'en\'eralis\'ee locale\fg ({\og Jacobienne locale, groupe de bivecteurs de Witt universel, et symbole mod\'er\'e\fg}, \textit{C.~R. Acad. Sci. Paris S\'er.~I Math.}  \textbf{318} (1994), \numero 8, p\ptbl 743--746). Par ailleurs, si dans le cas d'un morphisme projectif et lisse, les composantes connexes du sch\'ema de Picard sont propres, il n'en va plus de m\^eme dans le cas singulier. Se pose naturellement le probl\`eme de la compactification des sch\'emas de Picard: ce probl\`eme a \'et\'e \'etudi\'e en d\'etail, notamment dans (Altman~A.B. \& Kleiman~S., {\og Compactifying the Picard scheme\fg}, \emph{Adv. in Math.} \textbf{35} (1980), \numero 1, p\ptbl 50--112. et {\og Compactifying the Picard scheme. II\fg}, \emph{Amer. J.~Math.} \textbf{101} (1979), \numero 1, p\ptbl 10--41). Le cas des courbes avait \'et\'e \'etudi\'e auparavant (D'Souza~C., {\og Compactification of generalised Jacobians\fg}, \emph{Proc. Indian Acad. Sci. Sect. A Math. Sci.} \textbf{88} (1979), \numero 5, p\ptbl 419--457). On sait m\^eme exactement quand la jacobienne compactifi\'ee d'une courbe est irr\'eductible (Rego~C.J., {\og The compactified Jacobian\fg} \emph{Ann. Sci. \'Ec. Norm. Sup. (4)} \textbf{13} (1980), \numero 2, p\ptbl 211--223, celle-ci \'etant l'adh\'erence de la jacobienne (ordinaire) lorsque la courbe est g\'eom\'etriquement int\`egre et localement planaire; pour une construction en famille des jacobiennes compactifi\'ees, voir (Esteves~E., {\og Compactifying the relative Jacobian over families of reduced curves\fg}, \emph{Trans. Amer. Math. Soc.} \textbf{353} (2001), \numero 8, p\ptbl 3045--3095). Depuis, les r\'esultats d'existence du sch\'ema de Picard dans le cas propre ont progress\'e depuis l'\'edition originale de \SGA 2; \cf (Murre~J.P., {\og On contravariant functors from the category of pre-schemes over a field into the category of abelian groups (with an application to the Picard functor)\fg}, \emph{Publ. Math. Inst. Hautes \'Etudes Sci.} \textbf{23} (1964), p\ptbl 5-43) et surtout (Artin~M., {\og Algebraization of formal moduli. I\fg}, in \emph{Global Analysis (Papers in Honor of K\ptbl Kodaira)}, Univ. Tokyo Press, Tokyo, 1969, p\ptbl 21--71). Voir aussi (Raynaud~M., {\og Sp\'ecialisation du foncteur de Picard\fg}, \emph{Publ. Math. Inst. Hautes \'Etudes Sci.} \textbf{38} (1970), p\ptbl 27--76) dans le cas d'un sch\'ema propre sur un anneau de valuation discr\`ete, mais pas n\'ecessairement plat. Pour une discussion de r\'esultats plus r\'ecents, en particulier ceux d'Artin, pour le foncteur de Picard des sch\'emas propres et plats, en particulier dans le cas cohomologiquement plat en dimension $0$, voir le chapitre VIII de (Bosch~S., Lütkebohmert~W. \& Raynaud~M., \emph{N\'eron models}, Ergebnisse der Mathematik und ihrer Grenzgebiete (3), vol.~21, Springer-Verlag, Berlin, 1990) et les r\'ef\'erences cit\'ees. Beaucoup plus r\'ecemment, des r\'esultats tr\`es fins ont \'et\'e obtenus dans le cas de courbes relatives $f:X\to S$ \sisi{au dessus}{au-dessus} du spectre $S$ d'un anneau de valuation discr\`ete \`a corps r\'esiduel parfait. Plus pr\'ecis\'ement, on suppose que $f$ est propre et plat, $X$~r\'egulier et $f_*\OX=\Oo_S$. En revanche, on ne suppose pas $f$ cohomologiquement plat en dimension~$0$, \ie on ne suppose pas $H^1(X,\Oo)$ sans torsion. Le sch\'ema de Picard n'est alors pas repr\'esentable, que ce soit par un sch\'ema ou un espace alg\'ebrique, d\`es que la torsion en question est non nulle. Soit~$J$ le mod\`ele de N\'eron de la fibre g\'en\'erique de $f$: c'est un quotient du foncteur de Picard $P$. Alors, Raynaud a montr\'e que le noyau de l'application tangente $H^1(X,\Oo)=\mathrm{Lie}(P)\to \mathrm{Lie}(J)$ co\"incide avec le sous-groupe de torsion du $H^1$ et que le conoyau a m\^eme longueur (voir le th\'eor\`eme 3.1 de (Liu~Q., Lorenzini~D. \& Raynaud~M., {\og N\'eron models, Lie algebras, and reduction of curves of genus one\fg}, \emph{Invent. Math.} \textbf{157} (2004), p\ptbl 455--518). Ce r\'esultat permet aux auteurs pr\'ecit\'es d'\'etudier le lien entre les conjectures de Birch-Swinnerton-Dyer et Artin-Tate (voir th\ptbl 6.6 de {\loccit}). En ce qui concerne le sch\'ema de Picard local, voir la th\`ese de Boutot, cit\'ee note de l'\'editeur~\eqref{Boutot}~page~\pageref{Boutot}.}. De fa\c con g\'en\'erale, l'\'elimination d'hypoth\`eses de puret\'e dans divers th\'eor\`emes d'existence, notamment de repr\'esentabilit\'e de foncteurs comme les foncteurs de Hilbert, ou de Picard etc., \`a~l'aide des techniques d\'evelopp\'ees dans ce s\'eminaire, m\'erite une \'etude syst\'ematique.

2\textsuperscript{o}) On peut se proposer de donner des conditions n\'ecessaires et suffisantes maniables, en termes de profondeur, pour que la condition de finitude universelle envisag\'ee dans \Ref{XII.5.10} soit v\'erifi\'ee. Lorsque $S$ est le spectre d'un corps, il r\'esulte facilement de \EGA III 1.4.15 qu'il faut et il suffit que les $\R^ig_\ast(F)$ ($i\leq n$) soient coh\'erents, ce qui s'exprime bien en termes de profondeur gr\^ace \`a \Ref{VIII}~\Ref{VIII.2.3}. Dans le cas g\'en\'eral, on notera cependant qu'il ne suffit pas d'exiger que la condition pr\'ec\'edente soit v\'erifi\'ee pour toutes les fibres $U_s\subset X_s$ ($s\in S$), m\^eme dans le cas o\`u $n=0$. Prendre par exemple $X=S$, $S$ le spectre d'un anneau de valuation discr\`ete, $U$ l'ouvert r\'eduit au point g\'en\'erique, $F=\Oo_U$.

3\textsuperscript{o}) Voici cependant une condition \emph{suffisante} assurant qu'on est sous les conditions de l'hypoth\`ese de \Ref{XII.5.10}: Il suffit que $f$ soit plat, et que pour tout $s\in S$ et tout $x\in Y_s=X_s- U_s$, on ait
\[
\prof \Oo_{X_s, x}\geq n+2.
\]
En effet, compte tenu du lemme \Ref{XII.2.5} (\cf relation (\Ref{eq:XII.16}) apr\`es \Ref{XII.2.5}), il s'ensuit que l'on a alors $g_\ast(\Oo_U)\simeq\OX$ et $\R^ig_\ast(\Oo_U)=0$ pour $0<i\leq n$, et les m\^emes relations seront \'evidemment v\'erifi\'ees apr\`es tout changement de base $S'\to S$.
\end{remarques}

\chapter{Probl\`emes et conjectures} \label{XIII}

\section[Relations entre r\'esultats globaux et locaux]
{Relations entre r\'esultats globaux et locaux. Probl\`emes affines li\'es \`a la dualit\'e} \label{XIII.1}

Il\pageoriginale est bien connu que beaucoup d'\'enonc\'es concernant un sch\'ema
projectif $X$ peuvent se formuler en termes d'\'enonc\'es concernant
un certain anneau gradu\'e, ou mieux un anneau local complet,
\sisi{\ignorespaces}{\`a} savoir l'anneau de coordonn\'ees homog\`enes
de $X$ (\ie l'anneau affine du c\^{o}ne projetant $\tilde{X}$ de
$X$), ou son compl\'et\'e (\ie le compl\'et\'e de l'anneau local du sommet
de $\tilde{X}$). L'int\'er\^et de cette reformulation est qu'elle
permet souvent \`a partir de r\'esultats globaux connus, de
conjecturer, voire de d\'emontrer, des r\'esultats analogues pour des
anneaux locaux noeth\'eriens complets plus g\'en\'eraux que ceux qui
interviennent r\'eellement dans l'\'enonc\'e global, par exemple pour
des anneaux locaux qui ne sont pas n\'ecessairement d'\'egale
caract\'eristique. Ainsi, le th\'eor\`eme de dualit\'e de Serre pour
l'espace projectif \Ref{XII}~\Ref{XII.1.1} a sugg\'er\'e l'utile
th\'eor\`eme de dualit\'e locale \sisi{V.5.6}{\Ref{V}~\Ref{V.2.1}}. Le
th\'eor\`eme fondamental de Serre sur la cohomologie des
\sisi{m}{M}odules coh\'erents sur l'espace projectif (finitude,
comportement a\sisi{s}{}symptotique pour $n$ grand de $\H^i(X,
F(n))$, {\cf \EGA III~2.2.1}) se g\'en\'eralise en un th\'eor\`eme de
structure pour les invariants locaux $\H^i_{\mm}(M)$, voir
\sisi{V.5.7}{\Ref{V}~\Ref{V.3}}. De m\^eme, les th\'eor\`emes de
Lefschetz pour le groupe fondamental, et le groupe de Picard (\og
crit\`eres d'\'equivalence{\fg}), bien familiers dans le cas
classique et \'etendus par la suite \`a un corps de base quelconque,
ont sugg\'er\'e les th\'eor\`emes de Lefschetz \og locaux\fg des expos\'es
\Ref{X} et \Ref{XI}. Bien entendu, les th\'eor\`emes locaux \`a leur
tour sont des outils pr\'ecieux pour obtenir des \'enonc\'es globaux.
Par exemple la dualit\'e locale permet de formuler une propri\'et\'e
asymptotique globale \Ref{XII}~\Ref{XII.1.3} (i) par l'annulation
de certains invariants locaux $\H^i(F_x)$. De fa\c con plus
substantielle, les th\'eor\`emes Lefschetz locaux, impliquant par
exemple la \og puret\'e\fg ou la parafactorialit\'e de certains
anneaux locaux intersections compl\`etes (\Ref{X}~\Ref{X.3.4} et
\Ref{XI}~\Ref{XI.3.13}) permettent dans les th\'eor\`emes de
Lefschetz globaux de se d\'ebarrasser de certaines hypoth\`eses de
non singularit\'e, comme dans \Ref{X}~\Ref{X.3.5}, \Ref{X.3.6},
\Ref{X.3.7}.

Une\pageoriginale autre g\'en\'eralisation utile des th\'eor\`emes concernant les sch\'emas projectifs sur un corps $k$ consiste \`a remplacer $k$ par un sch\'ema de base g\'en\'eral. Ainsi, la suite de {\EGA III} donnera\nde{en fait, cette g\'en\'eralisation ne s'y trouve pas; voir note plus bas, et la note de l'\'editeur~\eqref{noteconrad} de la page~\pageref{noteconrad}} une g\'en\'eralisation dans ce sens de la dualit\'e de Serre\sfootnote{\cf S\'eminaire Hartshorne, cit\'e \`a la fin de {\Exp \Ref{IV}}.}; les th\'eor\`emes de finitude et de comportement asymptotique des $\H^i(X, F(n))$ ont \'et\'e \'enonc\'es dans \EGA III~2.2.1 sur un sch\'ema de base g\'en\'eral, enfin les th\'eor\`emes de Lefschetz peuvent \'egalement se d\'evelopper pour un morphisme projectif, comme on~a vu dans \Ref{XII}~\Ref{XII.4.9}, gr\^ace au th\'eor\`eme local \Ref{XII}~\Ref{XII.4.7}. Bien entendu, le fait de travailler sur un sch\'ema de base g\'en\'eral conduit aussi \`a des \'enonc\'es essentiellement nouveaux, tels le \og th\'eor\`eme de comparaison\fg \EGA III~4.15 et le th\'eor\`eme d'existence de faisceaux \EGA III~5.1.4 (qui, on l'a vu par ailleurs dans \Ref{IX}, rel\`event des m\^emes th\'eor\`emes-clefs de nature cohomologique que les th\'eor\`emes de Lefschetz pour $\pi_1$ et $\Pic$).

Il s'impose alors de d\'egager des th\'eor\`emes qui englobent simultan\'ement les deux g\'en\'eralisations que nous venons de signaler des \'enonc\'es concernant les sch\'emas projectifs sur un corps. Les objets naturels pour une telle g\'en\'eralisation commune sont les \emph{anneaux noeth\'eriens s\'epar\'es et complets pour une topologie $I$-adique}. Leur \'etude, \`a ce point de vue, n'a pas encore \'et\'e abord\'ee s\'erieusement, et me semble \`a l'heure actuelle le sujet le plus int\'eressant dans la th\'eorie locale des faisceaux coh\'erents. Voici un probl\`eme typique dans cette direction:

\refstepcounter{subsection}\label{XIII.1.1}
\refstepcounter{footnoteAG}\sfootnotetext{Cette conjecture, et la conjecture 1.2, ci-dessous, sont fausses, comme l'a montr\'e \sisi{R.~HARTSHORNE (non publi\'e)}{R\ptbl Hartshorne, {\og Affine duality and cofinite modules\fg}, \emph{Invent. Math.} \textbf{9} (1969/70), p.~145-164, section 3.}}
\begin{enonce*}{Conjecture 1.1}
[\og Deuxi\`eme th\'eor\`eme de finitude affine\fg$^{(**)}$\sisi{}{\setcounter{footnoteNDE}{1}}\ndemark]
\ndetext{toutefois, si $A$ est local complet (\resp r\'egulier de caract\'eristique positive) et $J$ est l'id\'eal maximal, l'\'enonc\'e est vrai pour $M$ de type fini (\resp $M=A$), \cf (Hartshorne~R., {\og Affine duality and cofinite modules\fg}, \emph{Invent. Math.} \textbf{9} (1969/70), p\ptbl 145-164, corollaire 1.4) (\resp (Huneke~C. \& Sharp~R., {\og Bass numbers of local cohomology modules\fg}, \emph{Trans. Amer. Math. Soc.} \textbf{339} (1993), \numero 2, p\ptbl 765--779), qui contient d'ailleurs des r\'esultats bien plus forts). Pour des m\'ethodes compl\`etement diff\'erentes ($D$-modules) permettant d'aborder la caract\'eristique nulle, voir (Lyubeznik~G., {\og Finiteness properties of local cohomology modules (an application of $D$-modules to commutative algebra)\fg}, \emph{Invent. Math.} \textbf{113} (1993), \numero 1, p\ptbl 41--55); voir aussi du m\^eme auteur {\og Finiteness properties of local cohomology modules for regular local rings of mixed characteristic: the unramified case\fg}, \emph{Comm. Algebra} \textbf{28} (2000), \numero 12, p\ptbl 5867--5882, Special issue in honor of Robin Hartshorne, et {\og Finiteness properties of local cohomology modules: a characteristic-free approach\fg}, \emph{J.~Pure Appl. Algebra} \textbf{151} (2000), \numero 1, p\ptbl 43--50. La notion de module cofini a \'evolu\'e depuis sous la houlette de Hartshorne. On dit que $M$ est $J$ cofini s'il est de support contenu dans $V(J)$ et si tous les $\Ext^i_A(A/J, \H^i_J(M))$ sont de type finis. Sur ce sujet, voir par exemple (Delfino~D. \& Marley~Th., {\og Cofinite modules and local cohomology\fg}, \emph{J. Pure Appl. Algebra} \textbf{121} (1997), \numero 1, p\ptbl 45--52).}
Soient $M$ un module de type fini sur un anneau noeth\'erien $A$ (qu'on supposera au besoin quotient d'un r\'egulier), $J$ un id\'eal de $A$, prouver que les modules $\H^i_J(M)$ sont \og $J$-cofinis\fg, \ie que les modules
$$
\Hom_A(A/J, \H^i_J(M))$$
sont de type fini.
\end{enonce*}

Rappelons qu'on d\'esigne par $\H^i_J (M)$ le module $\H^i_Y (X, \tilde{M})$ (o\`u $X=\Spec (A)$, $Y=V(J)$) de \Exp \Ref{I}, interpr\'et\'e dans \Ref{II} en termes de limite inductive de cohomologies\pageoriginale de complexes de Koszul, ou encore pour $i\ge 2$ le module $\H^{i-1} (X-Y, \tilde{M})$. \sisi{A}À vrai dire, \Ref{XIII.1.1} devrait \^etre une cons\'equence d'un \'enonc\'e plus pr\'ecis, impliquant que les $\H^i_J ({M})$ sont dans une sous-cat\'egorie ab\'elienne convenable ${\DDJ}$ de la cat\'egorie ${\CJ}$ des $A$-modules de support $\subset Y=V(J)$, telle que $H\in \Ob {\DDJ}$ implique que $H$ est $J$-cofini. (\sisi{NB}{N.B.} La cat\'egorie des modules $H$ de support contenu dans $V(J)$ et qui sont $J$-cofinis n'est malheureusement pas stable par passage au quotient!). Le probl\`eme essentiel consisterait alors \`a d\'efinir ${\DDJ}$. De fa\c con plus pr\'ecise, la solution du probl\`eme \Ref{XIII.1.1} devrait sortir (du moins si $A$ est quotient d'un r\'egulier) d'une th\'eorie de dualit\'e, g\'en\'eralisant \`a la fois la dualit\'e locale, et la th\'eorie de dualit\'e des morphismes projectifs \`a laquelle il a \'et\'e fait allusion plus haut, et qui serait du genre suivant:

\refstepcounter{subsection} \label{XIII.1.2}
\begin{enonce*}{Conjecture 1.2}[\og Dualit\'e affine\fg\sfootnotemark]
\sfootnotetext{Cette conjecture, fausse telle quelle, a cependant \'et\'e \'etablie sous une forme assez voisine par \sisi{R.~HARTSHORNE, \emph{Affine duality and cofinite modules} (\`a para\^itre).}{R.~Hartshorne, {\og Affine duality and cofinite modules\fg}, \emph{Invent. Math.} \textbf{9} (1969/70), p\ptbl 145-164.}} Supposons $A$ r\'egulier, s\'epar\'e et complet pour la topologie $J$-adique. Soit $C^{\boule}(A)$ une r\'esolution injective de $A$.
\begin{enumeratei}
\item
Prouver que le foncteur
$$D_J: L_{\boule}\mto \Hom_J(L_{\boule}, C^{\boule}(A))$$
de la cat\'egorie des complexes de $A$-modules, libres de type fini en toute dimension et \`a degr\'es limit\'es sup\'erieurement, (o\sisi{u}{\`u} les morphismes sont les homomorphismes de complexes \`a homotopie pr\`es) dans la cat\'egorie des complexes de $A$-modules $K^{\boule}$, injectifs en toute dimension et \`a degr\'es limit\'es sup\'erieurement (o\sisi{u}{\`u} les homomorphismes sont d\'efinis de la m\^eme fa\c con), est pleinement fid\`ele.

\item
Prouver que pour tout $K^{\boule}$ de la forme $D_J(L_{\boule})$, les $\H^i(K^{\boule})(= \Ext_Y^i (X; L_{\boule}, \OX))$ sont $J$-cofinis.

\item
De fa\c con plus pr\'ecise, prouver que les $K^{\boule}$ qui sont homotopes \`a un complexe de la forme $D_J(L_{\boule})$ peuvent se caract\'eriser par des propri\'et\'es de finitude des $\H^i(K^{\boule})$, plus fortes que celle envisag\'ee dans (ii), par exemple par\pageoriginale la propri\'et\'e $\H^i(K^{\boule})\in \Ob {\DDJ}$, o\`u ${\DDJ}$ est une cat\'egorie ab\'elienne convenable, comme envisag\'ee plus haut.
\end{enumeratei}
\end{enonce*}

Notons que le probl\`eme est r\'esolu par l'affirmative lorsque $A$ est local et que $J$ en est un id\'eal de d\'efinition (\cf \Exp \Ref{IV}), et aussi lorsque $J$ est l'id\'eal nul. Dans ces deux cas, exceptionnellement, on peut se borner \`a prendre pour ${\DDJ}$ la cat\'egorie des \sisi{m}{M}odules \`a support $V(J)$ qui sont $J$-cofinis, (ce qui dans le deuxi\`eme cas signifie simplement, qu'on prend la cat\'egorie des \sisi{m}{M}odules de type fini sur $A$). Une solution affirmative de la conjecture \Ref{XIII.1.2} en g\'en\'eral en donnerait une pour \Ref{XIII.1.1}, en prenant pour $L_{\boule}$ le dual d'une r\'esolution libre de type fini de $M$. D'autre part, une solution affirmative de \Ref{XIII.1.1} donnerait une r\'eponse affirmative \`a la premi\`ere partie de la conjecture suivante que nous formulons sous forme \og globale\fg:

\begin{conjecture} \label{XIII.1.3}
Soit $X\subset \PP^r_k$ un sous-sch\'ema ferm\'e du sch\'ema projectif type qui soit localement une intersection compl\`ete et dont toute composante irr\'eductible soit de codimension $\ge s$. Soit $U=\PP^r_k -X$.
\begin{enumeratei}
\item
Prouver que pour tout Module coh\'erent $F$ sur $U$, on a
$$
\dim \H^i(U, F)<+\infty{\rm \ pour\ } i\ge s.{}\sfootnote{La part\sisi{}{ie} (i) de cette conjecture est prouv\'ee par R.~\sisi{HARTSHORNE}{Hartshorne} lorsque $\sisi{Y}{X}$ est lisse \sisi{dans \emph{Ample Vector Bundles}, Pub. Math. \numero 29, (1966)}{\cf \emph{Ample subvarieties of algebraic varieties}, Notes written in collaboration with C\ptbl Musili, Lect. Notes in Math, vol.~156, Springer-Verlag, Berlin-New York, 1970, th\'eor\`eme III.5.2}. Cet auteur a \'egalement trouv\'e un exemple pour (ii), \cf \sisi{R.~HARTSHORNE, \emph{Cohomological dimension of algebraic varieties}, \`a para\^itre aux Ann. of Math.}{R.~Hartshorne, {\og Cohomological dimension of algebraic varieties\fg}, \emph{Ann. of Math. (2)} \textbf{88} (1968), p\ptbl 403--450, exemple page 449.}}
$$

\item
Donner un exemple, avec $X$ connexe et r\'egulier, o\`u on a
$$
\H^s(U, F)\neq 0.
$$
\end{enumeratei}
\end{conjecture}

Pour voir que (i) est un cas particulier de \Ref{XIII.1.1} on consid\`ere
$$
\H^i(U, F(\cdot))=\sisi{\coprod}{\tbigoplus}_n \H^i (U, F(n))= \H^i(E^{r+1} -\sisi{X}{\tilde X}, \tilde{F})$$
comme un module sur l'anneau affine $k[t_0, \ldots, t_r]$ du c\^{o}ne projetant $E^{r+1}$ de $\PP^r$. Ce module n'est autre que $\H_J^{i+1}(M)$, o\`u $J$ est l'id\'eal du c\^{o}ne projetant $\tilde X$ de $X$ dans $E^{r+1}$. D'autre part de l'hypoth\`ese faite sur $X$, qui implique que $\sisi{X}{\tilde X}$ est aussi une intersection compl\`ete de co\sisi{-}{}dimension $\ge s$ en tout\pageoriginale point de $E^{r+1}$ distinct de l'origine, r\'esulte que $\H^{i+1}_J(M)$ est nul en dehors de l'origine pour $i\ge s$. S'il est donc $J$-cofini comme le veut \Ref{XIII.1.1}, il est \afortiori $\mm$-cofini, ce qui implique facilement qu'il est de dimension finie en tout degr\'e\refstepcounter{toto}\nde{\label{noteHM}Hartshorne a prouv\'e (Hartshorne~R., {\og Cohomological dimension of algebraic varieties\fg}, \emph{Ann. of Math. (2)} \textbf{88} (1968), p\ptbl403--450) que la cohomologie $H^{n-1}(\PP^n_k-X,F)$ est nulle pour $F$ coh\'erent et $X$ de dimension positive ($k$ alg\'ebriquement clos). En fait, gr\^ace essentiellement \`a la dualit\'e de Serre et au th\'eor\`eme de Lichtenbaum --- nullit\'e de la cohomologie des faisceaux coh\'erents en dimension maximale des vari\'et\'es quasi-projectives irr\'eductibles non compl\`etes ---, on se ram\`ene \`a prouver que le compl\'et\'e formel $\hat \PP^n_k$ et $X$ ont m\^eme corps de fonctions rationnelles. C'est le point difficile (th\'eor\`eme 7.2 de {\loccit}), autrement dit $\PP^n_k$ est $G3$ dans la terminologie de (Hironaka~H. \& Matsumura~H., {\og Formal functions and formal embeddings\fg}, \emph{J.~Math. Soc. Japan} \textbf{20} (1968), p\ptbl 52--82). Ces auteurs ont prouv\'e ind\'ependamment les r\'esultats pr\'ec\'edents, et en fait beaucoup mieux. Ils ont prouv\'e que $X$ est universellement $G3$ et ont calcul\'e le corps des fonctions rationnelles du compl\'et\'e formel d'une vari\'et\'e ab\'elienne le long d'une sous-vari\'et\'e de dimension positive. C'est dans cet article qu'appara\^it pour la premi\`ere fois les conditions $G1,G2,G3$ d\'esormais classiques.} Noter d'ailleurs que la conjecture \Ref{XIII.1.3} se pose d\'ej\`a pour une courbe irr\'eductible non singuli\`ere $X$ dans $\PP^3$, on ignore si dans ce cas les $\H^2(\PP^3-X, \OX(n))$ sont de dimension finie, ou s'ils sont n\'ecessairement nuls\sfootnote{La question vient d'\^etre r\'esolue affirmativement par R\ptbl \sisi{HARTSHORNE}{Hartshorne}\sisi{}{ (Hartshorne~R., {\og Ample vector bundles\fg}, \emph{Publ. Math. Inst. Hautes \'Etudes Sci.} \textbf{29} (1966), p\ptbl63--94, th\'eor\`eme 8.1)} et H. \sisi{HIRONAKA}{Hironaka}.}. On ignore m\^eme s'il existe une courbe irr\'eductible de $\PP^3$ qui ne soit ensemblistement l'intersection de deux hypersurfaces\label{courbesurface}\nde{voir la conjecture~\Ref{XII.3.5} et la note correspondante.}.

\begin{probleme} \label{XIII.1.4}
Donner une variante affine du \og th\'eor\`eme de comparaison\fg \textup{\EGA III~4.1.5} comme un th\'eor\`eme de commutation des foncteurs $\H^i_J$ avec certaines limites projectives.
\end{probleme}

Enfin, dans le pr\'esent ordre d'id\'ees, j'avais pos\'e le probl\`eme suivant: \sisi{S}{s}oit $A$ un anneau local noeth\'erien r\'egulier \emph{complet}, $K$ son corps des fractions, prouver que $\Ext^i_A(K, A)=0$ pour tout $i$. Une r\'eponse affirmative a \'et\'e donn\'ee sur le champ par \sisi{Mr AUSLANDER}{M.~Auslander}, l'hypoth\`ese de r\'egularit\'e peut \^etre remplac\'ee par celle que $A$ est int\`egre et il est vrai en fait que $\Ext^i_A(K, M)=0$ pour tout $i$, d\`es que $M$ est de type fini sur $A$. Cela r\'esulte aussit\^{o}t de l'\'enonc\'e suivant, d\^u \`a \sisi{AUSLANDER}{Auslander}: si $A$ est un anneau local noeth\'erien complet, alors pour tout module de type fini $M$ sur $A$, les foncteurs $\Ext_A^i(., M)$ transforment limites inductives en limites projectives\nde{\'ecrire $K=\varinjlim_{a\neq 0}A[1/a]$ et observer que $a\neq 0$ est $A$-r\'egulier.}.

\section{Probl\`emes li\'es au $\pi_0$: th\'eor\`emes de Bertini locaux} \label{XIII.2}

Soit $A$ un anneau local noeth\'erien complet, $f$ un \'el\'ement de son
id\'eal maximal, $X=\Spec (A)$, $Y=\Spec(A/fA)$. L'utilisation de la
technique \og Lefschetz\fg locale permet de donner des crit\`eres
pour que $Y'=X'\cap Y$ (o\`u $X'=X-\{\mm\}$) soit connexe, en termes
d'hypoth\`eses sur $X'$. Ainsi, il suffit que l'on ait: a) $X'$\pageoriginale
connexe\sisi{}{,} b) $\prof \Oo_{X', x} \ge 2$ pour tout point ferm\'e $x$
de $X'$\sisi{}{,} c) $f$ est $A$-r\'egulier. On notera cependant que les
hypoth\`eses b) et c) ne sont pas de nature purement
topologiques, par exemple ne sont pas invariant\sisi{}{e}s en
rempla\c cant $X$ par $X_{\text{red}}$. Dans la situation analogue
pour un sch\'ema projectif $X'$sur un corps et une section
hyperplane $Y'$ de $X'$, l'utilisation du \og th\'eor\`eme de
Bertini\fg et du \og th\'eor\`eme de connexion\fg de Zariski permet
d'obtenir en fait des r\'esultats d'allure nettement plus
satisfaisants, qui m'avaient amen\'e dans le s\'eminaire oral \`a
\'enoncer une conjecture, que j'ai r\'esolue depuis par l'affirmative.
\'Enon\c cons donc ici:

\begin{theoreme} \label{XIII.2.1}
Soient $A$ un anneau local noeth\'erien complet, $X$ son spectre, $\aaa$ le point ferm\'e de $X$, $X'=X-\{ \aaa\}$. Supposons que $X$ satisfasse les conditions (o\`u $k$ d\'esigne un entier $\ge 1$):
\begin{itemize}
\item[a$_k$)] Les composantes irr\'eductibles de $X'$ sont de dimension $\ge k+1$.

\item[b$_k$)] $X'$ est connexe en dimension $\ge k$, \ie on ne peut d\'econnecter $X'$ par une partie ferm\'ee de dimension $<k$ (\cf \Ref{III}~\Ref{III.3.8}).
\end{itemize}

Soit $m$ un entier, $0\le m\le k$, et soient $f_1, \dotsc, f_m\in \rr(A)$, posons $B=A/\sum_i f_iA$, $Y=\Spec(B)=V(f_1)\cap\dotsb\cap V(f_m)$, $Y'=X'\cap Y=Y-\{\aaa\}$. Alors $Y$ satisfait les conditions a$_{k-m}$), b$_{k-m}$). En particulier, pour toute suite de $m\le k$ \'el\'ements $f_1, \dotsc, f_m$ de $\rr(A)$, $Y'=X'\cap V(f_1)\cap\dotsb\cap V(f_m)$ est connexe.
\end{theoreme}

Il est d'ailleurs facile de voir que si la derni\`ere conclusion est valable (il suffit \'evidemment d'y faire $m=k$), et en excluant le cas o\`u $X$ serait irr\'eductible de dimension~$0$ ou $1$, il s'ensuit que les composantes irr\'eductibles de $X'$ sont de dimension $\ge k+1$, et $X'$ est connexe en dimension $\geq k$, de sorte qu'en un sens, \Ref{XIII.2.1} est un r\'esultat \og le meilleur possible\fg.

\bigskip

Donnons le principe de la d\'emonstration de~\Ref{XIII.2.1}. Seule la condition b$_{k-m}$) pour $Y$ offre un probl\`eme. On est ramen\'e facilement pour $k$ donn\'e au cas o\`u $X$ est int\`egre, et m\^eme (en passant au normalis\'e, qui est fini sur $X$) au\pageoriginale cas o\`u $X$ est \emph{normal}. Si $k=1$, donc $\dim X' \ge 2$, alors $X'$ est de profondeur $\ge 2$ en ses points ferm\'es et on peut appliquer le r\'esultat rappel\'e au d\'ebut du \sisi{N${}^\circ$}{n$^\circ$}, qui montre que \hbox{$Y'=X'\cap V(f)$} est connexe. Dans le cas $k\ge 1$, on suppose le th\'eor\`eme d\'emontr\'e pour les \hbox{$k'\!<\!k$}. Par r\'ecurrence sur $m$, on est ramen\'e au cas o\`u $m=1$, \ie \`a v\'erifier que pour $f_1\in\nobreak \rr(A)$, $X'\cap V(f_1)$ est connexe en dimension $\ge k-1$. S'il ne l'\'etait pas, \ie s'il \'etait disconnect\'e par un $Z'$ de dimension $< k-1$, il existerait une suite $f_2, \dotsc, f_k$ telle que \hbox{$X'\cap V(f_1)\cap\dotsb\cap V(f_k)$} soit disconnexe, et dans cette suite on peut choisir arbitrairement $f_2\in \rr(A)$ soumis \`a la seule condition de ne s'annuler sur aucun point d'une certaine partie finie $F$ de $X'$ (\sisi{\ignorespaces}{\`a} savoir l'ensemble des points maximaux de $Z'$). D'ailleurs, on v\'erifie facilement, utilisant le fait que $X'$ est normal, donc satisfait la condition ({$S_2$}) de Serre\sfootnote{\Cf \EGA IV 5.7.2.} qu'il existe une partie finie $F'$ de $X'$ telle que $f\in \rr(A)$, $V(f)\cap F'=\emptyset$ implique que $V(f)\cap X'$ satisfait \'egalement \`a la condition ({$S_2$}). On peut alors choisir $f_2$ de telle sorte que $f_2$ ne s'annule ni sur $F$ ni sur $F'$, donc que $X'\cap V(f_2)$ satisfasse ({$S_2$}). Mais alors\sisi{}{,} en vertu du th\'eor\`eme de \sisi{HARTSHORNE}{Hartshorne} \Ref{III}~\Ref{III.3.6}\sisi{}{,} $X'\cap V(f_2)$ est connexe en codimension~$1$ donc (comme toute composante de $X'\cap V(f_2)$ est de dimension $\ge k$) il est connexe en dimension $\ge k-1$. Appliquant l'hypoth\`ese de r\'ecurrence \`a $V(f_2)=\Spec(A/f_2A)$, il s'ensuit que $X'\cap V(f_2)\cap V(f_1)\cap V(f_3)\cap \dotsb\cap V(f_k)$ est connexe, alors qu'on l'avait construit disconnexe, absurde.

Signalons quelques corollaires int\'eressants:

\begin{corollaire} \label{XIII.2.2}
Soit $f\colon X\to Y$ un morphisme propre de pr\'esch\'emas localement noeth\'erien, avec $Y$ int\`egre, $y_0\in Y$, $y_1$ le point g\'en\'erique de $Y$. On suppose
\begin{enumeratea}
\item
$Y$ est unibranche en $y_0$, toute composante irr\'eductible de $X$ domine $Y$.
\item
Les composantes irr\'eductibles de $X_1$ sont de dimension $\ge k+1$, et $X_1$ est connexe en dimension $\ge k$.
\end{enumeratea}

Alors\pageoriginale les composantes irr\'eductibles de $X_0$ sont de dimension $\ge k+1$, et $X_0$ est connexe en dimension $\ge k$.
\end{corollaire}

En effet, le th\'eor\`eme de connexion de Zariski (\cf\ \EGA III 4.3.1) implique que $X_0$ est connexe; pour montrer qu'il n'est pas disconnect\'e par une partie ferm\'ee de dimension $<k$, on est ramen\'e \`a montrer que ces anneaux locaux en des points $x\in X_0$ tels que $\dim \overline{x}<k$ ont un spectre non disconnect\'e par $x$. Or ceci est vrai sans supposer ni $f$ propre, ni $Y$ unibranche en $y_0$. On est ramen\'e pour le voir au cas o\`u $X$ est int\`egre dominant $Y$, et si l'on veut, $Y$ affine de type fini sur $\ZZ$, de sorte que l'on est sous les conditions de la formule des dimensions pour $\Oo_{X, x}$ sur $\Oo_{Y, y_0}$. Utilisant dans ce cas la finitude de la cl\^{o}ture normale, on peut m\^eme supposer $X$ \emph{normal}, donc en vertu d'un th\'eor\`eme de \sisi{NAGATA}{Nagata}\sfootnote{\Cf \EGA IV 7.8.3 (i) (ii) (v).}, le compl\'et\'e d'un anneau local $\Oo_{X, x}$ de $X'$ est encore normal, donc (si $\Oo_{X, x}$ est de dimension $N$) $\Spec (\hat{\Oo}_{X, x})$ est connexe en dimension $\ge N-1$. Soit $n=\dim\Oo_{\sisi{X}{Y}, y_0}$, alors $\degtr~k(x)/k(y)<k$ implique $\dim \Oo_{X, x}>n+(k+1)-k=n+1$, compte-tenu \sisi{que }{de} $\dim X_1\ge k+1$, et prenant un syst\`eme $f_1, \dotsc, f_n$ de param\`etres de $\Oo_{Y, y_0}$ qu'on rel\`eve en des \'el\'ements de $\Oo_{X, x}$, on voit par \Ref{XIII.2.1} que $\Spec (\hat{\Oo}_{X, x}/\sum f_i\hat{\Oo}_{X, x})$ est connexe en dimension $\ge 1$, \ie n'est pas disconnect\'e par son point ferm\'e, ou ce qui revient au m\^eme, $\Spec(\hat{\Oo}_{X_0, x})$ n'est pas disconnect\'e par son point ferm\'e; \afortiori il en est ainsi de $\Spec (\Oo_{X_0, x})$.

Comme dans le cas du th\'eor\`eme de connexion ordinaire, on peut varier \Ref{XIII.2.2} en prenant les fibres g\'eom\'etriques (sur les cl\^{o}tures alg\'ebriques des corps r\'esiduels), \`a condition de supposer $Y$ \emph{g\'eom\'etriquement} unibranche en $y_0$, ou (sans autre hypoth\`ese que $Y$ noeth\'erien) que $f$ est universellement ouvert. Appliquant ceci au cas o\`u $Y$ est le sch\'ema dual d'un sch\'ema projectif $\PP_k^r$ sur un corps, on retrouve une forme renforc\'ee du r\'esultat global qui avait inspir\'e \Ref{XIII.2.1}, savoir:

\refstepcounter{subsection}\label{XIII.2.3}
\begin{enonce*}{Corollaire 2.3\ndemark}\ndetext{pour une tr\`es belle preuve directe, voir
(Fulton~W. \& Lazarsfeld.~R., {\og Connectivity and its applications in algebraic geometry\fg},
in \emph{Algebraic geometry (Chicago, Ill., 1980)}, Lect. Notes in Math., vol.~862,
Springer, Berlin-New York, 1981, p\ptbl 26--92, th\'eor\`eme 2.1). \Cf aussi \cite{HL}, cit\'e dans la note de l'\'editeur~\eqref{HL}~page~\pageref{HL}.}
Soit $X$ un sous-sch\'ema ferm\'e de $\PP^r_k$ ($k$ un corps), on suppose les composantes irr\'eductibles de $X$ de dimension $\ge \sisi{k}{l}+1$, et $X$ g\'eom\'etriquement connexe\pageoriginaled en dimension $\ge \sisi{k}{l}$. Alors pour toute suite $H_1, \dotsc, H_m$ de $m$ hyperplans de $\PP^r_k$ ($0\le m\le \sisi{k}{l}-1$), $X\cap H_1\cap \dotsb \cap H_m$ satisfait la m\^eme condition avec $\sisi{k}{l}-m$, en particulier est g\'eom\'etriquement connexe en dimension $\ge \sisi{k}{l}-1$.
\end{enonce*}

On peut d'ailleurs modifier de fa\c con
\'evidente cet \'enonc\'e, pour le cas o\`u on s'est donn\'e un morphisme
propre $X \to \PP^r_k$, qui n'est pas n\'ecessairement une
immersion, et une extension analogue est possible pour
\Ref{XIII.2.1} (en consid\'erant un sch\'ema propre sur~$X'$). Ces
\'enonc\'es se d\'eduisent d'ailleurs formellement des \'enonc\'es donn\'es
ici, compte tenu du th\'eor\`eme de connexion ordinaire qui nous
ram\`ene au cas d'un morphisme fini.

\begin{corollaire} \label{XIII.2.4}
Soit $A$ un anneau local noeth\'erien normal complet de dimension $\ge k+2$. Soient $X=\Spec (A)$, $X'=X-\{ \aaa\}$, $f_1, \dotsc, f_k$ des \'el\'ements de $\rr(A)$, alors $Y'=X'\cap V(f_1) \cap \dotsb \cap V(f_k)$ est connexe, et $\pi_1(Y')\to \pi_1(X')$ est surjectif.
\end{corollaire}

On proc\`ede comme dans \SGA 1~X~2.11.

Dans tout ceci, il n'\'etait question que de questions de \emph{connexit\'e}. Or dans le cas global, des th\'eor\`emes bien connus affirment que pour une vari\'et\'e projective irr\'eductible $X\subset \PP^r_k$, $k$ alg\'ebriquement clos, son intersection avec un hyperplan $H$ assez \og g\'en\'eral\fg est irr\'eductible, (et non seulement connexe): c'est le \emph{th\'eor\`eme de \sisi{BERTINI}{Bertini}}, prouv\'e par \sisi{ZARISKI}{Zariski}, qui implique \`a son tour par le th\'eor\`eme de connexion de Zariski, que pour \emph{tout} $H$, $H\cap X$ est connexe (bien que non n\'ecessairement irr\'eductible). On peut d'ailleurs proc\'eder en sens inverse, en prouvant ce dernier r\'esultat par une technique \`a la Lefschetz, et en d\'eduisant le th\'eor\`eme de Bertini, en se ramenant au cas o\`u $X$ est normale, et utilisant le r\'esultat suivant: pour $H$ \og assez g\'en\'eral\fg, $X\cap H$ est \'egalement normal. Cela sugg\`ere la

\refstepcounter{subsection}\label{XIII.2.5}
\begin{enonce*}{Conjecture 2.5\ndemark}\ndetext{voir la note de l'\'editeur suivante.}
Soit $A$ un anneau local noeth\'erien complet normal. Montrer qu'il existe $f\in\rr(A)$ non nul tel que $Y'=X'\cap V(f)=Y-\{ \aaa\}$ (o\`u\pageoriginale $Y=\Spec (A/fA)$) soit normal (donc irr\'eductible par \Ref{XIII.2.1} si $\dim A\ge 3$).
\end{enonce*}

Pour bien faire, il faudrait montrer que dans un sens convenable, il existe m\^eme \og beaucoup\fg d'\'el\'ements $f$ ayant la propri\'et\'e en question, par exemple qu'on peut choisir $f$ dans une puissance arbitraire de l'id\'eal maximal. Utilisant le crit\`ere de normalit\'e de Serre, et la remarque faite plus haut pour la propri\'et\'e $(S_2)$ de Serre, on voit qu'on aurait une r\'eponse affirmative \`a \Ref{XIII.2.5} si on en avait une \`a la

\refstepcounter{subsection}\label{XIII.2.6}
\begin{enonce*}{Conjecture 2.6\ndemark}
\ndetext{on trouve maintenant une preuve de cette conjecture dans la litt\'erature, et donc la pr\'ec\'edente doit \'egalement \^etre consid\'er\'ee comme prouv\'ee comme indiqu\'e plus haut. On peut trouver \'egalement deux tentatives de preuves, publi\'ees ant\'erieurement mais h\'elas infructueuses, de Flenner et Trivedi. Voir Trivedi V., {\og Erratum: ``A local Bertini theorem in mixed characteristic"\fg}, \emph{Comm. Algebra} \textbf{25} (1997), \numero 5, p\ptbl1685-1686. Toutefois, l'\'editeur n'a pas v\'erifi\'e que la preuve est d\'esormais compl\`ete.}
Soient $A$ un anneau local noeth\'erien complet, $U$ une partie ouverte de son spectre $X$, $F$ une partie finie de $X'=X-\{ \aaa\}$. On suppose que $U$ est r\'egulier. Prouver qu'il existe un $f \in\rr(A)$ tel que $V(f)\cap U$ soit r\'egulier, et $V(f)\cap F=\emptyset$.
\end{enonce*}

Voir, pour un r\'esultat de type \og Bertini local\fg, \sisi{CHOW}{Chow}~\cite{XIII.2}

\section{Probl\`emes li\'es au $\pi_1$} \label{XIII.3} Ici encore, on~a de nombreuses questions, sugg\'er\'ees par les r\'esultats globaux ou les r\'esultats transcendants.

\refstepcounter{subsection}\label{XIII.3.1}
\begin{enonce*}{Conjecture 3.1\ndemark}\ndetext{l'\'enonc\'e analogue est vrai pour les sch\'emas (connexes) $X$ de type fini sur un corps s\'eparablement clos $k$ sous l'hypoth\`ese de d\'esingularisation forte pour tous les $\bar k$-sch\'emas (de type fini), en particulier si $k$ est de caract\'eristique nulle ou $X$ de dimension $\leq 2$. On se ram\`ene pour ce faire au cas des surfaces quasi-projectives par des techniques type Lefschetz d\'evelopp\'ees par Mme Raynaud, \cf notes \textit{supra}; voir \SGA 7.1, th\'eor\`eme II.2.3.1.} Soient $A$ un anneau local noeth\'erien complet \`a corps r\'esiduel alg\'ebriquement clos, $X=\Spec (A)$, $ X'= X- \{\aaa\}$, $\aaa$ le point ferm\'e. Supposons les composantes irr\'eductibles de $X$ de dimension $\ge 2$, $X'$ connexe.
\begin{enumeratei}
\item
Prouver que $\pi_1(X')$ est topologiquement \`a engendrement fini.
\item
Si $p$ est l'exposant caract\'eristique du corps r\'esiduel $k$ de $A$, prouver que le plus grand groupe topologique quotient de $\pi_1(X')$ qui est \og d'ordre premier \`a $p$\fg est de pr\'esentation finie.
\end{enumeratei}
\end{enonce*}

Pour la partie (i), utilisant la th\'eorie de la descente \SGA 1~IX~5.2 et le th\'eor\`eme\pageoriginale 2.4, on est ramen\'e au cas o\`u $A$ est normal de dimension $2$. Dans ce cas, une m\'ethode syst\'ematique pour \'etudier le groupe fondamental de $X'$, inaugur\'ee par \sisi{MUMFORD}{Mumford}~\cite{XIII.5} dans le cadre transcendant, consiste \`a \emph{d\'esingulariser} $X$ \ie \`a consid\'erer un morphisme projectif birationnel $Z\to X$, avec $Z$ int\`egre r\'egulier, induisant un isomorphisme $Z'=Z|X'\to X'$; il est plausible qu'un tel $Z$ existe toujours, c'est en \sisi{tous}{tout} cas ce que d\'emontre la m\'ethode d'\sisi{ABHYANKAR}{Abhyankar}~\cite{XIII.1} dans le cas d'\og \'egales caract\'eristiques{\fg}\sfootnote{La possibilit\'e de \og r\'esoudre\fg $A$ est \sisi{prouv\'e}{prouv\'ee} maintenant en toute g\'en\'eralit\'e par Abhyankar~\cite{XIII.8}.}. Soit~$C$ la fibre du point ferm\'e de $X$ par $Z\to X$, c'est une courbe alg\'ebrique sur le corps r\'esiduel $k$, connexe en vertu du th\'eor\`eme de connexion. La solution de \Ref{XIII.3.1} semble alors li\'ee au

\begin{probleme} \label{XIII.3.2}
Avec les notations pr\'ec\'edentes, mettre en relations $\pi_1(X')$ avec les invariants topologiques de $C$, en particulier son groupe fondamental, (pour faire appara\^itre la g\'en\'eration topologique finie de $\pi_1(X')$, en utilisant par exemple \SGA \sisi{}{\textup{I}, th\'eor\`eme~}\textup{X~2.6}).
\end{probleme}

Une autre m\'ethode serait de consid\'erer $A$ comme une alg\`ebre finie sur un anneau local \emph{r\'egulier} complet $B$ de dimension $2$, ramifi\'e suivant une courbe $C$ contenue dans $\Spec(B)=Y$. On est conduit ainsi au

\begin{probleme} \label{XIII.3.3}
Soient $A$ un anneau local r\'egulier complet de dimension $2$, \`a corps r\'esiduel alg\sisi{.}{\'ebriquement} clos $k$, $X$ son spectre, $C$ une partie ferm\'ee de $X$ de dimension~$1$. D\'efinir des invariants locaux de la courbe immerg\'ee $C$, ayant un sens ind\'ependant de la caract\'eristique r\'esiduelle, et dont la connaissance permette de calculer le groupe fondamental de $X-C$ par g\'en\'erateurs et relations lorsque $k$ est de caract\'eristique nulle. Prouver que lorsque $k$ est de caract\'eristique $p> 0$, le groupe fondamental \og tame\fg de $X-C$ est quotient du pr\'ec\'edent, et les deux groupes fondamentaux (en \sisi{car.}{caract\'eristique}~$0$, et en \sisi{car.}{caract\'eristique} $p> 0$) ont un m\^eme quotient maximal d'ordre premier \`a $p$.
\end{probleme}

Bien entendu, \Ref{XIII.3.3} nous montre que dans \Ref{XIII.3.1}, il y a lieu \'egalement de remplacer~$X'$ par un sch\'ema de la forme $X-Y$, o\`u $Y$ est une partie ferm\'ee de $X$ qui\pageoriginale est de codimension $\ge 2$ dans toute composante de $X$ la contenant. Lorsqu'on abandonne cette restriction sur $Y$, il doit exister encore un r\'esultat de finitude analogue, \`a condition de mettre des restrictions genre \og tame\fg sur la ramification en les points maximaux des composantes irr\'eductibles de $Y$ qui sont de codimension $1$.

\begin{probleme} \label{XIII.3.4}
Soit $A$ un anneau local noeth\'erien complet de dimension $2$, \`a corps r\'esiduel alg\'ebriquement clos. Soit encore $X=\Spec(A)$, $X'=X-\{\aaa\}$. Trouver des propri\'et\'es de structure particuli\`eres de $\pi_1(X')$ pour le cas o\`u $A$ est une intersection compl\`ete.
\end{probleme}

Une solution satisfaisante de ce probl\`eme permettrait peut-\^etre de r\'esoudre le vieux probl\`eme suivant:

\refstepcounter{subsection}\label{XIII.3.5}
\begin{enonce*}{Conjecture 3.5\ndemark}\ndetext{ce probl\`eme est, \`a l'automne 2004, encore ouvert.} Trouver une courbe irr\'eductible dans $\PP^3_k$ ($k$ corps alg\'ebriquement clos), de pr\'ef\'erence non singuli\`ere, qui ne soit pas ensemblistement l'intersection de deux hypersurfaces. \sisi{(KNESER~\cite{XIII.4} montre qu'on peut l'obtenir toujours comme intersection de trois hypersurfaces).}{\ignorespaces}
\end{enonce*}
\sisi{\ignorespaces}{\noindent(Kneser~\cite{XIII.4} montre qu'on peut l'obtenir toujours comme intersection de trois hypersurfaces).}

\section[Probl\`emes li\'es aux $\pi_i$ sup\'erieurs]{Probl\`emes li\'es aux $\pi_i$ sup\'erieurs: th\'eor\`emes de Lefschetz locaux et globaux pour les espaces analytiques complexes\protect\ndemark} \label{XIII.4}

Soit $X$ un sch\'ema\ndetext{les \'enonc\'es sont pr\'ecis\'es dans les Commentaires (section~\Ref{XIII.6}). Les conjectures qui y apparaissent sont devenues des th\'eor\`emes, \cf les notes de bas de page de la section~\Ref{XIII.6}.} localement de type fini sur le corps des complexes $\CC$, on peut lui associer un espace analytique $X^h$ sur $\CC$, d'o\`u des invariants d'homotopie et d'homologie $\pi_i(X^h)$, $\H_i(X^h)$, $\H^i(X^h)$ etc. \sisi{...}{} On sait d'ailleurs que $X$ est connexe si et seulement si $X^h$ l'est, donc que l'on a une bijection
$$
\pi_0(X^h)\to \pi_0(X)$$
De m\^eme, comme tout rev\^etement \'etale $X'$ de $X$ d\'efinit un rev\^etement \'etale $X^{\prime h}$ de~$X^h$\pageoriginale, on~a un homomorphisme canonique
$$
\pi_1(X^h)\to \pi_1(X),
$$
dont on sait, utilisant un th\'eor\`eme de Grauert-Remmert, qu'il identifie le deuxi\`eme groupe au compactifi\'e du premier pour la topologie des sous-groupes d'indice fini (ce qui exprime simplement le fait que $X'\sisi{\rightsquigarrow}{\mto} X^{\prime h}$ est une \'equivalence de la cat\'egorie des rev\^etements \'etales de $X$ avec la cat\'egorie des rev\^etements \'etales \emph{finis} de $X^h$). Il s'ensuit que les r\'esultats de ce s\'eminaire (par voie purement alg\'ebrique) sur $\pi_0(X)$ et $\pi_1(X)$ implique des r\'esultats pour $\pi_0(X^h)$ et $\pi_1(X^h)$ (qui sont de nature transcendante). Si d'ailleurs $X$ est propre, la suite exacte bien connue $0\to \ZZ \to \CC\to \CC^*\to 0$ permet de montrer que le groupe de N\'eron-Severi de $X$ (quotient de son groupe de Picard par la composante connexe de l'\'el\'ement neutre) est isomorphe \`a un sous-groupe de $\H^2(X^h, \ZZ)$; dans le cas non singulier k\"ahl\'erien, c'est le sous-groupe not\'e $\H^{(1, 1)}(X^h, \ZZ)$ (classes de type $(1, 1)$):
$$
\Pic(X)/\Pic^0(X)\subset \H^2(X, \ZZ).
$$

Par suite, les renseignements que nous avons obtenus sur les groupes de Picard impliquent des renseignements, tr\`es partiels il est vrai, sur les groupes $\H^2(X^h, \ZZ)$. Il est tentant de compl\'eter tous ces r\'esultats fragmentaires par des conjectures.

Des indications tr\`es pr\'ecises, allant dans le m\^eme sens que ceux qu'on vient de signaler, sont fournies par un classique th\'eor\`eme de \sisi{LEFSCHETZ}{Lefschetz}~\cite{XIII.7}. Il affirme que si~$X$ est un espace analytique projectif non singulier irr\'eductible de dimension $n$, et si~$Y$ est une section hyperplane non singuli\`ere, alors l'injection
$$Y^{n-1}\to X^n$$
induit\pageoriginale un homomorphisme
$$
\pi_i(Y^{n-1})\to \pi_i(X^n)$$
qui est un \emph{isomorphisme pour} $i\le n-2$, un \emph{\'epimorphisme pour} $i=n-1$. Il en r\'esulte l'\'enonc\'e analogue pour les homomorphismes
$$
\H_i(Y^{n-1})\to \H_i(X^n)$$
sur l'homologie (enti\`ere pour fixer les id\'ees), tandis qu'en cohomologie,
$$
\H^i(X^n)\to\H^i(Y^{n-1})
$$
est un isomorphisme en dimension $i\le n-2$, un monomorphisme en dimension $i=n-1$. Nous avons obtenu des variantes de ces r\'esultats dans le cadre des sch\'emas, pour $\pi_0$, $\pi_1$, $\Pic$, valables d'ailleurs sans hypoth\`eses de non singularit\'e dans une large mesure, \cf \Exp \Ref{XII}. \ De plus, dans l'\'elimination des hypoth\`eses de non singularit\'e, nous avons utilis\'e de fa\c con essentielle des variantes \og locales\fg de ces th\'eor\`emes de Lefschetz globaux. Tout ceci sugg\`ere les probl\`emes suivants, qui devront sans doute \^etre attaqu\'es simultan\'ement\sfootnote{Les formulations \Ref{XIII.4.1} \`a \Ref{XIII.4.3} qui suivent sont provisoires. Voir conjectures \ref{XIII.6.A} \`a \ref{XIII.6.D} plus bas, dans \og commentaires \`a l'\Exp \Ref{XIII}\fg pour des formulations plus satisfaisantes, ainsi que \Exp \Ref{XIV}.}.

\begin{probleme} \label{XIII.4.1}
Soient $X$ un espace analytique, $Y$ une partie analytique ferm\'ee de~$X$ (ou simplement une partie ferm\'ee?)\nde{la signification de cette question n'est pas claire; en effet, l'\'enonc\'e m\^eme du probl\`eme ne semble pas avoir de sens dans ce cas puisque la codimension de $Y$ dans $X$ n'est pas d\'efinie lorsque~$Y$ n'est plus suppos\'e analytique.} telle que pour tout $x\in Y$, l'anneau local $\OXx$ soit une intersection compl\`ete. Soit $n$ la codimension complexe de $Y$ dans $X$. L'homomorphisme canonique
$$
\pi_i(X-Y)\to \pi_i(X)$$
est-il un isomorphisme pour $i\le n-2$, et un \'epimorphisme pour $i=n-1$?
\end{probleme}

\enlargethispage{\baselineskip}%
Dans ce probl\`eme et les \sisi{pr\'ec\'edents}{suivants}, on suppose \'evidemment implicitement choisi un point-base pour d\'efinir les groupes d'homotopie. Pour \'enoncer le probl\`eme suivant, il faut d\'efinir, pour un espace analytique $X$ (plus g\'en\'eralement, pour un espace localement connexe par arc) et un $x\in X$, des invariants\pageoriginale locaux $\pi_i^x(X)$\sfootnote{Si $i\ge 2$. Pour le cas $i\le 1$, \cf \emph{Commentaires} dans \numero 6 ci-dessous, page~\pageref{p25y}.}.\refstepcounter{toto}\label{pagexpo} Pour ceci, on choisit une application non constante $f$ de l'intervalle $[0, 1]$ dans $X$, telle que $f(0)=x$ et $f(t)\neq x$ pour $t\neq 0$ (il en existe si $x$ n'est pas un point isol\'e). Alors pour tout voisinage $U$ de $x$, il existe un $\varepsilon>0$ tel que $0<t<\varepsilon$ implique $f(t)\in U$, et les groupes d'homotopies $\pi_i(U-x, f(t))$ sont essentiellement ind\'ependants de $t$ (ils sont, pour $t$ variable, reli\'es par un syst\`eme transitif d'isomorphismes), on peut les noter $\pi_i(U-x, f)$. On pose alors\refstepcounter{toto}\label{pip15}\label{pXIII.15}
$$
\pi_i^x(X)= \varprojlim_{U}\pi_{i-1}(U-x, f)
$$
la limite projective \'etant prise sur le syst\`eme des voisinages ouverts $U$ de $x$. En toute rigueur, cette limite d\'epend de $f$, et devrait \^etre not\'ee $\pi_i^x(X, f)$, mais on v\'erifie que pour $f$ variable, ces groupes sont isomorphes entre eux\sfootnote{Du moins si $x$ ne disconnecte pas $X$ au voisinage de $x$, \cf \emph{Commentaires} ci-dessous, page~\pageref{p25y}.}, de fa\c con pr\'ecise ils forment un syst\`eme local sur l'espace des chemins du type envisag\'e issus de $x$. Ces invariants sont la version homotopique des invariants de cohomologie locale $\H^i_x(F)$ pour un faisceau $F$ sur $X$, introduit dans \Ref{I}, et devraient jouer le r\^{o}le de \emph{groupes d'homotopie locale relatifs} de $X$ modulo $X-x$. Leur annulation pour $i\le n$ et pour tout $x\in Y$, o\`u~$Y$ est une partie ferm\'ee de $X$ de dimension topologique $\le d$, devrait entra\^iner que les homomorphismes
$$
\pi_i(X-Y) \to \pi_i(X)
$$
sont bijectifs pour $i<n-d$, et surjectif pour $i=n-d$\sfootnote{Pour une formulation corrig\'ee, \cf \emph{Commentaires} ci-dessous, page~\pageref{p25y}.}. De ce point de vue, \Ref{XIII.4.1} impliquerait (pour $Y$ r\'eduit \`a un point) une conjecture de nature purement locale, s'exprimant par \refstepcounter{toto}
$$
\pi_i^x(X)=0 \rm{\ pour\ } i\le n-1\label{p15l}
$$
lorsque $X$ est une intersection compl\`ete de dimension $n$ en $x$.

\sisi{A}À titre d'exemple d'invariants locaux $\pi_i^x(X)$, notons que si $x$ est un point non singulier de \sisi{\ignorespaces}{$X$ de} dimension complexe $n$, alors\pageoriginale
$$
\pi_i^x(X)=\pi_{i-1}(S^{2n-1}),
$$
o\`u $S^{2n-1}$ d\'esigne la sph\`ere de dimension $2n-1$. En particulier dans ce cas $\pi_{i}^x(X)=0$ pour $i\le 2n-1$, ce qui correspond au fait que si d'une \emph{vari\'et\'e} topologique $X$ on enl\`eve une partie ferm\'ee $Y$ de codimension $\ge m$, alors $\pi_i(X-Y)\to \pi_i(X)$ est un isomorphisme pour $i\le m-2$ et un \'epimorphisme pour $i=m-1$.

Ceci pos\'e:

\begin{probleme} \label{XIII.4.2}
Soient $X$ un espace analytique, $x\in X$, $t$ une section de $\OX$ \ s'annulant en $x$, $Y$ l'ensemble des z\'eros de $t$. Supposons les conditions suivantes satisfaites:
\begin{enumeratea}
\item
$t$ est r\'eguli\`ere en $x$ (\ie non diviseur de $0$ en $x$, hypoth\`ese peut-\^etre superflue, d'ailleurs).
\item
En les points $x'$ de $X-Y$ voisins de $x$, $\Oo_{X, x'}$ est une intersection compl\`ete (hypoth\`ese qui doit pouvoir se remplacer par la suivante plus g\'en\'erale si \Ref{XIII.4.1} est vraie: pour $x'$ comme ci-dessus, $\pi_i^{x'}(X)=0$ pour $i\le n-1$).
\item
En les points $y$ de $Y-\{ x\}$ voisins de $x$, on a
$$
\prof \Oo_{X, y}\ge n$$
(il suffit par exemple que l'on ait $\prof\Oo_{X, x}\ge n$).
\end{enumeratea}

Sous ces conditions, l'homomorphisme canonique
$$
\pi_i^x(Y)\to \pi_i^x(X)$$
est-il un isomorphisme pour $i\le n-2$, un \'epimorphisme pour $i=n-1$?
\end{probleme}

Voici enfin une variante globale de \Ref{XIII.4.2}, qui devrait s'en d\'eduire par consid\'eration du c\^{o}ne projetant en son origine, et qui g\'en\'eraliserait les th\'eor\`emes de Lefschetz classiques:

\begin{probleme} \label{XIII.4.3}
Soient\pageoriginale $X$ un espace analytique projectif, muni d'un Module inversible~$L$ ample, $t$ une section de $L$, $Y$ l'ensemble des z\'eros de $t$. Supposons:
\begin{enumeratea}
\item
$t$ est une section r\'eguli\`ere (\sisi{\ignorespaces}{hypoth\`ese} peut-\^etre superflue).
\item
Pour tout $x\in X-Y$, $\Oo_{X, x}$ est une intersection compl\`ete (devrait pouvoir se remplacer par $\pi_i^{x}(X)=0$ pour $i\le n-1$).
\item
Pour tout $x\in Y$, $\prof \Oo_{X, x}\ge n$.
\end{enumeratea}

Sous ces conditions, l'homomorphisme
$$
\pi_i(Y)\to \pi_i(X)$$
est-il un isomorphisme pour $i\le n-2$, un \'epimorphisme pour $i=n-1$?
\end{probleme}

Nous laisserons au lecteur le soin d'\'enoncer des conjectures analogues de nature \emph{cohomologique}\sfootnote{\Cf \Exp \Ref{XIV} les r\'esultats correspondants en th\'eorie des sch\'emas.}, les hypoth\`eses et conclusions portant alors sur les invariants cohomologiques locaux (\`a coefficients dans un groupe donn\'e). En tout \'etat de cause, le r\'esultat-clef semble devoir \^etre \Ref{XIII.4.2}, quand l'hypoth\`ese $b)$ y est prise sous forme resp\'ee, \ --- qu'on se place au point de vue de l'homologie, ou de l'homotopie.

Nous avons \'enonc\'e ces conjectures dans le cadre transcendant, dans
l'espoir d'y int\'eresser les topologues et de les convaincre que
les questions du type \og Lefschetz\fg sont loin d'\^etre closes.
Bien entendu, maintenant que nous sommes sur le point de disposer
d'une bonne th\'eorie de la cohomologie des sch\'emas (\`a coefficients
finis), gr\^ace aux travaux r\'ecents de M. \sisi{ARTIN}{Artin}, les
m\^emes questions se posent dans le cadre des sch\'emas, et il est
difficile de douter qu'elles ne re\c coivent une r\'eponse positive,
dans un avenir prochain\sfootnote{voir note pr\'ec\'edente.}.

\section{Probl\`emes li\'es aux groupes de Picard locaux} \label{XIII.5}

Un\pageoriginale premier probl\`eme fondamental, signal\'e pour la premi\`ere fois par \sisi{MUMFORD}{Mumford}~\cite{XIII.5} dans un cas particulier, est le suivant. Soient $A$ un anneau local complet de corps r\'esiduel $k$, $X=\Spec(A)$, $U=\Spec(A)-\{a\}$, o\`u $a$ est l'id\'eal maximal de $A$ \ie le point ferm\'e de $\Spec(A)$. On se propose de construire un syst\`eme projectif strict $G$ de groupes localement alg\'ebriques $G_i$ sur $k$, et un isomorphisme naturel \begin{equation*} \label{eq:XIII.5.+} \tag{$+$} {\Pic(U)\simeq G(k)} \end{equation*} o\`u on pose \'evidemment $G(k)=\varprojlim G_i(k)$. De fa\c con heuristique, on se propose de \og mettre une structure de groupe alg\'ebrique\fg (ou, du moins, pro-alg\'ebrique, en un sens convenable) sur le groupe $\Pic(\sisi{X'}{U})$.

Il est \'evident que tel quel, le probl\`eme n'est pas assez pr\'ecis, car la donn\'ee d'un isomorphisme~(\Ref{eq:XIII.5.+}) est loin de caract\'eriser le pro-objet $G$. Si $A$ contient un sous-corps not\'e encore $k$, qui soit un corps de repr\'esentants, on peut pr\'eciser le probl\`eme, en exigeant que pour une extension $k'$ de $k$ variable, on ait un isomorphisme, fonctoriel en $k'$: \begin{equation*} \label{eq:XIII.5.+'} \tag{$+'$} {\Pic(U')\simeq G(k')} \end{equation*} o\`u $U'$ est l'ouvert analogue \`a $U$ dans $\Spec(A')$, $A'=A\hat{\otimes}_k k'$. On peut proc\'eder de fa\c con analogue m\^eme si $A$ n'a pas de corps de repr\'esentants, pourvu que $k$ soit parfait, ce qui permet alors de construire fonctoriellement un $A'$ \og par extension r\'esiduelle $k'/k$\fg. D'ailleurs, lorsque $A$ admet un corps de repr\'esentants, la structure alg\'ebrique qu'on trouvera sur $\Pic(U)$ d\'ependra essentiellement du choix de ce corps de repr\'esentants (comme on voit d\'ej\`a sur le c\^{o}ne projetant d'une courbe elliptique), il semble donc qu'il faille partir d'un \og pro-anneau alg\'ebrique sur $k$\fg \`a la \sisi{GREENBERG}{Greenberg}~\cite{XIII.3}, pour arriver \`a d\'efinir le pro-objet $G$. Il est d'ailleurs concevable que dans le cas o\`u il n'y a pas de corps de repr\'esentants donn\'e, on ne trouve qu'un syst\`eme projectif de groupes\pageoriginaled \emph{quasi}-alg\'ebriques au sens de Serre, ou plut\^{o}t quasi-localement alg\'ebriques (les groupes $G_i$ obtenus ne seront pas en g\'en\'eral de type fini sur $k$, mais seulement localement de type fini sur $k$). Il est m\^eme possible qu'on ne trouvera en g\'en\'eral qu'une structure encore plus faible sur $\Pic(U)$, du genre de celles rencontr\'ees par \sisi{NERON}{N\'eron}~\cite{XIII.6} dans sa th\'eorie de d\'eg\'en\'erescence des vari\'et\'es ab\'eliennes d\'efinies sur des corps locaux.

Une m\'ethode pour attaquer le probl\`eme, \'egalement introduite par \sisi{MUMFORD}{Mumford}, consiste \`a d\'esingulariser $X$, \ie \`a consid\'erer un morphisme projectif birationnel $Y\to X$ avec $Y$ r\'egulier. Lorsque $U$ est r\'egulier (\ie $a$ est un point singulier isol\'e), on peut souvent trouver $Y$ de telle fa\c con que $Y\sisi{|U}{_|U}=V\to U$ soit un isomorphisme. Dans ce cas, on aura donc $$ \Pic(U)\simeq\Pic(V)\simeq\Pic(Y)/\Im \ZZ^I, $$ o\`u $I$ est l'ensemble des composantes irr\'eductibles de la fibre $Y_{a}$\sisi{, }{} (chacune de celles-ci d\'efinissant un \'el\'ement de $\Pic(Y)$, \'etant un diviseur localement principal, gr\^ace \`a $Y$ r\'egulier). D'autre part, utilisant la technique de G\'eom\'etrie Formelle \EGA III~4~et~5, notamment le th\'eor\`eme d'existence, on trouve $$ \Pic(Y)\simeq\varprojlim\Pic(Y_n), $$ o\`u $Y_n=Y\otimes_A A_n$, $A_n=A/\mm^{n+1}$. Lorsque $A$ admet un corps de repr\'esentants $k$, on dispose de la th\'eorie des sch\'emas de Picard des sch\'emas projectifs $Y_n$ sur $k$, donc on~a $$ \Pic(Y_n)\simeq \mathbf{Pic}_{Y_n/k}(k). $$

Cela fournit donc une construction d'un syst\`eme projectif de groupes localement alg\'ebriques $\mathbf{Pic}_{Y_n/k}/\Im \ZZ^I$, qui est le syst\`eme cherch\'e.\refstepcounter{toto}\nde{\label{Boutot}la question a \'et\'e grandement \'eclaircie par les r\'esultats de {Boutot} (Boutot~J.-F., \emph{Sch\'ema de {P}icard local}, Lect. Notes in Math., vol.~632, Springer, Berlin, 1978). En particulier, si $A$ est une $k$-alg\`ebre locale compl\`ete (noeth\'erienne) de profondeur $\geq 2$ telle que $H^2_\m(A)$ est de dimension finie sur $k$, le groupe de Picard local est un sch\'ema en groupes localement de type fini sur $k$, d'espace tangent \`a l'origine $H^2_\m(A)$. Si $A$ est de plus normal de dimension $\geq 3$, le crit\`ere de normalit\'e de Serre~\Ref{XI}~\Ref{XI.3.11} joint au corollaire~\Ref{V}~\Ref{V.3.6} assurent la finitude requise et, d\`es lors, l'existence du sch\'ema de Picard local. Voir aussi (Lipman~J., {\og The Picard group of a scheme over an Artin ring\fg}, \emph{Publ. Math. Inst. Hautes \'Etudes Sci.} \textbf{46} (1976), p\ptbl 15--86) pour une d\'emarche plus proche de celle de Grothendieck esquiss\'ee plus haut.} Dans le cas envisag\'e ici, on peut d'ailleurs voir (utilisant que $a$ est un point singulier isol\'e) que les composantes connexes des sous-groupes images universelles\pageoriginale dans ce syst\`eme projectif forment un syst\`eme projectif \emph{essentiellement constant}, donc en l'occurrence on trouve un groupe localement alg\'ebrique $G$ comme solution du probl\`eme. Si on suppose m\^eme $A$ normal de dimension $2$, alors une remarque de \sisi{MUMFORD}{Mumford} (disant que la matrice d'intersection des composantes de $Y_{a}$ dans $X$ est d\'efinie n\'egative\nde{Mumford~D., {\og The topology of normal singularities of an algebraic surface and a criterion for simplicity\fg}, \emph{Publ. Math. Inst. Hautes \'Etudes Sci.} \textbf{9} (1961), p\ptbl 5--22.}) implique que $G$ est m\^eme un groupe alg\'ebrique, \ie de type fini sur $k$ (le nombre de ses composantes connexes \'etant d'ailleurs \'egal au d\'eterminant de la matrice d'intersection envisag\'ee il y a un instant).

Si par contre $a$ n'est pas une singularit\'e isol\'ee, on se convainc sur des exemples (avec $A$ de dimension $2$) qu'on trouve un syst\`eme projectif de groupes alg\'ebriques, ne se r\'eduisant pas \`a un seul groupe alg\'ebrique.

\refstepcounter{toto}\label{pXIII20}%
Une fois qu'on disposerait d'une bonne notion de \og \emph{sch\'ema de Picard local}\fg, il y aurait lieu de renforcer la notion de parafactorialit\'e, en disant que $A$ est \og \emph{g\'eom\'etriquement parafactoriel}\fg, lorsque non seulement $A$ et m\^eme $\widehat{A}$ sont parafactoriels, mais que le sch\'ema de Picard local $G(\widehat{A})$ est le groupe trivial (ce qui est plus fort, lorsque le corps r\'esiduel n'est pas alg\'ebriquement clos, que de dire que $G$ n'a d'autre point rationnel sur $k$ que l'unit\'e). On se rend compte de la n\'ecessit\'e d'une notion renforc\'ee de parafactorialit\'e, en se rappelant qu'il existe des anneaux locaux complets normaux de dimension $2$ qui sont factoriels, mais qui admettent des alg\`ebres finies \'etales qui ne le sont pas\nde{les anneaux factoriels \`a hens\'elis\'e non factoriel arrivent naturellement lorsqu'on \'etudie les espaces de modules de fibr\'es vectoriels. Voir par exemple (Dr\'ezet~J.-M., {\og Groupe de Picard des vari\'et\'es de modules de faisceaux semi-stables sur $\mathbf{P}\sb 2$\fg}, in \emph{Singularities, representation of algebras, and vector bundles (Lambrecht, 1985)},
Lect. Notes in Math., vol.~1273, Springer, Berlin, 1987, p\ptbl 337--362). \textit{Stricto sensu}, Dr\'ezet montre que le compl\'et\'e n'est pas factoriel, mais en fait la preuve donne le r\'esultat pour l'hens\'elis\'e: le point est que le th\'eor\`eme de slice \'etale de Luna (Luna~D., {\og Slices \'etales\fg}, in \emph{Sur les groupes alg\'ebriques}, M\'em. Soc. math. France, vol.~33, Soci\'et\'e math\'ematique de France, Paris, 1973, p\ptbl 81--105) d\'ecrit l'anneau local d'un quotient au sens de la g\'eom\'etrie invariante pr\`es d'un point semi-stable localement pour la topologie \'etale.}. Un anneau local \og \emph{g\'eom\'etriquement factoriel}\fg serait alors un anneau~$A$ normal tel que tous \sisi{c}{l}es localis\'es de dimension $\geq 2$ soient g\'eom\'etriquement parafactoriels, ou mieux, tel que les localis\'es de $\widehat{A}$ soient parafactoriels\sfootnote{Pour une notion plus souple d'anneau local \og g\'eom\'etriquement factoriel\fg, \cf \emph{Commentaires}, page~\pageref{commpara}.}. Bien entendu, il serait int\'eressant de trouver une \og bonne\fg d\'efinition de ces notions, ind\'ependante de la th\'eorie, encore \`a faire, des sch\'emas de Picard locaux.\nde{voir page~\pageref{commpara}: un anneau local est g\'eom\'etriquement factoriel (\resp parafactoriel) si son hens\'elis\'e strict est factoriel (\resp parafactoriel).}

Il est en tous cas plausible qu'on aura besoin de ces notions si on d\'esire obtenir des \'enonc\'es du type suivant: Soit $A$ un \og bon anneau\fg (par exemple une alg\`ebre de type fini sur $\ZZ$, ou sur un anneau local complet, par exemple sur un corps). Soit~$U$ l'ensemble des $x\in X=\Spec(A)$ tels que $\OXx$ soit \og g\'eom\'etriquement\pageoriginale factoriel\fg, alors $U$ est ouvert\sisi{?.}{?} Ou encore: Soit $f\colon X\to Y$ un morphisme plat de type fini avec $Y$ localement noeth\'erien, soit $U$ l'ensemble des $x\in X$ tels que $\Oo_{X_{f(x)}, x}$ soit \og g\'eom\'etriquement factoriel\fg, alors $U$ est ouvert, du moins sous des conditions suppl\'ementaires sympathiques sur $f$\sisi{?.}{?} Je doute qu'avec la notion habituelle d'anneau factoriel, il existe des \'enonc\'es vrais de ce type.

Nous avons soulev\'e ici, dans un cas particulier, la question de l'\'etude des propri\'et\'es g\'eom\'etriques d'anneaux locaux \og variables\fg, par exemple les $\OXx$ pour $x$ parcourant un pr\'esch\'ema $X$. Lorsque $X$ est un sch\'ema de type fini sur un corps, par exemple, on sait\nde{voir \EGA IV.16.} qu'il existe sur $X$ un syst\`eme projectif d'alg\`ebres finies $P^n_{X/k}$ (obtenu en compl\'etant $X\times_k X$ le long de la diagonale), dont la fibre en tout point $x\in X$ rationnel sur $k$ est isomorphe au syst\`eme projectif des $\OXx/\mm_X^{n+1}$. Il est alors naturel de relier l'\'etude des compl\'et\'es des anneaux locaux $\OXx$, pour $x$ variable, \`a celle de la \og famille alg\'ebrique d'anneaux locaux complets\fg donn\'ee par les $P^n$, en notant que pour tout $x\in X$ (rationnel sur $k$ ou non), on obtient un anneau local complet
$$
P^{\infty}(x)=\varprojlim P^n(x)
$$
(o\`u $P^n(x)=$~fibre~r\'eduite~$P^n\otimes_{\OXx}k(x)$). Un int\'er\^et particulier s'attachera par exemple \`a l'anneau complet associ\'e ainsi au point g\'en\'erique, et on s'attendra \`a ce que ses propri\'et\'es alg\'ebrico-g\'eom\'etriques (s'exprimant par exemple par ses groupes de Picard, ou d'homotopie, ou d'homologie, locaux), seront essentiellement celles des compl\'et\'es $\widehat{\Oo}_{X, x}$ pour $x$ dans un ouvert dense convenable $U$.

On peut, de fa\c con g\'en\'erale, se proposer de faire l'\'etude simultan\'ee des anneaux locaux complets obtenus ainsi \`a partir d'un syst\`eme projectif adique $(P_n)$ d'alg\`ebres finies sur un sch\'ema donn\'e $X$. Il est plausible qu'on trouvera, moyennant certaines conditions de r\'egularit\'e (telle la platitude des $P_n$)\pageoriginale que les groupes d'homotopie locale proviennent d'un syst\`eme projectif de sch\'emas en groupes finis sur $X$, et qu'on aura des r\'esultats analogues pour les groupes de Picard locaux. En ce qui concerne ces \sisi{derni\`eres}{derniers}, un premier cas int\'eressant qui m\'erite d'\^etre investigu\'e est celui o\`u on part d'une surface alg\'ebrique $X$ ayant des courbes singuli\`eres, et qu'on se propose d'\'etudier les sch\'emas de Picard locaux en les points variables sur celles-ci, en termes d'un pro-sch\'ema en groupes convenable d\'efini sur le lieu singulier.

\section[Commentaires]{Commentaires\protect\sfootnotemark} \label{XIII.6}

Le\pageoriginale point de vue de\sfootnotetext{R\'edig\'es en Mars~1963.} la \og Cohomologie \'etale\fg des sch\'emas et des progr\`es r\'ecents dans cette th\'eorie, nous am\`ene\sisi{}{nt} \`a pr\'eciser et \`a \'elargir en m\^eme temps certains des probl\`emes pos\'es. Pour la notion de \og topologie\fg et de \og topologie \'etale d'un sch\'ema\fg, je renvoie \`a M\ptbl\sisi{ARTIN}{Artin}, \sisi{Grothendieck Topologies}{\emph{Grothendieck Topologies}}, Harvard University 1962 (notes mim\'eographi\'ees)\sfootnote{Ou de pr\'ef\'erence, \`a \SGA 4.}.

Cette th\'eorie, par une notion plus fine de localisation que celle que fournit la traditionnelle \og topologie de Zariski\fg, am\`ene \`a attacher un int\'er\^et particulier aux anneaux \emph{strictement locaux}\refstepcounter{toto}\label{pXIII.6}, \ie les anneaux locaux hens\'eliens \`a corps r\'esiduel s\'eparablement clos. Pour tout anneau local $A$ de corps r\'esiduel $k$, et toute cl\^{o}ture s\'eparable $k'$ de~$k$, on peut trouver un homomorphisme local de $A$ dans un anneau strictement local~$A'$, la \emph{cl\^{o}ture strictement locale} de $A$, de corps r\'esiduel $k'$, ayant une propri\'et\'e universelle \'evidente. $A'$ est hens\'elien, plat sur $A$, et $A'\otimes_A k\simeq k'$; il est noeth\'erien si et seulement si $A$ l'est. (\Cf \loccit Chap\ptbl III, section 4)\sfootnote{Ou \EGA IV~18.8.}. Si $X$ est un pr\'esch\'ema, et $x$ un point de~$X$, $x'$ un point au-dessus de $x$, spectre d'une cl\^{o}ture s\'eparable $k'$ de $k=k(x)$, on est amen\'e \`a d\'efinir l'\emph{anneau strictement local} de $X$ en $x'$, $\Oo'_{X, x'}$, comme la cl\^{o}ture strictement locale de l'anneau local habituel $\OXx$, relativement \`a l'extension r\'esiduelle $\sisi{K'/K}{k'/k}$. Ce sont les anneaux \emph{strictement} locaux des points \og g\'eom\'etriques\fg de $X$ qui, au point de vue de la topologie \'etale, sont sens\'es refl\'eter les propri\'et\'es locales du pr\'esch\'ema $X$. Ils jouent aussi, \`a bien des \'egards, le r\^{o}le qu'on faisait jouer aux \emph{compl\'et\'es} des anneaux locaux de $X$ (disons, en les points \`a corps r\'esiduel alg\'ebriquement clos), tout en restant \og plus proches\fg de $X$ et permettant un passage plus ais\'e aux \og points voisins\fg.

Il y a lieu alors de reprendre un bon nombre de questions, qu'on pose g\'en\'eralement pour les anneaux locaux complets (\'eventuellement restreints \`a avoir\pageoriginale un corps r\'esiduel alg\'ebriquement clos), pour les anneaux locaux noeth\'eriens hens\'eliens (\resp les anneaux strictement locaux noeth\'eriens). C'est ainsi que les probl\`emes topologiques soulev\'es dans les n$^{\textup{os}}$~\Ref{XIII.2} et \Ref{XIII.3}, se posent plus g\'en\'eralement pour les anneaux strictement locaux. On peut d'ailleurs \'enoncer \`a titre conjectural, pour les \og bons\fg anneaux strictement locaux, certaines propri\'et\'es de simple connexion et d'acyclicit\'e pour les fibres g\'eom\'etriques du morphisme canonique $\Spec(\widehat{A})\to\Spec(A)$, qui montreraient que pour beaucoup de propri\'et\'es de nature \og topologique\fg, il revient au m\^eme de les prouver pour l'anneau $A$, ou pour son compl\'et\'e $\widehat{A}$. Certains r\'esultats obtenus d\'ej\`a dans cette voie\sfootnote{\Cf M\ptbl\sisi{ARTIN}{Artin} dans \SGA 4~XIX} permettent d'esp\'erer qu'on disposera bient\^{o}t de r\'esultats complets dans cette direction.

\refstepcounter{toto}\label{commpara}%
La notion de localisation \'etale fournit une d\'efinition qui semble raisonnable de la notion d'anneau local \og \emph{g\'eom\'etriquement parafactoriel}\fg ou \og \emph{g\'eom\'etriquement factoriel}\fg\sisi{,}{} (dont le besoin a \'et\'e signal\'e dans \numero\Ref{XIII.5}, p\ptbl\sisi{20}{\pageref{pXIII20}}): on appellera ainsi un anneau local dont la cl\^{o}ture strictement locale est parafactorielle, \resp factorielle. Des hypoth\`eses de cette nature s'introduisent effectivement d'une fa\c con naturelle dans l'\'etude de la cohomologie \'etale des pr\'esch\'emas\nde{voir par exemple (Strano~R., {\og The Brauer group of a scheme\fg}, \emph{Ann. Mat. Pura Appl. (4)} \textbf{121} (1979), p\ptbl 157--169) o\`u l'hypoth\`ese de g\'eom\'etrique parafactorialit\'e des anneaux locaux d'un sch\'ema $X$ permet parfois de montrer la co\"incidence des groupes de Brauer de $X$ (calcul\'es en termes d'alg\`ebres d'Azumaya) et du groupe de Brauer cohomologique de $X$.}. Ainsi, si $X$ est un pr\'esch\'ema localement noeth\'erien dont les anneaux strictement locaux sont factoriels (\ie dont les anneaux locaux ordinaires sont \og g\'eom\'etriquement factoriels{\fg}), on montre que les $\H^i(X_{\et}, \GG_m)$ sont des groupes de torsion pour $i\geq 2$ (ce qui permet parfois d'exprimer ces groupes en termes des groupes de cohomologie \`a coefficients dans les groupes $\bbmu_n$ des racines $n$\sisi{.i{\^emes}}{-i\`emes} de l'unit\'e), et si $X$ est int\`egre de corps des fractions $K$, l'homomorphisme naturel $\H^2(X_{\et}, \GG_m)\to\H^2(K, \GG_m)=\Br(K)$ est injectif\nde{le lien entre groupe de Brauer et groupe de Picard est intime. Citons \`a ce propos les r\'esultats suivants de Saito (Saito~S., {\og Arithmetic on two-dimensional local rings\fg}, \emph{Invent. Math.} \textbf{85} (1986), \numero 2, p\ptbl 379--414) dans le cas des surfaces, le premier \'etant local l'autre global. Soit $A$ un anneau local excellent de dimension $2$, normal et hens\'elien \`a corps r\'esiduel \textit{fini} et $X$ le compl\'ementaire du point ferm\'e dans $\Spec(A)$. Alors, on~a une dualit\'e parfaite de groupes \emph{de torsion} $\Pic(X)\times \Br(X)\to \QQ/\ZZ$ --- par groupe de Brauer de $X$, on entend groupe de Brauer \emph{cohomologique} $\textrm{Br}(X)=H^2_\et(X,G_m)$. Dans le cas global, on~a la g\'en\'eralisation suivante d'un r\'esultat de Lichtenbaum (Lichtenbaum~S., {\og Duality theorems for curves over $p$-adic fields\fg}, \emph{Invent. Math.} \textbf{7} (1969), p\ptbl 120--136): soit $k$ le corps des fractions d'un anneau de valuation discr\`ete complet $\Oo$ \`a corps r\'esiduel fini et $X$ une courbe projective, lisse et g\'eom\'etriquement compl\`ete sur $k$. Le groupe $\Pic^0(X)$ est muni de la topologie induite de la topologie adique de $k$ et $\Pic(X)$ est le groupe topologique qui fait de $\Pic^0(X)$ un sous-groupe ouvert. Alors, on~a une dualit\'e parfaite de groupes topologiques $\Pic(X)\times \Br(X)\to \QQ/\ZZ$. Notons que cet \'enonc\'e, qui concerne les courbes, se d\'emontre bien entendu en consid\'erant un mod\`ele r\'egulier (propre et plat) de $X$ sur $\Oo$: c'est un r\'esultat sur les surfaces.} des exemples montrent que ces conclusions peuvent \^etre en d\'efaut, m\^eme pour $X$ local, si on suppose seulement $X$ factoriel au lieu de g\'eom\'etriquement factoriel\sfootnote{\Cf A\ptbl \sisi{GROTHENDIECK}{Grothendieck}, le groupe de Brauer~II (S\'eminaire Bourbaki \numero 297, Nov.~1965), notamment 1.8 et 1.11~b.}.

Concernant les probl\`emes de type Lefschetz local et global soulev\'es dans \Ref{XIII.3.4}, et leurs analogues en th\'eorie des sch\'emas, la version homologique de ces questions s'est consid\'erablement clarifi\'ee, tout r\'esultant formellement de trois th\'eor\`emes g\'en\'eraux, l'un concernant la dimension cohomologique de\pageoriginale certains sch\'emas affines (\resp des espaces de Stein), tels les sch\'emas affines $X$ de type fini sur un corps alg\'ebriquement clos: leur dimension cohomologique est $\leq\dim X$ (\og th\'eor\`eme de Lefschetz affine{\fg})\sfootnote{\Cf \SGA 4~XIV.}: l'autre \'etant un th\'eor\`eme de dualit\'e pour la cohomologie (\`a coefficients discrets) d'un morphisme projectif\sfootnote{\Cf \SGA 4~XVIII.}, enfin le dernier un th\'eor\`eme de \emph{dualit\'e locale} de nature analogue\sfootnote{\Cf \SGA 5~I.}. En G\'eom\'etrie Alg\'ebrique, seul ce dernier n'est pas d\'emontr\'e au moment d'\'ecrire ces lignes (il l'est cependant en \sisi{car.~$0$}{caract\'eristique $0$}, utilisant la r\'esolution des singularit\'es par \sisi{HIRONAKA}{Hironaka}). D'ailleurs, dans le cadre transcendant, on dispose d\`es \`a pr\'esent de la dualit\'e globale et locale, d\'emontr\'ees r\'ecemment par \sisi{VERDIER}{Verdier}\nde{voir Verdier~J.-L., {\og Dualit\'e dans la cohomologie des espaces localement compacts\fg}, in \emph{S\'eminaire Bourbaki}, vol.~9, Soci\'et\'e math\'ematique de France, Paris, 1995, \Exp 300, p\ptbl 337--349.}. Bornons-nous \`a indiquer que dans l'\'enonc\'e des versions homologiques des probl\`emes \Ref{XIII.4.2} et \Ref{XIII.4.3} (qui d\'esormais m\'eritent le nom de conjectures), les conditions \og \`a l'infini\fg a) et c) sont certainement superflues, seule \'etant importante la \emph{structure cohomologique locale} de $X-Y$, qu'on supposera par exemple localement intersection compl\`ete de dimension $\geq n$. De plus, dans \Ref{XIII.4.3} disons, le fait que $Y$ soit une section hyperplane ne devrait pas jouer, et doit pouvoir se remplacer par la seule hypoth\`ese que $X$ est compacte et $X-Y$ est Stein (\ie dans le cas de la G\'eom\'etrie Alg\'ebrique, $X$ est propre sur $k$ et $X-Y$ affine; comme nous le disions, la version homologique de cette conjecture est d\'emontr\'ee pour les espaces alg\'ebriques sur le corps $\CC$)\sfootnote{\Cf \Exp \Ref{XIV}}.

\refstepcounter{toto}\label{p25y}%
Dans la d\'efinition (p\ptbl\sisi{15}{\pageref{pip15}}) des $\pi_i^x(X)$, on doit supposer $i\geq 2$. Pour $i=0, 1$, il n'y a pas de d\'efinition raisonnable des $\pi_i^x(X)$\sisi{,}{;} il y a lieu de les remplacer par $\H_0^x(X)$ et $\H_1^x(X)$, d\'efinis respectivement comme le co\sisi{-}{}noyau et le noyau dans l'homomorphisme naturel
$$
\varprojlim \H_0(U-(x), \ZZ)\to\varprojlim\H_0(U, \ZZ).
$$

\sisi{A}À la rigueur et pour la commodit\'e des formulations, on pourra poser $\pi_i^x(X)=\H_i^x(X)$ pour $i\leq 1$, sinon il faut compl\'eter les assertions ult\'erieures concernant les $\pi_i^x$ par les assertions correspondantes pour $\H_0^x, \H_1^x$. Si $x$ est un point isol\'e de $X$, il convient de poser $\pi_i^x(X)=0$ pour $i\neq 0$, $\pi_0^x(X)=\H_0^x(X)=\ZZ$.

L'assertion\pageoriginale que les $\pi_i^x(U, f)$ soient isomorphes entre eux n'est vraie que lorsque $X$ n'est pas disconnect\'e par $x$ au voisinage de $x$, \ie si $\pi_i^x(X)=0$ pour $i=0, 1$. Dans le cas g\'en\'eral $\pi_i^x(X)$ ne peut d\'esigner qu'une \emph{famille} de groupes, pas n\'ecessairement isomorphes entre eux; cependant l'\'ecriture $\pi_i^x(X)=0$ garde un sens \'evident.

\refstepcounter{toto}\label{p26y}\sisi{Page~15, ligne~-~11}{Page~\pageref{p15l}}, o\`u je pr\'evois que l'annulation des invariants homotopiques locaux $\pi_i^x(X)$ pour $x\in Y$, $i\leq n$ doit entra\^iner la bijectivit\'e de $\pi_i(X-Y)\to\pi_i(X)$ pour $i<n-d$, la surjectivit\'e pour $i=n-d$, il convient d'\^etre prudent, faute de pouvoir disposer dans le contexte pr\'esent (comme en G\'eom\'etrie Alg\'ebrique) de points \og g\'en\'eraux\fg en lesquels les conditions locales devront aussi s'appliquer. Il sera sans doute n\'ecessaire, pour cette raison, de faire appel \`a des invariants homotopiques locaux relatifs
$$
\pi_i^Y(X, f)=\pi_i^Y(X, x)=\varprojlim_U \pi_{i-1}(U-U\cap Y, f(t))\qquad\text{pour }i\geq 2\quoi,
$$
(et d\'efinition {\sisi{ad-hoc}{\textit{ad hoc}}} comme ci-dessus pour $i=0, 1$), o\`u $Y$ est une partie ferm\'ee de~$X$; ou de suppl\'eer \`a l'absence de points g\'en\'eraux en exprimant les hypoth\`eses sur $X$ en termes de propri\'et\'es de nature topologique (pour la topologie \'etale) des spectres des anneaux locaux de $X$, ce qui permet de r\'ecup\'erer des points g\'en\'eraux. La m\^eme r\'eserve s'applique \`a la g\'en\'eralisation des conjectures \Ref{XIII.4.2} et \Ref{XIII.4.3} au cas o\`u $X-Y$ n'est pas suppos\'e localement une intersection compl\`ete, g\'en\'eralisation sugg\'er\'ee dans l'\'enonc\'e des conditions~b) de ces conjectures.

Pour formuler les versions expurg\'ees des conjectures \Ref{XIII.4.2} et \Ref{XIII.4.3} sugg\'er\'ees par les r\'esultats auxquels on~a fait allusion plus haut, il convient de poser la

\begin{supdefinition} \label{XIII.6.1} Soit $X$ un espace topologique, $Y$ une partie localement ferm\'ee de $X$, et $n$ un entier. On dit que $X$ est de \emph{profondeur homotopique} $\geq n$ le long de $Y$, et on \'ecrit $\mathrm{prof\:hpt}_Y(X)\geq n$, si pour tout $x\in Y$, on~a $\pi_i^Y(X, x)=0$ pour $i<n$.
\end{supdefinition}

Il doit\pageoriginale \^etre \'equivalent de dire que pour tout ouvert $X'$ de $X$, et tout $x\in X'\cap U=U'$ (o\`u $U=X-Y$), l'homomorphisme canonique
$$
\pi_i(U', x)\to\pi_i(X', x)$$
est un isomorphisme pour $i<n-1$, un monomorphisme pour $i=n-1$\nde{lorsque le couple $(X, Y)$ est de plus poly\'edral, cette \'equivalence est vraie; \cf (Eyral~C., {\og Profondeur homotopique et conjecture de Grothendieck\fg}, \emph{Ann. Sci. \'Ec. Norm. Sup. (4)} \textbf{33} (2000), \numero 6, p\ptbl 823--836).}.

\begin{supdefinition} \label{XIII.6.2}
Soit $X$ un espace analytique complexe, $n$ un entier, on dit que la \emph{profondeur homotopique rectifi\'ee} de $X$ est $n$\sfootnote{Dans la premi\`ere \'edition de ces notes, nous avions employ\'e le terme: \og vraie profondeur homotopique\fg. Dans la pr\'esente version, nous suivons \EGA IV~10.8.1.}, si pour toute partie analytique localement ferm\'ee $Y$ de $X$, on a
\begin{equation*} \label{eq:XIII.6.x} \tag{$x$} {\mathrm{prof\:htp}_Y(X)\geq n-\dim Y}
\end{equation*}
(o\`u, bien entendu, $\dim Y$ d\'esigne la dimension \emph{complexe} de $Y$).
\end{supdefinition}

Il doit \^etre \'equivalent de dire que pour toute partie analytique irr\'eductible $Y$ localement ferm\'ee dans $X$, il existe une partie analytique ferm\'ee $Z$ de $Y$, de dimension $<\dim Y$, telle que la relation (\Ref{eq:XIII.6.x}) soit valable pour $Y-Z$ au lieu de $Y$. Cela permettrait par exemple dans la d\'efinition \Ref{XIII.6.2} de se borner au cas o\`u $Y$ est non singuli\`ere.\refstepcounter{toto}\nde{\label{HL}toutes les conjectures qui suivent, convenablement rectifi\'ees si j'ose dire, sont devenues des th\'eor\`emes gr\^ace au travail de Hamm et L\^e D\~ung Tr\'ang (Hamm~H.A. \& L\^e D\~ung Tr\'ang, {\og Rectified homotopical depth and Grothendieck conjectures\fg}, in \emph{The Grothendieck Festschrift, Vol. II}, Progr. Math., vol.~87, Birkh\"auser, Boston, 1990, p\ptbl 311--351) cit\'e \cite{HL} dans ce qui suit. En ce qui concerne les deux d\'efinitions conjecturalement \'equivalentes de la profondeur rectifi\'ee, elles sont m\^emes \'equivalentes \`a une troisi\`eme, qui s'exprime en termes de stratification de Withney (\cf \loccit, th\'eor\`eme 1.4).}

\enlargethispage{\baselineskip}%
La conjecture suivante, de nature purement topologique, est dans la nature d'un \og \emph{th\'eor\`eme de Hurewicz local}\fg.

\refstepcounter{aconjecture} \label{XIII.6.A}
\begin{enonce*}{Conjecture A}[\og \sisi{th.}{th\'eor\`eme} de Hurewic\sisi{s}{z} local{\fg}\ndemark]\ndetext{comme observ\'e dans \cite{HL}, exemple 3.1.3, cette conjecture est fausse d\'ej\`a pour $X=\{\boldsymbol{z}\in\nobreak\CC^n\mid\nobreak z_1^2+z_2^3+\cdots+z_n^3=0\}, n\geq 4$ et $Y$ r\'eduit \`a l'origine. Mais, convenablement modifi\'ee, elle est vraie (th\'eor\`eme 3.1.4 de {\loccit}).}
Soit $X$ un espace topologique, $Y$~une partie localement ferm\'ee, soumis au besoin \`a des conditions de \og smoothness\fg genre triangulabilit\'e locale de la paire $(X, Y)$, $n$ un entier $\geq 3$. Pour qu'on ait $\mathrm{prof\:htp}_Y(X)\geq n$, il faut et il suffit que l'on ait\pageoriginale
$$
\H^i_Y(\ZZ_X)=0\quad\text{pour } i<n$$
(on dit alors que $X$ est de \emph{profondeur cohomologique} $\geq n$ le long de $Y$), et que les groupes fondamentaux locaux
$$
\sisi{\pi^2_Y}{\pi_2^Y}(X, x)=\varprojlim_{U\ni x} \pi_1(U-U\cap Y)$$
soient nuls (on dit alors que $X$ est \og pur\fg le long de $Y$).
\end{enonce*}

On notera que si $X$ est un espace analytique, $Y$ un sous-espace analytique, et si~$X$ est pur le long de $Y$, alors pour tout $x\in Y$, l'anneau local $\OXx$, ainsi que ses localis\'es par rapport \`a des id\'eaux premiers contenant l'id\'eal d\'efinissant le germe~$Y$ en~$x$ (\ie dans l'image inverse $Y_x$ de $Y$ par $\Spec(\OXx)=X_x\to X$) sont purs au sens de \Exp \Ref{X}; il semble plausible que la r\'eciproque soit \'egalement vraie. Des remarques analogues valent pour la profondeur cohomologique, \'etant entendu qu'on travaille avec la topologie \'etale sur les $\Spec(\OXx)$.

La conjecture \Ref{XIII.4.1} se g\'en\'eralise alors en la

\refstepcounter{aconjecture} \label{XIII.6.B}
\begin{enonce*}{Conjecture B}[\og Puret\'e{\fg}\ndemark]\ndetext{cette conjecture est prouv\'ee, m\^eme dans le cas o\`u $E$ est singulier, dans \cite{HL}: c'est le th\'eor\`eme 3.2.1.}
Soient $E$ un espace analytique, $X$ une partie analytique de $E$. On suppose que $E$ est non singulier de dimension $N$ en $x\sisi{}{\in X}$, et que $X$ peut se d\'ecrire par $p$ \'equations analytiques au voisinage de tout point. Alors la profondeur homotopique rectifi\'ee de $X$ est $\geq N-p$.
\end{enonce*}

En particulier, une intersection compl\`ete locale de dimension $n$ en tout point serait de profondeur homotopique rectifi\'ee $\geq n$, ce qui n'est autre que la conjecture \Ref{XIII.4.1}.

Les conjectures \Ref{XIII.4.2} et \Ref{XIII.4.3} se g\'en\'eralisent respectivement en:

\refstepcounter{aconjecture} \label{XIII.6.C}
\begin{enonce*}{Conjecture C}[\og Lefschetz local{\fg}\ndemark]\ndetext{cette conjecture est prouv\'ee dans \cite{HL}, ce m\^eme sous sa forme forte de la remarque qui suit, \cf th\'eor\`eme 3.3.1 de~\loccit}
Soient $X$ un espace analytique, $Y$ une partie analytique ferm\'ee, $x$ un point de $Y$, on suppose que $X-Y$ est Stein au voisinage de $x$ (par exemple $Y$ d\'efini par une \'equation en $x$), et que $X-Y$ est de profondeur homotopique rectifi\'ee $\geq n$ au voisinage de $x$ (par exemple, est en tout point de $X-Y$ voisin de $x$, une intersection compl\`ete de dimension $\geq n$, \cf conjecture \Ref{XIII.6.B}). Alors l'homomorphisme canonique
$$
\pi_i^x(Y)\to\pi_i^x(X)$$
est un isomorphisme pour $i<n-1$, un \'epimorphisme pour $i=n-1$.
\end{enonce*}

\refstepcounter{aconjecture} \label{XIII.6.D}
\begin{enonce*}{Conjecture D}[\og Lefschetz global{\fg}\ndemark]
\ndetext{cette conjecture est l\`a encore prouv\'ee dans \cite{HL}, ce m\^eme sous sa forme forte de la remarque qui suit, \cf th\'eor\`eme 3.4.1 de~\loccit}
Soient\pageoriginaled $X$ un espace analytique \emph{compact}, $Y$ un sous-espace analytique de $X$ tel que $U=X-Y$ soit Stein, et soit de profondeur homotopique rectifi\'ee $\geq n$ (par exemple intersection compl\`ete de dimension $\geq n$ en tout point). Alors l'homomorphisme canonique
$$
\pi_i(Y)\to\pi_i(X)$$
est un isomorphisme pour $i<n-1$, un \'epimorphisme pour $i=n-1$.
\end{enonce*}

\begin{remarquestar}
Lorsque, dans les \'enonc\'es \Ref{XIII.6.C} et \Ref{XIII.6.D}, on
remplace l'hypoth\`ese que $X-Y$ est Stein par l'hypoth\`ese que
$X-Y$ est r\'eunion de $c+1$ ouverts de Stein (qui jouera le
r\^{o}le d'une hypoth\`ese de \og concavit\'e\fg topologique), les
conclusions doivent \^etre modifi\'ees simplement en y rempla\c cant
$n$ par $n-c$.\nde{signalons enfin le r\'esultat suivant de Fulton,
\`a comparer avec le r\'esultat de Fulton-Hansen cit\'e note de
l'\'editeur~\eqref{FH}~page~\pageref{FH}: soient $X$ et $H$ sont des
sous-sch\'emas ferm\'es de $\PP^m_\CC$, $n$ la dimension de $X$ et $d$
la codimension de $H$. Alors, l'application
$$
\pi_i(X, X\cap H)\to\pi_i(\PP^n_\CC, H)$$
est un isomorphisme si $i\leq n-d$ et est surjective si $i=n-d-1$; voir (Fulton~W., {\og Connectivity and its applications in algebraic geometry\fg}, in \emph{Algebraic geometry (Chicago, Ill., 1980)}, Lect. Notes in Math., vol.~862, Springer, Berlin-New York, 1981, p\ptbl 26--92).}
\end{remarquestar}

Explicitons pour finir, dans le \og cas global\fg \Ref{XIII.6.D}, la conjecture concernant le groupe fondamental (obtenue en faisant $n=3$):

\refstepcounter{aconjecture}\label{XIII.6.D'}
\begin{enonce*}{Conjecture D$\boldsymbol'$}[Lefschetz global pour groupe fondamental\ndemark]
\ndetext{cette conjecture est d\'emontr\'ee dans \cite{HL}, \cf th\'eor\`eme 3.5.1 de~\loccit}
Soient $X$ un espace analytique compact sur le corps des complexes, $Y$ une partie analytique ferm\'ee telle que $U=X-Y$ soit Stein. Supposons de plus les conditions suivantes satisfaites:
\begin{enumeratei}
\item
Pour tout $x\in U$, le groupe fondamental local $\pi_2^x(X, x)$ est nul (\ie $X$ \emph{est \og pur en }$x${\fg}), ou seulement l'anneau local $\OXx$ est pur.
\item
Les anneaux locaux des points de $U$ sont \og connexes en dimension $\geq 2$\fg.
\item
Les anneaux locaux des points de $U$ sont de dimension $\geq 3$.
\end{enumeratei}

Sous ces conditions, pour tout $x\in Y$, l'homomorphisme
$$
\pi_1(Y, x)\to\pi_1(X, x)
$$
est\pageoriginale un isomorphisme (et $\pi_2(Y, x)\to\pi_2(X, x)$ un \'epimorphisme).
\end{enonce*}

On notera que les conditions locales (i)~(ii)~(iii) sur $U$ sont satisfaites si $U$ est localement intersection compl\`ete de dimension $\geq 3$. Du point de vue de la G\'eom\'etrie Alg\'ebrique, (lorsque $U$ provient d'un sch\'ema, encore not\'e $U$), les conditions (i) \`a (iii) correspondent \`a des hypoth\`eses sur les invariants locaux $\pi_i^x(U)$, \sisi{\ignorespaces}{\`a} savoir \hbox{$\pi_i^x(U)=0$} pour $i<3-\degtr k(x)/k$, pour les points $x$ tels que l'on ait respectivement $\degtr k(x)/k=0, 1, 2$. La condition globale sur $U$ ($U$ Stein) sera satisfaite si $X$ est projective et $Y$ une section hyperplane.

\chapterspace{-4}
\chapter[Profondeur et th\'eor\`emes de Lefschetz]
{Profondeur et th\'eor\`emes de Lefschetz en~cohomologie \'etale} \label{XIV}
\begin{center}
par Mme M. Raynaud\sfootnote{D'apr\`es des notes in\'edites de A.~Grothendieck.}
\end{center}

\vspace*{1cm}
Au\pageoriginale \numero \Ref{XIV.1}, nous d\'efinissons une notion de \og \textit{profondeur \'etale}\fg qui est l'analogue en cohomologie \'etale de la notion de profondeur \'etudi\'ee dans \Ref{III}, en cohomologie des faisceaux coh\'erents. Apr\`es une partie technique, nous d\'emontrons au \numero \Ref{XIV.4} des \og th\'eor\`emes de Lefschetz\fg, le th\'eor\`eme central \'etant \Ref{XIV.4.2}. Soit $X$ un sch\'ema, $Y$ une partie ferm\'ee de $X$, $U$ l'ouvert compl\'ementaire $X-Y$ et $F$ un faisceau ab\'elien sur $X$, pour la topologie \'etale; d'une mani\`ere g\'en\'erale le but des th\'eor\`emes de Lefschetz est de montrer que, si $F$ satisfait \`a certaines conditions locales sur $U$, exprimables en termes de profondeur \'etale au points de $U$, alors, sous certaines conditions suppl\'ementaires de nature globale sur $U$ (par exemple $U$ \textit{affine}), l'application naturelle des groupes de cohomologie \'etale
$$
\H^i(X, F)\to \H^i(Y, \sisi{F{|Y}}{F_{|Y}})$$
est un isomorphisme pour des valeurs $i<n$, o\`u $n$ est un certain entier explicit\'e. En prenant pour $F$ un faisceau constant, on obtient ainsi des conditions pour que $\pi_0(X)$ soit \'egal \`a $\pi_0(Y)$ et des conditions pour que les groupes fondamentaux rendus ab\'eliens de $X$ et $Y$ soient les m\^emes. Au \numero \Ref{XIV.5}, l'introduction d'une notion de \og profondeur g\'eom\'etrique\fg permet de donner des cas particuliers utiles des th\'eor\`emes de Lefschetz (\Ref{XIV.5.7}). Enfin au \numero \Ref{XIV.6}, nous signalons quelques conjectures, concernant notamment des variantes \og non commutatives\fg des th\'eor\`emes obtenus.

\section{Profondeur cohomologique et homotopique}\label{XIV.1}

\setcounter{subsection}{-1}

\subsection{}\label{XIV.1.0}
Fixons\pageoriginale les notations suivantes. Soient $X$ un sch\'ema\sfootnote{Conform\'ement \`a la nouvelle terminologie (\cf r\'e\'edition de \EGA I), nous appellerons ici \og sch\'ema\fg ce qui \'etait appel\'e pr\'ec\'edemment \og pr\'esch\'ema\fg et \og sch\'ema s\'epar\'e\fg ce qui \'etait appel\'e \og sch\'ema\fg.}, $Y$ une partie ferm\'ee de~$X$, $U$ l'ouvert compl\'ementaire et $i:Y=X-U\to X$ l'immersion canonique. Soient~$\underline{\Gamma}_Y$ le foncteur qui, \`a un faisceau ab\'elien sur $X$, fait correspondre le \og faisceau des sections \`a support dans $Y$\fg, c'est-\`a-dire $\underline{\Gamma}_Y=i_*i^!$ (\cf \SGA 4 IV 3.8 et VIII 6.6) et $\Gamma_Y$ le foncteur $\Gamma\cdot \underline{\Gamma}_Y$ (o\`u $\Gamma$ est le foncteur \og sections globales{\fg}). Consid\'erons la cat\'egorie d\'eriv\'ee $D^+(X)$ et le foncteur d\'eriv\'e $\R\underline{\Gamma}_Y$ (\resp $\R\Gamma_Y$) de $\underline{\Gamma}_Y$ (\resp de $\Gamma_Y$) (\cf \cite{XIV.3}). \'Etant donn\'e un complexe de faisceaux ab\'eliens $F$ sur $X$, \`a degr\'es born\'es inf\'erieurement, on peut le consid\'erer comme un \'el\'ement de $D^+(X)$; nous noterons alors ${\h}_Y^p(F)$ le $p$-i\`eme faisceau de cohomologie de $\R\underline{\Gamma}_Y(F)$ et $\H_Y^p(X, F)$ le $p$-i\`eme groupe de cohomologie de $\R\Gamma_Y(F)$. Les r\'esultats de (\SGA 4 V 4.3 et 4.4) s'\'etendent trivialement \`a ${\h}_Y^p(F)$ et $\H_Y^p(X, F)$.

\enlargethispage{\baselineskip}%
\begin{proposition} \label{XIV.1.1}
Soient $X$ un sch\'ema, $Y$ une partie ferm\'ee de $X$, $U$ l'ouvert compl\'ementaire et $i:U\to X$ l'immersion canonique. D\'esignons par $F$, soit un faisceau d'ensembles sur $X$, soit un faisceau en groupes sur $X$, soit un complexe de faisceaux ab\'eliens sur $X$, \`a degr\'es born\'es inf\'erieurement. Fixons les notations suivantes: si $X'\to X$ est un morphisme, $U'$ et $F'$ d\'esignent les images inverses de $U$ et $F$ sur $X'$; par ailleurs, si $\bar{y}$ est un point g\'eom\'etrique de $X$, $\bar{X}$ d\'esigne le localis\'e strict de $X$ en $\bar{y}$ et $\bar{U}$ et $\bar{F}$ les images inverses de $U$ et $F$ dans $\bar{X}$.

$1^\circ)$ Soient\pageoriginale $F$ un faisceau d'ensembles sur $X$ et $n$ un entier $\leq 2$; alors les conditions suivantes sont \'equivalentes:

\textup{(i)} Le morphisme canonique
$$F\to i_*i^*F$$
est injectif si $n\geq 1$, bijectif si $n\geq 2$.

\textup{(ii)} Pour tout sch\'ema $X'$ \'etale au-dessus de $X$, le morphisme canonique
$$
\H^0(X', F')\to \H^0(U', F')
$$
est injectif si $n\geq 1$, bijectif si $n\geq 2$.

Supposons de plus l'ouvert $U$ r\'etrocompact dans $X$; alors les conditions pr\'ec\'edentes sont \'equivalentes \`a la suivante:

\textup{(iii)} Pour tout point g\'eom\'etrique $\bar{y}$ de $Y$, le morphisme canonique
$$
\H^0(\bar{X}, \bar{F})\to \H^0(\bar{U}, \bar{F})$$
est injectif si $n\geq 1$, et bijectif si $n\geq 2$.

$2^\circ)$ Soient $F$ un faisceau en groupes sur $X$ et $n$ un entier $\leq 3$; alors les conditions suivantes sont \'equivalentes:

\textup{(i)} Le morphisme canonique
$$F\to i_*i^*F$$
est injectif si $n\geq 1$, bijectif si $n\geq 2$, et si $n\geq 3$, en plus des conditions pr\'ec\'edentes, le faisceau d'ensemble point\'es $\R^1i_*(i^*F)$ est nul.

\textup{(ii)} Pour tout sch\'ema $X'$ \'etale au-dessus de $X$, le morphisme canonique
$$
\H^0(X', F')\to \H^0(U', F')
$$
est\pageoriginale injectif si $n\geq 1$, bijectif si $n\geq 2$, et de plus le morphisme canonique
$$
\H^1(X', F')\to \H^1(U', F')
$$
est injectif si $n\geq 2$, bijectif si $n\geq 3$.

\textup{(ii bis)} Identique \`a \textup{(ii)}, sauf dans le cas $n=2$ o\`u l'on suppose seulement $\H^0(X', F')\to \H^0(U', F')$ bijectif.

Supposons de plus $U$ r\'etrocompact dans $X$; alors les conditions pr\'ec\'edentes sont aussi \'equivalentes \`a la suivante:

\textup{(iii)} Pour tout point g\'eom\'etrique $\bar{y}$ de $Y$, le morphisme canonique
$$
\H^0(\bar{X}, \bar{F})\to \H^0(\bar{U}, \bar{F})
$$
est injectif si $n\geq 1$, bijectif si $n\geq 2$; enfin si $n\geq 3$, en plus des conditions pr\'ec\'edentes, $\H^1(\bar{U}, \bar{F})$ est nul.

$3^\circ)$ Soit $F$ un complexe de faisceaux ab\'eliens, \`a degr\'es born\'es inf\'erieurement et $n$~un entier; alors les conditions suivantes sont \'equivalentes:

\textup{(i)} On a ${\h}_Y^p(F)=0$ pour $p<n$ (\cf 1.0).

\textup{(ii)} Pour tout sch\'ema $X'$ \'etale au-dessus de $X$, le morphisme canonique
$$
\H^p(X', F')\to \H^p(U', F')
$$
est bijectif pour $p<n-1$, injectif pour $p=n-1$.

Supposons $U$ r\'etrocompact dans $X$, alors les conditions qui pr\'ec\'edent sont aussi \'equivalentes \`a la suivante:

\textup{(iii)} Pour tout point g\'eom\'etrique $\bar{y}$ de $Y$, le morphisme canonique
$$
\H^p(\bar{X}, \bar{F})\to \H^p(\bar{U}, \bar{F})$$
est bijectif pour $p<n-1$, injectif pour $p=n-1$.

Dans\pageoriginale le cas o\`u $F$ est un faisceau ab\'elien et $n\geq 2$, les conditions \textup{(i)} et \textup{(ii)} sont aussi \'equivalentes \`a la suivante:

\textup{(ii bis)} Pour tout sch\'ema $X'$ \'etale au-dessus de $X$, le morphisme canonique
$$
\H^p(X', F')\to \H^p(U', F')$$
est bijectif pour $p<n-1$.
\end{proposition}

\skpt
\begin{proof}
$1^\circ)$ Il est clair que \textup{(i)}\SSI \textup{(ii)}.
Montrons que, si $U$ est r\'etrocompact dans $X$, \textup{(i)}\SSI
\textup{(iii)}. En effet, \textup{(i)} \'equivaut \`a dire que, pour
tout point g\'eom\'etrique $\bar{y}$ de $X$, le morphisme
$F_{\bar{y}}\to (i_*i^*F)_{\bar{y}}$ est injectif si $n\leq 1$ et
bijectif si $n\leq 2$ (\SGA 4 VIII 3.6). Comme ce morphisme est de
toutes fa\c cons bijectif lorsque $\bar{y}$ est un point
g\'eom\'etrique de $U$, on peut se borner aux points g\'eom\'etriques
$\bar{y}$ de $Y$. Or il r\'esulte du fait que $i$ est quasi-compact
et de (\SGA 4 VIII 5.3) que le morphisme
$$F_{\bar{y}}\to (i_*i^*F)_{\bar{y}}$$
s'identifie canoniquement au morphisme
$$
\H^0(\bar{X}, \bar{F})\to \H^0(\bar{U}, \bar{F}),
$$
d'o\`u l'\'equivalence de \textup{(i)} et \textup{(iii)}.

$2^\circ)$ \textup{(i)}\ALORS \textup{(ii)}. Les assertions sur le $\H^0$ r\'esultent de $1^\circ)$. Soit alors $i'$ l'immersion canonique de $U'$ dans $X'$; les assertions sur le $\H^1$ r\'esultent de la suite exacte (\SGA 4 XII 3.2)
$$0\to \H^1(X', i'_*i'^*F')\to \H^1(U', F')\to \H^0(X', \R^1i'_*(i'^*F')).
$$

(ii bis)\ALORS \textup{(i)}\pageoriginale. D'apr\`es $1^\circ)$, il suffit de montrer que, pour $n\geq 3$, on~a $\R^1i_*(i^*F)=0$. Or $\R^1i_*(i^*F)$ est le faisceau associ\'e au pr\'efaisceau $X'\mto \H^1(U', F')$, c'est-\`a-dire, par hypoth\`ese, le faisceau associ\'e au pr\'efaisceau $X'\mto \H^1(X', F')$, lequel est nul.

\textup{(i)}\SSI \textup{(iii)}. Tenant compte de $1^\circ)$, la seule chose qui reste \`a voir est que la relation $\R^1i_*(i^*F)=0$ \'equivaut au fait que $\H^1(\bar U, \bar F)=0$ pour tout point g\'eom\'etrique $\bar y$ de $Y$. Comme $\R^1i_*(i^*F))$ est nul hors de $Y$, il revient au m\^eme de dire que $\R^1i_*(i^*F)=0$ ou que $(\R^1i_*(i^*F))_{\bar y}=0$ pour tout point g\'eom\'etrique $\bar y$ de $Y$. Il suffit alors de noter que, $i$ \'etant quasi-compact, on~a $(\R^1i_*(i^*F))_{\bar y}=\H^1(\bar U, \bar F)$ (\SGA 4 VIII 5.3).

$3^\circ)$ \textup{(i)}\ALORS \textup{(ii)}. Soit $X'$ un sch\'ema \'etale au-dessus de $X$; on~a la suite exacte (\SGA 4 V 4.5)
\begin{equation*} \label{eq:XIV.1.1.*} \tag{$*$} {\to \H^p_{Y'}(X', F')\to H^p(X', F')\to H^p(U', F')\to};
\end{equation*}
donc \textup{(ii)} \'equivaut \`a $H^p_{Y'}(X', F')=0$ pour $p<n$ et pour tout sch\'ema $X'$ \'etale au-dessus de $X$. Consid\'erons alors la suite spectrale
$$
\E_2^{pq}=\H^p(X', {\h}_Y^q(F))\To \H_{Y'}^*(X', F');$$
par hypoth\`ese, ${\h}_Y^q(F)=0$ pour $q<n$, d'o\`u $\E_2^{pq}=0$ pour $p+q<n$ et par suite $\H_{Y'}^p(X', F')=0$ pour $p<n$.

\textup{(ii)}\ALORS \textup{(i)}. Le faisceau ${\h}_Y^p(F)$ est associ\'e au pr\'efaisceau $X'\mto \H_{Y'}^p(X', F')$; comme on~a d\'ej\`a remarqu\'e que \textup{(ii)} \'equivaut \`a la relation $\H_{Y'}^p(X', F')=0$ pour $p<n$ et pour tout sch\'ema $X'$ \'etale au-dessus de $X$, on~a bien ${\h}_Y^p(F)=0$ pour $p<n$.

\textup{(i)}\SSI \textup{(iii)}. Les faisceaux ${\h}_Y^p(F)$ sont concentr\'es sur $Y$; par suite il revient au m\^eme de dire que ${\h}_Y^p(F)=0$ ou de dire que, pour tout point\pageoriginale g\'eom\'etrique $\bar y$ de $Y$, la fibre $({\h}_Y^p(F))_{\bar y}$ est nulle. Or, $i$ \'etant quasi-compact, on d\'eduit de (\SGA 4 VIII 5.2) que l'on a $({\h}_Y^p(F))_{\bar y}={\h}_{\bar Y}^p(\bar X, \bar F)$. L'\'equivalence de \textup{(i)} et \textup{(iii)} en r\'esulte, compte tenu de l'analogue sur $\bar X$ de la suite exacte (\Ref{eq:XIV.1.1.*}).

(ii bis)\ALORS \textup{(ii)} dans le cas o\`u $F$ est un faisceau ab\'elien. La seule chose qui reste \`a montrer est que ${\h}_Y^{n-1}(F)=0$. Or, pour $n>2$, le faisceau ${\h}_Y^{n-1}(F)$ est associ\'e au pr\'efaisceau $X'\mto \H^{n-2}(U', F')=\H^{n-2}(X', F')$ donc est bien nul. Le cas $n=2$ r\'esulte du fait que ${\h}_Y^1(F)$ est le conoyau du morphisme $F\to i_*i^*F$.
\skipqed
\end{proof}

\begin{definition} \label{XIV.1.2}
Les notations sont celles de \Ref{XIV.1.1}. On dit que $F$ est de $Y$-profondeur \'etale $\geq n$ et on \'ecrit
$$
\prof_Y(F)\geq n$$
si $F$ satisfait aux conditions \'equivalentes \textup{(i)} et \textup{(ii)} de \Ref{XIV.1.1}. Si $F$ est un complexe de faisceaux ab\'eliens, on appelle $Y$-profondeur \'etale de $F$ la borne sup\'erieure des~$n$ pour lesquels $\prof_Y(F)\geq n$; on utilisera la m\^eme notation si $F$ est un faisceau d'ensembles \resp de groupes non n\'ecessairement commutatifs (de sorte qu'on a alors $0\leq \prof_Y(F)\leq 2$\sisi{}{,} \resp $0\leq \prof_Y(F)\leq 3$, lorsque le contexte ne permet pas de confusion sur celle des trois variantes envisag\'ees ici qu'on utilise).

Si $\LL$ est un ensemble de nombres premiers, on dit que \sisi{$X$ est de $Y$-profondeur \'etale pour $\LL \geq n$}{la $Y$-profondeur \'etale pour~$\LL$ de $X$ est $\geq n$} et on \'ecrit
$$
\prof_Y^{\LL}(X)\geq n
$$
si, pour tout faisceau constant de la forme $\ZZ/\ell\ZZ$ avec $\ell\in \LL$, on~a $\prof_Y(\ZZ/\ell\ZZ)\geq n$\pageoriginale. On d\'efinit de fa\c con \'evidente la $Y$-profondeur \'etale pour $\LL$ de $X$. Si $\LL=\PP$, ensemble de tous les nombres premiers, et s'il n'y a pas de risque de confusion avec la notation de \textup{(\EGA IV 5.7.1)} (relative au cas o\`u $F=\OX$), on omet $\LL$ dans la notation; sinon on \'ecrit $\prof \et(X)$.

Enfin on dit que $X$ est de $Y$-profondeur homotopique \sisi{\ignorespaces}{$\geq 3$} pour $\LL \sisi{\geq 3}{}$ et on \'ecrit
$$
\prof \hop_Y^{\LL}(X)\geq 3$$
si, pour tout faisceau de $\LL$-groupes constant fini $F$ sur $X$, on~a $\prof_Y(F)\geq 3$. Si $\LL=\PP$, on omet $\LL$ dans la notation.
\end{definition}

\begin{corollaire} \label{XIV.1.3}
Sous les conditions de \Ref{XIV.1.1}, si $\prof_Y(F)\geq n$, alors, pour tout ferm\'e~$Z$ de $Y$, on a
$$
\prof_Z(F)\geq n.
$$
\end{corollaire}

Faisons par exemple le raisonnement dans le cas o\`u $F$ est un complexe de faisceaux ab\'eliens \`a degr\'es born\'es inf\'erieurement. On utilise \Ref{XIV.1.1} $3^\circ)$ \textup{(ii)}. Soit $V=X-Z$ et consid\'erons, pour tout entier $p$, la suite de morphismes
$$
\H^p(X, F)\lto{f} \H^p(V, F)\lto{g} \H^p(U, F).
$$

Par hypoth\`ese $g$ et $f\circ g$ sont bijectifs pour $p<n-1$ et injectifs pour $p=n-1$; il en est donc de m\^eme de $f$. Comme le raisonnement est valable quand on remplace $X$ par un sch\'ema $X'$ \'etale au-dessus de $X$, ceci d\'emontre \Ref{XIV.1.3}.

\begin{corollaire} \label{XIV.1.4} Les notations sont celles de \Ref{XIV.1.1} $2^\circ)$. Si $X'$ est un sch\'ema sur\pageoriginale $X$, notons $\Phi'$ le foncteur qui fait correspondre \`a un rev\^etement \'etale de $X'$ sa restriction \`a~$U'$, et $\Phi'_{F'}$, le foncteur qui fait correspondre \`a un torseur\sfootnote{\ie un \og fibr\'e principal homog\`ene\fg dans une ancienne terminologie.} sous $F'$ sa restriction \`a~$U'$. Alors les conditions suivantes sont \'equivalentes:

\textup{(i)} On a $\prof_Y(F)\geq 1$ (\resp $\prof_Y(F)\geq 2$, \resp $\prof_Y(F)\geq 3$).

\textup{(ii)} Pour tout sch\'ema $X'$ \'etale au-dessus de $X$, le foncteur $\Phi'_{F'}$ est fid\`ele (\resp pleinement fid\`ele, \resp une \'equivalence de cat\'egories).

En particulier, pour que l'on ait $\prof_Y({X})\geq 1$ (resp.\ $\prof_Y({X})\geq 2$, resp.\ $\prof \hop_Y(X)\geq 3$), il faut et il suffit que le foncteur $\Phi'$ soit fid\`ele (resp.\ pleinement fid\`ele, \resp une \'equivalence de cat\'egories).
\end{corollaire}

Cela r\'esulte en effet de \Ref{XIV.1.1} $2^\circ)$ \textup{(ii)}, compte tenu de l'interpr\'etation de $\H^1(X', F')$ comme l'ensemble des classes (mod. isomorphismes) de torseurs sous $F'$ (\SGA 4 VII~2), et des rev\^etements \'etales $Z$ de degr\'e $n$ d'un sch\'ema comme associ\'es \`a des rev\^etements principaux galoisiens de groupe le groupe sym\'etrique $\mathfrak{S}_n$, \`a $Z$ \'etant associ\'e le rev\^etement $\underline{\Isom}_X(Z_0, Z)$, o\`u $Z_0$ est le rev\^etement trivial de $X$ de degr\'e $n$.

\begin{corollaire} \label{XIV.1.5} Sous les conditions de \Ref{XIV.1.1} $3^\circ)$, supposons que l'on ait $\prof_Y(F)\geq n$; alors on a
$$
\H_Y^n(X, F)\simeq \H^0(X, {\h}_Y^n(F)).
$$
\end{corollaire}

Le corollaire r\'esulte de la suite spectrale
$$
\E_2^{pq}=H^p(X, {\h}_Y^q(F))\To H_Y^*(X, F).
$$

En effet on~a par hypoth\`ese $\E_2^{pq}=0$ pour $q<n$, d'o\`u il r\'esulte que
$$
\H_Y^n(X, F)=\E_2^{0n}=\H^0(X, {\h}_Y^n(F)).
$$

\skpt
\begin{remarques} \label{XIV.1.6}
a) La notion\pageoriginale de $Y$-profondeur, sous la forme des conditions \'equivalentes \textup{(i)} et \textup{(ii)} de \Ref{XIV.1.1}, a un sens pour n'importe quel site. Dans le cas particulier o\`u $X$ est un sch\'ema localement noeth\'erien, muni de la topologie de Zariski, et $F$ un faisceau de $\OX$-modules coh\'erents, on trouve la notion usuelle de $Y$-profondeur comme borne inf\'erieure des profondeurs aux points de $Y$ (\Ref{III}).

b) Pour $n\leq 2$, la notion de $Y$-profondeur \'etale de $X$ est ind\'ependante de $\LL$. Pour $n=1$, elle signifie simplement que $U$ est dense dans $X$. En effet cette condition est n\'ecessaire pour que l'on ait $\prof_Y(F)\geq 1$, et elle est aussi suffisante car on peut supposer $X$ r\'eduit, cas dans lequel la condition $U$ dense dans $X$ se conserve par changement de base \'etale \textup{(\EGA IV 11.10.5) \textup{(ii)} b))}. Si $U$ est r\'etrocompact dans $X$, la relation $\prof_Y(X)\geq 1$ \'equivaut aussi \`a dire $Y$ ne contient aucun point maximal de~$X$ (\EGA I 6.6.5). Pour $n=2$ et $U$ r\'etrocompact dans $X$, la condition $\prof_Y(X)\geq 2$ \'equivaut au fait que, pour tout point g\'eom\'etrique $\bar y$ de $Y$, $\bar U$ est connexe non vide, c'est-\`a-dire que \og $Y$ ne disconnecte pas $X$, localement pour la topologie \'etale\fg.

c) Si $X$ est de $Y$-profondeur $\geq n$ pour $\LL$ et $U$ r\'etrocompact dans $X$, alors pour tout faisceau ab\'elien localement constant de $\LL$-torsion $F$ sur $X$, on~a $\prof_Y(X)\geq n$. En effet, la propri\'et\'e $\prof_Y(F)\geq n$ \'etant locale pour la topologie \'etale, on peut supposer $F$ constant; alors $F$ est limite inductive filtrante de faisceaux sommes finies de faisceaux de la forme $\ZZ/p^m\ZZ$, o\`u $m$ est un entier $>0$ et $p\in \LL$. En utilisant \Ref{XIV.1.1} \textup{(iii)} et \textup{(\SGA 4 VII 3.3)}, on voit que l'on peut se ramener au cas o\`u $F=\ZZ/p^m\ZZ$, puis\pageoriginale, par r\'ecurrence sur $m$, au cas o\`u $F=\ZZ/p\ZZ$ pour lequel l'assertion r\'esulte de la d\'efinition.

d) D'apr\`es \Ref{XIV.1.4}, si l'on a $\prof_Y(X)\geq 3$, le couple $(X, Y)$ est \textit{pur} au sens de \Ref{X}~\Ref{X.3.1}. En fait les couples purs que l'on rencontre dans la pratique (\cf \Ref{X}~\Ref{X.3.4}) satisfont \`a la condition plus forte de profondeur homotopique $\geq 3$, et cette notion peut donc se substituer avantageusement \`a celle de couple pur.

e) Soient $F$ un complexe de faisceaux ab\'eliens et $T(F)$ le complexe obtenu en appliquant \`a $F$ le foncteur translation (\cite{XIV.3}); alors on~a \'evidemment:
$$
\prof_Y(T(F))=\prof_Y(F)-1.
$$

f) Signalons que les travaux r\'ecents d'Artin-Mazur (\cite{XIV.1}) permettent de d\'efinir la notion de profondeur homotopique $\geq n$, quel que soit l'entier $n$ (pas seulement si $n\geq 3$).

g) Sous les conditions de \Ref{XIV.1.1} 3), pour que l'on ait $\prof_Y(F)=\infty$, il faut et il suffit que l'on ait $F\isomto Ri_*(i^*F)$ dans $D^+(X)$. En effet, les ${\h}_Y^p(F)$ sont les faisceaux de cohomologie du c\^{o}ne (= mapping cylinder) du morphisme canonique $F\to Ri_*i^*(F)$.
\end{remarques}

\begin{definition} \label{XIV.1.7}
Soient $X$ un sch\'ema, $x$ un point de $X$, $\bar x$ un point g\'eom\'etrique au-dessus de $x$ et $\bar X$ le localis\'e strict de $X$ en $\bar x$. Comme pr\'ec\'edemment $F$ d\'esigne, soit un faisceau d'ensembles sur $X$, soit un faisceau en groupes sur $X$, soit un complexe de faisceaux ab\'eliens sur $X$ \`a degr\'es born\'es inf\'erieurement, $\bar F$ son image inverse sur~$\bar X$ et $\LL$ un ensemble de nombres premiers. On dit que $F$ est de \emph{profondeur \'etale $\geq n$ au point $x$}\pageoriginale (\resp que la \emph{profondeur \'etale pour $\LL$ de $X$ en $x$ est $\geq n$}, \resp que la \emph{profondeur homotopique pour $\LL$ de $X$ en $x$ est $\geq 3$}) et l'on \'ecrit $\prof_x(F)\geq n$ (\resp $\prof_x^{\LL}(X)\geq n$, \resp $\prof \hop_x^{\LL}(X)\geq \sisi{3}{3}$) si l'on a $\prof_{\bar x}(\bar F)>n$ (\resp $\prof_{\bar x}^{\LL}(\bar X)>\nobreak n$, \resp $\prof \hop_{\bar x}^{\LL}(\bar X)\geq \sisi{3}{3}$). On d\'efinit de fa\c con \'evidente l'entier $\prof_x^{\LL}(X)$ et, si $F$ est un complexe de faisceaux ab\'eliens, l'entier $\prof_x(F)$. Si $\LL$ est l'ensemble de tous les nombres premiers, on omet $\LL$ dans la notation $\prof_x^{\LL}(X)$, sauf s'il y a un risque de confusion avec la notation de \textup{(\EGA IV 5.7.1)}, auquel cas on \'ecrit $\prof \et_x(X)$.
\end{definition}

On a alors la \textit{caract\'erisation ponctuelle} suivante de la profondeur:

\begin{theoreme} \label{XIV.1.8} Soient $X$ un sch\'ema, $Y$ une partie ferm\'ee de $X$ telle que l'ouvert $U=X-Y$ soit r\'etrocompact dans $X$. Si $F$ est, soit un faisceau d'ensembles sur $X$, soit un faisceau en groupes sur $X$, soit un complexe de faisceaux ab\'eliens sur $X$ \`a degr\'es born\'es inf\'erieurement, alors on a
$$
\prof_Y(F)=\inf_{y\in Y} \prof_y(F).
$$
\end{theoreme}

\subsubsection{} \label{XIV.1.8.1} Montrons d'abord que, pour tout point $y$ de $Y$, on~a l'in\'egalit\'e $\prof_y(F)\geq \prof_Y(F)$. Soient en effet $\bar y$ un point g\'eom\'etrique au-dessus de $y$, $\bar X$ le localis\'e strict de $X$ en $\bar y$, $\bar F$ et $\bar Y$ les images inverses de $F$ et $Y$ sur $\bar X$. On a d'apr\`es \Ref{XIV.1.7} et \Ref{XIV.1.3}
$$
\prof_y(F)=\prof_{\bar y}(\bar F)\geq \prof_{\bar Y}(\bar F)\geq \prof_Y(F),
$$
la derni\`ere in\'egalit\'e utilisant l'hypoth\`ese \og $U$ r\'etrocompact\fg, via les conditions \textup{(iii)} dans \Ref{XIV.1.1} et la transitivit\'e dans la formation des localis\'es stricts.

\subsubsection{} \label{XIV.1.8.2}
Inversement supposons que, pour tout point $y$ de $Y$, on ait $\prof_y(\sisi{Y}{F})\geq n$ ($n$~entier) et montrons que l'on a alors $\prof_Y(F)\geq n$.

Rappelons\pageoriginale d'abord les r\'esultats bien connus suivant \textup{(\SGA 4 VIII)}:

\begin{subsublemme} \label{XIV.1.8.2.1}
Soient $X$ un sch\'ema, $F$ un faisceau d'ensembles sur $X$ (\resp $G\to F$ un monomorphisme de faisceaux d'ensembles sur $X$). Alors, pour que deux sections $s$ et $s'$ de $F$ co\"incident (\resp pour qu'une section $s$ de~$F$ provienne d'une section de~$G$), il faut et il suffit qu'il en soit ainsi localement. En particulier, si~$s$ et~$s'$ sont deux sections de $F$, il existe un plus grand ouvert $V$ de $X$ sur lequel elles co\"incident (\resp si $s$ est une section de $F$ sur $X$, il existe un plus grand ouvert $V$ de $X$ tel que $s|V$ provienne d'une section de $G$ sur $V$). Cet ouvert est aussi l'ensemble des points~$x$ de $X$ tels que, d\'esignant par $\bar x$ un point g\'eom\'etrique au-dessus de~$x$, les sections~$s$ et~$s'$ aient m\^eme image dans la fibre $F_{\bar x}$ (\resp que l'image de $s$ dans $F_{\bar x}$ provienne d'un \'el\'ement de $G_{\bar x}$).
\end{subsublemme}

Revenons \`a la d\'emonstration de \Ref{XIV.1.8}.

$1^\circ)$ Cas o\`u $F$ est un faisceau d'ensembles. Si $n=1$, il suffit de montrer que le morphisme canonique
$$
\H^0(X, F)\to \H^0(U, F)$$
est injectif, le r\'esultat s'appliquant encore quand on remplace $X$ par un sch\'ema \'etale au-dessus de $X$. Soient $s$ et $s'$ deux sections de $F$ au-dessus de $X$ qui ont m\^eme image dans $\H^0(U, F)$ et soit $V$ le plus grand ouvert au-dessus duquel elles sont \'egales; on~a \'evidemment $V\supset U$. Supposons $V\neq X$ et soit $y$ un point maximal de $X-V$, $\bar y$ un point g\'eom\'etrique au-dessus\pageoriginale de~$y$, $\bar X$ le localis\'e strict de $X$ en $\bar y$ et $\bar V$ et $\bar F$ les images inverses de $V$ et $F$ sur $\bar X$. D'apr\`es le choix de $y$, on~a $\bar X-\bar y=\bar V$ et par hypoth\`ese, le morphisme
$$
\H^0(\bar X, \bar F)\to \H^0(\bar X-\bar y, \bar F)=\H^0(\bar V, \bar F)$$
est injectif. Il en r\'esulte que $s$ et $s'$ co\"incident au point $y$, ce qui est absurde. Si $n=2$, il suffit de montrer, compte tenu de ce qui pr\'ec\`ede, que le morphisme
$$
\H^0(X, F)\to \H^0(U, F)=\H^0(X, i_*i^*F)$$
est surjectif (o\`u $i$ est l'immersion canonique de $U$ dans $X$). Soient $s$ une section de $i_*i^*F$ au-dessus de $X$ et $V$ le plus grand ouvert au-dessus duquel elle provient d'une section de $F$. Supposons $V\neq X$ et soit $y$ un point maximal de $X-V$; avec les notations pr\'ec\'edentes, il r\'esulte de l'hypoth\`ese que le morphisme canonique
$$
\H^0(\bar X, \bar F)\to \H^0(\bar X-\bar y, \bar F)=\H^0(\bar V, \bar F)
$$
est bijectif; par suite $s|\bar V$ se prolonge au point $\bar y$, ce qui est absurde et ach\`eve la d\'emonstration dans le cas $1^\circ)$.

$2^\circ)$ Cas o\`u $F$ est un faisceau en groupes. Compte tenu de $1^\circ)$, la seule chose qui reste \`a montrer est que, dans le cas $n=3$, le morphisme
$$
\H^1(X, F)=\H^1(X, i_*i^*F)\to \H^1(U, F)$$
est bijectif. On sait d\'ej\`a qu'il est injectif gr\^ace \`a $1^\circ)$ et \`a \Ref{XIV.1.1} $2^\circ)$ (ii bis). Pour la surjectivit\'e, on utilise la suite exacte (\SGA 4 XII 3.2)
$$
0\to \H^1(X, i_*i^*F)\to \H^1(U, F)\lto{d}\H^0(X, \R^1i_*(i^*F)).
$$
Soient $s\in \H^1(U, F)$ et $V\supset U$ le plus grand ouvert au-dessus duquel $d(s)=0$; c'est\pageoriginale aussi le plus grand ouvert tel que $s$ provienne d'un \'el\'ement de $\H^1(V, F)$. Supposons $V\neq X$ et soit $y$ un point maximal $X-V$; si $\bar X$ est le localis\'e strict de $X$ en un point g\'eom\'etrique $\bar y$ au-dessus de~$y$, on a, avec des notations \'evidentes, la suite exacte
$$
0\to \H^1(\bar X, \bar i_*(\bar i^*\bar F))\to \H^1(\bar U, \bar F)\lto{\bar d}\H^0(\bar X, \R^1\bar i_*(\bar i^*\bar F)).
$$
Comme $i:U\to X$ est un quasi-compact, $\R^1\bar i_*(\bar i^*\bar F)$ est l'image inverse de $\R^1i_*(i^*F)$ par le morphisme $\bar X\to X$, d'o\`u $\H^0(\bar X, \R^1\bar i_*(\bar i^*\bar F))=(\R^1i_*(i^*F))_{\bar y}$. Par hypoth\`ese et vu que $y\in Y$, le morphisme
$$
\H^1(\bar X, \bar F)\to \H^1(\bar V, \bar F)
$$
est bijectif. L'image $\bar s$ de $s$ dans $\H^1(\bar U, \bar F)$, qui se prolonge \`a $\bar V$ par d\'efinition de $V$, se prolonge donc aussi \`a $\bar X$; il en r\'esulte que $\bar d(\bar s)=0$ donc l'image de $d(s)$ dans la fibre g\'eom\'etrique $(\R^1i_*(i^*F))_{\bar y}$ est nulle; mais ceci contredit la d\'efinition de $V$, d'o\`u le cas~$2^\circ)$.

$3^\circ)$ Cas o\`u $F$ est un complexe de faisceaux ab\'eliens, \`a degr\'es born\'es inf\'erieurement. On raisonne par r\'ecurrence sur $n$. La conclusion est satisfaite pour $n$ assez petit, puisque $F$ est \`a degr\'es born\'es inf\'erieurement. Supposons donc que $\prof_Y(F)\geq n-1$ et montrons que $\prof_Y(F)\geq n$, sachant que, pour tout point $y$ de $Y$, on~a $\prof_y(F)\geq n$. Il suffit de voir que le morphisme canonique
\begin{equation*} \label{eq:XIV.1.8.*} \tag{$*$} {\H^{n-2}(X, F)\to \H^{n-2}(U, F)}
\end{equation*}
est surjectif et que
\begin{equation*} \label{eq:XIV.1.8.**} \tag{$**$} {\H^{n-1}(X, F)\to \H^{n-1}(U, F)}
\end{equation*}
est injectif (le r\'esultat s'appliquant quand on remplace $X$ par un sch\'ema \'etale au-dessus\pageoriginale de $X$).

a) Surjectivit\'e de (\Ref{eq:XIV.1.8.*}). La d\'emonstration est analogue \`a celle de $2^\circ)$. Compte tenu de \Ref{XIV.1.5} et de (\SGA 4 V 4.5), on~a la suite exacte
$$
\H^{n-2}(X, F)\to \H^{n-2}(U, F)\lto{d}\H^{n-1}_Y(X, F)=\H^0(X, {\h}_Y^{n-1}(F)).
$$
Soit $s\in \H^{n-2}(U, F)$ et $V\supset U$ le plus grand ouvert au-dessus duquel $d(s)=0$, lequel est aussi le plus grand ouvert tel que $s$ se prolonge \`a $\H^{n-2}(V, F)$. Supposons $V\neq X$ et soient $y$ un point maximal de $X-V$ et $\bar X$ le localis\'e strict de $X$ en un point g\'eom\'etrique $\bar y$ au-dessus de~$y$. Comme $i:U\to X$ est quasi-compact, la formation de ${\h}_Y^{n-1}(F)$ commute au changement de base $\bar X\to X$ et l'on a donc (avec des notations \'evidentes) la suite exacte
$$
\H^{n-2}(\bar X, \bar F)\to \H^{n-2}(\bar U, \bar F)\lto{\bar d}\H^{n-1}_{\bar Y}(\bar X, \bar F)=({\h}_Y^{n-1}(F))_{\bar y},
$$
la derni\`ere \'egalit\'e r\'esultant de l'hypoth\`ese de r\'etrocompacit\'e sur $U$.

Or on~a par hypoth\`ese l'isomorphisme
$$
\H^{n-2}(\bar X, \bar F)\isomto \H^{n-2}(\bar X-\bar y, \bar F)=\H^{n-2}(\bar V, \bar F);
$$
par suite l'image $\bar s$ de $s$ dans $\H^{n-2}(\bar U, \bar F)$, qui se prolonge (par d\'efinition de $V$) \`a $\H^{n-2}(\bar V, \bar F)$, se prolonge aussi \`a $\H^{n-2}(\bar X, \bar F)$; mais ceci montre que $\bar d(\bar s)=0$, c'est-\`a-dire que $d(s)$ est nulle en $\bar y$, ce qui est absurde.

b) Injectivit\'e de (\Ref{eq:XIV.1.8.**}). En utilisant la surjectivit\'e de (\Ref{eq:XIV.1.8.*}), on obtient la suite exacte
$$0\to \H^0(X, {\h}_Y^{n-1}(F))\to \H^{n-1}(X, F)\to \H^{n-1}(U, F)$$
et il faut montrer que tout \'el\'ement $s\in \H^0(X, {\h}_Y^{n-1}(F))$ est nul. Soit $V$\pageoriginale le plus grand ouvert au-dessus duquel $s=0$. Supposons que l'on ait $V\neq X$ et soient $y$ un point maximal de $X-V$, $\bar X$ un localis\'e strict de $X$ en un point g\'eom\'etrique $\bar y$ au-dessus de~$y$. On a, par hypoth\`ese de r\'ecurrence et par \Ref{XIV.1.8.1}, la relation $\prof_{\bar Y}(\bar F)\geq \prof_Y(F)\geq n-1$, d'o\`u le fait que l'application $e$ du diagramme qui suit est injective:
$$
\H^0(\bar X, {\h}_Y^{n-1}(\bar F))= ({\h}_Y^{n-1}(F))_{\bar y}\lto{e} \H^{n-1}(\bar X, \bar F) \lto{f} \H^{n-1}(\bar V, \bar F).
$$
Il en est de m\^eme de $f$ en vertu de l'hypoth\`ese; l'\'egalit\'e de gauche r\'esulte de l'hypoth\`ese de r\'etrocompacit\'e sur $U$. Soit $\bar s$ l'image de $s$ dans $({\h}_Y^{n-1}(F))_{\bar y}$; comme $s$ s'annule au-dessus de $V$, on~a $f\cdot e(\bar s)=0$, d'o\`u $\bar s=0$, ce qui contredit le choix de $y$ et ach\`eve la d\'emonstration.

\begin{remarque} \label{XIV.1.9}
Un r\'esultat analogue \`a \Ref{XIV.1.8} est sans doute valable dans le cas o\`u l'on remplace le topos \'etale d'un sch\'ema $X$ par un \og topos localement de type fini \fg, c'est-\`a-dire d\'efinissable par un site localement de type fini (\SGA 4 VI 1.1). Pour le voir, il faut utiliser un r\'esultat de P. Deligne (\SGA \sisi{4 XVII}{4 VI.9}), affirmant qu'il y a \og suffisamment de foncteurs fibres\fg dans un tel topos.
\end{remarque}

Nous allons d\'eduire de \Ref{XIV.1.8} des cas importants o\`u l'on peut d\'eterminer la profondeur \'etale.

\refstepcounter{subsection}\label{XIV.1.10}
\begin{enonce*}{Th\'eor\`eme 1.10}[Th\'eor\`eme de semi-puret\'e cohomologique\ndemark]\ndetext{\label{gabpur} Gabber a prouv\'e depuis --- en 1994 --- la conjecture de puret\'e cohomologique absolue de Grothendieck: si $Y$ est un sous-sch\'ema ferm\'e de sch\'emas noeth\'eriens absolus de codimension pure $c$ et $n$ un entier inversible sur $X$, alors $\SheafH_Y^q(\Lambda)$ est nul si $q\neq 2c$ et vaut $\Lambda_Y(-c)$ (twist de Tate) sinon, o\`u on~a pos\'e $\Lambda=\ZZ/n\ZZ$. Voir (Fujiwara~K., {\og A Proof of the Absolute Purity Conjecture (after Gabber)\fg}, in \emph{Algebraic geometry 2000, Azumino (Hotaka)}, Adv. Stud. in Pure Math., vol.~36, 2002, p\ptbl 153-183). Pour des applications \`a l'existence du complexe dualisant, voir (\SGA 5, Lect. Notes in Math., vol.~589, Springer-Verlag, 1977, p\ptbl 1672), expos\'e 1 et \loccit, \S8. Cette conjecture avait \'et\'e prouv\'ee dans le cas $n=\ell^\nu$ avec $\ell$ premier inversible sur $X$ assez grand en utilisant de fa\c con cruciale la K-th\'eorie (Thomason~R.W., {\og Absolute cohomological purity\fg}, \emph{Bull. Soc. Math. France} \textbf{112} (1984), \numero 3, p\ptbl 397--406). On retrouve la $K$-th\'eorie dans le preuve de Gabber par le biais de la suite spectrale de Atiyah-Hirzebruch-Thomason reliant cohomologie \'etale et $K$-th\'eorie, m\'ethode d\'ej\`a utilis\'ee dans l'approche de Thomason. Outre ce r\'esultat, l'autre argument fondamental est la g\'en\'eralisation du th\'eor\`eme de Lefschetz cit\'e note~\eqref{gabaf}, page~\pageref{gabaf}.}
D\'esignons par~$X$, soit un sch\'ema lisse sur un corps $k$, soit un sch\'ema r\'egulier excellent \textup{(\EGA IV 7.8.2)} de caract\'eristique nulle (N.B. si l'on admet la r\'esolution des singularit\'es au sens de \textup{(\SGA 4 XIX)}, il suffit de supposer, plus g\'en\'eralement, que $X$ est un sch\'ema r\'egulier excellent d'\'egale caract\'eristique).\pageoriginaled Soient $Y$ une partie ferm\'ee de $X$ et $\LL$ l'ensemble des nombres premiers distincts de la caract\'eristique de $X$. Alors on a
$$
\prof_Y^{\LL}(X)=2\codim(Y, X).
$$
\end{enonce*}

\begin{proof}
Il r\'esulte de \Ref{XIV.1.8} que l'on a
$$
\prof_Y^{\LL}(X)=\inf_{y\in Y} \prof_y(X).
$$

Comme d'autre part $\codim(Y, X)=\inf_{y\in Y} \dim \Oo_{X, y}$, on est ramen\'e \`a montrer que
$$
\prof_y^{\LL}(X)=2 \dim \Oo_{X, y},
$$
ce qui r\'esulte de (\SGA 4 XVI 3.7 et XIX 3.2).
\skipqed
\end{proof}

\begin{theoreme}[Th\'eor\`eme de puret\'e homotopique\ndemark] \label{XIV.1.11}
Si $X$ est un sch\'ema localement noeth\'erien qui est r\'egulier (\resp dont les anneaux locaux sont des intersections compl\`etes), $Y$ une partie ferm\'ee de $X$ telle que $\codim(Y, X)\geq 2$ (resp.\ \hbox{$\codim(Y, X)\geq3$}), alors on a\ndetext{\label{jongoort} r\'ecemment, de Jong et Oort ont obtenu l'\'enonc\'e de puret\'e suivant: soit $\tilde S\to S$ une r\'esolution des singularit\'es du spectre $S$ d'un anneau local noeth\'erien normal de dimension $2$ et soit~$U$ le compl\'ementaire du point ferm\'e $s$ dans $S$. Supposons de plus que $k(s)$ soit alg\'ebriquement clos. Alors, pour tout nombre premier $p$, en particulier si $S$ est de caract\'eristique $p$, le morphisme de restriction $H^1_\et(\tilde S,\QQ_p)\to H^1_\et(U,\QQ_p)$ est bijectif (de Jong A.J. \& Oort~F, {\og Purity of the stratification by Newton polygons\fg}, \emph{J.~Amer. Math. Soc.} \textbf{13} (2000), \numero 1, p\ptbl 209--241, th\'eor\`eme 3.2). Si $k=\CC$ et $A$ est le compl\'et\'e d'une singularit\'e de surface, ce r\'esultat est d\^u \`a Mumford (voir page~\pageref{XIII.7}, \cite{XIII.5}).}
$$
\prof \hop_Y(X)\geq 3.
$$
\end{theoreme}

Il r\'esulte en effet de \Ref{XIV.1.8} que l'on a $\prof \hop_Y(X)=\inf_{y\in Y} \prof \hop_y(X)$. Or les anneaux strictement locaux de $X$ aux diff\'erents points de $Y$ sont des anneaux r\'eguliers de dimension $\geq 2$ (\resp d'intersection compl\`ete et de dimension $\geq 3$). Il r\'esulte alors du th\'eor\`eme de puret\'e \Ref{X}~\Ref{X.3.4} que $\prof \hop_y(X)\geq 3$, ce qui d\'emontre le th\'eor\`eme.

\begin{exemple} \label{XIV.1.12}
Soient $X$ un sch\'ema localement noeth\'erien, $Y$ une partie ferm\'ee de~$X$ et $n=1$ ou $2$. Alors, si l'on a $\prof_Y(\OX)\geq n$ ($\prof_Y(\OX)$ d\'esignant la $Y$-profondeur au sens des faisceaux coh\'erents (\cf 1.6 a)), on~a aussi $\prof_Y(X)\geq n$; c'est \'evident pour $n=1$ et, pour $n=2$, cela n'est\pageoriginale autre que le th\'eor\`eme de Hartshorne (\Ref{III}~\Ref{III.1}). Par contre l'assertion analogue est fausse pour $n\geq 3$. Prenons par exemple un espace affine de dimension $\geq 3$ sur un corps de caract\'eristique $\neq 2$ et faisons op\'erer le groupe $\ZZ/2\ZZ$ par sym\'etrie par rapport \`a l'origine. Soient $X$ le quotient et $Y=\{x\}$ l'image de l'origine dans $X$. Alors $\Oo_{X, x}$ est un anneau de Cohen-Macaulay, donc on~a $\prof_x(\OX)\geq 3$; mais l'espace affine priv\'e de l'origine est un rev\^etement \'etale de $X-\{x\}$ qui ne se prolonge pas en un rev\^etement \'etale de $X$; donc on~a d'apr\`es \Ref{XIV.1.4} $\prof_Y(X)=2$.
\end{exemple}

Le th\'eor\`eme suivant est l'analogue de (\EGA IV 6.3.1):

\begin{theoreme} \label{XIV.1.13}
Soient $f:X\to S$ un morphisme de sch\'emas, $Y$ une partie ferm\'ee de $X$, $Z$ une partie ferm\'ee de $S$, telles que $f(Y)\subset Z$. On suppose que les anneaux locaux de $X$ au diff\'erents points de $Y$ sont noeth\'eriens et que les ouverts $X-Y$ et $S-Z$ sont r\'etrocompacts dans $X$ et $S$ respectivement. Soient $p$, $q$, $r$ des entiers tels que $p\geq -r$, $q\geq 0$, $\LL$ un ensemble de nombres premiers et $F$ un complexe de faisceaux ab\'eliens de $\LL$-torsion sur $S$, tel que les faisceaux de cohomologie $\H^i(F)$ soient nuls pour $i<-r$. On suppose que
\begin{enumeratea}
\item
Le morphisme $f$ est localement ($p+q+r-2$)-acyclique pour $\LL$ \textup{(\SGA 4 XV 1.11)}.
\item
On a
$$
\prof_Z(F)\geq p.
$$
\item
Pour tout point $s$ de $Z$, on a
$$
\prof_{Y_s}^{\LL}(X_s)\geq q.
$$
\end{enumeratea}
\noindent
Alors on a\pageoriginale
$$
\prof_Y(f^*F)\geq p+q.
$$
\end{theoreme}

Nous aurons besoin du lemme suivant:

\begin{sublemme} \label{XIV.1.13.1}
Soient $\LL$ un ensemble de nombres premiers, $n$ et $r$ des entiers, $f:X\to S$ un morphisme localement $n$-acyclique pour $\LL$. Soient $F$ un complexe de faisceaux ab\'eliens, \`a faisceaux de cohomologie de $\LL$-torsion, tel que $\H^i(F)=0$ pour $i<-r$, $Z$ une partie ferm\'ee de $S$ telle que $S-Z$ soit r\'etrocompact dans $S$ et $T=f^{-1}(Z)$. Alors le morphisme canonique
$$f^*({\h}_Z^i(F))\to {\h}_T^i(f^*F)$$
est bijectif pour $i<n-r+2$ et injectif pour $i=n-r+2$.
\end{sublemme}

Posons $U=S-Z$ et $V=X-T$, de sorte que l'on a le carr\'e cart\'esien
$$
\xymatrix{
V \ar[r]^-g\ar[d]_k& U\ar[d]^j\\
X \ar[r]^-f&S
}$$
Consid\'erons le diagramme commutatif suivant dont les lignes sont exactes
$$
\xymatrix{ \ar[r]&f^*({\h}_Z^i(F))\ar[r]\ar[d]& f^*({\h}^i(F))\ar[r]\ar[d]^{\wr} &f^*({\h}^i(\R j_*(j^*F)))\ar[r]\ar[d]& \\
\ar[r] &{\h}_T^i(f^*F)\ar[r]& {\h}^i(f^*F) \ar[r]&{\h}^i(\R k_*(k^*f^*F)) \ar[r]&\text{;}}$$
il en r\'esulte que l'on est ramen\'e \`a montrer que le morphisme
$$f^*({\h}^i(\R j_*(j^*F)))\to {\h}^i(\R k_*(k^*f^*F))$$
est\pageoriginale bijectif pour $i<n-r+1$ et injectif pour $i=n-r+1$. Or un tel morphisme provient du morphisme suivant entre suites spectrales d'hypercohomologie
$$
\xymatrix{f^*E_2^{p, q}=f^*(\R^pj_*({\h}^q(j^*F)))\ar@{=>} [r]\ar[d]&f^*({\h}^*(\R j_*(j^*F)))\ar[d]\\
E'^{p, q}_2=\R^pk_*({\h}^q(k^*f^*F))\ar@{=>}[r]&{\h}^*(\R k_*(k^*f^*F)).}
$$
Comme $j$ est quasi-compact, il r\'esulte de (\SGA 4 XV 1.10) que le morphisme $f^*(E^{p, q}_2)\to E'^{p, q}_2$ est bijectif pour $p\leq n$ et injectif pour $p=n+1$; en particulier il est bijectif pour $p+q\leq n-r$ et injectif pour $p+q=n-r+1$. La conclusion en r\'esulte aussit\^{o}t.

Revenons \`a la d\'emonstration de \Ref{XIV.1.13}. Soit $T=f^{-1}(Z)$. D'apr\`es \Ref{XIV.1.13.1} et la condition a), le morphisme canonique $f^*(\h_Z^i(F))\to\h_T^i(f^*F)$ est un isomorphisme pour $i\leq p+q$. Il r\'esulte donc de b) que $\h^i_T(f^*F)=0$ pour $i<p$ et, pour $i<p+q, \ \h^i_T(f^*F)$ restreint \`a $T$ est l'image inverse d'un faisceau $G^i$ sur $Z$. Soit
$$f_T:T\to Z$$
la restriction de $f$ \`a $T$. Il r\'esulte alors de c) et du corollaire qui suit que $\h^j_Y(f^*_T(G^i))=0$ pour $j<q$. On en conclut que
$$
\h^j_Y(\h^i_T(f^*F))=0\text{ pour } i+j<p+q,
$$
car l'in\'egalit\'e $i+j<p+q$ entra\^ine ou bien $i<p$ et alors $\h^i_T(f^*F)=0$, ou bien $j<q$ et alors $\h_Y^j(\h_T^i(f^*F))=0$. \'Etant donn\'e que l'on a, avec les notations de \Ref{XIV.1.0}, $\underline{\Gamma}_Y=\underline{\Gamma}_Y.\underline{\Gamma}_T$, on~a la suite spectrale
\begin{equation*} \label{eq:XIV.1.13.2} \tag{1.13.2}
E_2^{i, j}=\h_Y^j(\h_T^i(f^*F))\To\h^*_Y(f^*F);
\end{equation*}
comme\pageoriginale $E_2^{i, j}=0$ pour $i+j<p+q$, on voit que $\h_Y^k(f^*F)=0$ pour $k<p+q$. Le th\'eor\`eme sera donc d\'emontr\'e si l'on prouve le corollaire suivant (qui est le cas particulier de 1.13 obtenu en y faisant $Z=S, \ r=p=0$ et $F$ r\'eduit au degr\'e $0$).

\begin{corollaire}\label{XIV.1.14}
Soient $f:X\to S$ un morphisme, $Y$ une partie ferm\'ee telle que l'ouvert compl\'ementaire $X- Y$ soit r\'etrocompact dans $X$ et que les anneaux locaux de~$X$ aux diff\'erents points de~$Y$ soient noeth\'eriens. Soient $\LL$ un ensemble de nombres \hbox{premiers}, $q$ un entier et $F$ un faisceau ab\'elien de $\LL$-torsion sur $S$. Supposons que~$f$ soit localement $(q-2)$-acyclique pour $\LL$ et que, pour tout point $s$ de $S$, on ait \hbox{$\prof_{Y_s}^\LL(X_s)\!\geq\! q$}. Alors, on~a $\prof_Y(f^*F)\geq q$.
\end{corollaire}

$1^\circ)$ R\'eduction au cas o\`u $X$ et $S$ sont des sch\'emas strictement locaux, $f$ un morphisme local et $Y$ r\'eduit \`a un point ferm\'e de $X$.

D'apr\`es \Ref{XIV.1.8}, pour \'etablir \Ref{XIV.1.14}, il faut montrer que l'on a pour tout point $y$ de $Y$:$$
\prof_y(f^*F)\geq q.
$$

Soient $s=f(y)$, $\bar s$ un point g\'eom\'etrique au-dessus de $s$, $\bar y$ un point g\'eom\'etrique au-dessus de~$y$ et de $\bar s$, $\bar X$ et $\bar S$ les localis\'es stricts de $X$ et $S$ en $\bar y$ et $\bar s$ respectivement, $\bar f:\bar X\to\bar S$ le morphisme canonique et $\bar F$ l'image inverse de $F$ sur $\bar S$. Comme on~a la relation $\prof_y(f^*F)=\prof_{\bar y}(\bar f^*\bar F)$, il suffit de montrer que les hypoth\`eses de \Ref{XIV.1.14} se conservent quand on remplace $f$ (\resp $Y$, \resp $F$) par $\bar f$ (\resp $\{\bar y\}$, \resp $\bar F$). La condition de r\'etrocompacit\'e r\'esulte de l'hypoth\`ese\pageoriginale noeth\'erienne sur $\Oo_{X, x}$, impliquant que $\bar X$ est noeth\'erien. D'apr\`es (\SGA 4 XV~1.10 \textup{(i)}), $\bar f$ est encore localement $(q-2)$-acyclique pour $\LL$. D'autre part la fibre $(\bar X)_{\bar s}$ de $\bar X$ au-dessus de $\bar s$ s'identifie au localis\'e strict de $X_s$ en $\bar y$, donc satisfait \`a la relation $\prof_{\bar y}((\bar X)_{\bar s})\geq q$. Comme une relation analogue est trivialement v\'erifi\'ee pour les fibres du $\bar S$-sch\'ema $\bar X$, autres que la fibre ferm\'ee, ceci ach\`eve la r\'eduction.

$2^\circ)$ Cas o\`u $X$ et $S$ sont strictement locaux, $f$ un homomorphisme local et $Y$ r\'eduit au point ferm\'e de $X$. Soient
$$g:Y=X-\{y\}\to S$$
le morphisme structural de $U$. On doit montrer que le morphisme canonique
$$u_i:H^i(X, f^*F)\to H^i(U, f^*F)$$
est bijectif pour $i\leq q-2$ et injectif pour $i=q-1$. Consid\'erons le diagramme commutatif
$$
\xymatrix{H^i(X, f^*F)\ar[rr]^{u_i}&& H^i(U, f^*F)\\&H^i(S, F)\ar[lu]^{v_i}\ar[ru]_{w_i}}$$
Le morphisme $v_i$ est \'evidemment bijectif pour tout $i$. D'autre part $g$ est localement $(q-2)$-acyclique pour $\LL$; de plus ses fibres sont $(q-2)$-acycliques pour $\LL$, comme il r\'esulte du fait que $\prof_y(X_s)\geq q$ et que les fibres de $f$ sont $(q-2)$-acycliques pour~$\LL$; comme $g$ est quasi-compact puisque $X$ est noeth\'erien, il r\'esulte de (\SGA 4 XV 1.16) que $g$ est $(q-2)$-acyclique pour $\LL$. Par suite $w_i$, donc aussi $u_i$, est bijectif pour $i\leq q-2$ et injectif pour $i=q-1$, ce qui ach\`eve la d\'emonstration\pageoriginale de \Ref{XIV.1.14}.

\begin{corollaire} \label{XIV.1.15}
Soient $f:X\to S$ un morphisme de sch\'emas, $\LL$ un ensemble de nombres premiers, $m$ et $r$ des entiers et $F$ un complexe de faisceaux ab\'eliens de $\LL$-torsion sur $S$ tel que $\h^i(F)=0$ pour $i<-r$. Soit $x$ un point de $X,\ s=f(x)$ et supposons que l'anneau local $\Oo_{X, x}$ soit noeth\'erien. Alors, si $f$ est localement $m$-acyclique pour $\LL$, on~a la relation
\begin{equation*} \label{eq:XIV.115..*} \tag{$*$} {\prof_x(f^*F)\geq\inf(\prof_s(F)+\prof^\LL_xX_s), n)\text{ o\`u }n=m-r+2.}
\end{equation*}
En particulier, si $n\geq\prof_s(F)+\prof_x^\LL(X_s)$, par exemple si $f$ est localement acyclique pour $\LL$, on a
\begin{equation*} \label{eq:XIV.15.**} \tag{$**$} {\prof_x(f^*F)\geq\prof_s(F)+\prof_x^\LL(X_s).}
\end{equation*}
Si $\LL$ est r\'eduit \`a un \'el\'ement $\ell$ et si l'on a $n\geq \prof_s(F)+\prof_x^\LL(X_s)$, l'in\'egalit\'e pr\'ec\'edente est une \'egalit\'e.
\end{corollaire}

On se ram\`ene au cas o\`u $s$ et $x$ sont des points ferm\'es, en prenant les localis\'es stricts de $S$ et $x$ en des points g\'eom\'etriques $\bar s$ au-dessus de $s$ et $\bar x$ au-dessus de $x$ et de $\bar s$. Si on~a l'in\'egalit\'e $n\geq\prof_s(F)+\prof_x^\LL(X_s)$, alors ($*$) s'obtient \`a partir de \Ref{XIV.1.13} en y faisant $p=\prof_s(F)$ et $q=\prof_x^\LL(X_s)$ (l'hypoth\`ese que $S-\{s\}$ est r\'etrocompact dans $S$ r\'esulte du fait que $X-X_s$ est r\'etrocompact dans $X$ et que $f$ est surjectif puisqu'il est $(-1)$-acyclique (sauf peut-\^etre si la conclusion de 1.15 est vide)). Si l'on~a $n<\prof_s(F)+\prof_x^\LL(X_s)$, l'in\'egalit\'e \eqref{eq:XIV.115..*} s'obtient encore \`a partir de \Ref{XIV.1.13} en y faisant par exemple $p=\prof_s(F)$ et $q=n-p$. Il reste \`a d\'emontrer\pageoriginale la derni\`ere assertion. Soient $p=\prof_s(F)$ et $q=\prof_x(X_s)$; il r\'esulte de \eqref{eq:XIV.1.13.2} que l'on a
$$
\h_x^{p+q}(f^*(F))\simeq \h_x^q(f^*(\h^p_s(F))).
$$
Comme $\prof_s(F)=p$, le faisceau $\h^p_s(F)$ est un faisceau de $\ell$-torsion, constant sur~$s$, diff\'erent de z\'ero. Par suite le faisceau $G=f^*(\h_s^p(F))$ est un faisceau de $\ell$-torsion, constant sur $X_s$, non nul, donc contient un sous-faisceau isomorphe \`a $\ZZ/\ell\ZZ$; comme $\h_x^q(\ZZ/\ell\ZZ)$ est diff\'erent de z\'ero, on~a bien $\h_x^q(G)\neq 0$.

\begin{corollaire} \label{XIV.1.16}
Soient $f:X\to S$ un morphisme r\'egulier de sch\'emas excellents (\textup{\EGA IV 7.8.2}) de caract\'eristique nulle, $\ell$ un nombre premier et $F$ un complexe de faisceaux de $\ell$-torsion sur $S$. Soient $x\in X$, $s=f(x)$; alors on a
$$
\prof_x(f^*F)=\prof_s(F)+2\dim(\Oo_{X, x}).
$$
\end{corollaire}

En effet $f$ \sisi{est r\'egulier par d\'efinition de \og excellent\fg donc}{est} localement acyclique (\SGA 4 XIX 4.1). Il r\'esulte alors de \Ref{XIV.1.15} que l'on a
$$
\prof_x(f^*F)=\prof_s(F)+\prof_x(X_s).
$$
Or on~a d'apr\`es \Ref{XIV.1.10}
$$
\prof_x(X_s)=2\dim\Oo_{X_s, x},
$$
d'o\`u le r\'esultat.

\begin{remarque} \label{XIV.1.17}
Il r\'esulte de \Ref{XIV.1.15} que \Ref{XIV.1.13} reste vrai quand on remplace b) et c) par\pageoriginale les conditions:

b$'$) Pour tout point $s\in f(Y)$, on~a $\prof_s(F)\geq p$.

c$'$) Pour tout point $x\in Y$, si $s=f(x)$, on~a $\prof^\LL_x(X_s)\geq q$.
\end{remarque}

Dans le cas d'un faisceau d'ensembles ou de groupes, on~a le th\'eor\`eme suivant analogue \`a \Ref{XIV.1.13}.

\begin{theoreme} \label{XIV.1.18}
Soient $f:X\to S$ un morphisme de sch\'emas, $Y$ une partie ferm\'ee de $X$ telle que $X-Y$ soit r\'etrocompact dans $X$ et que, pour tout point $x$ de $Y$, l'anneau local $\Oo_{X, x}$ soit noeth\'erien.

$1^\circ)$ Soient $F$ un faisceau d'ensembles sur $S$ et $n$ un entier \'egal \`a $1$ ou $2$. Supposons que $f$ soit localement $(n-2)$-acyclique et que, pour tout point $s$ de $f(Y)$, on ait:
$$
\prof_{Y_s}(X_s)+\prof_s(F)\geq n.
$$
Alors on a:$$
\prof_Y(f^*F)\geq n.
$$

$2^\circ)$ Soient $\LL$ un ensemble de nombres premiers et $F$ un faisceau de ind-$\LL$-groupes. Supposons que $f$ soit localement $1$-asph\'erique pour $\LL$ (\textup{\SGA 4 XV 1.11}) et que, pour tout point $s$ de $f(Y)$, on ait:
$$
\prof\hop^\LL_{Y_s}(X_s)+\prof_s(F)\geq 3.
$$
Alors, on a:$$
\prof_Y(f^*F)\geq 3.
$$
\end{theoreme}

\begin{proof}
On\pageoriginale se ram\`ene, comme dans \Ref{XIV.1.14} et \Ref{XIV.1.15}, au cas o\`u $X$ et $S$ sont des sch\'emas strictement locaux, $f$ un homomorphisme local et $Y$ le point ferm\'e $x$ de $X$. Soit $s=f(x)$ le point ferm\'e de $S$; on~a le diagramme commutatif:
$$
\xymatrix{ X-X_s\ar[r]^i\ar[d]_h&X-\{x\}\ar[rd]_g\ar[rr]^j&&X\ar[ld]^f\\
S-\{s\}\ar[rr]_k&&S. }$$

$1^\circ)$ a) \textit{Cas} $n=1$.

Si l'on a $\prof_s(F)\geq 1$, alors le morphisme $F\to k_*k^*F$ est injectif, donc le morphisme $f^*F\to f^*(k_*k^*F)$ est aussi injectif. D'autre part il r\'esulte du fait que $f$ est localement $(-1)$-acyclique, que le morphisme $f^*(k_*k^*F)\to (j.i)*(f^*\sisi{F{|X-X_s}}{F_{|X-X_s}}))$ est injectif. Finalement, le morphisme compos\'e $f^*F\to (j.i)*(f^*\sisi{F{|X-X_s}}{F_{|X-X_s}}))$ est injectif, ce qui montre que l'on a $\prof_{X_s}(f^*F)\geq 1$, donc aussi $\prof_x(f^*F)\geq 1$.

Si l'on a $\prof_x(X_s)\geq 1$, on consid\`ere le diagramme commutatif
\begin{equation*} \label{eq:XIV.18.*} \tag{$*$}
\begin{array}{c}
\xymatrix{
H^0(X, f^*F)\ar[r]^-{v}\ar[d]_\wr&H^0(X-\{x\}, f^*F)\ar[d]\\
H^0(X_s, f^*F)\ar[r]^-{v'}&H^0(X-\{x\}, f^*F);
}
\end{array}
\end{equation*}
Par hypoth\`ese, $v'$ est injectif donc il en est de m\^eme de $v$.

b) \textit{Cas} $n=2$. On consid\`ere le diagramme commutatif
\begin{equation*} \label{eq:XIV.1.18.**} \tag{$**$}
\begin{array}{c}
\xymatrix{ H^0(S, F)\ar[d]_{m}^{\wr}\ar[rr]^u&&H^0(S-\{s\}, F)\ar[d]^{\wr}_{ n}\\
H^0(X, f^*F)\ar[r]^-v&H^0(X-\{x\}, f^*F)\ar[r]^w&H^0(X-X_s, f^*F);
}
\end{array}
\end{equation*}
on doit\pageoriginaled montrer que $v$ est bijectif. Le morphisme $m$ est \'evidement bijectif, et, comme $f$ est $0$-acyclique, $n$ est aussi bijectif.

Si l'on a $\prof_s(F)\geq 2$, $u$ est bijectif. Comme on l'a vu dans a), la seule hypoth\`ese $\prof_s(F)\geq 1$ entra\^ine la relation $\prof_{X_s}(f^*F)\geq 1$; par suite $v$ et $w$ sont injectifs; il r\'esulte alors de \eqref{eq:XIV.1.18.**} que $v$ est bijectif.

Si l'on a $\prof_x(X_s)\geq 2$, alors $g$ est $0$-acyclique (car il est localement $0$-acyclique et ses fibres sont $0$-acycliques). Il en r\'esulte que $v\cdot m$ est bijectif, donc $v$ est bijectif.

Si l'on a $\prof_s(F)\geq 1$ et $\prof_x(X_s)\geq 1$, alors on sait d\'ej\`a que $v$ et $w$ sont injectifs. Soient $z$ un point maximal de $X_s-\{x\}$ (un tel point existe d'apr\`es l'hypoth\`ese $\prof_x(X_s)\geq 1)$, $\bar Z$ le localis\'e strict de $X$ en un point g\'eom\'etrique au-dessus de $z$ et $\overline{f^*F}$ l'image inverse de $f^*F$ sur $\bar Z$. Consid\'erons le diagramme commutatif
\begin{equation*}
\xymatrix@C=0mm{
&H^0(S, F)\ar[rr]\ar[ld]_{m}^{\sim}&&H^0(S-\{s\}, F)\ar[rd]_{\sim}^{n}\\
H^0(X, f^*F)\ar[rr]^-{v}\ar[rd]_{m'}^{\sim}&&H^0(X-\{x\}, f^*F)\ar[rr]^-{w}\ar[ld]^r&& H^0(X-X_s, f^*F)\ar[ld]^{n'}_{\sim}\\
&H^0(\bar Z, f^*F)\ar[rr]&&H^0(\bar Z-\{\bar z\}, \overline{f^*F})}
\end{equation*}
le morphisme $m'\cdot m$ est \'evidemment bijectif et il r\'esulte du fait que $f$ est localement $0$-acyclique que $n'\cdot n$ est bijectif; par suite $m'$ et $n'$ sont aussi bijectifs. Comme $w$ est injectif, $r$ est aussi injectif et par suite $v$ est bijectif.

$2^\circ)$ Compte tenu de b), on sait d\'ej\`a que $\prof_x(f^*F)\geq 2$.

Si\pageoriginale l'on a $\prof_s(F)\geq 3$, alors $\R^1k_*(k^*F)=1$\nde{le torseur trivial est successivement not\'e $0$ ou $1$ par la suite; on~a laiss\'e cette double notation, qui, r\'eflexion faite, n'apporte aucune ambigu\"it\'e.}. Comme $f$ est localement $1$-asph\'erique, on~a $\R^1(j\cdot i)_*(f^*\sisi{F{|X-X_s}}{F_{|X-X_s}})\!=\!f^*(\R^1k_*(k^*F))\!=\!1$. On a donc \hbox{$\prof_{X_s}(f^*F)\!\geq\! 3$} et par suite on~a $\prof_x(f^*F)\geq 3$.

Si l'on a $\prof_x(X_s)\geq 3$, alors $g$ est $1$-asph\'erique (car $g$ est localement $1$-asph\'erique et ses fibres sont $1$-asph\'eriques). On a donc $H^1(X-\{x\}, f^*F)=H^1(S, F)=1$ et par suite $\prof_x(f^*F)\geq 3$.

Si l'on a $\prof_s(F)\geq 2$ et $\prof_x(X_s)\geq 1$, on utilise la suite exacte (\SGA 4 XII 3.2):
$$1\to\R^1j_*(i_*(f^*\sisi{F{|X-X_s}}{F_{|X-X_s}}))\to\R^1(j.i)_*(f^*\sisi{F{|X-X_s}}{F_{|X-X_s}}) \to j_*(\R^1i_*(f*\sisi{F{|X-X_s}}{F_{|X-X_s}})).
$$
Comme $f$ et $g$ sont localement $1$-asph\'eriques, on a
\begin{align*}
\R^1(j\cdot i)_*(f^*\sisi{F{|X-X_s}}{F_{|X-X_s}})&\simeq f^*(\R^1k_*(k^*F)) \\
\R^1i_*(f^*\sisi{F{|X-X_s}}{F_{|X-X_s}})&\simeq g^*(\R^1k_*(k^*F));
\end{align*}
la suite exacte pr\'ec\'edente s'\'ecrit alors sous la forme
\enlargethispage{\baselineskip}%
\begin{equation} \label{eq:XIV.18.***}\tag{${*}{*}{*}$}
1\to\R^1j_*(i_*(f^*\sisi{F{|X-X_s}}{F_{|X-X_s}}))\to f^*(\R^1k_*(k^*F))\lto{a} j_*(j^*(f^*(\R^1k_*(k^*F))).
\end{equation}
L'hypoth\`ese $\prof_s(F)\geq 2$ montre que le morphisme $F\to k_*k^*F$ est bijectif; en appliquant $g^*$, on trouve compte tenu du fait que $g$ est localement $0$-acyclique, $f^*\sisi{F{|X-\{x\}}}{F_{|X-\{x\}}}=i_*(f^*\sisi{F{|X-X_s}}{F_{|X-X_s}})$. L'hypoth\`ese $\prof_x(X_s)\geq 1$ montre que le morphisme~$a$ est injectif (noter que $f^*(\R^1k_*(k^*F))$ est un faisceau \'egal \`a $1$ en dehors de $X_s$ et constant sur $X_s$). Il r\'esulte alors de \eqref{eq:XIV.18.***} que l'on a $\R^1j_*(f^*\sisi{F{|X-\{x\}}}{F_{|X-\{x\}}})=1$, donc $\prof_x(f^*F)\geq 3$\pageoriginale.

\enlargethispage{.8\baselineskip}%
Si l'on a $\prof_s(F)\geq 1$ et $\prof_x(f^*F)\geq 2$, on consid\`ere le faisceau en espaces homog\`enes $G$ d\'efini par la suite exacte
$$
1\to F\to k_*k*F\to G\to 1.
$$
En appliquant \`a cette suite exacte le foncteur exact $g^*$ et en utilisant (\SGA 4 XII 3.1), on obtient le diagramme commutatif suivant dont les lignes sont exactes:
$$
\xymatrix@R=5mm{f^*(k_*k^*F)\ar[r]\ar[d]&f^*G\ar[r]\ar[d]^b&1\\
j_*(g^*(k_*k^*F))\sisi{\ar[r]}{\ar[r]^-{u}}&j_*(g^*G)\ar[r]&\R^1j_*(g^*F) \ar[r]&\R^1j_*(g^*(k_*k^*F)).}
$$
Comme $\prof_x(X_s)\geq 2$, le morphisme $b$ est bijectif donc $u$ est surjectif et l'on a ainsi une application \`a noyau r\'eduit \`a l'\'el\'ement neutre:
$$
1\to\R^1j_*(g^*F)\to\R^1j_*(g^*(k_*k^*F))=R.
$$
Comme $g^*(k_*k^*F)\simeq i_*(f^*\sisi{F{|X-X_s}}{F_{|X-X_s}})$ (car $g$ est localement $0$-acyclique), $R$ s'identifie au premier terme de la suite exacte \eqref{eq:XIV.18.***}; or on~a vu dans le cas pr\'ec\'edent que $R=1$ d\`es que l'on a $\prof_x(X_s)\geq 1$, ce qui d\'emontre que $\prof_x(f^*F)\geq 3$ et ach\`eve la d\'emonstration de \Ref{XIV.1.18}.
\skipqed
\end{proof}

Les corollaires qui suivent sont des g\'en\'eralisations de (\SGA 4 XVI 3.2 et 3.3).

\begin{corollaire} \label{XIV.1.19}
Soient $f:X\to S$ un morphisme plat, \`a fibres s\'eparables, de sch\'emas localement noeth\'eriens et $Y$ une partie ferm\'ee de $X$. Supposons\pageoriginale que pour tout point $s\in f(Y)$, la fibre $Y_s$ soit rare\nde{\og rare\fg = \og d'int\'erieur vide\fg, \cf Bourbaki TG~IX.52.} dans $X_s$ et que l'une des deux conditions suivantes soit v\'erifi\'ee:
\begin{enumeratea}
\item
l'adh\'erence de $f(Y)$ est rare dans $S$.
\item
$X_s$ est g\'eom\'etriquement unibranche aux points de $Y_s$.
\par\noindent
Alors on a
$$
\prof_Y(X)\geq 2.
$$
\end{enumeratea}
\end{corollaire}

Il r\'esulte en effet de l'hypoth\`ese faite sur $f$ que $f$ est localement $0$-acyclique (\SGA 4 XV 4.1). On applique alors \Ref{XIV.1.13}. L'hypoth\`ese $Y_s$ rare dans $X_s$ (\resp $\overline{f(Y)}$ rare dans~$S$) \'equivaut d'apr\`es \Ref{XIV.1.6} b) \`a la relation $\prof_{Y_s}(X_s)\geq 1$ (\resp $\prof_{\overline{f(Y)}}(S)\geq 1$). L'hypoth\`ese $X_s$ g\'eom\'etriquement unibranche en chaque point de $Y_s$ \'equivaut \`a dire que le localis\'e strict de $X_s$ en un point g\'eom\'etrique de $Y_s$ est irr\'eductible; sachant que $Y_s$ est rare dans $X_s$, cela entra\^ine \'evidemment $\prof_{Y_s}(X_s)\geq 2$, gr\^ace \`a \Ref{XIV.1.8}. Dans l'un et l'autre cas \Ref{XIV.1.13} donne bien $\prof_Y(X)\geq 2$.

\begin{corollaire} \label{XIV.1.20}
Soient $f:X\to S$ un morphisme r\'egulier (\textup{\EGA IV 6.8.1}) de sch\'emas localement noeth\'eriens, $Y$ une partie ferm\'ee de $X$. Supposons que, pour tout point $s\in f(Y)$, l'une des conditions suivantes soit r\'ealis\'ee:
\begin{enumeratea}
\item
On a $\codim(Y_s, X_s)\geq 2$.
\item
On a $\codim(Y_s, X_s)\geq 1$ et $\prof_s(S)\geq 1$.
\item
On a $\prof\hop_s(S)\geq 3$.
\end{enumeratea}
\noindent
Alors on a\pageoriginale
$$
\prof\hop_Y(X)\geq 3.
$$
\end{corollaire}

Cela r\'esulte en effet de \Ref{XIV.1.18}, \'etant donn\'e que l'hypoth\`ese a) implique $\prof\hop_{Y_s}(X_s)\geq 3$ (\cf \Ref{XIV.1.11}), et que la condition $\codim(Y_s, X_s)\geq 1$ implique \'evidemment $\prof_Y(X)\geq 2$.

\section{Lemmes techniques} \label{XIV.2}

\subsection{} \label{XIV.2.1}
Soient $S$ un sch\'ema localement noeth\'erien, $f:X\to S$ un morphisme localement de type fini, $t$ un point de $S$. Si $x\in X$ est tel que $s=f(x)\in\Spec\Oo_{S, t}$, on pose
$$
\delta_t(x)=\degtr k(x)/k(s)+\dim(\overline{\{s\}}),
$$
o\`u $\overline{\{s\}}$ d\'esigne l'adh\'erence de $s$ dans $\Spec\Oo_{S, t}$, $k(x)$ et $k(s)$ les corps r\'esiduels de $x$ et $s$ respectivement. Si $S$ est un anneau local de point ferm\'e $t$, on \'ecrit aussi $\delta(x)$ au lieu de $\delta_t(x)$ (\cf \SGA 4 XIV 2.2).
\begin{sublemme} \label{XIV.2.1.1} Soit un carr\'e cart\'esien
$$
\xymatrix{X'\ar[r]^{h}\ar[d]_{f'}&X\ar[d]^f\\S'\ar[r]^g&S},
$$
o\`u $S$ et $S'$ sont des anneaux locaux noeth\'eriens, de points ferm\'es $t$ et $t'$ respectivement, $g$ un morphisme fid\`element plat tel que $g^{-1}(t)=t'$\pageoriginale, $f$ un morphisme localement de type fini. Soient $x'\in X',\ x=h(x'),\ s=f(x),\ s'=f'(x')$; alors on a
$$
\delta(x')\leq\delta(x).
$$
De plus l'in\'egalit\'e pr\'ec\'edente est une in\'egalit\'e si et seulement si l'on a:
$$
\degtr k(x)/k(s)=\degtr k(x')/k(s')\text{ et }\dim(\overline{\{s\}})=\dim(\overline{\{s'\}}).
$$
En particulier, \'etant donn\'e $x\in X$, on peut trouver $x'$ tel que l'on ait $\delta(x)=\delta(x')$.
\end{sublemme}

On a en effet (\EGA IV 6.11)
$$
\dim(\overline{\{s\}})=\dim g^{-1}(\overline{\{s\}}).
$$
Il en r\'esulte que, pour tout point $s'$ de $g^{-1}(s)$, on~a la relation $\dim(\overline{\{s'\}})\leq \dim(\overline{\{s\}})$, et que, $s$ \'etant donn\'e, on peut trouver $s'\in g^{-1}(s)$, tel qu'on ait l'\'egalit\'e. D\'esignons alors par $Z$ l'adh\'erence sch\'ematique de $x$ dans la fibre $X_s$ de $X$ en~$s$, et soit $Z'=Z\times_{\Spec k(s)}\Spec k(s')$. Alors, $Z'$ est \'equidimensionnel de dimension $\degtr k(x)/k(s)$; on~a donc, pour tout point $x'\in Z'_{x}$,
$$
\degtr k(x')/k(s')\leq\degtr k(x)/k(s)\text{, et on~a l'\'egalit\'e}$$
lorsque $x'$ est un point maximal de $Z'_x$. D'o\`u aussit\^{o}t la conclusion annonc\'ee.

\subsection{} \label{XIV.2.2} Soient $f:X\to S$ un morphisme localement de type fini et $T$ une partie ferm\'ee de $S$. Soient $x\in X, s=f(x)$; nous poserons
$$
\delta_T(x)=\degtr k(x)/k(s)+\codim(\overline{\{s\}}\cap T, \overline{\{s\}})=\inf_{t\in T\cap\overline{\{s\}}}\delta_t(x).
$$

\begin{sublemme} \label{XIV.2.2.1}
Soit\pageoriginale un carr\'e cart\'esien
$$
\xymatrix{X'\ar[r]^{h}\ar[d]_{f'}&X\ar[d]^f\\S'\ar[r]^g&S},
$$
o\`u les sch\'emas $S$ et $S'$ sont localement noeth\'eriens, cat\'enaires, le morphisme $f$ localement de type fini et $g$ fid\`element plat. Soient $T$ une partie ferm\'ee de $S$, $T'$ une partie ferm\'ee de $S'$, telles que $g(T')\subset T$, $x'$ un \'el\'ement de $X'$, $x=h(x')$ et
$$h_{x'}:\Spec\Oo_{X', x'}\to\Spec\Oo_{X, x}$$
le morphisme induit par $h$. Alors on a:
$$
\delta_T(x)-\delta_{T'}(x')\leq\dim h_{x'}^{-1}(x).
$$
\end{sublemme}

Soient $s'=f'(x'), s=f(x)$. On a, par d\'efinition:
\begin{multline*}
\delta_T(x)-\delta_{T'}(x')=\degtr k(x)/k(s)-\degtr k(x')/k(s')\\
+\codim(\overline{\{s\}}\cap T, \overline{\{s\}})-\codim(\overline{\{s'\}}\cap T', \overline{\{s'\}}).
\end{multline*}
Puisque $g$ est fid\`element plat, il r\'esulte de (\EGA IV 6.1.4) que l'on a
\begin{align*} \label{eq:XIV.2.2.*} \tag{$*$}
\codim(\overline{\{s\}}\cap T, \overline{\{s\}})&=\codim(g^{-1}(\overline{\{s\}})\cap g^{-1}(T), g^{-1}(\overline{\{s\}}))\\&\leq\codim(g^{-1}(\overline{\{s\}})\cap T', g^{-1}(\overline{\{s\}}));
\end{align*}
comme $S'$ est cat\'enaire, on a, d'apr\`es ($\EGA 0_{\textup{IV}}$ 14.3.2 b)):
\begin{multline*}
\codim(\overline{\{s'\}}\cap T', g^{-1}(\overline{\{s\}}))= \codim(\overline{\{s'\}}\cap T', \overline{\{s'\}})+ \codim(\overline{\{s'\}}, g^{-1}(\overline{\{s\}}))\\
= \codim(\overline{\{s'\}}\cap T', g^{-1}(\overline{\{s\}})\cap T')+ \codim(g^{-1}(\overline{\{s\}})\cap T', g^{-1}(\overline{\{s\}})).
\end{multline*}
On\pageoriginale d\'eduit de cette relation et de \eqref{eq:XIV.2.2.*}
$$
\delta_T(x)-\delta_{T'}(x')\leq\degtr k(x)/k(s)-\degtr k(x')/k(s')+\codim(\overline{\{s'\}}, g^{-1}(\overline{\{s\}})).
$$
Calculons $\codim(\overline{\{s'\}}, g^{-1}(\overline{\{s\}}))=\dim\Oo_{S'_s, s'}$ (o\`u $S'_s$ est le fibre de $S'$ en $s$). Soit $Z$ l'image ferm\'ee de $x$ dans $X_s$ et $Z'\subset X'_s$ le sch\'ema d\'efini par le carr\'e cart\'esien
$$
\xymatrix{Z'\ar[r]\ar[d]&Z\ar[d]\\S'_s\ar[r]&\Spec k(s)}.
$$
Le morphisme $Z\to \Spec k(s)$ est plat, localement de type fini et l'on a $\dim Z=\degtr k(x)/k(s)$. Il r\'esulte alors de (\EGA IV 6.1.2) que
$$
\dim(\Oo_{Z', x'})=\dim(\Oo_{S'_s, s'})+\degtr k(x)/k(s)-\degtr k(x')/k(s');$$
compte tenu du fait que $Z'_{s'}\simeq Z\otimes_{k(s)}k(s')$, on obtient alors:
$$
\delta_T(x)-\delta_{T'}(x')\leq\dim(\Oo_{Z', x'}).
$$
Or $\Spec(\Oo_{Z', x'})$ s'identifie \`a la fibre en $x$ du morphisme
$$(\Spec(\Oo_{X', x'}))_s\to(\Spec(\Oo_{X, x}))_s,
$$
donc aussi \`a la fibre en $x$ de $h_{x'}$, ce qui d\'emontre le th\'eor\`eme.

\subsection{} \label{XIV.2.3}
Les d\'emonstrations des th\'eor\`emes du \numero \Ref{XIV.4} sont bas\'ees sur la th\'eorie de la dualit\'e; elles utilisent les lemmes qui suivent. Soit $m$ un\pageoriginale entier puissance d'un nombre premier $\ell$; si $X$ est un sch\'ema, tous les faisceaux consid\'er\'es sur $X$ sont des faisceaux de $\ZZ/m\ZZ$-modules; on~a alors la notion de complexe dualisant sur $X$ (\SGA 5 I 1.7). Supposons qu'il existe un tel complexe $K$ sur $X$; alors, pour chaque point g\'eom\'etrique~$\bar x$ \sisi{au dessus}{au-dessus} d'un point $x$ de $X$, on d\'eduit de $K$ (\cf \SGA 5 I 4.5) un complexe dualisant~$K_{\bar x}$ sur $\Spec k(\bar x)$, de sorte que l'on a $K_{\bar x}\simeq\ZZ/m\ZZ[n]$ (le crochet d\'esignant le foncteur de translation) pour un certain entier $n$ ne d\'ependant que de $x$. Nous poserons
$$
\delta_K(x)=n.
$$
Si $K$ est normalis\'e au point $x$ (\SGA 5 I 4.5), on~a donc $n=0$.

\begin{sublemme} \label{XIV.2.3.1} Soit $X$ un sch\'ema localement noeth\'erien, muni d'un complexe dualisant $K$. Si $x$ et $x'$ sont deux points de $X$, tels que $x$ soit une sp\'ecialisation de $x'$ et que l'on ait $\codim(\overline{\{x\}}, \overline{\{x'\}})=1$, alors on a
$$
\delta^K(x)=\delta^{K}(x')-2.
$$
\end{sublemme}

On peut d'abord se ramener au cas o\`u $X$ est un sch\'ema strictement local. Soient en effet $\bar X$ le localis\'e strict de $X$ en un point g\'eom\'etrique $\bar x$ au-dessus de $x$, $i:\bar X\to X$ le morphisme canonique, $\bar x'$ un point g\'eom\'etrique de $\bar X$ au-dessus de $x'$. Alors $i^*K$ est un complexe dualisant sur $\bar X$ et l'on a (\SGA 5 I 4.5)
$$
(i^*K)_{\bar x}\simeq K_{\bar x}\text{ et } (i^*K)_{\bar x'}\simeq K_{\bar x'},
$$
ce\pageoriginale qui ach\`eve la r\'eduction au cas strictement local.

Si $j:\overline{\{x'\}}\to X$ d\'esigne l'immersion du sous-sch\'ema ferm\'e r\'eduit de $X$, d'espace sous-jacent $\overline{\{x'\}}$, alors $\R^!j(K)$ est un complexe dualisant sur $\overline{\{x'\}}$ et on voit tout de suite, en utilisant (\SGA 5 I 4.5) qu'il suffit de d\'emontrer le lemme pour $\overline{\{x'\}}$. On est ainsi ramen\'e au cas o\`u $X$ est un sch\'ema strictement local int\`egre de dimension $1$.

Soient alors $X'$ le normalis\'e de $X$ et $f:X'\to X$ le morphisme canonique; $f$ est un morphisme entier, surjectif, radiciel, et il en r\'esulte que $f^*K$ est un complexe dualisant sur $X'$ et qu'il suffit de d\'emontrer le lemme pour $X'$ et pour les points \sisi{au dessus}{au-dessus} de~$x$ et~$x'$. On est ainsi ramen\'e au cas o\`u $X$ est un sch\'ema local, int\`egre, r\'egulier de dimension $1$, mais on sait (\cf \SGA 5 I 4.6.2 et 5.1) qu'alors $\mu_m[2]$ et $\ZZ/m\ZZ$ sont des complexes dualisants, normalis\'es respectivement aux points $x$ et $x'$; le lemme en r\'esulte aussit\^{o}t.

\begin{sublemme} \label{XIV.2.3.2} Soient $S$ un sch\'ema local noeth\'erien, $f:X\to S$ un morphisme de type fini. Si $K$ est un complexe dualisant sur $S$, normalis\'e au point ferm\'e $t$ de $S$ et si $\R^!f(K)=K'$ est un complexe dualisant sur $X$ (\cf \SGA 5 I 3.4.3), on a, pour tout point $x$ de $X$:$$
\delta^{\sisi{K}{{K'}}}(x)=2\delta(x).
$$
\end{sublemme}

En effet soient $s=f(x)$ et $x'$ un point ferm\'e de la fibre $X_s$; alors on~a $\delta^{K'}(x')=\delta^K(s)$ et d'apr\`es \Ref{XIV.2.3.1}
$$
\delta^K(s)=2\codim(\overline{\{t\}}, \overline{\{s\}})=2\dim(\overline{\{s\}}).
$$
Comme\pageoriginaled on peut choisir pour $x'$ une sp\'ecialisation de $x$, on~a d'apr\`es \Ref{XIV.2.3.1}
$$
\delta^{K'}(x)=\delta^{K'}(x')+2\codim(\overline{\{x\}}, \overline{\{x'\}})=\delta^{K'}(x')+2\degtr k(x)/k(s);$$
le lemme en r\'esulte aussit\^{o}t.

Le lemme suivant servira seulement pour la r\'eciproque du th\'eor\`eme de Lefschetz, dans le \numero\Ref{XIV.4}:

\begin{sublemme} \label{XIV.2.3.3}
Soit un carr\'e cart\'esien
$$
\xymatrix{X'\ar[r]^{h}\ar[d]_{f'}&X\ar[d]^f\\S'\ar[r]^g&S},
$$
o\`u $S$ est un sch\'ema strictement local excellent de caract\'eristique nulle, $S'$ le compl\'et\'e de $S$ et $f$ un morphisme de type fini. Soient $\ell$ un nombre premier, $x\in X$, $Z$ l'adh\'erence sch\'ematique de $X'_x$ dans $X'$, et $i:X'_x\to Z$, $j:Z\to X$ les morphismes canoniques. Alors, si $k:X'\to R$ est une immersion ferm\'ee de $X'$ dans un sch\'ema $R$ r\'egulier excellent, de caract\'eristique nulle, le complexe
$$K'=i^*(R^!(k.j)(\ZZ/\ell\ZZ))$$
est un complexe dualisant sur $X'_x$ constant (c'est-\`a-dire ayant un seul faisceau de cohomologie non nul, isomorphe \`a $\ZZ/\ell\ZZ$).
\end{sublemme}

\enlargethispage{\baselineskip}%
Compte tenu de (\SGA 5 I 3.4.3), la seule chose \`a d\'emontrer est que $K'$ est constant. Or, comme $Z$ est excellent, l'ensemble des points de $Z$ dont\pageoriginale les anneaux locaux sont r\'eguliers est un ensemble ouvert $U$ (\EGA IV 7.8.3 \textup{(iv)}), et $U$ contient \'evidemment $X'_x$ qui est r\'egulier. Soit alors
$$u:U\to R$$
l'immersion canonique de $U$ dans $R$; il r\'esulte du th\'eor\`eme de puret\'e (\SGA 4 XIX 3.2 et 3.4) et de l'isomorphisme
$$(\mu_l)_S\simeq(\ZZ/\ell\ZZ)_S$$
($S$ strictement local) que l'on a
$$
\R^!u(\ZZ/\ell\ZZ)\simeq\ZZ/\ell\ZZ[2c],
$$
o\`u $c$ est une fonction localement constante sur $U$, n\'ecessairement constante au voisinage de $X'_x$, car les fibres de $g$ sont g\'eom\'etriquement int\`egres d'apr\`es (\EGA IV 18.9.1) donc $X'_x$ int\`egre. Le lemme en r\'esulte aussit\^{o}t.

\section{R\'eciproque du th\'eor\`eme de Lefschetz affine} \label{XIV.3}

Le pr\'esent num\'ero sera utilis\'e au \numero \Ref{XIV.4} pour prouver une r\'eciproque au \og th\'eor\`eme de Lefschetz\fg; un lecteur qui n'est int\'eress\'e que par la partie directe dudit th\'eor\`eme peut donc omettre la lecture du pr\'esent num\'ero.

\subsection{} \label{XIV.3.1}

Rappelons l'\'enonc\'e du th\'eor\`eme de Lefschetz affine\nde{\label{gabaf}Gabber a prouv\'e la g\'en\'eralisation suivante. Soit $Y$ un sch\'ema strictement local de type arithm\'etique sur un sch\'ema r\'egulier noeth\'erien $S$ de dimension $\leq 1$. Soit $f:X\to Y$ un morphisme affine de type fini, $\Lambda=\ZZ/n\ZZ$ avec $n$ inversible sur $X$ et $F$ un $\Lambda$-faisceau. Alors, $H^q(X,F)=0$ si $q>\delta(F)$. On en d\'eduit le th\'eor\`eme de Lefschetz local suivant. Soit $\Oo$ est strictement local de type arithm\'etique sur $S$. Pour tout $f\in \Oo$ non diviseur de z\'ero et tout $\Lambda$-faisceau $F$ sur $\Spec(\Oo[f^{-1}])$, on~a $H^q(\Spec(\Oo[f^{-1}]),F)=0$ pour $q>\dim(\Oo)$. \Cf (Fujiwara~K., {\og A Proof of the Absolute Purity Conjecture (after Gabber)\fg}, in \emph{Algebraic geometry 2000, Azumino (Hotaka)}, Adv. Stud. in Pure Math., vol.~36, 2002, p\ptbl 153-183, $\S$5) et surtout l'article d'Illusie (Illusie~L., {\og Perversit\'e et variation\fg}, \emph{Manuscripta Math.} \textbf{112} (2003), p\ptbl 271-295). Ce r\'esultat est un des points cruciaux utilis\'es par Gabber dans sa d\'emonstration du th\'eor\`eme de puret\'e de Grothendieck (\cf la note~\eqref{gabpur}, page~\pageref{gabpur}).} (\SGA 4~XIX~6.1~bis):

Soient $S$ un sch\'ema strictement local excellent de caract\'eristique nulle, $f:X\to S$ un morphisme affine de type fini et $F$ un faisceau de torsion\pageoriginale sur $X$. Alors, si l'on pose
$$
\delta(F)=\sup\{\delta(x)|x\in X \text{ et } F_{\bar{x}} \neq 0\},
$$
on a
$$
H^q(X, F)=0 \text{ pour } q>\delta(F).
$$

Avant d'\'enoncer la r\'eciproque, prouvons quelques lemmes.

\begin{lemme} \label{XIV.3.2} Soient $K$ un corps, $\ell$ un nombre premier distinct de la caract\'eristique de $K$ et $F$ un faisceau de $\ell$-torsion sur $K$, constructible, non nul. Supposons que la $\ell$-dimension cohomologique de $K$ (\textup{\SGA 4 X~1}) soit \'egal \`a $n$ (ceci est r\'ealis\'e par exemple si $K$ est le corps des fractions d'un anneau strictement local excellent int\`egre, de caract\'eristique nulle, de dimension $n$ (\textup{\SGA 4~XIX~6.3}), ou si $K$ est une extension de type fini de degr\'e de transcendance $n$ d'un corps s\'eparablement clos (\textup{\SGA 4~X~2.1})). Alors on peut trouver une extension s\'eparable finie $L$ de $K$, telle que l'on ait:
$$H^n(L, \sisi{F{|L}}{F_{|L}}) \neq 0.
$$
\end{lemme}

On peut trouver une extension finie $K'$ de $K$, telle que les restrictions de $F$ et de $\mu_\ell$ \`a $\Spec K'$ soient des faisceaux constants. On a alors $\cd_{\ell}(K') = \cd_{\ell}(K) = n$ (\SGA 4~X~2.1), et il r\'esulte de (\cite{XIV.2}~II~\S 3~Prop\ptbl4~\textup{(iii)}) que l'on peut trouver une extension finie $L$ de $K'$ telle que l'on ait\pageoriginale
$$
H^n(L, \mu_{\ell})\neq 0, \quad\text {\ie } H^n(L, \ZZ/\ell\ZZ) \neq 0.
$$
Or le foncteur $H^n(L,\cdot)$ est exact \`a droite sur la cat\'egorie des faisceaux de $\ell$-torsion, puisque $\cd_{\ell}(L)=n$; comme $F$ admet un quotient isomorphe \`a $\ZZ/\ell\ZZ$, on~a aussi $H^n(L, \sisi{F{|L}}{F_{|L}}) \neq 0$.

\begin{corollaire} \label{XIV.3.3}
Soient $k$ un corps, $K$ une extension de type fini de degr\'e de transcendance $n$ de $k$, $F$ un faisceau de $\ell$-torsion sur $K$ constructible, non nul, avec $\ell$ premier \`a la caract\'eristique de $k$. Alors on peut trouver une extension finie s\'eparable~$L$ de~$K$ telle que, si $u: \Spec L \to \Spec k$ d\'esigne le morphisme canonique, on ait
$$R^n u_*(\sisi{F{|\Spec L}}{F_{|\Spec L}}) \neq 0.
$$
\end{corollaire}

Lorsque le corps $k$ est s\'eparable\sisi{}{ment} clos, le corollaire est un cas particulier de \Ref{XIV.3.2}. Dans le cas g\'en\'eral, on peut trouver une extension s\'eparable finie $k_1$ de $k\sisi{'}{}$ telle que les composantes irr\'eductibles de $K \otimes_k k_1$ soient g\'eom\'etriquement irr\'eductibles (\EGA IV~4.5.11); soit $K_1$ l'une d'elles. Si $k'$ est une cl\^{o}ture s\'eparable de $k_1$, alors $K' = K_1 \otimes_{k_1} k'$ est un corps, et l'on a d'apr\`es (\EGA IV~4.2)
$$
\deg \tr K'/k' = \deg \tr K/k = n.
$$
Il r\'esulte alors de \Ref{XIV.3.2} que l'on peut trouver une extension finie s\'eparable $L'$ de $K'$ telle que l'on ait $H^n(L', \sisi{F{|L'}}{F_{|L'}}) \neq 0$. Mais on~a $k' = \varinjlim_i k_i$, o\`u $k_i$ parcourt les extensions finies de $k_1$ contenues dans $k'$, et par\pageoriginale suite $K' = \varinjlim_i (k_i \otimes_{k_1} K_1)$. Il en r\'esulte que l'on peut trouver un indice $i$ et une extension finie s\'eparable $L$ de $k_i \otimes_{k_1} K_1 =K_i$, telle que l'on ait $L' \simeq L \otimes_{K_i} K'$. L'extension $L$ de $K$ r\'epond \`a la question; en effet il r\'esulte du diagramme commutatif
$$
\xymatrix{
&{\Spec L}\ar[dl]_v\ar[dr]^u&\\
{\Spec k_i} \ar[rr]^w &{}&{\Spec k}\,,
}
$$
avec $w$ fini donc $R^q w_* = 0$ si $q>0$, que l'on a
$$R^n u_*(\sisi{F{|\Spec L}}{F_{|\Spec L}}) \simeq w_*(R^n v_*(\sisi{F{|\Spec L}}{F_{|\Spec L}})).
$$
Or $R^nv_*(\sisi{F{|\Spec L}}{F_{|\Spec L}}) \neq 0$, puisque $H^n(L', \sisi{F{|L'}}{F_{|L'}}) \neq 0$; on~a donc aussi $R^nu_*(\sisi{F{|\Spec L}}{F_{|\Spec L}}) \neq 0$.

Rappelons le lemme connu suivant (\cf $\EGA 0_{\textup{III}}$~10.3.1.2 et \EGA IV~18.2.3):

\begin{lemme} \label{XIV.3.4}
Soient $X$ un sch\'ema, $x$ un point de $X$, $K$ une extension s\'eparable finie de $k(x)$. Alors il existe un sch\'ema $X_1$ \'etale au-dessus de $X$, affine, et un point $x_1 \in X_1$ au-dessus de $x$, tels que $k(x_1)$ soit $k(x)$-isomorphe \`a $K$.
\end{lemme}

Nous utiliserons au \numero \Ref{XIV.4} la forme technique qui suit de la r\'eciproque de \Ref{XIV.3.1}.

\begin{proposition} \label{XIV.3.5}
Soit\pageoriginale un carr\'e cart\'esien
$$
\xymatrix{
{X'}\ar[r]^h\ar[d]_{f'} &{X}\ar[d]^f\\
{S'}\ar[r]^g &{S}\,,
}
$$
o\`u les sch\'emas $S$ et $S'$ sont strictement locaux excellents de caract\'eristique nulle, le morphisme $f$ localement de type fini, $g$ r\'egulier (\textup{\EGA IV~6.8.1}) surjectif, la fibre ferm\'ee de $g$ r\'eduite au point ferm\'e de $S'$. \'Etant donn\'e un $S$-sch\'ema $X_1$ (\resp un $S$-morphisme $f_1$, etc.), nous noterons $X_1'$ (\resp $f_1'$, etc.) le sch\'ema $X_1 \times_S S'$ (\resp le morphisme $(f_1)_{(S')}$, etc.). Soit $F$ un faisceau de $\ZZ/m\ZZ$-modules sur $X'$ ($m$ puissance d'un nombre premier $\ell$), constructible, satisfaisant aux conditions suivantes:
\begin{enumeratei}
\item
Pour tout point $x \in X$, on peut trouver une extension finie s\'eparable~$K$ de $k(x)$ telle que la restriction de $F$ \`a la fibre $(X')_{(\Spec K)}$ provienne par image r\'eciproque d'un faisceau constructible sur $\Spec K$.
\item
Pour tout morphisme $f_1: X_1 \to S$, avec $X_1$ \'etale au-dessus de $X$, affine, pour tout point $s \in S$ et pour tout entier $q>0$, on peut trouver une extension finie s\'eparable~$K$ de $k(s)$ telle que la restriction de $R^q f'_{1*}(F|X'_1)$ \`a la fibre $S'_{(\Spec K)}$ provienne par image r\'eciproque d'un faisceau constructible sur $\Spec K$.
\end{enumeratei}

Soit $n$ un entier, et supposons que pour tout sch\'ema $X_1$ \'etale au-dessus de $X$, affine, on ait\pageoriginale
$$
H^i(X_1', F)=0 \text{ pour } i>n.
$$
Alors, si $\bar{x}'$ est un point g\'eom\'etrique au-dessus du point $x' \in X'$, tel que $F_{\bar{x}'} \neq 0$, on~a
$$
\delta(x') \leq n.
$$
\end{proposition}

Soit $Z'$ l'ensemble des points $x'$ de $X'$ tels que l'on ait $F_{\bar{x}'} = 0$. Alors, si $Z=h(Z')$, on~a d'apr\`es \textup{(i)} $Z' = h^{-1}(Z)$; soient $x' \in X'$, $x=h(x')$, $s'=f'(x')$, $s=f(x)$. Il r\'esulte de \Ref{XIV.2.1.1} et du fait que la fonction $\delta$ diminue par sp\'ecialisation, qu'il suffit de d\'emontrer l'in\'egalit\'e $\delta(x') \leq n$ lorsque $x$ est un point maximal de $Z$ et $x'$ tel que l'on ait
$$r=\deg \tr k(x)/k(s) = \deg \tr k(x')/k(s') \text{ et }d=\dim \overline{\{s\}} = \dim\overline{\{s'\}}$$
Soit $x'$ un tel point; il suffit de montrer que l'on peut trouver un sch\'ema $X_1$ \'etale sur~$X$, affine, tel que l'on ait
$$
H^{d+r}(X_1', F) \neq 0.
$$

L'ensemble $Z'$ est constructible (\SGA 4~IX~2.4), donc il en est de m\^eme de $Z$ (\EGA IV~1.9.12); on peut alors supposer, quitte \`a restreindre $X$ \`a un voisinage de~$x$, que $Z$ est un ferm\'e irr\'eductible de point g\'en\'erique $x$. Soit $T=f(Z)$; $T$ est un ensemble constructible contenu dans $\overline{\{s\}}$; on peut donc trouver un ouvert affine $U$ de~$S$, tel que $s \in U$ et que $T \cap U = T_U$ soit un ferm\'e irr\'eductible de $U$ de point g\'en\'erique~$s$.

Soit\pageoriginale alors $V$ un sch\'ema \'etale sur $X$, affine, dont l'image dans $X$ contienne $x$ et dont l'image dans $S$ soit contenue dans $U$; soient $Z_V$ l'image inverse de $Z$ dans $V$ et $u: Z_V \to T_U$ le morphisme canonique. Soi\sisi{en}{}t $W$ un sch\'ema \'etale sur $U$, affine, on note alors $T_W$ l'image inverse de $T_U$ dans $W$ et soit $X_1 = W\times_U V$. Comme $F$ est nul en dehors de $Z'$, on~a la suite spectrale
$$
E^{pq}_2 = H^p((T_W)', R^qu_*'(F\sisi{}{_{|(Z_V)'}})) \To H^*(X'_1, F).
$$
Nous allons montrer que l'on peut choisir $V$ et $W$ de telle sorte que l'on ait
\begin{enumeratea}
\item
$E^{pq}_2=0$ pour $p>d$ et pour $q>r$.
\item
$E^{dr}_2 \neq 0$.
\end{enumeratea}
\noindent
Il r\'esultera alors de la suite spectrale la relation $H^{d+r}(X_1', F) \neq 0$.
\begin{enumerate}
\item[$1^\circ)$]
Posons $G_q = R^q u'_*(F\sisi{}{_{|(Z_V)'}})$; alors on a:
$$(G_q)_{\bar{s}'}=H^q((Z_V)'_{\bar{s}'}, F\sisi{}{{_{|(Z_V)'}}}_{\bar{s}'}),
$$
car $s'$ est un point maximal de $(T_U)'$. Comme la fibre $(Z_V)'_{\bar{s}'}$ est un sch\'ema affine de type fini de dimension $r$ sur un corps s\'eparablement clos, il r\'esulte de \Ref{XIV.3.1} que l'on a
$$(G_q)_{\bar{s}'}=0 \text{ pour } q>r.
$$

Pour $q>r$, soit $Y'_q$ l'ensemble des points de $(T_U)'$ o\`u la fibre g\'eom\'etrique\pageoriginale de $G_q$ est $\neq 0$ et $Y_q=g(Y'_q)$; alors on~a $Y'_q = g^{-1}(Y_q)$ d'apr\`es \textup{(ii)}, donc $Y_q$ est un sous-ensemble constructible de $T_U$ (\SGA 4~XIX~5.1 et \EGA 1.9.12) qui ne contient pas $s$; quitte \`a restreindre $U$ \`a voisinage ouvert de $s$, on peut supposer que l'on a $G_q=0$ pour $q>r$, donc $E_2^{pq}=0$ pour $q>r$.

Par ailleurs, comme $(T_W)'$ est un sch\'ema affine de type fini au-dessus de $g^{-1}(\overline{\{s\}})$, on~a quel que soit $q$ (\cf \Ref{XIV.3.1}):
$$H^p((T_W)', G_q)=0 \text{ pour } p>\dim g^{-1}(\overline{\{s\}})=d,
$$
d'o\`u la condition a).

\item[$2^\circ)$]
Montrons que l'on peut choisir $V$ de telle sorte que l'on ait $(G_r)_{\bar{s}'} \neq 0$. D'apr\`es~\textup{(i)}, il existe un faisceau constructible $I$, d\'efini sur une extension finie s\'eparable $K$ de $k(x)$, dont l'image inverse sur $(X')_{(\Spec K)}$ soit isomorphe \`a $\sisi{F{|(X')_{(\Spec K)}}}{F_{|(X')_{(\Spec K)}}}$. D'apr\`es \Ref{XIV.3.3}, on trouve une extension finie s\'eparable $L$ de $K$ telle que, si $v:\Spec L \to \Spec k(s)$ d\'esigne le morphisme canonique, on ait $R^r v_*(I) \neq 0$. Comme le morphisme $S'_s \to \Spec k(s)$ est r\'egulier, on~a d'apr\`es (\SGA 4~XIX~4.2):
$$R^rv_*'(\sisi{F{|(\Spec L)'}}{F_{|(\Spec L)'}}) \simeq (R^rv_*(I))' \neq 0.
$$
D'apr\`es le lemme \Ref{XIV.3.4}, on peut trouver un sch\'ema $X_2$ \'etale sur $X$, affine, et un point~$x_2$ de $X_2$ au-dessus de $x$, tel que $L$ soit $k(x)$-isomorphe \`a $k(x_2)$ et l'on peut supposer $X_2$ au-dessus de $U$. Comme $x$ est un point maximal de $Z$, on a\pageoriginale
$$
\Spec L \simeq \varprojlim_V Z_V,
$$
o\`u $V$ parcourt les voisinages ouverts affines de $x_2$. On en d\'eduit par passage \`a la limite (\SGA 4~VII~5.8), apr\`es restriction \`a la fibre g\'eom\'etrique en $\bar{s}'$:
$$
(R^rv_*'(\sisi{F{|\Spec I}}{F_{|\Spec I}})')_{\bar{s}'} = \varinjlim_V (R^ru_*'(F\sisi{}{_{|(Z_V)'}})_{\bar{s}'},
$$
ce qui montre bien que l'on peut trouver $V$ tel que l'on ait $(G_r)_{\bar{s}'} \neq 0$.

\item[$3^\circ)$]
Le sch\'ema $V$ ayant \'et\'e choisi dans $2^\circ)$, montrons que l'on peut choisir le sch\'ema~$W$ de telle sorte que l'on ait
$$
E^{dr}_2 = H^d((T_W)', G_r) \neq 0.
$$
D'apr\`es \textup{(ii)}, il existe un faisceau constructible $J$, d\'efini sur une extension finie s\'eparable $K$ de $k(s)$, dont l'image inverse sur $(S')_{(\Spec K)}$ soit isomorphe \`a $\sisi{{G_r}{|(S')_{(\Spec K)}}}{{G_r}_{|(S')_{(\Spec K)}}}$. D'apr\`es le lemme \Ref{XIV.3.2}, on peut trouver une extension finie s\'eparable $L$ de $K$, telle que l'on ait $H^d(\Spec L, J) \neq 0$. Comme le morphisme $(S')_{(\Spec L)} \to \Spec L$ est acyclique (\SGA 4~XIX~4.1 et XV~1.10 et 1.16), on a
$$
H^d((\Spec L)', \sisi{{G_r}{|(S')_{(\Spec L)'}}}{{G_r}_{|(S')_{(\Spec L)'}}}) = H^d(\Spec L, J) \neq 0.
$$
D'apr\`es \Ref{XIV.3.4}, on peut trouver un sch\'ema $U_1$ \'etale sur $U$, affine, et un point~$s_1$ au-dessus de $s$, tels que $k(s_1)$ soit $k(s)$-isomorphe \`a $L$. Or, $s$ \'etant un point maximal de~$T_U$, on~a
$$
\Spec L \simeq \varprojlim_W T_W,
$$
o\`u\pageoriginaled $W$ parcourt les voisinages ouverts affines de $s_1$. On en d\'eduit que $(\Spec L)' \simeq \varprojlim_W(T_W)'$, et par passage \`a la limite (\SGA 4~VII~5.8):
$$
H^d((\Spec L)', \sisi{{G_r}{|(\Spec L)'}}{{G_r}_{|(\Spec L)'}}) \simeq \varinjlim_W H^d((T_W)', \sisi{{G_r}{|(T_W)'}}{{G_r}_{|(T_W)'}});
$$
par suite on peut trouver $W$ tel que l'on ait
$$
H^d((T_W)', \sisi{{G_r}{|(T_W)'}}{{G_r}_{|(T_W)'}})\neq 0,
$$
ce qui ach\`eve la d\'emonstration du th\'eor\`eme.
\end{enumerate}

\begin{corollaire} \label{XIV.3.6}
Les hypoth\`eses concernant $S, S', f, f', m$ sont celles de \Ref{XIV.3.5}. D\'esignons maintenant par $F$ un complexe de faisceaux de $\ZZ/m\ZZ$-modules sur $X'$, \`a degr\'es born\'es inf\'erieurement et \`a cohomologie constructible, et dont les faisceaux de cohomologie satisfont aux conditions \textup{(i)} et \textup{(ii)} de \Ref{XIV.3.5}. Soit $n$ un entier, et supposons que, pour tout sch\'ema $X_1$ \'etale sur $X$, affine, on ait
$$
H^i(X'_1, F)=0\text{ pour } i>n.
$$
Alors, si $\bar{x}'$ est un point g\'eom\'etrique au-dessus d'un point $x'$ de $X'$, tel que l'on ait, pour un entier $j$, $(\underline{H}^j(F))_{\bar{x}'} \neq 0$, on a
$$
\delta(x')\leq n-j.
$$
\end{corollaire}

Soit $T'$ l'ensemble des points de $X'$ o\`u la conclusion de \Ref{XIV.3.7} est en d\'efaut et supposons $T' \neq \emptyset$; soit $T=f(T')$, $x$ un point maximal de $T$ et $x'$ un point de $X'$ au-dessus de~$x$. Soit $j$ le plus grand\pageoriginale entier tel que l'on ait $(\underline{H}^j(F))_{\bar{x}'} \neq 0$; on~a donc $r=\delta(x)>n-j$. Soit $Z'_q$ l'ensemble des points o\`u la fibre g\'eom\'etrique de $\underline{H}^q(F)$ est $=0$ et $Z_q =h(Z'_q)$; on voit comme dans la d\'emonstration de \Ref{XIV.3.5} que $Z_q$ est constructible. On a \'evidemment $Z'_q=\emptyset$ pour $q>n$ et pour $q$ suffisamment petit. Les autres valeurs de $q$ se r\'epartissent en trois sous-ensembles. Soit
$$
Q_1=\{q\mid x\in Z_q\text{ et une g\'en\'erisation de }x, \text{ distincte de }x, \not\in Z_q\}.
$$
On a $j \in Q_1$ et on peut trouver un voisinage ouvert affine $U_1$ de $x$, tel que, pour tout $q \in Q_1, U_1 \cap Z_q$ soit un ferm\'e irr\'eductible de point g\'en\'erique $x$. Si $q \in Q_1$, on a
\begin{equation*} \label{eq:XIV.3.6.*} \tag{$*$} \delta(\sisi{\underline{H}^q(F){|U_1'}}{\underline{H}^q(F)_{|U_1'}})=\delta(x)\quad\text{(pour la d\'efinition de $\delta(\underline{H}^q(F))$ \cf \Ref{XIV.3.1})}.
\end{equation*}
Soit
$$
Q_2=\{q\mid\text{aucune g\'en\'erisation de $x$ n'appartient \`a }Z_q\}.
$$
Alors, si $j<q\leq n$, on~a $q\in Q_2$, et l'on peut trouver un voisinage ouvert affine $U_2$ de $x$, tel que, pour tout $q \in Q_2$, on ait $Z_q \cap U_2 = \emptyset$; on~a ainsi
\begin{equation*} \label{eq:XIV.3.6.**} \tag{$**$} \sisi{\underline{H}^q(F){|U'_2}}{\underline{H}^q(F)_{|U'_2}}=0\text{ pour }q\in Q_2.
\end{equation*}
Soit enfin
$$
Q_3= \{q\mid Z_q \text{ contient des g\'en\'erisations strictes de }x\}.
$$
Alors on peut trouver un voisinage ouvert affine de $x$, $U_3$, tel que, pour tout $q \in Q_3$, tous les points maximaux de $Z_q \cap U_3$ soient des g\'en\'erisations de $x$. Si $q \in Q_3$, on a
\begin{equation*} \label{eq:XIV.3.6.***} \tag{${*}{*}{*}$} \delta(\sisi{\underline{H}^q(F){|U'_3}}{\underline{H}^q(F)_{|U'_3}})\leq n-q.
\end{equation*}

Pour\pageoriginale tout sch\'ema $X_1$ \'etale au-dessus de $U_1 \cap U_2 \cap \sisi{{}}{U}_3$, affine, consid\'erons la suite spectrale d'\sisi{}{h}ypercohomologie
$$
E^{pq}_2= H^p(X_1', \underline{H}^q(F)) \To H^*(X'_1, F).
$$
On a $E^{pq}_2=0$ pour $q\sisi{>}{\in Q_2}$ d'apr\`es \eqref{eq:XIV.3.6.**}. On a $E^{pq}_2=0$ pour $p+q \geq r+j$ sauf peut-\^etre pour $p=r, q=j$. En effet c'est clair si $q\in Q_2$; si $q\in Q_1$, on~a alors $p>r$ sauf si $p=r, q=j$ et cela r\'esulte de \Ref{XIV.3.1} compte tenu de (\Ref{eq:XIV.3.6.*}); enfin si $q\in Q_3$, comme $r>n-j$, on~a $p>n-q$ et l'assertion r\'esulte de \Ref{XIV.3.1} compte tenu de (\Ref{eq:XIV.3.6.***}). Vu que $H^{r+j}(X_1', F)=0$, il r\'esulte de la suite spectrale que l'on a
$$H^r(X_1', \underline{H}^j(F))=0;$$
or ceci entra\^ine, d'apr\`es \Ref{XIV.3.5}, $\delta(x)<r$, ce qui est absurde.

\begin{corollaire} \label{XIV.3.7}
Soient $S$ un sch\'ema strictement local excellent de caract\'eristique nulle, $f:X \to S$ un morphisme localement de type fini, $m$ une puissance d'un nombre premier, $F$ un complexe de faisceaux de $\ZZ/m\ZZ$-modules sur $X$, born\'e\sisi{s}{} inf\'erieurement, \`a cohomologie constructible, et $n$ un entier. Alors les conditions suivantes sont \'equivalentes:
\begin{enumeratei}
\item
Pour tout sch\'ema $X_1$ \'etale au-dessus de $X$, affine, on a
$$
H^i(X_1, F)=0\text{ pour }i>n.
$$

\item
Pour\pageoriginale tout point g\'eom\'etrique $\bar{x}$ au-dessus du point $x$ de $X$, et pour tout entier~$j$ tel que l'on ait $(\underline{H}^j(F))_{\bar{x}} \neq 0$, on a
$$
\delta(x) \leq n-j.
$$
\end{enumeratei}
\end{corollaire}

\textup{(i)} \ALORS \textup{(ii)} est le cas particulier de \Ref{XIV.3.6} obtenu en faisant $S=S'$.

\textup{(ii)} \ALORS \textup{(i)} r\'esulte imm\'ediatement de \Ref{XIV.3.1}, en utilisant la suite spectrale d'hypercohomologie
$$
H^p(X_1, \underline{H}^q(F)) \To H^*(X_1, F).
$$

\section{Th\'eor\`eme principal et variantes} \label{XIV.4}

\setcounter{subsection}{-1}

\subsection{} \label{XIV.4.0}
Soient $g\colon X\to S$ un morphisme s\'epar\'e de type fini, $T$ une partie ferm\'ee de $S$, \hbox{$Z=g^{-1}(T)$} et $F$ un complexe de faisceaux ab\'eliens sur $X$, \`a degr\'es born\'es inf\'erieurement. Nous appelons \emph{$i$-i\`eme groupe de cohomologie de $F$, \`a support propre, \`a support dans~$Z$} le groupe
$$
\H_{Z!}^{i}(X/S, F)=\H_{T}^{i}(S, \R_{!}g(F)),
$$
o\`u $\R_{!}g$ d\'esigne \og l'image directe \`a support propre\fg (\SGA 4~XVII). Dans le cas particulier o\`u $g$ est propre, on~a simplement
$$
\H_{Z!}^{i}(X/S, F)=\H_{Z}^{i}(X, F).
$$

\begin{proposition} \label{XIV.4.1} Soient $f\colon U\to S$ un morphisme de type fini, $F$ un complexe de faisceaux ab\'eliens sur $U$, \`a degr\'es born\'es inf\'erieurement. Supposons\pageoriginale que l'on ait une factorisation de $f$:
$$
\xymatrix@C=15pt{ U\ar[rr]^-{i}\ar[rd]_-{f} & & X\ar[ld]^-{g}\\
& S & }$$
o\`u $i$ est une immersion ouverte et $g$ un morphisme s\'epar\'e de type fini, et d\'esignons par $G$ un complexe de faisceaux ab\'eliens sur $X$, \`a degr\'es born\'es inf\'erieurement, qui prolonge $F$. Soient $Y$ un sous-sch\'ema ferm\'e de $X$ d'espace sous-jacent $X-U$, de sorte que l'on a un diagramme commutatif:
$$
\xymatrix@C=15pt{ Y\ar[rr]^-{j}\ar[rd]_-{h} & & X\ar[ld]^-{g}\\
& S & }$$
Soient enfin $n$ un entier et $T$ une partie ferm\'ee de $S$. Alors les conditions suivantes sont \'equivalentes:
\begin{enumeratei}
\item
On a $\prof_{T}(\R_{!}f(F))\geq n$.
\item
Le morphisme canonique
$$
\underline{H}_{T}^{i}(\R_{!}g(G))\to\underline{H}_{T}^{i}(\R_{!}h(j^{\ast}G))$$
est bijectif pour $i<n-1$, injectif pour $i=n-1$.
\item
Pour tout sch\'ema $S^{\prime}$ \'etale au-dessus de $S$, si l'on d\'esigne par\pageoriginale $X^{\prime}$ (\resp $f^{\prime}$, \resp etc.) le sch\'ema $X\times_{S}S^{\prime}$ (\resp le morphisme $f_{(S^{\prime})}$, \resp etc.), le morphisme canonique
$$
\H_{g^{\prime -1}(T^{\prime})!}^{i}(X^{\prime}/S^{\prime}, G^{\prime})\to\H_{h^{\prime -1}(T^{\prime})!}^{i}(Y^{\prime}/S^{\prime}, j^{\prime\ast}G^{\prime})$$
est bijectif pour $i<n-1$, injectif pour $i=n-1$.
\end{enumeratei}
\end{proposition}

Consid\'erons dans la cat\'egorie d\'eriv\'ee $\D^{+}(X)$ (\cf \cite{XIV.3}) le triangle distingu\'e
$$
\xymatrix@C=15pt{ & j_{\ast}j^{\ast}G\ar[ld] & \\
i_{!}F\ar[rr] & & G.\ar[lu]}$$
En appliquant \`a ce triangle le foncteur $\R_{!}g$, on obtient le triangle
\begin{equation*} \label{eq:XIV.4.1.*} \tag{$*$}
\begin{array}{c}
\xymatrix@C=0pt{
& \R_{!}h(j^{\ast}G)\ar[ld] & \\
\R_{!}f(F)\ar[rr] & & \R_{!}g(G).\ar[lu]
}
\end{array}
\end{equation*}

Montrons que \textup{(i)}\SSI\textup{(ii)}. En effet, d'apr\`es la d\'efinition \Ref{XIV.1.2}, \textup{(i)} \'equivaut \`a la relation
$$
\underline{H}_{T}^{i}(\R_{!}f(F))=0\quad\textrm{ pour }i<n\textrm{;}$$
or on d\'eduit de \eqref{eq:XIV.4.1.*} la suite exacte de faisceaux
$$
\to\underline{H}_{T}^{i}(\R_{!}f(F))\to\underline{H}_{T}^{i}(\R_{!}g(G))\to\underline{H}_{T}^{i}(\R_{!}h(j^{\ast}G))\to,
$$
d'o\`u\pageoriginale l'\'equivalence de \textup{(i)} et \textup{(ii)}.

\textup{(i)}\SSI\textup{(iii)}. En effet \textup{(i)} \'equivaut \`a dire que, pour tout sch\'ema $S^{\prime}$ \'etale au-dessus de $S$, on~a la relation
\begin{equation*} \label{eq:XIV.4.1.**} \tag{$**$} {\H_{T^{\prime}}^{i}(S^{\prime}, \R_{!}f^{\prime}(F^{\prime}))=0\quad\textrm{ pour }i<n.}
\end{equation*}
Or on d\'eduit de \eqref{eq:XIV.4.1.*} la suite exacte de groupes ab\'eliens
$$
\to\H_{T^{\prime}}^{i}(S^{\prime}, \R_{!}f^{\prime}(F^{\prime}))\to\H_{T^{\prime}}^{i}(S^{\prime}, \R_{!}g^{\prime}(G^{\prime}))\to \H_{T^{\prime}}^{i}(S^{\prime}, \R_{!}h^{\prime}(j^{\prime\ast}G^{\prime}))\to\;;$$
compte tenu de \Ref{XIV.4.0}, cette suite exacte s'\'ecrit sous la forme
$$
\to\H_{T^{\prime}}^{i}(S^{\prime}, \R_{!}f^{\prime}(F^{\prime}))\to\H_{g^{\prime -1}(T^{\prime})!}^{i}(X^{\prime}/S^{\prime}, G^{\prime})\to\H_{h^{\prime -1}(T^{\prime})!}^{i}(Y^{\prime}/S^{\prime}, j^{\prime\ast}G^{\prime})\sisi{}{\to.}$$
L'\'equivalence de \textup{(i)} et \textup{(iii)} en r\'esulte, compte tenu de la forme \eqref{eq:XIV.4.1.**} de \textup{(i)}.

\setcounter{subsubsection}{-1} \setcounter{subsection}{2}

\subsubsection{} \label{XIV.4.2.0}
Lorsque $f\colon U\to S$ est affine, nous allons donner des conditions \emph{locales} sur $F$ pour que les conditions \textup{(i)} \`a \textup{(iii)} de \Ref{XIV.4.1} soient v\'erifi\'ees. Dans la suite, les sch\'emas consid\'er\'es sont des sch\'emas excellents de caract\'eristique nulle, les faisceaux sont des faisceaux de $\ZZ/m\ZZ$-modules, o\`u $m$ est une puissance d'un nombre premier. Si l'on disposait de la r\'esolution des singularit\'es au sens de (\SGA 4~XIX), les r\'esultats \'enonc\'es, ainsi que leurs d\'emonstrations, seraient encore valables pour des sch\'emas excellents d'\'egale caract\'eristique, avec $m$ premier \`a la caract\'eristique.

\setcounter{subsection}{1}

\begin{theoreme} \label{XIV.4.2}
Soient\pageoriginale $S$ un sch\'ema excellent de caract\'eristique nulle et \hbox{$f\colon U\!\to\! S$} un morphisme s\'epar\'e de type fini. Soient $F$ un complexe de faisceaux de $\ZZ/m\ZZ$-modules sur $U$, \`a degr\'es born\'es inf\'erieurement et \`a cohomologie constructible, $n$ un entier et~$T$ une partie ferm\'ee de $S$. Alors les conditions suivantes sont \'equivalentes:
\begin{enumeratei}
\item
Pour tout sch\'ema $U_{1}$ \'etale au-dessus de $U$, affine sur $S$, on a, en d\'esignant par~$f_{1}$ le morphisme structural de $U_{1}$ et par $F_{1}$ la restriction de $F\sisi{_{1}}{}$ \`a $U_{1}$:
$$
\prof_{T}(\R_{!}f_{1}(F_{1}))\geq n$$
(\cf prop\ptbl \Ref{XIV.4.1} sur la signification de cette relation).
\item
Pour tout point $u$ de $U$, on a:
$$
\prof_{u}(F)\geq n-\delta_{T}(u),
$$
o\`u l'on pose (\cf \Ref{XIV.2.2}): $\delta_{T}(u)=\degtr(k(x)/k(s))+\codim\big(\overline{\{ s\} }\cap T, \overline{\{ s\} } \big)$.
\end{enumeratei}
\end{theoreme}

\begin{proof}
$1^\circ)$ Soient $t$ un point de $T$, $\overline{S}$ le localis\'e strict de $S$ en un point g\'eom\'etrique au-dessus de $t$ et $S^{\prime}$ le compl\'et\'e de $\overline{S}$, de point ferm\'e $t^{\prime}$; alors $\overline{S}$ est excellent d'apr\`es (\EGA IV~7.9.5), donc $S^{\prime}$ est un sch\'ema strictement local complet excellent. \'Etant donn\'e un sch\'ema $U$ sur $X$ (\resp un $S$-morphisme $f$, \resp etc.), nous d\'esignerons par $U^{\prime}$ (\resp $f^{\prime}$, \resp etc.) le sch\'ema $U\times_{S}S^{\prime}$ (\resp le morphisme $f_{(S^{\prime})}$, \resp etc.). On a le carr\'e cart\'esien
$$
\xymatrix{ U^{\prime}\ar[r]^-{h}\ar[d]_-{f^{\prime}} & U\ar[d]^-{f}\\
S^{\prime}\ar[r]^-{g} & S, }
$$
dans\pageoriginale lequel le morphisme $g$ est r\'egulier (\EGA IV~7.8.2). Montrons qu'il suffit de prouver que (pour tout point $t\in T$) les deux propri\'et\'es suivantes sont \'equivalentes:
\begin{itemize}
\item[$\textup{(i)}_{t}$] Pour tout sch\'ema $U_{1}$ \'etale sur $U$, affine sur $S$, posant $f_{1}\colon U_{1}\to S$, on a
$$
\prof_{t^{\prime}}(\R_{!}f_{1}^{\prime}(F_{1}^{\prime}))\geq n.
$$

\item[$\textup{(ii)}_{t}$] Pour tout point $u^{\prime}$ de $U^{\prime}$, on a
$$
\prof_{u^{\prime}}(F^{\prime})\geq n-\delta_{t^{\prime}}(u^{\prime}).
$$
\end{itemize}
Il suffit de d\'emontrer le lemme suivant:

\begin{sublemme} \label{XIV.4.2.1}
On a \textup{(i)}\SSI$\textup{(i)}_{t}$ pour tout $t\in T$ et \textup{(ii)}\SSI$\textup{(ii)}_{t}$ pour tout $t\in T$.
\end{sublemme}

\textup{(i)}\SSI$\textup{(i)}_{t}$ pour tout $t\in T$. En effet \textup{(i)} \'equivaut \`a dire que, pour tout sch\'ema $U_{1}$ \'etale au-dessus de $U$, affine sur $S$, on a
$$
\prof_{T}(\R_{!}f_{1}(F_{1}))\geq n\;;$$
or d'apr\`es \Ref{XIV.1.8}
$$
\prof_{T}(\R_{!}f_{1}(F_{1}))=\inf_{t\in T}\prof_{t}(\R_{!}f_{1}(F_{1})).
$$
Comme $g^{\ast}(\R_{!}f_{1}(F_{1}))\simeq\R_{!}f_{1}^{\prime}(F_{1}^{\prime})$ (\SGA 4~XVII), on~a d'apr\`es \Ref{XIV.1.16}
$$
\prof_{t}(\R_{!}f_{1}(F_{1}))=\prof_{t^{\prime}}(\R_{!}f_{1}^{\prime}(F_{1}^{\prime})),
$$
donc\pageoriginale \textup{(i)} \'equivaut \`a dire que l'on a, pour tout $t\in T$, $\prof_{t^{\prime}}(\R_{!}f_{1}^{\prime}(F_{1}^{\prime}))\geq n$, ce qui n'est autre que $\textup{(i)}_{t}$.

\textup{(ii)}$_{t}$ pour tout $t\in T$\ALORS\textup{(ii)}. En effet soit $u\in U$; on doit montrer la relation
$$
\prof_{u}(F)\geq n-\delta_{T}(u),
$$
o\`u $\delta_{T}(u)=\inf_{t\in T\cap\overline{\{ s\} }}\delta_{t}(u)$ (\cf \Ref{XIV.2.2}); on est donc ramen\'e \`a montrer que l'on a, pour tout $t\in T\cap\overline{\{ s\} }$
$$
\prof_{u}(F)\geq n-\delta_{t}(u).
$$
Soient $u^{\prime}$ un point de $U^{\prime}$ tel que l'on ait $h(u^{\prime})=u$ et $\delta_{t^{\prime}}(u^{\prime})=\delta_{t}(u)$ (\cf \Ref{XIV.2.1.1}). Comme $h$ est localement acyclique (\SGA 4~XIX~4.1), il r\'esulte de \Ref{XIV.1.16} et du fait que~$u^{\prime}$ est un point g\'en\'erique de $U_{u}^{\prime}$ que l'on a
$$
\prof_{u^{\prime}}(F^{\prime})=\prof_{u}(F).
$$
Mais on~a d'apr\`es $\textup{(ii)}_{t}$ $\prof_{u}(F)=\prof_{u^{\prime}}(F^{\prime})\geq n-\delta_{t}(u)$, ce qui d\'emontre \textup{(ii)}.

\textup{(ii)}\ALORS$\textup{(ii)}_{t}$ pour tout $t$. Avec les notations de \Ref{XIV.2.2.1}, pour tout point $u^{\prime}$ de $U^{\prime}$, on~a gr\^ace \`a \Ref{XIV.1.16}
$$
\prof_{u^{\prime}}(F^{\prime})\geq\prof_{u}(F)+2\dim h_{u^{\prime}}^{-1}(u)\geq\prof_{u}(F)+\dim h_{u^{\prime}}^{-1}(u).
$$
Compte tenu de \Ref{XIV.2.2.1} et \textup{(ii)}, on obtient
$$
\prof_{u^{\prime}}(F^{\prime})\geq n-\delta_{T}(u)+\dim h_{u^{\prime}}^{-1}(u)\geq n-\delta_{t^{\prime}}(u^{\prime}),
$$
ce\pageoriginaled qui n'est autre que $\textup{(ii)}_{t}$.

$2^\circ)$ $\textup{(ii)}_{t}$\SSI$\textup{(i)}_{t}$. On se ram\`ene imm\'ediatement au cas o\`u $F$ est \`a degr\'es born\'es, en tronquant $F$ \`a un rang suffisamment \'elev\'e. On peut r\'ealiser $S^{\prime}$ comme ferm\'e d'un sch\'ema local r\'egulier complet, donc excellent; il r\'esulte alors de (\SGA 5~I~3.4.3) qu'il existe un complexe dualisant $K$ sur $S^{\prime}$ et que $\R^{!}f^{\prime}(K)=K^{\prime}$ est un complexe dualisant sur $U^{\prime}$. Nous choisirons $K$ de telle sorte que l'on ait $\delta^{K}(t^{\prime})=0$ (pour la d\'efinition de $\delta^{K}(t^{\prime})$, \cf \Ref{XIV.2.3}), et noterons $\D F^{\prime}$ le dual de $F^{\prime}$ par rapport \`a $K^{\prime}$. On peut reformuler l'hypoth\`ese $\textup{(ii)}_{t}$ de la fa\c con suivante:

\begin{sublemme} \label{XIV.4.2.2}
Soit $u^{\prime}$ un point de $U^{\prime}$; alors les conditions suivantes sont \'equivalentes:
\begin{enumeratei}
\item
On a $\prof_{u^{\prime}}(F^{\prime})\geq n-\delta_{t^{\prime}}(u^{\prime})$.
\item
On a $\big(\underline{H}^{q}(\D F^{\prime}) \big)_{\overline{u}^{\prime}}=0$ pour $q>-n-\delta_{t^{\prime}}(u^{\prime})$ ($\overline{u}^{\prime}$ point g\'eom\'etrique au-dessus de $u^{\prime}$).
\end{enumeratei}
\end{sublemme}

Soient $\overline{U}^{\prime}$ le localis\'e strict de $U^{\prime}$ en $\overline{u}^{\prime}$ et $\overline{F}^{\prime}$ l'image r\'eciproque de $F$ par le morphisme $\overline{U}^{\prime}\to U^{\prime}$. La relation $\prof_{u^{\prime}}(F^{\prime})\geq n-\delta_{t^{\prime}}(u^{\prime})$ \'equivaut par d\'efinition \`a la suivante:
\begin{equation*} \label{eq:XIV.4.2.*} \tag{$*$} {\H_{\overline{u}^{\prime}}^{i}(\overline{F}^{\prime})=0\quad\textrm{ pour }i>n-\delta_{t^{\prime}}(u^{\prime}).}
\end{equation*}
Soit $\D\big(\H_{\overline{u}^{\prime}}^{i}(\overline{F}^{\prime})\big)$ le dual du groupe ab\'elien $\H_{\overline{u}^{\prime}}^{i}(\overline{F}^{\prime})$ par rapport \`a $\ZZ/m\ZZ$. D'apr\`es \Ref{XIV.2.3.2}, $K^{\prime}[-2\delta_{t^{\prime}}(u^{\prime})]=K^{\prime\prime}$ satisfait \`a $\delta^{K^{\prime\prime}}(u^{\prime})=0$; comme $F^{\prime}$ est \`a cohomologie constructible, on~a $\overline{\D F^{\prime}}=\D(\overline{F}^{\prime})$ et le th\'eor\`eme de\pageoriginale dualit\'e locale (\SGA 5~I~4.5.3) montre alors que l'on a
$$
\D\big(\H_{\overline{u}^{\prime}}^{i}(\overline{F}^{\prime})\big) \simeq\big(\H^{-i-2\delta_{t^{\prime}}(u^{\prime})}(\D F^{\prime})\big)_{\overline{u}^{\prime}}.
$$
Donc \eqref{eq:XIV.4.2.*} \'equivaut \`a la relation
\begin{equation*} \label{eq:XIV.4.2.**} \tag{$**$}
\big(\underline{H}^{q}(\D F^{\prime})\big)_{\overline{u}^{\prime}}=0\quad\text{pour }q>-n-\delta_{t^{\prime}}(u^{\prime}).
\end{equation*}

Nous sommes maintenant en mesure de d\'emontrer le th\'eor\`eme. La relation \textup{(ii)}$_{t}$ \'equivaut \`a la relation \eqref{eq:XIV.4.2.**}. Soit $G^{q}=\underline{H}^{q}(\D F^{\prime})$; le th\'eor\`eme de Lefschetz affine (\Ref{XIV.3.1}) entra\^ine en particulier que, pour tout sch\'ema $U_{1}$ \'etale sur $U$, affine sur $S$, on a
$$
\H^{p}(U_{1}^{\prime}, G^{q})=0\quad\textrm{ pour }p>\delta(G^{q}),
$$
o\`u $\delta(G^{q})$ est la borne sup\'erieure des $\delta_{t^{\prime}}(u^{\prime})$ pour les $u^{\prime}$ tels que l'on ait $G_{\overline{u}^{\prime}}^{q}\neq 0$; d'apr\`es \eqref{eq:XIV.4.2.**} on~a $\delta(G^{q})\leq -n-q$, donc $\textup{(ii)}_{t}$ entra\^ine la relation
$$
\H^{p}(U_{1}^{\prime}, \underline{H}^{q}(\D F^{\prime}))=0\quad\textrm{ pour }p>-q-n.
$$
Compte tenu de la suite spectrale d'hypercohomologie du foncteur \og sections sur $U_{1}^{\prime}$\fg par rapport au complexe $\D F^{\prime}$:
$$
\E_{2}^{pq}=\H^{p}(U_{1}^{\prime}, \underline{H}^{q}(\D F^{\prime}))\To\H^{\ast}(U_{1}^{\prime}, \D F^{\prime}),
$$
on obtient la relation
\begin{equation*} \label{eq:XIV.4.2.***} \tag{${*}{*}{*}$}
\H^{i}(U_{1}^{\prime}, \D F^{\prime})=0\quad\textrm{ pour }i>-n.
\end{equation*}

Inversement\pageoriginale supposons v\'erifi\'ee la relation pr\'ec\'edente, pour tout $U_{1}$ \'etale sur $U$, affine sur $S$. Appliquons la proposition \Ref{XIV.3.6} en y rempla\c cant $S$ par $\overline{S}^{\prime}$; les hypoth\`eses de \Ref{XIV.3.6} concernant $\sisi{S}{\overline{S}}$ sont satisfaites, car, pour tout sch\'ema $\overline{U}_{1}$ \'etale sur $\overline{U}$, affine, on peut trouver un sch\'ema au-dessus de $\overline{U}_{1}$ qui provienne par image r\'eciproque d'un sch\'ema \'etale au-dessus de $U$, affine sur $S$; quant aux hypoth\`eses concernant $F$, elles sont satisfaites gr\^ace \`a \Ref{XIV.2.3.3}. On a ainsi, pour tout point $u^{\prime}$ de $U^{\prime}$ tel que $\big(\underline{H}^{q}(\D F^{\prime})\big)_{\overline{u}^{\prime}}\neq 0$:
$$
\delta_{t^{\prime}}(u^{\prime})\leq -n-q,
$$
ce qui n'est autre que la relation \eqref{eq:XIV.4.2.**}; on~a donc prouv\'e l'\'equivalence
$$
\text{\textup{(ii)}$_{t}$\SSI\eqref{eq:XIV.4.2.***}}.
$$
Nous allons transformer la relation \eqref{eq:XIV.4.2.***}; on~a d'abord
$$
\H^{i}(U_{1}^{\prime}, \D F^{\prime})=\big(\underline{H}^{i}(\R f_{1\ast}^{\prime}(\D F_{1}^{\prime}))\big)_{t^{\prime}}\;;$$
mais d'apr\`es (\SGA 5~I~1.12), il existe un isomorphisme canonique
$$
\R f_{1\ast}^{\prime}(\D F_{1}^{\prime})\simeq\D\big(\R_{!}f_{1}^{\prime}(F_{1}^{\prime})\big),
$$
o\`u $\D\big(\R_{!}f_{1}^{\prime}(F_{1}^{\prime})\big)$ d\'esigne le dual de $\R_{!}f_{1}^{\prime}(F_{1}^{\prime})$ par rapport \`a $K$. On voit ainsi que $\textup{(ii)}_{t}$ \'equivaut \`a
$$
\big(\underline{H}^{i}(\D(\R_{!}f_{1}^{\prime}(F_{1}^{\prime})))\big)_{t^{\prime}}=0\quad\textrm{ pour }i>-n.
$$
Appliquant\pageoriginale de nouveau le th\'eor\`eme de dualit\'e locale (\SGA 5~I~4.5.3), mais cette fois-ci au point $t^{\prime}$, on trouve que
$$
\big(\underline{H}^{i}(\D(\R_{!}f_{1}^{\prime}(F_{1}^{\prime})))\big)_{t^{\prime}}\simeq\D\big(\H_{t^{\prime}}^{-i}(\R_{!}f_{1}^{\prime}(F_{1}^{\prime}))\big),
$$
et finalement $\textup{(ii)}_{t}$ \'equivaut \`a la relation
$$
\H_{t^{\prime}}^{i}(\R_{!}f_{1}^{\prime}(F_{1}^{\prime}))=0\quad\textrm{ pour }i<n,
$$
c'est-\`a-dire $\prof_{t^{\prime}}(\R_{!}f_{1}^{\prime}(F_{1}^{\prime}))\geq n$, ce qui ach\`eve la d\'emonstration du th\'eor\`eme.
\skipqed
\end{proof}

\begin{subremarque} \label{XIV.4.2.3}
Le raisonnement se simplifie assez consid\'erablement lorsqu'on suppose que $S$ admet (du moins localement) un complexe dualisant (par exemple est localement immergeable dans un sch\'ema r\'egulier). Cela \'evite le recours \`a un compl\'et\'e (le passage au cas $S$ strictement local \'etant imm\'ediat), \`a \Ref{XIV.2.3.3} et \`a l'\'enonc\'e technique peu plaisant \Ref{XIV.3.6}, qu'on peut alors remplacer par la r\'ef\'erence plus sympathique \Ref{XIV.3.7}.
\end{subremarque}

\begin{corollaire} \label{XIV.4.3}
Soient $S$ un sch\'ema excellent de caract\'eristique nulle et $f\colon U\to S$ un morphisme s\'epar\'e de type fini, tel que $U$ soit \emph{r\'eunion de $c+1$ ouverts, affines sur~$S$}. Soient $F$ un complexe de faisceaux de $\ZZ/m\ZZ$-modules, \`a degr\'es born\'es inf\'erieurement et \`a cohomologie constructible, $n$ un entier et $T$ une partie ferm\'ee de $S$. Supposons que, pour tout point $u\in U$, on ait
$$
\prof_{u}(F)\geq n-\delta_{T}(u).
$$
Alors on a\pageoriginale
$$
\prof_{T}(\R_{!}f(F))\geq n-c.
$$
\end{corollaire}

Soit en effet $U_{\sisi{i}{j}}$, $0\leq \sisi{i}{j}\leq c$, un recouvrement de $U$ par des ouverts $U_{\sisi{i}{j}}$, affines sur~$S$. Reprenant les notations de la d\'emonstration de \Ref{XIV.4.2}, on a, pour tout $\sisi{i}{j}$,
$$
\H^{i}(U_{\sisi{i}{j}}^{\prime}, \underline{H}^{q}(\D F^{\prime}))=0\quad\textrm{ pour }i>-n.
$$
En utilisant la suite spectrale qui relie la cohomologie de $U$ \`a celle du recouvrement form\'e par les $U_{\sisi{i}{j}}$ (\SGA 4~V~2.4), la relation pr\'ec\'edente montre que l'on a
$$
\H^{i}(U^{\prime}, \underline{H}^{q}(\D F^{\prime}))=0\quad\textrm{ pour }i>-n+c.
$$
Le corollaire r\'esulte alors de la fin de la d\'emonstration de \Ref{XIV.4.2}.

\begin{corollaire} \label{XIV.4.4}
Soient $S$ un sch\'ema excellent de caract\'eristique nulle, $g\colon X\to S$ un morphisme, $U$ un ouvert de $X$, r\'eunion de $c+1$ ouverts affines sur $S$, $Y$ un sous-sch\'ema ferm\'e d'espace sous-jacent $X-U$ et $j\colon Y\to X$ le morphisme naturel. Soient~$F$ un complexe de faisceaux de $\ZZ/m\ZZ$-modules sur $X$, \`a degr\'es born\'es inf\'erieurement et \`a cohomologie constructible, $T$ une partie ferm\'ee de $S$ et $n$ un entier. Supposons que, pour tout point $u$ de $U$, on ait
$$
\prof_{u}(F)\geq n-\delta_{T}(u).
$$
Alors\pageoriginale le morphisme canonique
$$
\H_{g^{-1}(T)!}^{i}(X/S, F)\to\H_{(g^{-1}(T)\cap Y)!}^{i}(Y/S, j^{\ast}F)$$
est bijectif pour $i<n-c-1$, injectif pour $i=n-c-1$.
\end{corollaire}

Cela r\'esulte imm\'ediatement de \Ref{XIV.4.1} et \Ref{XIV.4.3}.

\begin{corollaire}[Th\'eor\`eme de Lefschetz local] \label{XIV.4.5}
Soient $S$ un sch\'ema local hens\'elien excellent de caract\'eristique nulle, $t$ le point ferm\'e de $S$, $X$ un sch\'ema propre sur $S^{\prime}=S-\{ t\}$ et $U$ un ouvert de $X$, r\'eunion de $c+1$ ouverts affines. Soient $Y$ un sous-sch\'ema ferm\'e de $X$, d'espace sous-jacent $X-U$, $j\colon Y\to X$ le morphisme canonique, $F$ un complexe de faisceaux de $\ZZ/m\ZZ$-modules sur $X$, \`a degr\'es born\'es inf\'erieurement et \`a cohomologie constructible, et $n$ un entier. Supposons que, pour tout point $u$ de~$U$, on ait
$$
\prof_{u}(F)\geq n-\delta_{t}^{\prime}(u), \quad\text{o\`u } \delta_{t}^{\prime}(u)=\delta_{t}(u)-1.
$$
Alors le morphisme canonique
$$
\H^{i}(X, F)\to\H^{i}(Y, j^{\ast}F)$$
est bijectif pour $i<n-c-1$, injectif pour $i=n-c-1$.
\end{corollaire}

Soit $f\colon U\to S$ le morphisme canonique; il r\'esulte de
\Ref{XIV.4.2}, appliqu\'e en rempla\c cant $n$ par $n+1$, que l'on a
$$
\prof_{t}(\R_{!}f(F\sisi{|U}{_{|U}}))\geq n+1-c.
$$
La relation\pageoriginale pr\'ec\'edente montre que le morphisme canonique
$$
\H^{i}(S, \R_{!}f(F\sisi{|U}{_{|U}}))\to\H^{i}(S^{\prime}, \R_{!}f(F\sisi{|U}{_{|U}}))$$
est bijectif pour $i<n-c$, injectif pour $i=n-c$. Comme $\R_{!}f(F\sisi{|U}{_{|U}})$ est nul en dehors de $S^{\prime}$, on~a $\H^{i}(S, \R_{!}f(F\sisi{|U}{_{|U}}))\simeq\big(\underline{H}^{i}(\R_{!}f(F\sisi{|U}{_{|U}}))\big)_{t}=0$, et par suite
\begin{equation*} \label{eq:XIV.4.5.*} \tag{$*$}
\H^{i}(S^{\prime}, \R_{!}f(F\sisi{|U}{_{|U}}))=0\quad\text{pour }i<n-c.
\end{equation*}
Soient $g\colon X\to S^{\prime}$, $h\colon Y\to S^{\prime}$, $f^{\prime}\colon U\to S^{\prime}$ les morphismes canoniques. Il r\'esulte du triangle distingu\'e
$$
\xymatrix@C=0pt{ & \R h_{\ast}(j^{\ast}F)\ar[ld] & \\
\R_{!}f^{\prime}(F\sisi{|U}{_{|U}})\ar[rr] & & \R g_{\ast}(F)\ar[lu] }$$
que la condition \eqref{eq:XIV.4.5.*} \'equivaut au fait que le morphisme
$$
\H^{i}(S^{\prime}, \R g_{\ast}(F))\to\H^{i}(S^{\prime}, \R h_{\ast}(j^{\ast}F))$$
est bijectif pour $i<n-c-1$, injectif pour $i=n-c-1$. Comme ce morphisme s'identifie canoniquement au morphisme
$$
\H^{i}(X, F)\to\H^{i}(Y, j^{\ast}F),
$$
la conclusion en r\'esulte aussit\^{o}t.

\begin{corollaire}[Th\'eor\`eme de Lefschetz global] \label{XIV.4.6}
Soient\pageoriginale $S$ le spectre d'un corps, $X$ un sch\'ema propre sur $S$ et $U$ un ouvert de $X$ r\'eunion de $c+1$ ouverts affines. Soient~$Y$ un sous-sch\'ema ferm\'e de $X$, d'espace sous-jacent $X-U$, $j\colon Y\to X$ le morphisme canonique, $F$ un complexe de faisceaux de $\ZZ/m\ZZ$-modules sur $X$, \`a degr\'es born\'es inf\'erieurement et \`a cohomologie constructible et $n$ un entier. Supposons que, pour tout point $u$ de $U$, on ait
$$
\prof_{u}(F)\geq n-\dim\big(\overline{\{ u\} }\big).
$$
Alors le morphisme canonique
$$
\H^{i}(X, F)\to\H^{i}(Y, j^{\ast}F)$$
est bijectif pour $i<n-c-1$, injectif pour $i=n-c-1$.

Plus g\'en\'eralement, si $g\colon X\to S$ est un morphisme s\'epar\'e de type fini, les hypoth\`eses sur $S$, $U$, $Y$, $F$ \'etant les m\^emes que pr\'ec\'edemment, alors le morphisme canonique
$$
\H_{!}^{i}(X/S, F)\to\H_{!}^{i}(Y/S, j^{\ast}F)$$
(o\`u $\H_{!}^{i}$ d\'esigne la cohomologie \`a support propre, c'est-\`a-dire $\H_{!}^{i}(X\sisi{}{/S}, F)=\H^{i}(S, \R_{!}g(F))$) est bijectif pour $i<n-c-1$, injectif pour $i=n-c-1$.
\end{corollaire}

Le corollaire est un cas particulier de \Ref{XIV.4.4}, avec $T=S$.

Voici une r\'eciproque partielle \`a \Ref{XIV.4.3}:

\begin{proposition} \label{XIV.4.7}
Soient\pageoriginale $S$ un sch\'ema noeth\'erien, $f\colon U\to S$ un morphisme de type fini. Supposons qu'il existe un complexe dualisant $K$ sur $S$ et que $\R^{!}f(K)$ soit un complexe dualisant sur $U$. Soient $T$ une partie ferm\'ee de $S$ et $c$ un entier. Alors les conditions suivantes sont \'equivalentes:
\begin{enumeratei}
\item
Pour tout complexe de faisceaux de $\ZZ/m\ZZ$-modules $F$ sur $U$, \`a degr\'es born\'es inf\'erieurement et \`a cohomologie constructible, et pour tout entier $n$ tel que l'on ait, pour tout point $u$ de $U$,
$$
\prof_{u}(F)\geq n-\delta_{T}(u),
$$
on a
$$
\prof_{T}(\R_{!}f(F))\geq n-c.
$$

\item
Pour tout faisceau de $\ZZ/m\ZZ$-modules $G$ sur $U$, constructible, et pour tout point $t\in T$, on a
$$
\big(\R^{p}f_{\ast}(G)\big)_{\overline{t}}=0\quad\textrm{ pour }p>\delta(G, f, t)+c$$
(rappelons d'apr\`es (\textup{\SGA \sisi{X~}{}4~XIX~6.0}) que $\delta(G, f, t)=\sup\{ \delta_{t}(u)|t\in\overline{\{ u\} }\textrm{ et }G_{\overline{u}}\neq 0\}$).
\end{enumeratei}
\end{proposition}

N.B. La condition \textup{(ii)} est satisfaite en vertu de \Ref{XIV.3.1} si $f$ est s\'epar\'e et si $U$ est, localement sur $S$ pour la topologie \'etale, r\'eunion de $c+1$ ouverts affines sur $S$, donc \Ref{XIV.4.7} contient \Ref{XIV.4.3}\sfootnote{Du moins dans le cas \sisi{ou}{o\`u} $S$ admet localement un complexe dualisant, \sisi{p. ex.}{par exemple} $S$ immergeable localement dans un sch\'ema r\'egulier.}.

On peut \'evidemment supposer que $S$ est local et que $T$ est le point ferm\'e $t$ de~$S$. La d\'emonstration de \textup{(ii)}\ALORS\textup{(i)} est essentiellement identique\pageoriginale \`a la partie $2^\circ)$ de la d\'emonstration de \Ref{XIV.4.2}. Montrons rapidement que \textup{(i)}\ALORS \textup{(ii)}. Le th\'eor\`eme de dualit\'e locale (\SGA 5~I~4.3.2) appliqu\'e \`a $\D G$ montre que
$$
\D\big(\H_{\overline{u}}^{i}(\D G)\big) \simeq\big(\underline{H}^{-i-2\delta_{t}(u)}(G)\big)_{\overline{u}}.
$$
Comme $G$ est r\'eduit au degr\'e $0$, on~a donc $\H_{\overline{u}}^{i}(\D G)=0$ sauf peut-\^etre pour $i=-2\delta_{t}(u)$; plus pr\'ecis\'ement
\sisi{\begin{align*}
\prof_{u}(\D G)=-2\delta_{t}(u) &\;\textrm{ si }\; G_{\overline{u}}\neq 0, \\
\prof_{u}(\D G)=\infty &\;\textrm{ si }\; G_{\overline{u}}=0.
\end{align*}}
{\[
\prof_{u}(\D G)=
\begin{cases}
-2\delta_{t}(u) &\text{si } G_{\overline{u}}\neq 0, \\
\infty &\text{si } G_{\overline{u}}=0.
\end{cases}
\]}
Il en r\'esulte que l'on a, quel que soit $u\in U$:
$$
\prof_{u}(\D G)\geq -n-\delta_{t}(u).
$$
Il r\'esulte alors de l'hypoth\`ese \textup{(i)} que l'on a $\prof_{t}(\R_{!}f(\D G))\geq -n-c$. On transforme cette relation en utilisant l'isomorphisme $\R_{!}f(\D G)\simeq\D(\R f_{\ast}(G))$ (\SGA 5 I~1.12) et en appliquant le th\'eor\`eme de dualit\'e locale au point $t$; on obtient ainsi
$$
\big(\underline{H}^{i}(\R f_{\ast}(G)\big)_{\overline{t}}=0\quad\textrm{ pour }i>n+c,
$$
ce qui n'est autre que \textup{(ii)}.

\subsection{} \label{XIV.4.8} Les hypoth\`eses \'etant celles de \Ref{XIV.4.4} avec $g$ propre (\resp \Ref{XIV.4.5}, \resp \Ref{XIV.4.6} avec $g$ propre), si $V$ est un voisinage ouvert de $Y$ dans $X$, le\pageoriginaled morphisme
$$
\H^{i}(V, F)\to\H^{i}(Y, j^{\ast}F)$$
est bijectif pour $i<n-c-1$, injectif pour $i=n-c-1$. Si $\sisi{i}{\iota}\colon V\to X$ est le morphisme canonique, il suffit en effet pour le voir d'appliquer \Ref{XIV.4.4} (\resp \Ref{XIV.4.5}, \resp \Ref{XIV.4.6}) au complexe $\R \sisi{i}{\iota}_{\ast}(\sisi{F{|V}}{F_{|V}})$. On peut se poser la question de savoir si le morphisme pr\'ec\'edent est bijectif si $i=n-c-1$, injectif pour $i=n-c$. Il suffit \'evidemment que les hypoth\`eses soient v\'erifi\'ees quand on remplace $n$ par $n+1$; la proposition qui suit montre qu'il suffit d'un peu moins.

\begin{proposition} \label{XIV.4.9}
Soient $S$ un sch\'ema local excellent de caract\'eristique nulle, de point ferm\'e $t$ (\resp en plus des conditions pr\'ec\'edentes on suppose $S$ hens\'elien), $f\colon X\to S$ un sch\'ema propre sur $S$ (\resp propre sur $S-\{ \sisi{s}{t}\}$) et $U$ un ouvert de $X$ r\'eunion de $c+1$ ouverts affines. Soient $Y$ un sous-sch\'ema ferm\'e de $X$, d'espace sous-jacent $X-U$, $j\colon Y\to X$ le morphisme canonique, $F$ un complexe de faisceaux de $\ZZ/m\ZZ$-modules sur $X$, \`a degr\'es born\'es inf\'erieurement et \`a cohomologie constructible, et $n$ un entier. On suppose que l'on a, pour tout point $u$ de $U$,
$$
\prof_{u}(F)\geq\inf(n-1, n-\delta_{t}(u))\quad(\textrm{\resp }\;\prof_{u}(F)\geq\inf(n-1, n+1-\delta_{t}(u))).
$$
Alors pour tout voisinage ouvert $V$ de $Y$ dans $X$, le morphisme canonique\pageoriginale
$$
\H_{f^{-1}(t)}^{i}(V, F)\to\H_{f^{-1}(t)\cap Y}^{i}(Y, j^{\ast}F)\quad(\textrm{\resp }\;\H^{i}(V, F)\to\H^{i}(Y, j^{\ast}F))$$
est bijectif pour $i<n-c-2$ et injectif pour $i=n-c-2$. De plus, il existe un voisinage ouvert $V_{0}$ de $Y$ dans $X$, tel que, pour tout autre tel $V$ avec $V\subset V_{0}$, le morphisme canonique
$$
\H_{f^{-1}(t)\cap V}^{i}(V, F)\to\H_{f^{-1}(t)\cap Y}^{i}(Y, j^{\ast}F)\quad(\textrm{\resp }\;\H^{i}(V, F)\to\H^{i}(Y, j^{\ast}F))$$
soit bijectif pour $i<n-c-1$, injectif pour $i=n-c-1$.
\end{proposition}

\begin{proof}
Posons pour simplifier $\delta_{t}^{\prime}(u)=\delta_{t}(u)$ (\resp $\delta_{t}^{\prime}(u)=\delta_{t}(u)-1$). On d\'eduit de \Ref{XIV.4.8} la premi\`ere assertion de \Ref{XIV.4.9}, car les hypoth\`eses de \Ref{XIV.4.4} (\resp \Ref{XIV.4.5}) sont v\'erifi\'ees quand on y remplace $n$ par $n-1$. Elles le sont aussi pour $n$ lui-m\^eme, sauf aux points $u$ tels que $\delta_{t}^{\prime}(u)=0$. Or, pour un $u\in U$, dire que l'on a $\delta_{t}^{\prime}(u)=0$ \'equivaut \`a dire que $u$ est un point \emph{ferm\'e de $U_{t}$} (\resp un point \emph{ferm\'e de $X$}). Soit $E$ l'ensemble des points de $U$ tels que $\delta_{t}^{\prime}(u)=0$; montrons que, \emph{pour tous les points $u\in E$, sauf un nombre fini}, on~a $\prof_{u}(F)\geq n$. Soient $\overline{S}$ le localis\'e strict de $S$ en $t$, $S^{\prime}$ le compl\'et\'e de~$\overline{S}$, de point ferm\'e $t^{\prime}$, et consid\'erons le carr\'e cart\'esien
$$
\xymatrix{ U^{\prime}\ar[r]^-{h}\ar[d]_-{f^{\prime}} & U\ar[d]^-{f}\\
S^{\prime}\ar[r]^-{g} & S.}$$
Les hypoth\`eses de profondeur aux points de $U$ se conservent quand on remplace\pageoriginale $U$ par $U^{\prime}$ et $F$ par l'image inverse $F^{\prime}$ de $F$ sur $U^{\prime}$. Soient en effet $u^{\prime}\in U^{\prime}$ et $u=h(u^{\prime})$. Si $u\not\in E$, on~a la relation $\prof_{u}(F)\geq n-\delta_{t}^{\prime}(u)$, et il r\'esulte de \Ref{XIV.4.2.1} que ceci entra\^ine la relation $\prof_{u^{\prime}}(F^{\prime})\geq n-\delta_{t}^{\prime}(u^{\prime})$. Si $u\in E$, $u^{\prime}$ est un point ferm\'e de $U_{t}^{\prime}$ (\resp un point ferm\'e de $X^{\prime}=X\times_{S}S^{\prime}$), et, comme la fibre $U_{u}^{\prime}$ de $h$ en $u$ est de dimension z\'ero et $h$ r\'egulier, il r\'esulte de \Ref{XIV.1.16} que l'on a $\prof_{u^{\prime}}(F^{\prime})=\prof_{u}(F)\geq n-1$.

Soient alors $K$ un complexe dualisant sur $S^{\prime}$, normalis\'e au point ferm\'e $t^{\prime}$, et $\D F^{\prime}$ le dual de $F^{\prime}$ par rapport \`a $\R^{!}f(K)$. D'apr\`es \Ref{XIV.4.2.2}, les hypoth\`eses de profondeur \'etale aux points de $U^{\prime}$ se traduisent par les relations:
$$
\big(\underline{H}^{q}(\D F^{\prime})\big)_{\overline{u}^{\prime}}=0\quad\textrm{ pour }q>-n-\delta_{t}^{\prime}(u^{\prime})\quad(\textrm{\resp }q>-n-2-\delta_{t^{\prime}}^{\prime}(u^{\prime})),
$$
si $u^{\prime}$ n'est pas un point de $E^{\prime}=h^{-1}(E)$,
$$
\big(\underline{H}^{q}(\D F^{\prime})\big)_{\overline{u}^{\prime}}=0\quad\textrm{ pour }q>-(n-1)\quad(\textrm{\resp }q>-n-1), \textrm{ si }u^{\prime}\in E^{\prime}.
$$

Soit $G=\underline{H}^{-(n-1)}(\D F^{\prime})$ (\resp $G=\underline{H}^{-n-1}(\D F^{\prime})$); comme $G$ est un faisceau constructible, l'ensemble des points en lesquels la fibre g\'eom\'etrique est non nulle est un ensemble constructible (\SGA 4~IX~2.4~\textup{(iv)}); or par hypoth\`ese cet ensemble est contenu dans l'ensemble $E^{\prime}$ des points ferm\'es de $U_{t^{\prime}}^{\prime}$ (\resp des points de $U^{\prime}$ ferm\'es dans $X^{\prime}$); il r\'esulte donc de \Ref{XIV.4.9.1} ci-dessous que cet ensemble est r\'eduit \`a un nombre fini de points. Appliquant \Ref{XIV.4.2.2}, on voit que, pour tous les points de $E^{\prime}$, sauf un nombre fini, on~a $\prof_{u^{\prime}}(F^{\prime})\geq n$. Il en r\'esulte bien par \Ref{XIV.1.16} que, pour tous les\pageoriginale points de $E^{\sisi{\prime}{}}$ sauf un nombre fini, on a
$$
\prof_{u}(F)\geq n.
$$

Soit $V$ un voisinage ouvert de $Y$ dans $X$, contenu dans le compl\'ementaire dans $X$ de l'ensemble fini de\sisi{}{s} points $u$ de $E$ pour lesquels on~a $\prof_{u}(F)=n-1$. Si $\sisi{i}{\iota}\colon V\to X$ est l'immersion canonique, soit
$$F_{1}=\R \sisi{i}{\iota}_{\ast}(\sisi{F{|V}}{F_{|V}})\;;$$
alors $F_{1}$ est un complexe de faisceaux sur $X$, \`a cohomologie constructible (\SGA 4~XIX~5.1) et \`a degr\'es born\'es inf\'erieurement. Nous allons voir que, pour tout point $u$ de $U$, le complexe $F_{1}$ v\'erifie la relation
\begin{equation*} \label{eq:XIV.4.9.*} \tag{$*$}
\prof_{u}(F_{1})\geq n-\delta_{t}^{\prime}(u).
\end{equation*}
Si $u\in U\cap V$, on~a $\prof_{u}(F_{1})=\prof_{u}(F)$, et la
relation \eqref{eq:XIV.4.9.*} est v\'erifi\'ee par hypoth\`ese sur les
points de $U$ qui n'appartiennent pas \`a $E$; pour ces derniers,
elle est aussi v\'erifi\'ee d'apr\`es le choix de $V$. Enfin, si $u\in
U$ et $u\not\in V$, on~a d'apr\`es \Ref{XIV.1.6} g)
$\prof_{u}(F_{1})=\infty$. On applique alors \Ref{XIV.4.4} (\resp
\Ref{XIV.4.5}) en rempla\c cant $F$ par $F_{1}$; on obtient le
r\'esultat annonc\'e, compte tenu du fait que l'on a, quel que soit
$i$:
$$
\H_{f^{-1}(t)}^{i}(X, \R \sisi{i}{\iota}_{\ast}(\sisi{F{|V}}{F_{|V}}))\simeq\H_{f^{-1}(t)\cap V}^{i}(V, F)\quad(\textrm{\resp }\H^{i}(X, \R \sisi{i}{\iota}_{\ast}(\sisi{F{|V}}{F_{|V}}))\simeq\H^{i}(V, F).
$$
\skipqed
\end{proof}

\begin{sublemme} \label{XIV.4.9.1} Un ensemble constructible $E$ contenu dans l'ensemble des points ferm\'es d'un sch\'ema $X$ noeth\'erien est r\'eduit \`a un nombre fini de points\pageoriginale.
\end{sublemme}

En effet, $E$ est r\'eunion fini\sisi{}{e} d'ensemble\sisi{}{s} de la forme $U\cap\complement V$, o\`u $U$ et $V$ sont des ouverts de $X$; par hypoth\`ese tous les points de $U\cap\complement V$ sont des points maximaux de cet ensemble, donc ils sont en nombre fini.

\section{Profondeur g\'eom\'etrique} \label{XIV.5}

Pour appliquer en pratique \Ref{XIV.4.2} et ses corollaires, il faut disposer d'un crit\`ere commode qui permette de v\'erifier les hypoth\`eses de profondeur \'etale aux points de $U$. Nous allons donner un tel crit\`ere, en utilisant le th\'eor\`eme de Lefschetz local~\sisi{(\Ref{XIV.4.5})}{\Ref{XIV.4.5}}.

\subsection{} \label{XIV.5.1} \sisi{Soient}{Soit} $A$ un anneau local noeth\'erien; quand nous parlerons de la profondeur \'etale de $A$, il s'agira de la profondeur au point ferm\'e. Nous allons introduire une notion de \og profondeur g\'eom\'etrique de $A$\fg, et utiliser \Ref{XIV.4.5} pour la comparer \`a la profondeur \'etale $\profet(A)$.

\begin{proposition} \label{XIV.5.2} Soit $A$ un anneau local noeth\'erien; supposons que $A$ soit isomorphe \`a un quotient d'un anneau local r\'egulier $B$ par un id\'eal $I$ (c'est vrai par exemple lorsque $A$ est complet, en vertu du th\'eor\`eme de Cohen (\textup{$\EGA 0_{\textup{IV}}$~19.8.8})). Soit $q$ le nombre minimal de g\'en\'erateurs de $I$; alors le nombre $\dim(B)-q$ est ind\'ependant du choix de $B$.
\end{proposition}

Le nombre minimal de g\'en\'erateurs de $I$ est aussi \'egal au rang du $k$-espace vectoriel $I\otimes_{B}k$, o\`u $k$ d\'esigne le corps r\'esiduel de $A$. On se ram\`ene tout de suite au cas o\`u $A$ est complet, car on~a $\widehat{A}\simeq\widehat{B}/\widehat{I}$ avec $\dim \widehat{B}=\dim B$ et $\rg_{k}(I\otimes_{B}k)=\rg_{k}(\widehat{I}\otimes_{\widehat{B}}k)$; pour la m\^eme raison on\pageoriginale peut supposer que l\sisi{es}{'} anneau\sisi{x}{} $B$ \sisi{son}{es}t complet\sisi{s}{}. Soient $B$ et $B^{\prime}$ deux anneaux locaux r\'eguliers complets, $f\colon B\to A$, $f^{\prime}\colon B^{\prime}\to A$ deux homomorphismes surjectifs et $I=\Ker(f)$, $I^{\prime}=\Ker(f^{\prime})$. On doit montrer que
$$
\dim B-\rg_{k}(I\otimes_{B}k)=\dim B^{\prime}-\rg_{k}(I^{\prime}\otimes_{B^{\prime}}k).
$$

Pla\c cons-nous d'abord dans le cas o\`u l'on a une factorisation de
la forme
$$
\xymatrix@C=15pt{ B\ar[rr]^-{f}\ar[rd]_-{g} & & A\\
& B^{\prime}\ar[ru]_-{f^{\prime}} & }$$
avec $g$ surjectif. Soit $J=\Ker(g)$; alors $J\subset I$ et $I/J=I^{\prime}$. Puisque $B^{\prime}$ est r\'egulier, $\dim(B^{\prime})=\dim(B)-\rg_{k}(J\otimes_{B}k)$ et $J$ est engendr\'e par des \'el\'ements faisant partie d'un syst\`eme r\'egulier de param\`etres de $B$. Il en r\'esulte que l'on a la suite exacte
$$0\to J\otimes_{B}k\to I\otimes_{B}k\to J/I\otimes_{B^{\prime}}k\to 0,
$$
et par suite
$$
\dim B-\rg_{k}(I\otimes_{B}k)=\dim B-\rg_{k}(J\otimes_{B}k)-\rg_{k}(J/I\otimes_{B^{\prime}}k)=\dim B^{\prime}-\rg_{k}(I^{\prime}\otimes_{B^{\prime}}k).
$$

Le cas g\'en\'eral se ram\`ene au pr\'ec\'edent; pour le voir, il suffit de montrer que l'on peut trouver un anneau local r\'egulier complet $B^{\prime\prime}$ et des homomorphismes surjectifs $g\colon B^{\prime\prime}\to B$ et $g^{\prime}\colon B^{\prime\prime}\to B^{\prime}$, rendant commutatif le diagramme\pageoriginale
\begin{equation*} \label{eq:XIV.5.2.*} \tag{$*$}
\begin{array}{c}
\xymatrix@C=10pt@R=15pt{ & B\ar[rd]^-{f} & \\
B^{\prime\prime}\ar[ru]^-{g}\ar[rd]_-{g^{\prime}} & & A\\
& B^{\prime}.\ar[ru]_-{f^{\prime}} &
}
\end{array}
\end{equation*}
Or, si $W$ est un anneau de Cohen de corps r\'esiduel $k$, on~a un morphisme local $W\to A$ qui se rel\`eve \`a $B$ et $B^{\prime}$ (\EGA IV~19.8.6), de sorte que l'on a le diagramme commutatif
$$
\xymatrix@R=15pt@C=10pt{ & B\ar[rd] & \\
W\ar[ru]\ar[rd] & & A\\
& B^{\prime}.\ar[ru] & }$$
On peut trouver des entiers $n$ et $n^{\prime}$ et des morphismes surjectifs \hbox{$h\colon W\llbracket T_{1},\dots, T_{n}\rrbracket\!\to\! B$} et $h^{\prime}\colon W\llbracket T_{1}^{\prime},\dots, T_{n^{\prime}}^{\prime}\rrbracket\to B^{\prime}$ qui soient des morphismes de $W$-alg\`ebres ($\EGA 0_{\textup{IV}}$ 19.8.8); si l'on pose alors $B^{\prime\prime}=W\llbracket T_{1},\dots, T_{n}, T_{1}^{\prime},\dots, T_{n^{\prime}}^{\prime}\rrbracket$ et si l'on d\'efinit $g$ et~$g^{\prime}$ comme des morphismes de $W$-alg\`ebres tels que
$$g(T_{i})=h(T_{i}), \quad g(T_{i}^{\prime})=b_{i}, \quad g^{\prime}(T_{i})=b_{i}^{\prime}, \quad g^{\prime}(T_{i}^{\prime})=h^{\prime}(T_{i}^{\prime}),
$$
o\`u $b_{i}$ (\resp $b_{i}^{\prime}$) est un \'el\'ement de $B$ (\resp de $B^{\prime}$) qui rel\`eve \sisi{$h^{\prime}.f^{\prime}(T_{i}^{\prime})$ (\resp $h.f(T_{i})$)}{$(f^{\prime}\circ h^{\prime})(T_{i}^{\prime})$ (\resp $(f\circ\nobreak h)(T_{i})$)}, le diagramme \eqref{eq:XIV.5.2.*} est bien commutatif.

La proposition \Ref{XIV.5.2} justifie la d\'efinition suivante:

\begin{definition} \label{XIV.5.3}
Soient\pageoriginale $A$ un anneau local noeth\'erien, $\widehat{A}$ son compl\'et\'e, qui est donc isomorphe au quotient d'un anneau local r\'egulier complet $B$ par un id\'eal $I$; si $q$ est le nombre minimal de g\'en\'erateurs de $I$, on appelle \emph{profondeur g\'eom\'etrique} de $A$ le nombre
$$
\profgeom(A)=\dim B-q.
$$
\end{definition}

\begin{proposition} \label{XIV.5.4} Soit $A$ un anneau local noeth\'erien. Alors on a
$$
\profgeom(A)\leq\dim A,
$$
et on~a l'\'egalit\'e si et seulement si $A$ est un\sisi{}{e} intersection compl\`ete.
\end{proposition}

On peut supposer $A$ complet. Soit alors $A=B/I$, o\`u $B$ est un anneau local r\'egulier complet et $I$ un id\'eal de $B$. Si $(x_{1},\dots, x_{q})$ est un syst\`eme minimal de g\'en\'erateurs de $I$, on~a $\dim A\geq\dim B-q$, et dire que $\dim A=\dim B-q$ \'equivaut \`a dire que $(x_{1},\dots, x_{q})$ fait partie d'un syst\`eme de param\`etres de $B$ ($\EGA 0_{\textup{IV}}$~16.3.7); la proposition en r\'esulte imm\'ediatement.

\begin{proposition} \label{XIV.5.5} Soient $A$ et $A^{\prime}$ des anneaux locaux noeth\'eriens, $f\colon A\to A^{\prime}$ un homomorphisme local. Supposons que $f$ soit plat et que, d\'esignant par $k$ le corps r\'esiduel de $A$, $A^{\prime}\otimes_{A}k$ soit un corps, extension s\'eparable de $\sisi{K}{k}$. Alors on a
$$
\profgeom(A)=\profgeom(A^{\prime}).
$$
\end{proposition}

Quitte\pageoriginale \`a remplacer $A$ et $A^{\prime}$ par leurs compl\'et\'es, on peut supposer $A$ et $A^{\prime}$ complets (il r\'esulte de ($\EGA 0_{\textup{III}}$~10.2.1) que l'hypoth\`ese de platitude est conserv\'ee et cela est \'evident pour les autres hypoth\`eses). Soit alors $A=B/I$, o\`u $B$ est un anneau local r\'egulier et $I$ un id\'eal de $B$. Comme $A^{\prime}$ est formellement lisse sur $A$ ($\EGA 0_{\textup{IV}}$~19.8.2), il r\'esulte de ($\EGA 0_{\textup{IV}}$~19.7.2) que l'on peut trouver un anneau local noeth\'erien complet $B^{\prime}$ et un homomorphisme local $B\to B^{\prime}$, tel que $B^{\prime}$ soit un $B$-module plat et que l'on ait $B^{\prime}\otimes_{B}A\simeq A^{\prime}$. On a donc $A^{\prime}\simeq B^{\prime}/IB^{\prime}$; de plus l'anneau $B^{\prime}$ est r\'egulier; en effet, si $\mm$ est l'id\'eal maximal de $B$, $\mm B^{\prime}$ est l'id\'eal maximal de $B^{\prime}$, et, puisque $\mm$ est engendr\'e par une suite r\'eguli\`ere par d\'efinition de \og r\'egulier\fg, \sisi{$B^{\prime}$}{$\mm B^{\prime}$} est engendr\'e par une suite $B^{\prime}$-r\'eguli\`ere ($\EGA 0_{\textup{IV}}$~15.1.14). Comme on~a \'evidemment $\dim B=\dim B^{\prime}$, et comme $I$ et $IB^{\prime}$ ont m\^eme nombre minimal de g\'en\'erateurs, l'assertion en r\'esulte.

\refstepcounter{subsection}\label{XIV.5.6}
\begin{enonce*}{Th\'eor\`eme 5.6\ndemark}\ndetext{Illusie a montr\'e depuis l'in\'egalit\'e $\prof_x(\ZZ^/\ell^\nu\ZZ)\geq \profgeom_x(X/S)-\delta(x)+1$ pour $x$ point de $X$ un sch\'ema de type fini sur un trait $S$ de caract\'eristique r\'esiduelle premi\`ere \`a $\ell$, et $\nu\geq 1$. Si $S$ est de caract\'eristique nulle, c'est une cons\'equence du th\'eor\`eme~\Ref{XIV.5.6}; voir (Illusie~L., {\og Perversit\'e et variation\fg}, \emph{Manuscripta Math.} \textbf{112} (2003), p\ptbl 271-295).}
Soit $A$ un anneau local excellent de caract\'eristique nulle. Alors on~a
$$
\profet(A)\geq\profgeom(A).
$$
\end{enonce*}

On peut supposer $A$ strictement local complet, puisque la profondeur g\'eom\'etrique et la profondeur \'etale se conservent par passage \`a l'hens\'elis\'e strict et au compl\'et\'e d'apr\`es \Ref{XIV.5.5} et \Ref{XIV.1.16}. Soit $A\simeq B/I$, o\`u $B$ est un anneau local r\'egulier complet, et soit $(f_{1},\dots, f_{q})$ un syst\`eme minimal de g\'en\'erateurs de l'id\'eal $I$. On a donc\pageoriginale
$$
\pi=\profgeom(A)=\dim B-q.
$$
Consid\'erons l'immersion ferm\'ee
$$Y=\Spec A\to X=\Spec B,
$$
et \sisi{soient}{soit} $U=X-Y=\bigcup_{1\leq i\leq q}X_{f_{i}}$. Si $a$ d\'esigne le point ferm\'e de $X$, on doit montrer que, pour tout nombre premier $p$, on a
$$
\H_{a}^{i}(Y, \ZZ/p\ZZ)=0\quad\textrm{ pour }i<\pi.
$$
Comme $B$ est r\'egulier excellent, on~a $\profet(B)=2\dim X$ (\cf \Ref{XIV.1.10}) et par suite $\H_{a}^{i}(X, \ZZ/p\ZZ)=0$ pour $i<2\dim X$. Il suffit donc, pour d\'emontrer le th\'eor\`eme, de prouver que le morphisme
\begin{equation*} \label{eq:XIV.5.6.*} \tag{$*$}
\H_{a}^{i}(X, \ZZ/p\ZZ)\to\H_{a}^{i}(Y, \ZZ/p\ZZ)
\end{equation*}
est bijectif pour $i<\pi$. On applique pour cela le th\'eor\`eme de Lefschetz local \sisi{(\Ref{XIV.4.5})}{\Ref{XIV.4.5}} avec $n=\pi+q-1$, $c=q-1$ donc $n-c=\pi$. Notons que $U=X-Y$ est r\'eunion des $q$ ouverts affines $X_{f_{i}}$. Montrons que l'on a, pour tout point $u$ de $U$:
$$
\profet_{u}(X)\geq\pi+q-1-\dim\big(\overline{\{ u\} }\big)=\dim\Oo_{X, u}$$
(o\`u $\overline{\{ u\} }$ d\'esigne l'adh\'erence de $u$ \emph{dans} $X-\overline{\{ a\} }$). En effet il r\'esulte de \Ref{XIV.1.10} que l'on a
$$
\profet_{u}(X)=2\dim\Oo_{X, u}\geq\dim\Oo_{X, u}.
$$
En\pageoriginaled utilisant \Ref{XIV.4.5}, on voit que \eqref{eq:XIV.5.6.*} est bijectif pour $i<\pi$, ce qui ach\`eve la d\'emonstration du th\'eor\`eme.

\begin{corollaire} \label{XIV.5.7} Soient $S$ le spectre d'un corps de caract\'eristique nulle (\resp un sch\'ema local hens\'elien excellent de caract\'eristique nulle), $f\colon X\to S$ un sch\'ema propre sur $S$ (\resp sur $S-\{ s\}$). Soient $U$ une r\'eunion de $c+1$ ouverts de $X$, \emph{affines}, $Y$ un sous-sch\'ema ferm\'e d'espace sous-jacent $X-U$, $n$ et $m$ des entiers $>0$. On suppose que, pour tout point $u$ de $U$, on a
$$
\profgeom(\Oo_{X, u})\geq n-\dim\big(\overline{\{ u\} }\big)$$
($\overline{\{ u\} }$ adh\'erence de $u$ dans $X$). Alors le morphisme canonique
$$
\H^{i}(X, \ZZ/m\ZZ)\to\H^{i}(Y, \ZZ/m\ZZ)$$
est bijectif pour $i<n-c-1$, injectif pour $i=n-c-1$.
\end{corollaire}

On applique \Ref{XIV.4.5} et \Ref{XIV.4.6}. Les hypoth\`eses de profondeur \'etale aux points de $U$ sont v\'erifi\'ees car on~a d'apr\`es \Ref{XIV.5.6} \sisi{$$
\profet_{u}(X)\geq\profgeom\Oo_{X, u}\geq n-\dim\big(\overline{\{ u\} }\big).
$$
}{$$
\profet_{u}(X)\geq\profgeom(\Oo_{X, u})\geq n-\dim\big(\overline{\{ u\} }\big).
$$
}
\enlargethispage{\baselineskip}%
\section{Questions ouvertes} \label{XIV.6}

\subsection{} \label{XIV.6.1} On peut se demander si l'implication \textup{(ii)} \ALORS \textup{(i)} de \Ref{XIV.4.2} est valable plus g\'en\'eralement pour des faisceaux $F$ de torsion, pas n\'ecessairement annul\'es par un entier $m$ donn\'e et pas n\'ecessairement constructibles\pageoriginale. Dans le cas o\`u $S$ n'est pas de caract\'eristique nulle, il semble possible que cette implication reste valable, m\^eme pour les faisceaux de $p$-torsion ($p$ la caract\'eristique r\'esiduelle). Enfin il n'est pas clair non plus que l'hypoth\`ese $S$ excellent ne puisse \^etre lev\'ee.

\subsection{} \label{XIV.6.2}
Soient $X$ un sch\'ema propre sur un corps $k$ ou bien le compl\'ementaire du point ferm\'e d'un sch\'ema local hens\'elien, et $j:Y \to X$ un sous-sch\'ema ferm\'e de $X$, dont le compl\'ementaire $U$ est affine. Alors, si $F$ est un faisceau d'ensembles sur $X$ ou un faisceau en groupes non n\'ecessairement commutatifs, les \'enonc\'es \Ref{XIV.4.5} ou \Ref{XIV.4.6} et \Ref{XIV.4.9} ont encore un sens pour un tel $F$, \`a condition de se borner \`a des petites valeurs de $n$. Si $u$ est un point de $U$, on d\'esigne par $\bar{u}$ un point g\'eom\'etrique au-dessus de~$u$, par $X(\bar{u})$ le localis\'e strict de $X$ en $\bar{u}$ et par $F_{\bar{u}}$ la fibre de $F$ en $\bar{u}$. Alors, en faisant \'eventuellement certaines hypoth\`eses sur $X$ et sur $F$, par exemple en supposant $X$ excellent (\'eventuellement de caract\'eristique nulle, ou d'\'egale caract\'eristique en utilisant la r\'esolution des singularit\'es) et $F$ ind-fini (ou au besoin m\^eme $\LL$-ind-fini avec $\LL$ premier \`a la caract\'eristique de $X$), on aimerait d\'emontrer les \'enonc\'es suivants:
\begin{enumeratea}
\item
Soit $F$ un faisceau d'ensembles (\resp un faisceau en groupes) et supposons que, pour tout point $u$ de $U$, on ait
$$F_{\bar{u}} \to H^0(X(\bar{u})-\bar{u}, F)\text{ injectif si } \dim(\overline{\{u\}}) \leq 1$$
(c'est-\`a-dire, pour un tel $u$, on~a $\prof_u(F) \geq 1$). Alors, quand $V$ parcourt l'ensemble des voisinages ouverts de $Y$, le morphisme canonique\pageoriginale
$$
\varinjlim_V H^0(V, F) \to H^0(Y, j^*F)
$$
est bijectif (\resp on~a la conclusion pr\'ec\'edente et de plus le morphisme $\varinjlim_V H^1(V, F) \to H^1(Y, j^* F)$ est injectif). Si $F$ est constructible, on peut remplacer les $\varinjlim$ par la cohomologie de $V$ pour $V$ \og assez petit\fg.
\item
Soit $F$ un faisceau d'ensembles (\resp un faisceau en groupes) et supposons que, pour tout point $u$ de $U$, on ait $\prof_u(F) \geq 2 -\dim(\overline{\{u\}})$, ce qui se traduit aussi par les relations
\begin{gather*}
F_{\bar{u}}\to H^0(X(\bar{u})-\bar{u}, F)\text{ est bijectif si } \dim(\overline{\{u\}})=0\\
F_{\bar{u}}\to H^0(X(\bar{u})-\bar{u}, F)\text{ est injectif si } \dim(\overline{\{u\}})=1.
\end{gather*}
Alors le morphisme canonique
$$H^0(X, F) \to H^0(Y, j^*F)$$
est bijectif (\resp on~a la conclusion pr\'ec\'edente et de plus le morphisme $H^1(X, F) \to H^1(Y, j^*F)$ est injectif).
\item
Soit $F$ un faisceau en groupes ind-fini. Supposons que, pour tout point $u$ de $U$, on ait
\begin{gather*}
F_{\bar{u}}\to H^0(X(\bar{u})-\bar{u}, F)\text{ bijectif si } \dim(\overline{\{u\}})=0\text{ ou }1,\\
F_{\bar{u}}\to H^0(X(\bar{u})-\bar{u}, F)\text{ injectif si } \dim(\overline{\{u\}})=2.
\end{gather*}
Alors\pageoriginale, quand $V$ parcourt l'ensembles des voisinages ouverts de $Y$, les morphismes canoniques
$$
\varinjlim_V H^0(V, F) \to H^0(Y, j^*F) \text{ et } \varinjlim_V H^1(V, F) \to H^1(Y, j^*F)
$$
sont bijectifs. Si $F$ est constructible, on peut remplacer les $\varinjlim$ par la cohomologie de~$V$ pour $V$ assez petit.
\item
Soit $F$ un faisceau en groupes. Supposons que, pour tout point $u$ de $U$, on ait $\prof_u(F) \geq 3-\dim(\overline{\{u\}})$, ce qui se traduit aussi par les conditions
\begin{gather*}
F_{\bar{u}}\to H^0(X(\bar{u})-\bar{u}, F)\text{ bijectif, et } H^1(X(\bar{u})-\bar{u}, F)=0 \text{ si } \dim(\overline{\{u\}})=0,\\
F_{\bar{u}}\to H^0(X(\bar{u})-\bar{u}, F)\text{ bijectif si } \dim(\overline{\{u\}})=1,\\
F_{\bar{u}}\to H^0(X(\bar{u})-\bar{u}, F)\text{ injectif si } \dim(\overline{\{u\}})=2.
\end{gather*}
Alors les morphismes canoniques
$$H^0(X, F) \to H^0(Y, j^*F)\text{ et } H^1(X, F) \to H^1(Y, j^*F)$$
sont bijectifs.
\end{enumeratea}

Comme indication en faveur de ces \'enonc\'es\refstepcounter{toto}\nde{\label{MiR}tous les \'enonc\'es de~\Ref{XIV.6.2}, mises \`a part les variantes constructibles, ont \'et\'e prouv\'es par Mme Raynaud; voir (Raynaud~M., {\og Th\'eor\`emes de Lefschetz en cohomologie des faisceaux coh\'erents et en cohomologie \'etale. Application au groupe fondamental\fg}, \emph{Ann. Sci. \'Ec. Norm. Sup. (4)} \textbf{7} (1974), p\ptbl 29--52, corollaire III.1.3).}, signalons \Ref{XIII}~\Ref{XIII.2.1}, \Ref{X}~\Ref{X.3.4} et \Ref{XII}~\Ref{XII.3.5}. Signalons que, gr\^ace \`a l'argument de \Ref{XIV.4.8} et \Ref{XIV.4.9}, l'\'enonc\'e a) (\resp c)) r\'esulterait de b) (\resp d)).

\subsection{} \label{XIV.6.3}
Il r\'esulterait de d) l'\'enonc\'e suivant analogue \`a \Ref{XIV.5.6}: si $A$ est\pageoriginale un anneau local noeth\'erien (\'eventuellement excellent) et si l'on a $\prof \geom (A) \geq 3$, alors on~a $\prof \hop (A) \geq 3$. Pour voir ceci, on r\'ealise $Y'=\Spec A$ comme ferm\'e d'un sch\'ema local r\'egulier $X'=\Spec B$, dont le compl\'ementaire est r\'eunion de $q$ ouverts affines, avec la relation $\dim B -q=\prof \geom(A)$. On a, pour tout point $x$ de $X'$, si $n=\dim B$, $\prof \hop_x(X) \geq \inf(3, n-\dim(\overline{\{x\}}))$ (\cf \Ref{XIV.1.11}) et l'on d\'eduit de d) que ceci entra\^ine $\prof \hop_y(Y') \geq \inf(3, n-q-\dim(\overline{\{y\}}))$, pour tout point $y$ de $Y'$. Le r\'esultat s'obtient alors en prenant pour $y$ le point ferm\'e de $Y'$.

\subsection{} \label{XIV.6.4}
Un variante de \Ref{XIV.4.2}, tout au moins de l'implication \textup{(ii)} \ALORS \textup{(i)}, doit encore \^etre valable dans le cas analytique complexe, \`a condition de travailler avec des faisceaux \og analytiquement constructibles\fg (\cf \Ref{XIII}); la d\'emonstration serait analogue \`a celle de \Ref{XIV.4.2}, en utilisant la th\'eorie de la dualit\'e de J.-L. Verdier. Notons par contre que, pour l'analogue analytique complexe des variantes non commutatives signal\'ees dans \Ref{XIV.6.2}, on ne dispose pas m\^eme d'une m\'ethode d'attaque pour les \'enonc\'es concernant le groupe fondamental sugg\'er\'es par les r\'esultats des expos\'es \Ref{X}, \Ref{XII}, \Ref{XIII}, rappel\'es \`a fin de \Ref{XIV.6.2}. Les m\'ethodes du S\'eminaire semblent en effet irr\'em\'ediablement li\'ees au cas des rev\^etements \emph{finis} (qui peuvent \^etre \'etudi\'es en termes de faisceaux coh\'erents d'\sisi{Alg\`ebres}{alg\`ebres}).

\chapter*{Index des notations}

\begin{tabular}{p{9cm}p{3cm}}
$\Gamma_{Z}, \ \underline{\Gamma}_{Z}$\ \dotfill & \Ref{I}~\Ref{I.1}\\
$i_{!}, \ i^{!}$\ \dotfill & \Ref{I}~\Ref{I.1}\\
$\ZZ_{Z, X}$\ \dotfill & \Ref{I}~\Ref{I.1.6}\\
$H^{i}_{Z}(X, F), \underline{H}^{i}_{Z}(F)$\ \dotfill & \Ref{I}~\Ref{I.2}\\
$H^{i}_{J}(M), H^{i}((f), M), H^{i}(M)$\ \dotfill & \Ref{IV}~\Ref{IV.5.4}, \Ref{V}~\Ref{V.2} (p\ptbl\sisi{4}{\pageref{pV2.4}})\\
$\Ass M$, $\Supp M$, $\Ann \ M\sisi{}{,} \rr(\aaa)$\ \dotfill & \Ref{III}~\Ref{III.1.1}\\
$\prof_{I}(M), \prof_{Y}(F)$\ \dotfill & \Ref{III}~\Ref{III.2.3}, \Ref{III.2.8}\\
$\Ab$\ \dotfill & \Ref{IV}~\Ref{IV.1} (p\ptbl\sisi{1}{\pageref{pIV.1}})\\
$\Hom^{.}(F^{.}, G^{.})$\ \dotfill & \Ref{V}~\Ref{V.1.1}\\
$\Ext^{i}_{Z}(X;F, G), \underline{\Ext}^{i}_{Z}(F, G)$\ \dotfill & \Ref{VI}~\Ref{VI.1.1}\\
${\Et}(X), \LL(Z)$\ \dotfill & \Ref{X} (p\ptbl\sisi{1}{\pageref{pX.1}})\\
$\Lef(X, Y), \Leff(X, Y)$\ \dotfill & \Ref{X}~\Ref{X.2} (p\ptbl\sisi{2}{\pageref{pX.2}})\\
$\underline{P}(Z), \Pic(Z)$\ \dotfill & \Ref{XI} (p\ptbl\sisi{1}{\pageref{pXI.1}})\\
$\Pi_{i}^{x}(X)$\ \dotfill & \Ref{XIII} (p\ptbl\sisi{15}{\pageref{pXIII.15}})\\
$D^{+}(X)$\ \dotfill & \Ref{XIV}~\Ref{XIV.1.0}\\
$\prof_{Y}(F) \geq n$\ \dotfill & \Ref{XIV}~\Ref{XIV.1.2}\\
$\prof_{Y}^{\LL}(X) \geq n$\ \dotfill & \Ref{XIV}~\Ref{XIV.1.2}\\
$\profet_{Y}(X)$\ \dotfill & \Ref{XIV}~\Ref{XIV.1.2}\\
$\LL$\ \dotfill & \Ref{XIV}~\Ref{XIV.1.2}\\
$\prof_{x}(F) \geq n$\ \dotfill & \Ref{XIV}~\Ref{XIV.1.7}\\
$\prof_{x}^{\LL}(X)\geq n$\ \dotfill & \Ref{XIV}~\Ref{XIV.1.7}\\
$\prof\ \hop^{\LL}_{x}(X) \geq 3$\ \dotfill & \Ref{XIV}~\Ref{XIV.1.7}\\
$\profet_{x}(X)$\ \dotfill& \Ref{XIV}~\Ref{XIV.1.7}\\
$\delta_{t}(x), \ \delta (n)$\ \dotfill & \Ref{XIV}~\Ref{XIV.2.1}\\
$H^{i}_{Z!}(X/S, F)$\ \dotfill & \Ref{XIV}~\Ref{XIV.4.0}
\end{tabular}

\chapter*{Index terminologique}

\begin{footnotesize}
\noindent
\begin{tabular}{lll}
Auslander-Buchsbaum (th\'eor\`eme de -) & & \Ref{XI}~\Ref{XI.3.13}\\
comparaison (th\'eor\`eme de -) & comparison theorem & \Ref{IX}~\Ref{IX.1.1}, \Ref{XII}~\Ref{XII.2.1}\\
dualit\'e locale (th\'eor\`eme de -) & local duality theorem & \Ref{V}~\Ref{V.2.1}\\
dualit\'e projective (th\'eor\`eme de -) & projective duality theorem & \Ref{XII}~\Ref{XII.1.1} \\
enveloppe injective & injective envelope & \Ref{IV} (p\ptbl\sisi{15}{\pageref{pIV.13}})\\
excision (th\'eor\`eme d' -) & excision theorem & \Ref{I}~\Ref{I.2.2}, \Ref{VI}~\Ref{VI.1.3}\\
existence (th\'eor\`eme d' -) & existence theorem & \Ref{IX}~\Ref{IX.2.1}, \Ref{XII}~\Ref{XII.3.1}\\
extension essentielle & essential extension & \Ref{IV} (p\ptbl\sisi{15}{\pageref{pIV.13}})\\
finitude (th\'eor\`eme de -) & finiteness theorem & \Ref{VIII}~\Ref{VIII.2.1}, \Ref{XII}~\Ref{XII.1.5}\\
foncteur dualisant & dualizing functor & \Ref{IV}~\Ref{IV.4.1}, \Ref{IV.4.2}\\
forme r\'esidu & residue form & \Ref{IV}~\Ref{IV.5.5}\\
g\'eom\'etriquement factoriel, \resp & geometrically factorial (\resp &\\
g\'eom. parafactoriel (anneau local -) & geom. parafactorial) local ring & \Ref{XIII} (p\ptbl20, 24)\\
groupes d'homotopie locale & local homotopy groups & \Ref{XIII} (p\ptbl\sisi{15}{\pageref{pXIII.15}})\\
Gysin (homomorphisme de -) & Gysin homomorphism & \Ref{I} (p\ptbl\sisi{13}{\pageref{hgysin}})\\
Hartshorne (th\'eor\`eme de -) & & \Ref{III}~\Ref{III.3.6}\\
Lefschetz (condition de -) & Lefschetz condition & \Ref{X}~\Ref{X.2} (p\ptbl\sisi{2}{\pageref{pX.2}})\\
Lefschetz (condition de - effective) & effective Lefschetz condition & \Ref{X}~\Ref{X.2} (p\ptbl\sisi{2}{\pageref{pX.2}})\\
Lefschetz (th\'eor\`eme de - affine) & affine Lefschetz theorem & \Ref{XIV}~\Ref{XIV.3.1}\\
module dualisant & dualizing module & \Ref{IV}~\Ref{IV.4.1}\\
\sisi{module orthogonal d'un module}{orthogonal d'un sous-module} & \sisi{module orthogonal of a module}{orthogonal of a submodule} & \Ref{IV}~\Ref{IV.5}\\
parafactoriel (couple - de pr\'e\-sch\'emas) & parafactorial pair of preschemes & \Ref{XI}~\Ref{XI.3.1}\\
parafactoriel (anneau local -) & parafactorial local ring & \Ref{XI}~\Ref{XI.3.2}\\
profondeur & depth & \Ref{III}~\Ref{III.2.3}\\
profondeur \'etale & etale depth & \Ref{XIV}~\Ref{XIV.1.2}\\
profondeur g\'eom\'etrique d'un anneau & geometrical depth of a noetherian &\\
local noeth\'erien & local ring & \Ref{XIV}~\Ref{XIV.5.3}\\
profondeur homotopique & homotopical depth & \Ref{XIII}~\Ref{XIII.6} d\'ef\ptbl1 (p\ptbl\sisi{4}{\pageref{XIII.6.1}})\\
profondeur homotopique (pour $\LL$) & homotopical depth (with respect to $\LL$) & \Ref{XIV}~\Ref{XIV.1.2}\\
profondeur homotopique rectifi\'ee & rectified homotopical depth & \Ref{XIII} d\'ef\ptbl2 (p\ptbl\sisi{5}{\pageref{XIII.6.2}})\\
\end{tabular}
\newpage
\begin{tabular}{lll}
pur (anneau local -) & pure local ring & \Ref{X}~\Ref{X.3.2} \\
pure (couple - de pr\'esch\'emas) & pure pair of preschemes & \Ref{X}~\Ref{X.3.1}\\
puret\'e (th\'eor\`eme de - de Zariski-Nagata) & purity theorem of Zariski-Nagata & \Ref{X}~\Ref{X.3.4}\\
puret\'e cohomologique (th\'eor\`eme de -) & cohomological purity theorem & \Ref{XIV}~\Ref{XIV.1.11}\\
r\'egulier ($M$-r\'egulier) & $M$-regular & \Ref{III.2.1}\\
Samuel (conjecture de -) & & \Ref{XI}~\Ref{XI.3.14}\\
semi-puret\'e cohomologique (th\'eor\`eme de -) & cohomological semi-purity theorem & \Ref{XIV}~\Ref{XIV.1.10}\\
strictement local (anneau -) & strictly local ring & \Ref{XIII}~\Ref{XIII.6} (p\ptbl\sisi{23}{\pageref{pXIII.6}})\\
Zariski-Nagata (th\'eor\`eme de puret\'e de -) & purity theorem of Zariski-Nagata & \Ref{X}~\Ref{X.3.4}
\end{tabular}
\end{footnotesize}

\backmatter
\end{document}